\documentclass{amsproc}
\usepackage{enumitem}

\usepackage{empheq}

\usepackage[linesnumbered,ruled,vlined]{algorithm2e}
\SetCommentSty{mycommfont}
\SetKwInput{KwInput}{Input}         
\SetKwInput{KwOutput}{Output}              %

\usepackage{subcaption} %

\usepackage{hyperref}

\usepackage{mathtools,cases} %

\usepackage{stmaryrd} %

\newtheorem{theorem}{Theorem}[section]

\theoremstyle{definition}

\newtheorem{example}[theorem]{Example}

\theoremstyle{remark}
\newtheorem{remark}[theorem]{Remark}

\numberwithin{equation}{section}

\newcommand{\EE}{{\mathbb E}}
\newcommand{\RR}{{\mathbb R}}
\newcommand{\TT}{{\mathbb T}}
\newcommand{\ZZ}{{\mathbb Z}}

\newcommand{\bF}{{\mathbf F}}
\newcommand{\bG}{{\mathbf G}}
\newcommand{\bJ}{{\mathbf J}}
\newcommand{\bN}{{\mathbf N}}

\newcommand{\cA}{{\mathcal A}}
\newcommand{\cB}{{\mathcal B}}
\newcommand{\cE}{{\mathcal E}}
\newcommand{\cF}{{\mathcal F}}
\newcommand{\cG}{{\mathcal G}}

\newcommand{\cL}{{\mathcal L}}
\newcommand{\cM}{{\mathcal M}}
\newcommand{\cP}{{\mathcal P}}
\newcommand{\cU}{{\mathcal U}}
\newcommand{\cW}{{\mathcal W}}
\newcommand{\cZ}{{\mathcal Z}}

\newcommand{\Id}{\mathrm{Id}}

\newcommand{\matM}{\mathbf{M}}
\newcommand{\matB}{\mathbf{B}}

\newcommand{\ctrl}{\alpha}
\newcommand{\ctrldim}{k}
\newcommand{\dom}{\mathcal{Q}}
\newcommand{\domT}{\mathcal{Q}_{T}}

\DeclareMathOperator*{\argmax}{arg\,max}
\DeclareMathOperator*{\argmin}{arg\,min}

\DeclareMathOperator*{\arginf}{arg\,inf}
\newcommand{\prox}{\mathrm{prox}}
\newcommand{\diver}{\mathrm{div}}
\newcommand{\grad}{\nabla}
\newcommand{\indic}{\mathbf{1}}

\makeindex

\begin{document}

\title[Numerical methods for MFG and MFC]{Numerical Methods for Mean Field Games \\ and Mean Field Type Control}

\author{Mathieu Lauri{\`e}re}
\address{Princeton University, Operations Research and Financial Engineering (ORFE) department, Sherrerd Hall, Princeton, New Jersey, U.S.A.}
\curraddr{}
\email{lauriere@princeton.edu}
\thanks{The author was partially supported by NSF DMS-1716673 and ARO W911NF-17-1-0578.}


\date{}

\begin{abstract}
Mean Field Games (MFG) have been introduced to tackle games with a large number of competing players. Considering the limit when the number of players is infinite, Nash equilibria are studied by considering the interaction of a typical player with the population's distribution. The situation in which the players cooperate corresponds to Mean Field Control (MFC) problems, which can also be viewed as optimal control problems driven by a McKean-Vlasov dynamics. These two types of problems have found a wide range of potential applications, for which numerical methods play a key role since most models do not have analytical solutions. In these notes, we review several aspects of numerical methods for MFG and MFC. We start by presenting some heuristics in a basic linear-quadratic setting. We then discuss numerical schemes for forward-backward systems of partial differential equations (PDEs), optimization techniques for variational problems driven by a Kolmogorov-Fokker-Planck PDE, an approach based on a monotone operator viewpoint, and stochastic methods relying on machine learning tools. 
\end{abstract}

\maketitle

\setcounter{tocdepth}{2}
\tableofcontents

These lecture notes are to supplement the AMS short course on mean field games on January 13--14, 2020, and are meant to be an introduction to numerical methods for mean field games. See also the other proceedings of the AMS Short Course on Mean Field Games~\cite{carmona-AMS-tmp,delarue-AMS-tmp,gravesmalhame-AMS-tmp,lacker-AMS-tmp,ramanan-AMS-tmp}.

\section{\bf Introduction}

Mean field games (MFGs for short) have been introduced to study differential
games with an infinite number of players, under indistinguishability and symmetry assumptions. This theory has been introduced by J.-M. Lasry and P-L. Lions~\cite{MR2269875,MR2271747,MR2295621}, and independently by Caines, Huang, and Malham{\'e} under the name of Nash Certainty Equivalence principle~\cite{MR2346927,MR2352434}. While MFGs correspond to non-cooperative games and focus on Nash equilibria, the cooperative counterpart has been studied under the name mean field type control (or mean field control, or MFC for short) or optimal control of McKean-Vlasov dynamics~\cite{MR3395471,MR3134900}. For more details, we refer to Lions' lectures at Coll\`ege de France~\cite{PLL-CDF} and the notes by Cardaliaguet~\cite{Cardaliaguet-2013-notes}, as well as the books by Bensoussan, Frehse and Yam~\cite{MR3134900}, Gomes, Pimentel and Voskanyan~\cite{MR3559742}, Carmona and Delarue~\cite{MR3752669,MR3753660}, and the surveys by Gomes and Sa{\'u}de~\cite{MR3195844} and Caines, Huang and Malham{\'e}~\cite{CHM-2017-survey}. Mean field games and mean field control problems have attracted a surge of interest and have found numerous applications, from economics to the study of crowd motion or epidemic models. 
For concrete applications, it is often important to obtain reliable quantitative results on the solution. Since very few mean field problems have explicit or semi-explicit solutions, the introduction of numerical methods and their analysis are crucial steps in the development of this field. 

The goal of these notes is to provide an introduction to numerical aspects of MFGs and MFC problems. We discuss both numerical schemes and associated algorithms, and illustrate them with numerical examples. In order to alleviate the presentation while  conveying the main ideas, the mathematical formalism is not carried out in detail. The interested reader is referred to the references provided in the text for a more rigorous treatment (see in particular the lecture notes~\cite{MR3135339,achdoulauriere-2020-mfg-numerical}).

\subsection{Outline}

The rest of these notes is organized as follows. In the remainder of this section, we introduce a general framework for the type of problems we will consider.  Linear-quadratic problems are used in Section~\ref{AMS-num-sec:LQ} as a testbed to introduce several numerical strategies that can be applied to more general problems, such as fixed point iterations, fictitious play iterations and Newton iterations. Beyond the linear-quadratic structure, optimality conditions in terms of PDE systems are discussed in Section~\ref{AMS-num-sec:PDE-scheme} and two numerical schemes to discretize these PDEs are presented: a finite-difference scheme and a semi-Lagrangian scheme. In Section~\ref{AMS-num-sec:sol-strat}, we present some strategies to solve the discrete schemes and provide numerical illustrations, including to crowd motion with congestion effects. Section~\ref{AMS-num-sec:optim-variational} focuses on MFC problems and MFGs with a variational structure. For such problems, optimization techniques can be used and we present two of them: the Alternating Direction Method of Multipliers and a primal-dual algorithm. For MFGs satisfying a monotonicity condition, a different technique, based on a monotonic flow, is discussed in Section~\ref{AMS-num-sec:monotone-op}. The last two sections present techniques which borrow tools from machine learning. Section~\ref{AMS-num-sec:NNapprox} presents methods relying on neural network approximation to solve stochastic optimal control problems or PDEs. We apply these methods to MFC and to the (finite state MFG) master equation respectively. Section~\ref{AMS-num-sec:modelfree} discusses model-free methods, \textit{i.e.}, methods which learn the solution by trial and error, without knowing the full model. We conclude in Section~\ref{sec:conclusion} by mentioning other methods and research directions.

\subsection{Definition of the problems and notation}

Let $T$ be a fixed time horizon, let $d$ and $\ctrldim$ be integers and let $\dom \subseteq \RR^d$ be a spatial domain, which will typically be either the whole space $\RR^d$ or the unit torus $\TT^d$. We will use the notation $\domT = [0,T] \times \dom$, and $\langle \cdot, \cdot \rangle$ for the inner product of two vectors of compatible sizes.

Let $f: \dom \times L^2(\dom) \times \RR^\ctrldim \to \RR, (x,m,\ctrl) \mapsto f(x,m,\ctrl)$ and $g: \dom \times L^2(\dom) \to \RR, (x,m) \mapsto g(x,m)$ be respectively a running cost and a terminal cost. Let $\sigma \ge 0$ be a constant parameter for the volatility of the state's evolution.  Let $b:  \dom \times L^2(\dom) \times \RR^\ctrldim \to \RR^d, (x,m,\ctrl) \mapsto b(x,m,\ctrl)$ be a drift function.  These functions could be allowed to also depend on time at the expense of heavier notation. Here, $x,m$ and $\ctrl$ play respectively the role of the state of the agent, the mean-field term, (\textit{i.e.}, the population's distribution), and the control used by the agent. In general, the mean-field term is a probability measure, but here, for simplicity, we will assume that this probability measure has a density which is in $L^2(\dom)$. 

We consider the following \index{mean field game}mean field game. A (mean field) Nash equilibrium consists in a flow of probability densities $\hat{m}: \domT \to \RR$ and a feedback control $\hat{\ctrl}: \domT\to \RR^\ctrldim$ satisfying the following two conditions:
\begin{enumerate}
	\item $\hat{\ctrl}$ minimizes $J^{MFG}_{\hat{m}}$ where, for $m \in L^2(\domT)$,
\begin{align}
\label{AMS-num-eq:def-J-MFG}
	J^{MFG}_{m}: \ctrl \mapsto  \EE \left[\int_0^T f(X_t^{m, \ctrl}, m(t, \cdot), \ctrl(t,X_t^{m, \ctrl}) ) dt + g(X_T^{m, \ctrl}, m(T,\cdot)) \right]
\end{align}
under the constraint that the  process $X^{m, \ctrl} = (X_t^{m, \ctrl})_{t \ge 0}$ solves the stochastic differential equation (SDE)
\begin{equation}
\label{AMS-num-eq:dyn-X-general-MFG}
	d X_t^{m, \ctrl} = b(X_t^{m, \ctrl}, m(t,\cdot), \ctrl(t, X_t^{m, \ctrl})) dt + \sigma d W_t, \qquad t \ge 0,
\end{equation}
where $W$ is a standard $d$-dimensional Brownian motion, and $X_0^{m, \ctrl}$ has distribution with density $m_0$ given;
	\item For all $t \in [0,T]$, $\hat{m}(t,\cdot)$ is the probability density of the law of $X_t^{\hat{m}, \hat{\ctrl}}$.
\end{enumerate}
In the definition~\eqref{AMS-num-eq:def-J-MFG} of the cost function, the subscript $m$ is used to emphasize the dependence on the mean-field flow, which is fixed when an infinitesimal agent performs their optimization. The second condition ensures that if all the players use the control $\hat{\ctrl}$ computed in the first point, the law of their individual states is indeed $\hat{m}$.

Using the same running cost, terminal costs, drift function and volatility, we can look at the corresponding \index{mean field control}mean field control (MFC for short) problem. This problem corresponds to a social optimum and is phrased as an optimal control problem. It can be interpreted as a situation in which all the agents cooperate to minimize the average cost. The goal is to find a feedback control $\ctrl^*: \domT\to \RR^\ctrldim$ minimizing
\begin{align}
\label{AMS-num-eq:def-J-MFC}
	J^{MFC}: \ctrl \mapsto \EE \left[\int_0^T f(X_t^{\ctrl}, m^{\ctrl}(t,\cdot), \ctrl(t,X_t^{\ctrl}) ) dt + g(X_T^{\ctrl}, m^{\ctrl}(T,\cdot)) \right]
\end{align}
where $m^{\ctrl}(\cdot, t)$ is the probability density of the law of $X_t^{\ctrl}$, under the constraint that the process $X^{\ctrl} = (X_t^{\ctrl})_{t \ge 0}$ solves the SDE
\begin{equation}
\label{AMS-num-eq:dyn-X-general-MFC}
	d X_t^{\ctrl} = b(X_t^{\ctrl}, m^{\ctrl}(t,\cdot), \ctrl(t, X_t^{\ctrl})) dt + \sigma d W_t, \qquad t \ge 0,
\end{equation}
and $X_0^{\ctrl}$ has distribution with density $m_0$ given. 

Besides the interpretation as a social optimum for a large population of cooperative players, MFC problems also arise in risk management~\cite{MR2784835} or in optimal control with a cost involving a conditional expectation~\cite{MR4133380,achdoulaurierelions2020optimal}.

Let $m^* = m^{\ctrl^*}$ denote the optimal mean-field flow. Then, we have, for any MFG equilibrium $(\hat\ctrl, \hat m )$:
$$
	J^{MFC}(\ctrl^*) = J^{MFG}_{m^*}(\ctrl^*) \leq J^{MFG}_{\hat{m}}(\hat{\ctrl}).
$$
In general the inequality is strict, which leads to the notion of price of anarchy as we will illustrate below.

While we will mostly focus on the problems as introduced above, the ergodic setting is discussed in Section~\ref{AMS-num-sec:monotone-op} and finite state MFGs are considered in~\S~\ref{AMS-num-sec:NN-finite-master}.

\clearpage

\section{\bf Warm-up: Linear-Quadratic problems}
\label{AMS-num-sec:LQ}

We start by considering the important subclass of problems in which, on the one hand, the dynamics is linear in the state, the control and the mean of the state (and possibly the mean of the control), and, on the other hand, the cost is quadratic in these variables. Then, the MFG and MFC optimal controls $\hat\ctrl$ and $\ctrl^*$ can be characterized by a system of backward ODEs, which is in general coupled with a forward ODE for the evolution of the mean of the state.\footnote{In some cases, a suitable parametrization leads to a decoupled ODE system; see \textit{e.g.}~\cite{gravesmalhame-AMS-tmp} for an example.}

\subsection{Problem definition and MFG solution}
\label{AMS-num-sec:LQ-pb}
We borrow the following example from~\cite[Chapter 6]{MR3134900}. For ease of presentation, we restrict our attention to the one-dimensional case, \textit{i.e.}, $d=\ctrldim=1$. Consider
\begin{align*}
	f(x,m,\ctrl)	&= \frac{1}{2} \left[ Q x^2 + \bar{Q} \left(x - S \int_{\dom} \xi m(\xi) d\xi \right)^2 + C \ctrl^2 \right]
	\\
	g(x,m) 	&= \frac{1}{2} \left[ Q_T x^2 + \bar{Q}_T \left(x - S_T \int_{\dom} \xi m(\xi) d\xi \right)^2 \right]
	\\
	b(x,m,\ctrl) 	&= Ax + \bar{A} \int_{\dom} \xi m(\xi) d\xi + B\ctrl \, ,
\end{align*}
where $Q, C, \bar{Q}, Q_T, \bar{Q}_T$ are non-negative constants, and $A, \bar{A}, S, S_T$ and $B$ are constants. Let $\nu = \tfrac{1}{2} \sigma^2$.  We consider that the initial distribution is a normal $\mathcal{N}(\bar x_0, \sigma_0^2)$ for some $\bar{x}_0 \in \RR$ and $\sigma_0>0$.

Under suitable conditions on these coefficients, the MFG for the above model has a unique solution $(\hat{\ctrl}, \hat{m})$ which satisfies the following. The proof relies on dynamic programming and on a suitable ansatz for the value function. See \textit{e.g.}~\cite[Chapter 6]{MR3134900} for more details. The solution is given by:
 \begin{subequations}
     \begin{empheq}[left=\empheqlbrace]{align*}
		\int_\dom \xi \hat{m}(t,\xi) d \xi &= z_t,
		\\
		\hat{\ctrl}(t,x) &= -B(p_t x + r_t) / C,
		\\
		J^{MFG}_{\hat{m}}(\hat{\ctrl}) &= \int_\dom u(0,\xi) m_0(\xi) d\xi = \frac{1}{2} p_0(\sigma_0^2 + \bar x_0^2) + r_0 \bar x_0 + s_0,
		\\
		u(t,x) &= \frac{1}{2}p_t x^2 + r_t x + s_t,
     \end{empheq}
   \end{subequations}
where $(z, p, r, s)$ solve the following system of ordinary differential equations (ODEs):
    \begin{subequations}
     \begin{empheq}[left=\empheqlbrace\,\,]{align}
     \label{AMS-num-eq:LQ-ODE-z}
     \frac{d z}{dt} &= (A+\bar{A} - B^2 C^{-1}p_t) z_t - B^2 C^{-1} r_t, 			&& z_0 = \int_{\dom} \xi m_0(\xi) d\xi,
     \\
     \label{AMS-num-eq:LQ-ODE-P}
     -\frac{d p}{dt} &= 2Ap_t - B^2C^{-1}p_t^2 + Q + \bar{Q}, 				&& p_T = Q_T + \bar{Q}_T,
     \\
     \label{AMS-num-eq:LQ-ODE-r}
     - \frac{d r}{dt} &= (A - B^2 C^{-1}p_t )r_t + (p_t \bar{A} - \bar{Q} S) z_t, 	&& r_T = - \bar{Q}_T S_T z_T,
     \\
     \label{AMS-num-eq:LQ-ODE-s}
     - \frac{d s}{dt} &= \nu p_t - \frac{1}{2} B^2 C^{-1} r_t^2 + r_t \bar{A} z_t + \frac{1}{2} S^2 \bar{Q} z_t^2, && s_T = \frac{1}{2} \bar{Q}_T S_T^2 z_T^2.
     \end{empheq}
   \end{subequations} In this system, $z$ represents the mean of the population's distribution whereas $r$ (together with $p$) characterizes the best response.   In fact, $u(t,x) = \frac{1}{2}p_t x^2 + r_t x + s_t$ is the value function of an infinitesimal player when the population is in the Nash equilibrium. 
The last equation admits an explicit solution, in terms of $z,p$ and $r$. The second equation can be solved independently of the other ones. It is a Riccati equation and, under suitable conditions, admits a unique positive (symmetric if $d>1$) solution. The first and the third equations are coupled. Note that the equation for $z$ is forward in time whereas the equation for $r$ is backward in time. This forward-backward structure prevents the use of a simple time-marching method to solve numerically the system. This obstacle is at the heart of numerical methods for MFG and MFC problems. Using this LQ example to provide the main ideas, we present below a few strategies to tackle forward-backward systems. In the rest of this section, we assume that $p$ solving~\eqref{AMS-num-eq:LQ-ODE-P} is given and we focus on the system~\eqref{AMS-num-eq:LQ-ODE-z}--\eqref{AMS-num-eq:LQ-ODE-r}.

\subsection{Time discretization}
\label{AMS-num-sec:LQ-ODE-time-discr}

In order to solve this ODE system numerically, we first partition the interval $[0,T]$ into $N_T$ subintervals, where $N_T$ is a positive integer. Let $\Delta t = 1/N_T$ and $t_i = i \times \Delta t$ for $i=0,\dots,N_T$. The functions of time $(z_t)_{t \in [0,T]}$, $(p_t)_{t \in [0,T]}$, $(r_t)_{t \in [0,T]}$  and $(s_t)_{t \in [0,T]}$ are approximated respectively by vectors $Z = (Z^n)_{n=0,\dots,N_T}$, $P = (P^n)_{n=0,\dots,N_T}$, $R = (R^n)_{n=0,\dots,N_T}$ and $S = (S^n)_{n=0,\dots,N_T}$. The ODE system~\eqref{AMS-num-eq:LQ-ODE-z}--\eqref{AMS-num-eq:LQ-ODE-r} is then replaced by a finite-difference system. Focusing on~\eqref{AMS-num-eq:LQ-ODE-z} and~\eqref{AMS-num-eq:LQ-ODE-r}, we consider the following system: for $n \in \{0,\dots,N_T-1\},$
\begin{subequations}
     \begin{empheq}[left=\empheqlbrace\,\,]{align}
     \label{AMS-num-eq:LQ-ODE-z-discrete}
     &\frac{Z^{n+1} - Z^{n}}{\Delta t} = (A+\bar{A} - B^2 C^{-1}P^{n}) Z^{n+1} - B^2 C^{-1} R^{n}, 			
     \\
     \notag
    &Z^0 = \bar x_0,
     \\
     \label{AMS-num-eq:LQ-ODE-r-discrete}
     &- \frac{R^{n+1}-R^{n}}{\Delta t} = (A - B^2 C^{-1} P^n ) R^{n} + (P^n \bar{A} - \bar{Q} S) Z^{n+1}, 	
     \\
     \notag
      &\hfill R^{N_T} = - \bar{Q}_T S_T Z^{N_T}.
     \end{empheq}
   \end{subequations}

Note that for each equation (considered separately) the scheme is semi-implicit, since the first one is forward in time and the second in backward in time. This finite-difference system can be rewritten in a matrix form: 
\begin{equation}
\label{AMS-num-eq:LQ-ODE-matrix-form}
	\matM 	
	\begin{pmatrix}
		Z
		\\
		R
	\end{pmatrix}
	+ \matB 
	=
	0\, ,
\end{equation}
where $\matM \in \RR^{(N_T+1) \times (N_T+1)}$ and $\matB \in \RR^{(N_T+1)}$ are defined in order to take into account the dynamics as well as the initial and terminal conditions. One can thus obtain $(Z,R)$ directly by solving this linear system. However, this is specific to the LQ setting. As we will see in the next section, forward-backward PDE systems appearing in MFG are generally not linear. Hence we present in the rest of this section several solution strategies which do not exploit the linear structure and will be useful in a more general setting. The LQ problem is simply used as an illustration to provide the main ideas behind these methods.

\subsection{Picard (fixed point) iterations}

\begin{figure}[h]
	\begin{subfigure}{.45\columnwidth}
		\centering
		\includegraphics[width=\columnwidth]{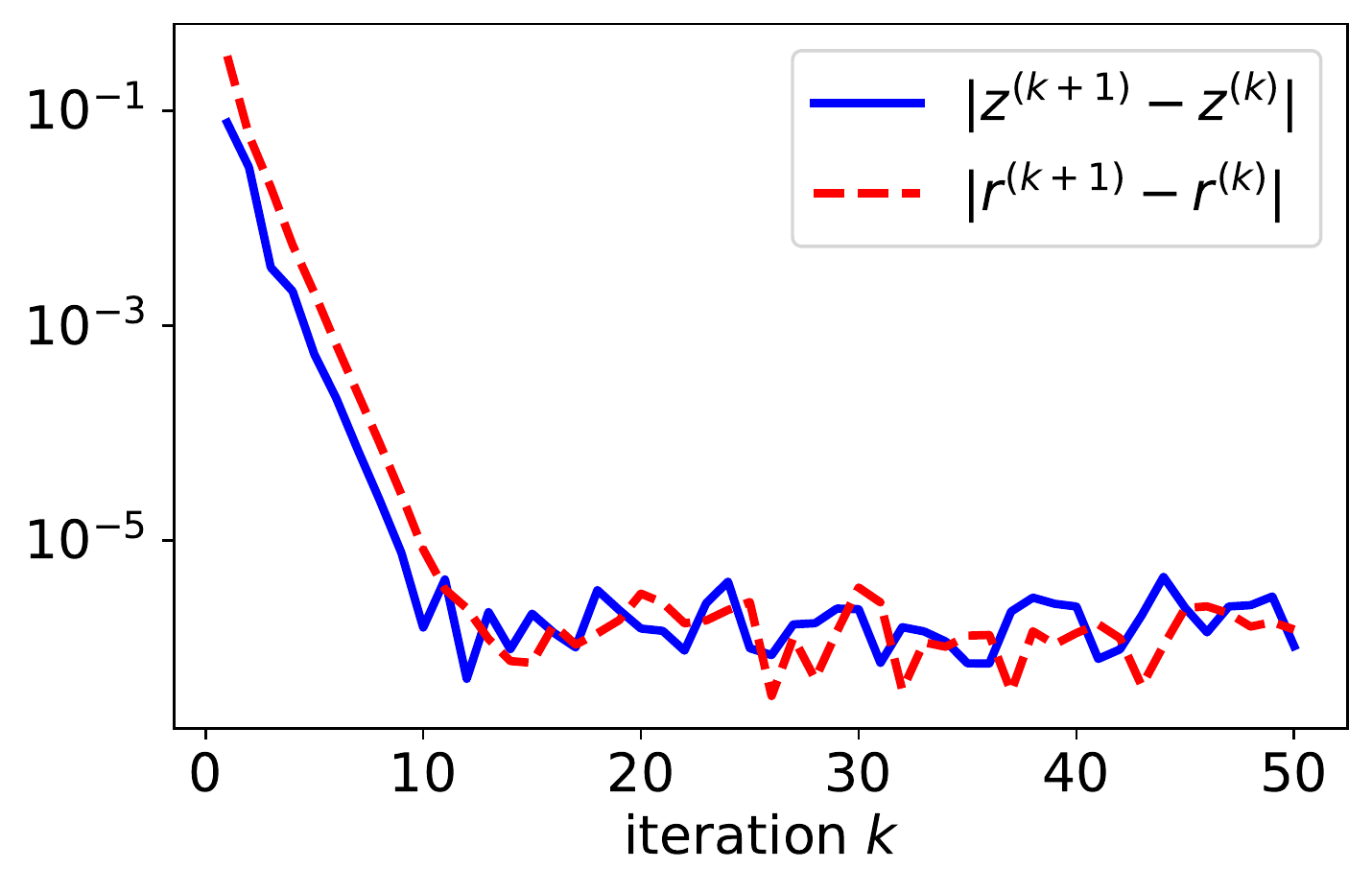}
	\end{subfigure}%
	\begin{subfigure}{.45\columnwidth}
		\centering 
		\includegraphics[width=\columnwidth]{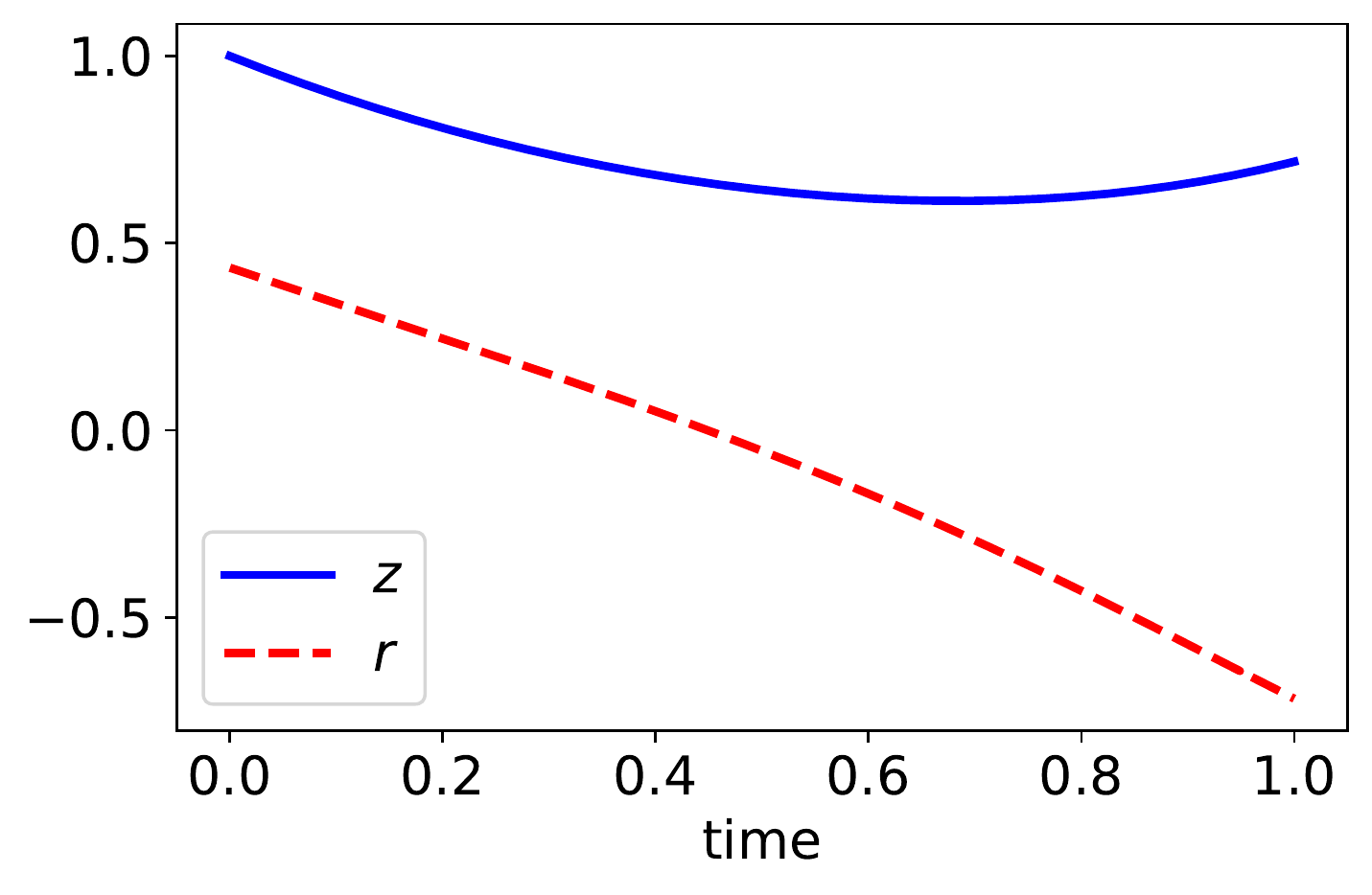}
	\end{subfigure}
	\caption{Picard iterations without damping for Test case~1 (see Table~\ref{AMS-num-tab:params-LQ-base}): $L^2$ difference between two successive iterates (left)  and $z^{(50)}$ and $r^{(50)}$ (right). The curves related to $z$ and $r$ are respectively in blue (full line) and in red (dashed line). } 
 	\label{AMS-num-fig:LQ-ODE-Picard-converge}
\end{figure}

To alleviate the presentation, we keep the ODE notation, \textit{i.e.},~\eqref{AMS-num-eq:LQ-ODE-z}\&\eqref{AMS-num-eq:LQ-ODE-r}, although for the implementation, we use their finite-difference counterpart, \textit{i.e.},~\eqref{AMS-num-eq:LQ-ODE-z-discrete}\&\eqref{AMS-num-eq:LQ-ODE-r-discrete}. A first approach consists in solving alternatively each ODE. Starting from an initial guess $z^{(0)}$, solve~\eqref{AMS-num-eq:LQ-ODE-r} in which $z$ is replaced by $z^{(0)}$. Denoting $r^{(1)}$ the solution, solve~\eqref{AMS-num-eq:LQ-ODE-z} in which $r$ is replaced by $r^{(1)}$. Then repeat these steps, plugging $z^{(1)}$ in~\eqref{AMS-num-eq:LQ-ODE-r} and so on. This procedure defines a map $\varphi$ such that $\varphi(z^{(0)}) = z^{(1)}$. If this map is a strict contraction, then the above iterations converge. A pseudo-code is given in Algorithm~\ref{AMS-num-algo:LQ-ODE-Picard} where (plain) Picard iterations correspond to the case $\delta(\cdot) \equiv 0$, and an example is displayed in Figure~\ref{AMS-num-fig:LQ-ODE-Picard-converge}. Here we used the values of parameters described in Table~\ref{AMS-num-tab:params-LQ-base} for Test case~1. Instead of fixing a priori a number of iterations, we can also consider a stopping criterion of the form:
$$
	\|r^{(\mathtt{k}+1)} - r^{(\mathtt{k})}\| < \varepsilon, \quad \hbox{ and } \quad \|z^{(\mathtt{k}+1)} - z^{(\mathtt{k})}\| < \varepsilon,
$$
for some threshold $\varepsilon>0$.

{\renewcommand{\arraystretch}{1.2}
\begin{table}[!ht]
\begin{center}
	\begin{tabular}{c|c|c|c|c|c|c}
				\hline
			Parameters & $Q, R ,S ,S_T, A, B,\sigma, \bar x_0$ & $\bar Q$  & $\bar Q_T$   & $Q_T$ & $\bar A$ & $\sigma_0$	\\
			\hline
			 Test case 1 & $1$ & $1$ & $1$ & $1$ & $1$ & $0.2$ \\
			\hline
			 Test case 2 & $1$ & $1$ & $2.45$ & $1$ & $1$  & $0.2$ \\
			\hline
			 Test case 3 & $1$ & $0$ & $0$ & $1$ & $[0,20]$ & $0.2$ \\
			 \hline
			 Test case 4 & $1$ & $0$ & $[0,20]$ & $1$ & $1$ &  $0.2$ \\
			 \hline
			 Test case 5 & $1$ & $0$ & $1$ & $[0,20]$ & $1$ &  $0.2$ \\
				 \hline 
				\end{tabular}
\caption{Parameters values for the LQ model introduced in \S~\ref{AMS-num-sec:LQ-pb}. The test cases are discussed in the text.}
\label{AMS-num-tab:params-LQ-base}
	\end{center}
\end{table}
}

\begin{figure}[h]
	\begin{subfigure}{.33\columnwidth}
		\centering
		\includegraphics[width=\columnwidth]{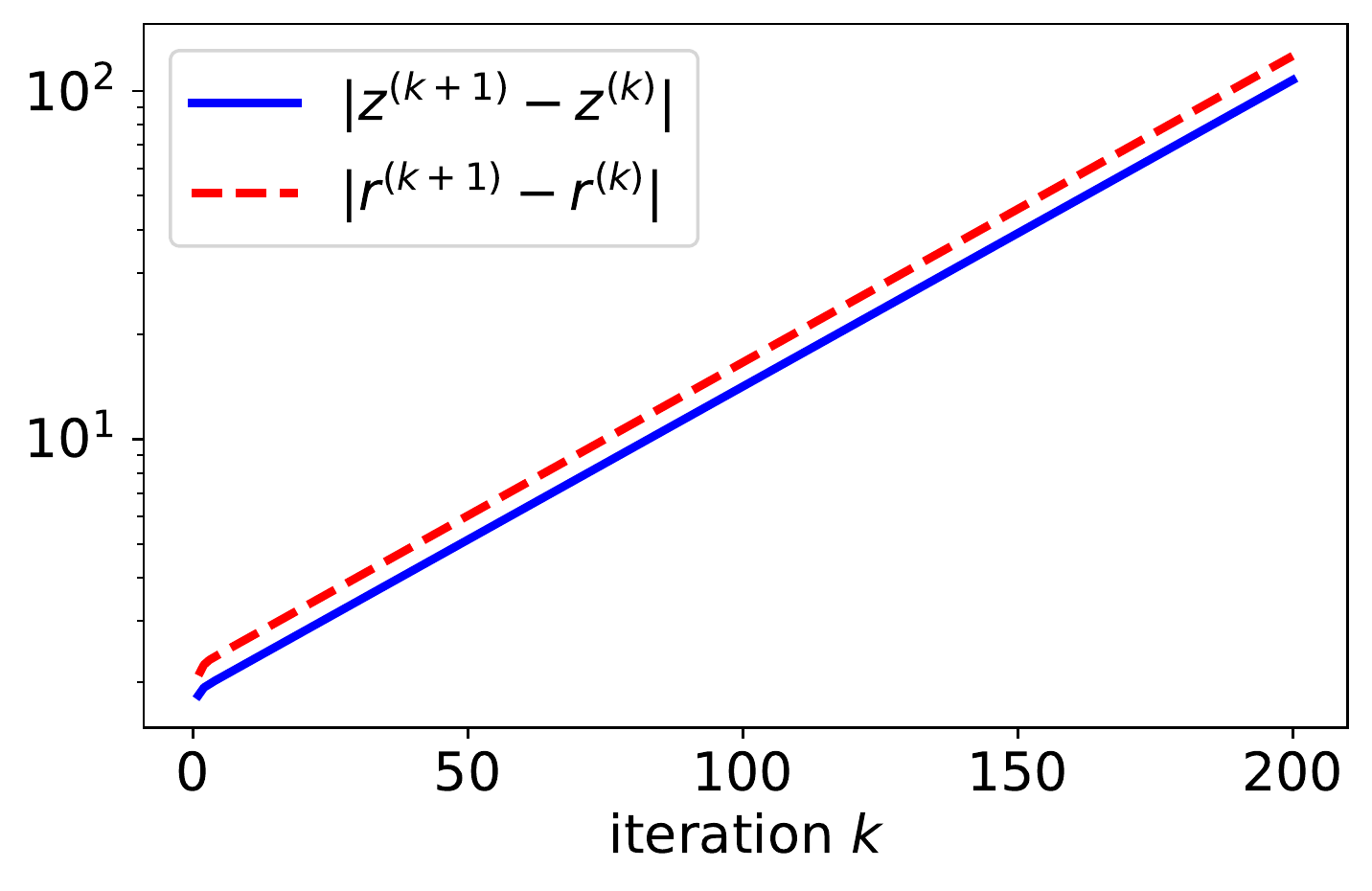}
	\end{subfigure}%
	\begin{subfigure}{.33\columnwidth}
		\centering 
		\includegraphics[width=\columnwidth]{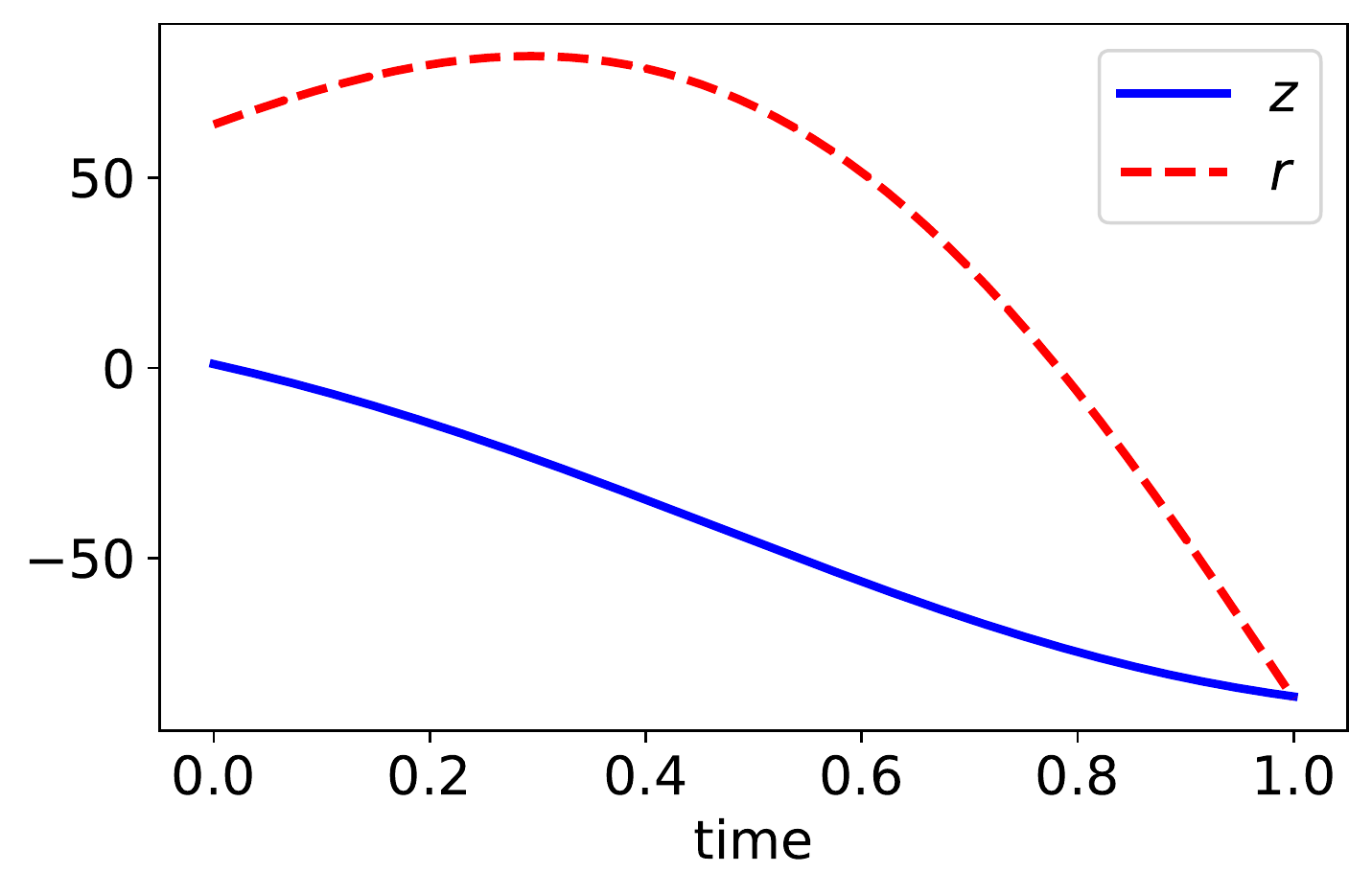}
	\end{subfigure}
	\begin{subfigure}{.33\columnwidth}
		\centering 
		\includegraphics[width=\columnwidth]{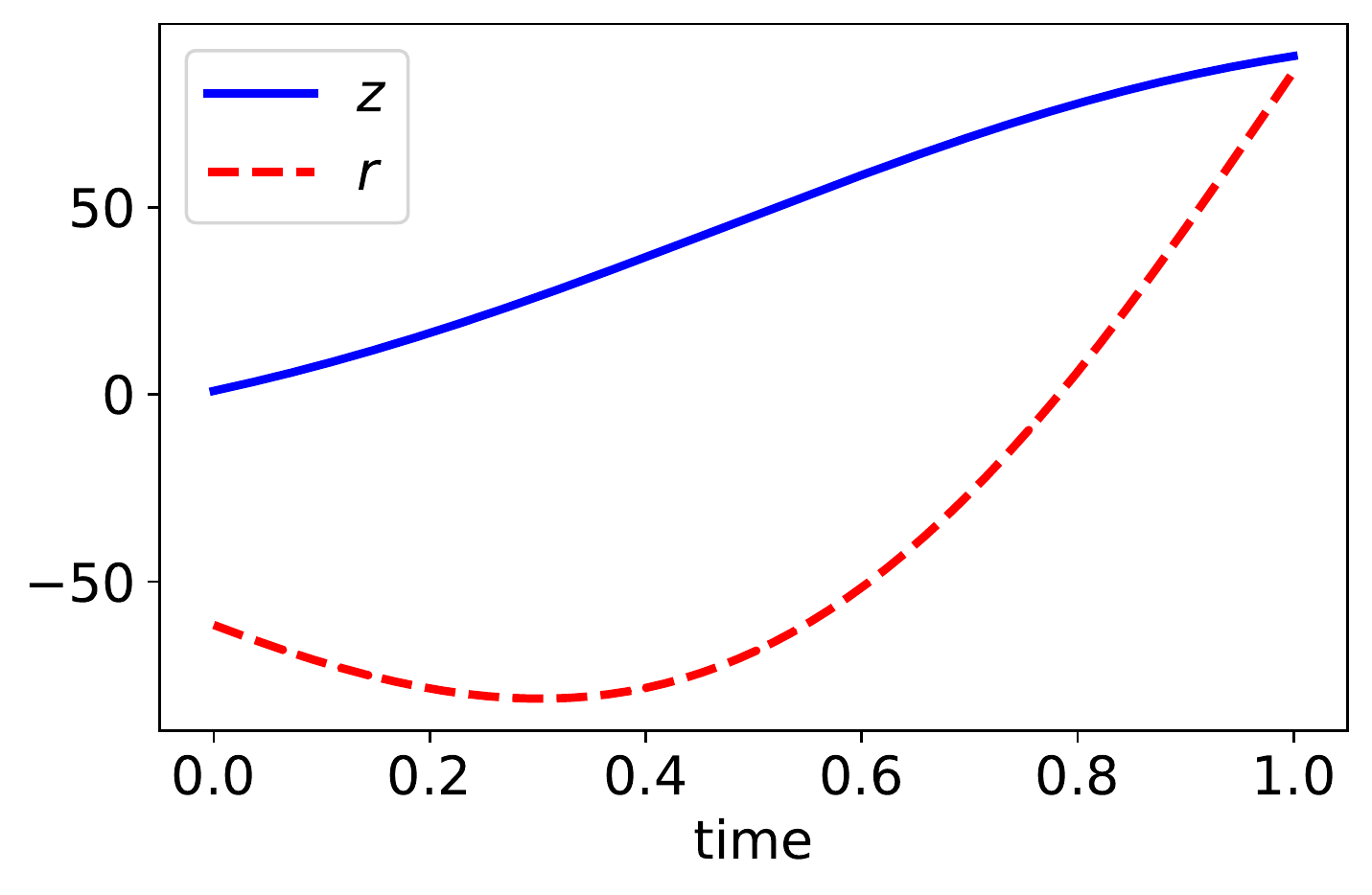}
	\end{subfigure}
	\caption{Picard iterations without damping for Test case~2 (see Table~\ref{AMS-num-tab:params-LQ-base}): $L^2$ difference between two successive iterates (left), $z^{(199)}$ and $r^{(199)}$ (middle) and $z^{(200)}$ and $r^{(200)}$ (right). The curves related to $z$ and $r$ are respectively in blue (full line) and in red (dashed line). } 
 	\label{AMS-num-fig:LQ-ODE-Picard-fail}
\end{figure}

However, this procedure fails on many examples. Figure~\ref{AMS-num-fig:LQ-ODE-Picard-fail} based on Test case~2 in Table~\ref{AMS-num-tab:params-LQ-base} provides an illustration: here we see that $z^{(\mathtt{k})}$ and $r^{(\mathtt{k})}$ end up reacting to each other from an iteration to the next one, and the overall process diverge. In the framework of MFG, a possible interpretation for this phenomenon is the following (recall that $z$ corresponds to the equilibrium mean and $r$ is part of the equilibrium control). Given a trajectory for the mean-field term, each agent compute their best response to this crowd behavior. Applying this best response leads to a new trajectory for the mean-field term. Computing once again the best response leads back to the first trajectory. This type of phenomenon is also well-known in two-player games.

\begin{figure}[h]
	\begin{subfigure}{.45\columnwidth}
		\centering
		\includegraphics[width=\columnwidth]{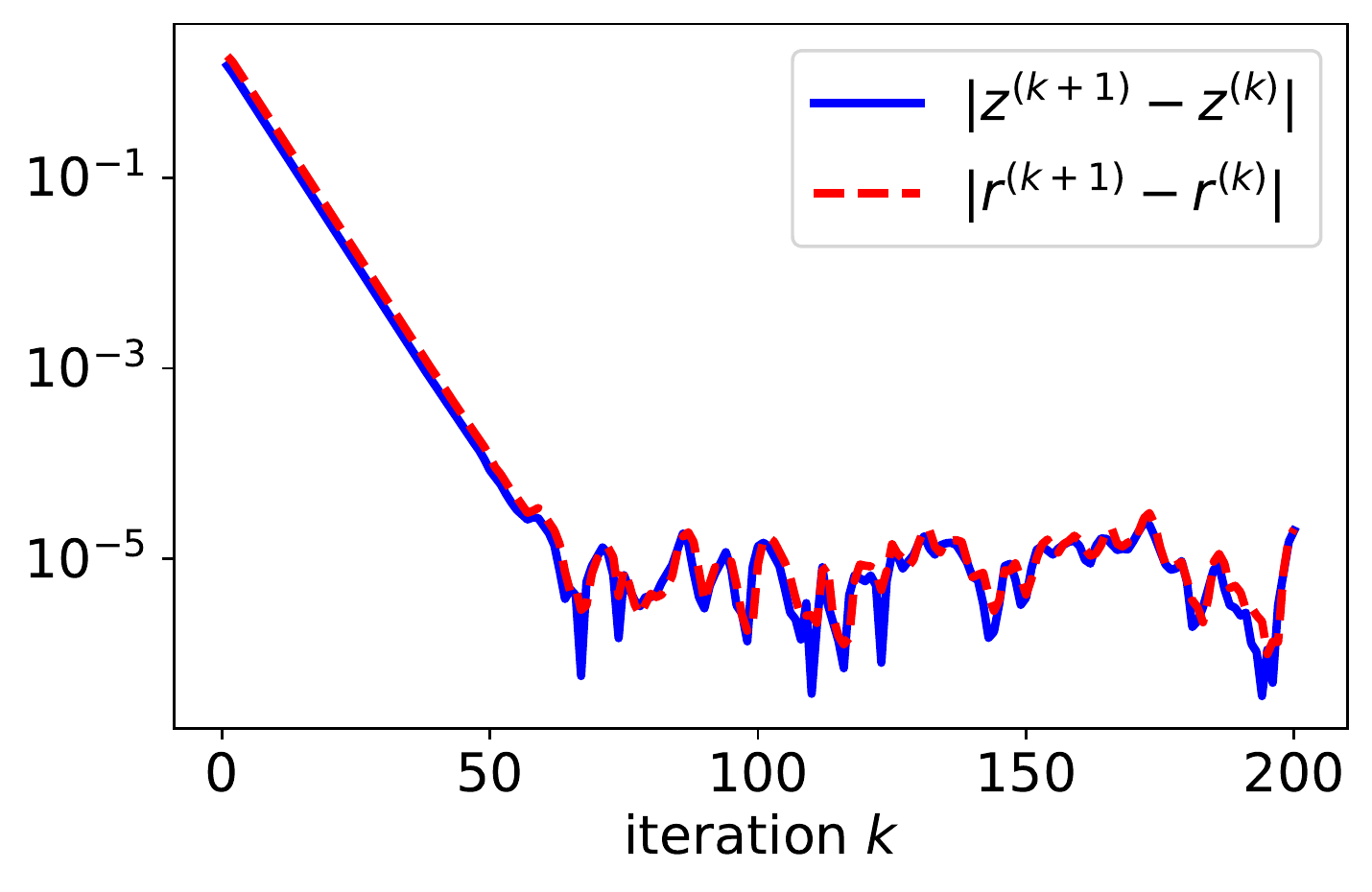}
	\end{subfigure}%
	\begin{subfigure}{.45\columnwidth}
		\centering 
		\includegraphics[width=\columnwidth]{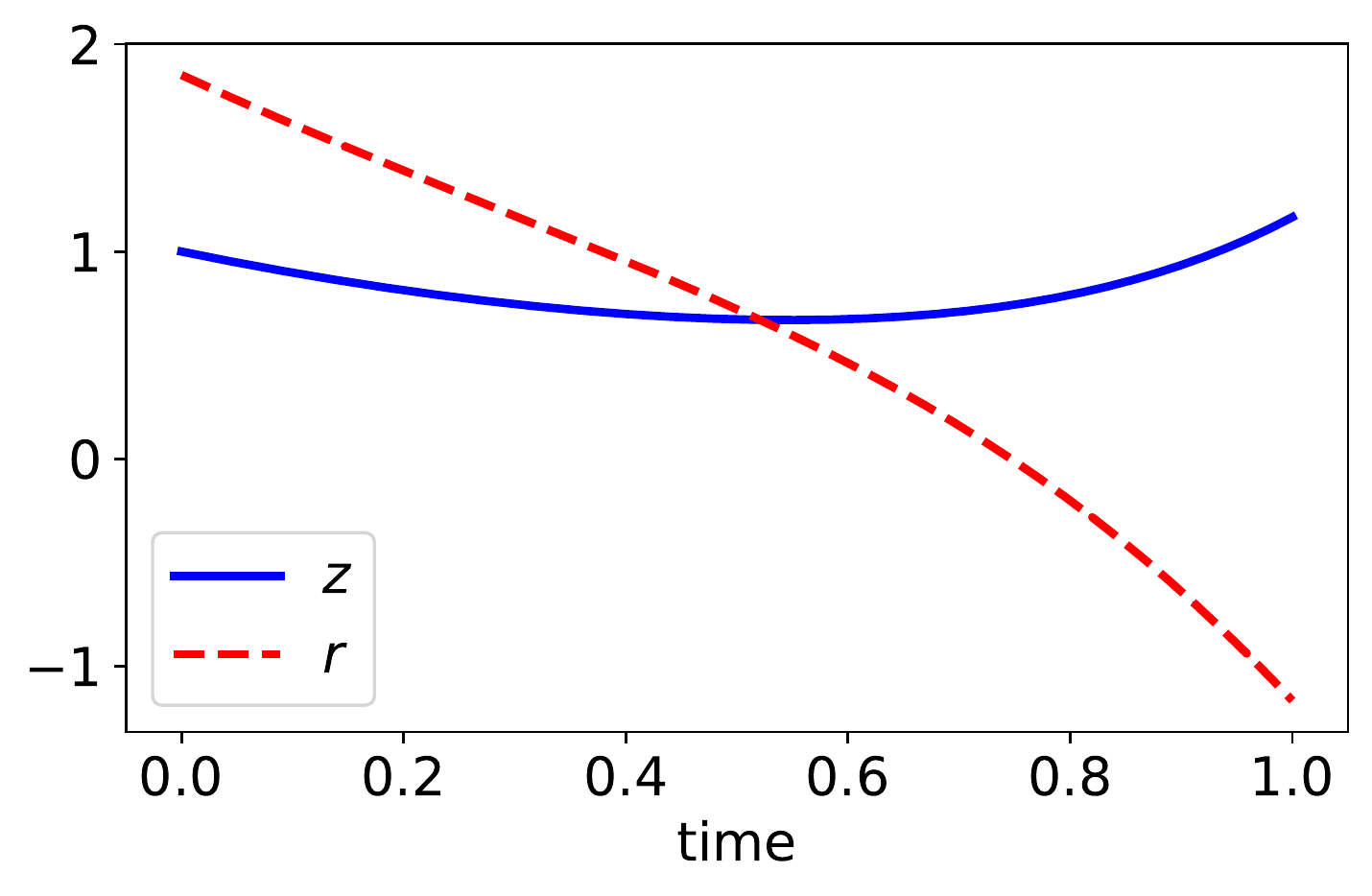}
	\end{subfigure}
	\caption{Picard iterations with constant damping $\omega = 0.1$ for Test case~2 (see Table~\ref{AMS-num-tab:params-LQ-base}): $L^2$ difference between two successive iterates (left) and $z^{(200)}$ and $r^{(200)}$ (right). The curves related to $z$ and $r$ are respectively in blue (full line) and in red (dashed line). } 
 	\label{AMS-num-fig:LQ-ODE-Picard-damped-converge}
\end{figure}

\begin{figure}[h]
	\begin{subfigure}{.33\columnwidth}
		\centering
		\includegraphics[width=\columnwidth]{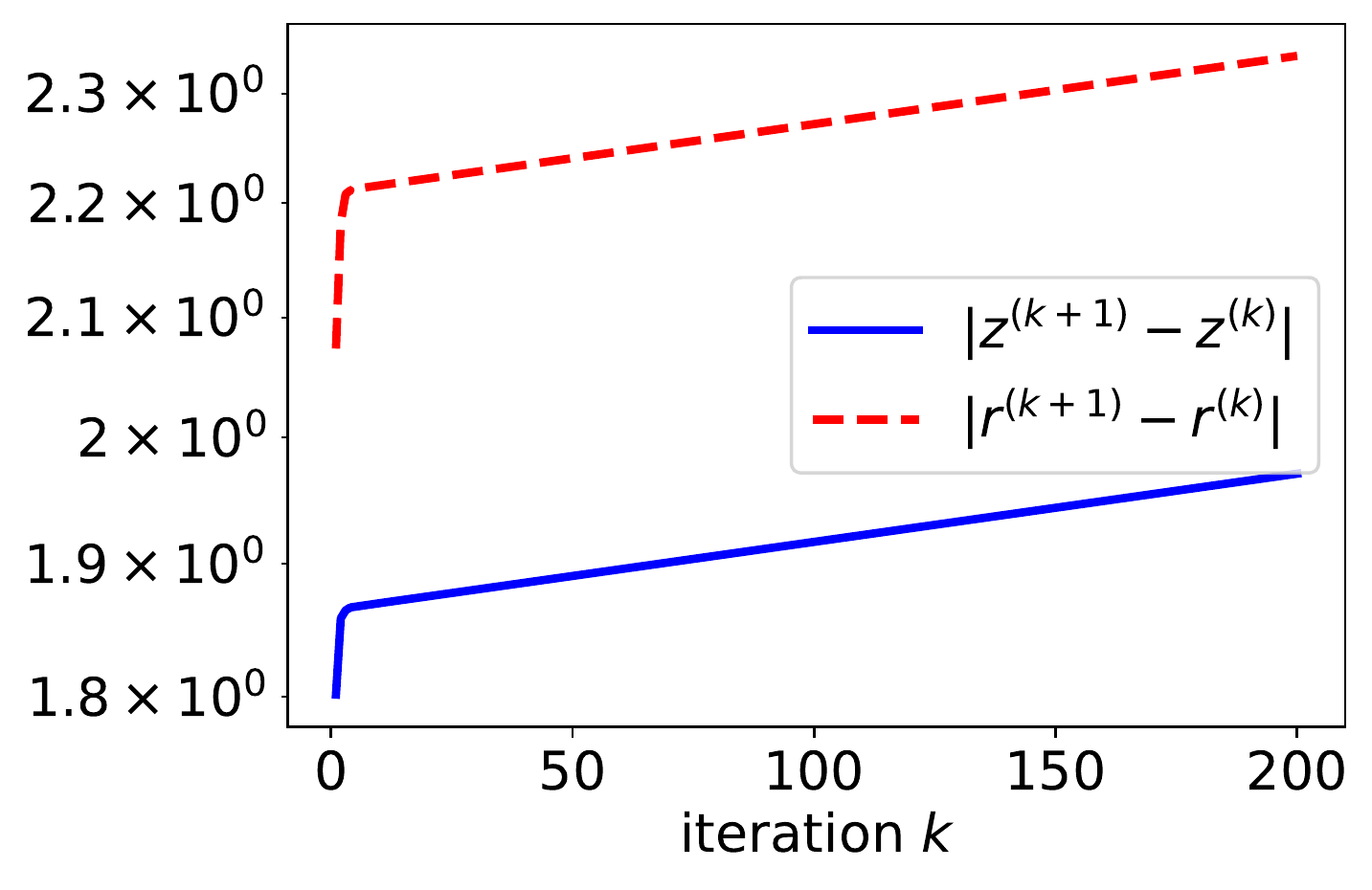}
	\end{subfigure}%
	\begin{subfigure}{.33\columnwidth}
		\centering 
		\includegraphics[width=\columnwidth]{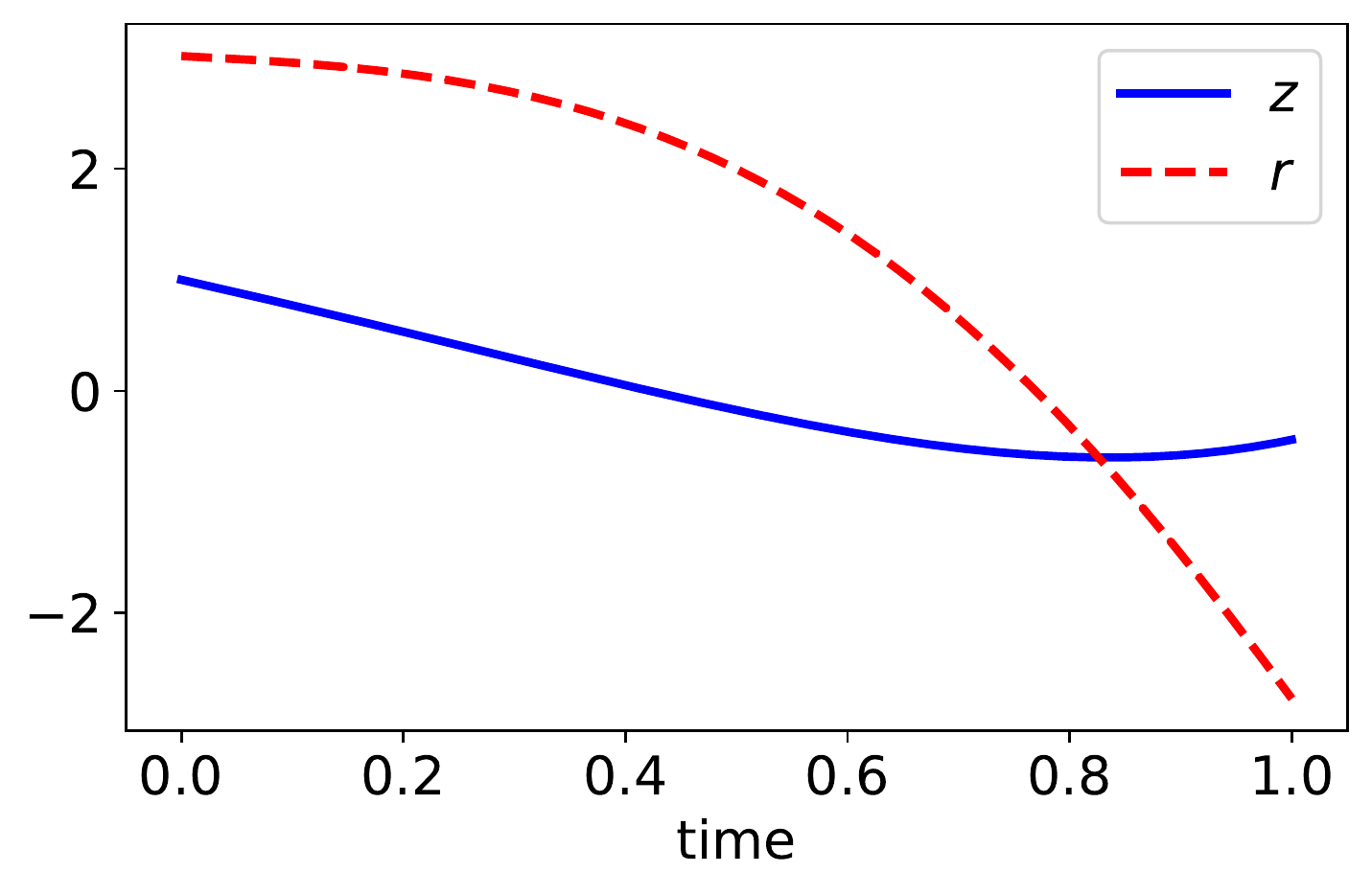}
	\end{subfigure}
	\begin{subfigure}{.33\columnwidth}
		\centering 
		\includegraphics[width=\columnwidth]{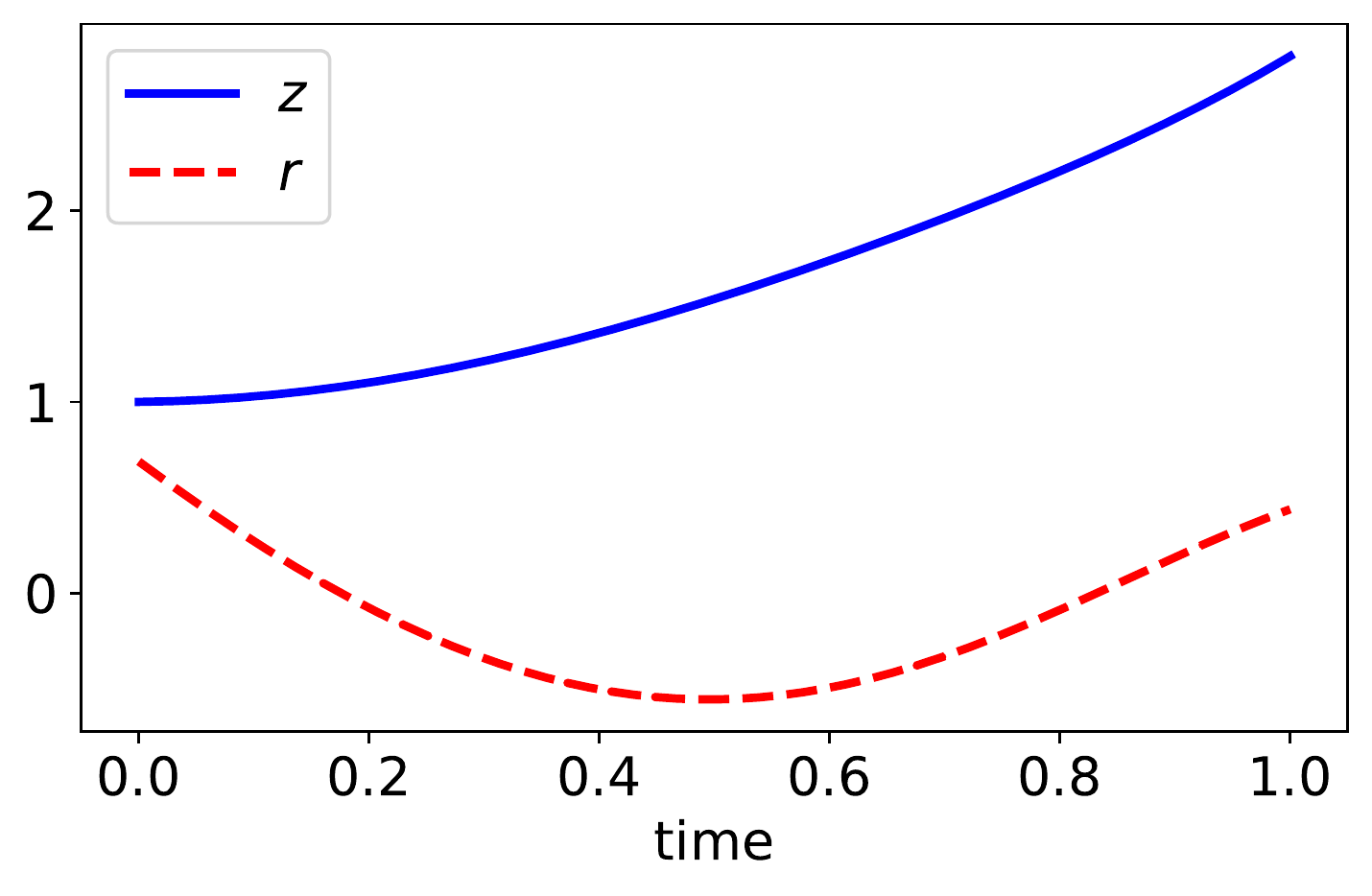}
	\end{subfigure}
	\caption{Picard iterations with constant damping $\omega = 0.01$ for Test case~2 (see Table~\ref{AMS-num-tab:params-LQ-base}): $L^2$ difference between two successive iterates (left), $z^{(199)}$ and $r^{(199)}$ (middle) and $z^{(200)}$ and $r^{(200)}$ (right). The curves related to $z$ and $r$ are respectively in blue (full line) and in red (dashed line). } 
 	\label{AMS-num-fig:LQ-ODE-Picard-damped-fail}
\end{figure}

A simple approach to try and fix the above method is to add damping to the iterations. Let $\omega \in [0,1)$ be a damping parameter. At iteration $\mathtt{k} \ge 0$, letting $z^{(\mathtt{k}+1)}$ be the solution to~\eqref{AMS-num-eq:LQ-ODE-z} with $r^{(\mathtt{k})}$, instead of directly plugging this new mean-field trajectory in~\eqref{AMS-num-eq:LQ-ODE-r}, we define
$$
	\tilde{z}^{(\mathtt{k}+1)} = \omega \tilde{z}^{(\mathtt{k})} + (1-\omega) z^{(\mathtt{k}+1)}
$$
and plug this trajectory in~\eqref{AMS-num-eq:LQ-ODE-r} to compute $r^{(\mathtt{k}+1)}$. This corresponds to taking  $\delta(\cdot) \equiv \omega$ in Algorithm~\ref{AMS-num-algo:LQ-ODE-Picard}. Such a modification of the plain Picard iterations typically helps to ensure convergence of the iterations but the choice of the damping parameter is not obvious: if it is too small, convergence may still fail, whereas if it is too large, convergence will be very slow. See Figure~\ref{AMS-num-fig:LQ-ODE-Picard-damped-converge} and Figure~\ref{AMS-num-fig:LQ-ODE-Picard-damped-fail}.

\begin{algorithm}[H]
\DontPrintSemicolon
  
  \KwInput{Initial guess $(\tilde z, \tilde r)$; damping $\delta(\cdot)$; number of iterations $\mathtt{K}$}
  \KwOutput{Approximation of $(\hat z, \hat r)$ solving~\eqref{AMS-num-eq:LQ-ODE-z}\&\eqref{AMS-num-eq:LQ-ODE-r}}
  Initialize $z^{(0)} = \tilde z^{(0)} = \tilde z, r^{(0)} = \tilde r$ \;
        \For{$\mathtt{k}=0,1,2,\dots, \mathtt{K}-1$}    
        { 
        	Let $r^{(\mathtt{k}+1)}$ be the solution to:
	$$
		 - \frac{d r}{dt} = (A - P_t B^2 C^{-1})r_t + (P_t \bar{A} - \bar{Q} S) \tilde{z}^{(\mathtt{k})}_t, 	\qquad r_T = - \bar{Q}_T S_T \tilde{z}^{(\mathtt{k})}_T
	$$\;
	Let $z^{(\mathtt{k}+1)}$ be the solution to:
	$$
		\frac{d z}{dt} = (A+\bar{A} - B^2 C^{-1}) z_t - B^2 C^{-1} r^{(\mathtt{k}+1)}_t, 			\qquad  z_0 = \bar x_0 
	$$\;
	Let $\tilde{z}^{(\mathtt{k}+1)} = \delta(\mathtt{k}) \tilde{z}^{(\mathtt{k})} + (1-\delta(\mathtt{k})) z^{(\mathtt{k}+1)}$
        }
        \Return{$(z^{(\mathtt{K})}, r^{(\mathtt{K})})$} 

\caption{Fixed-point Iterations \label{AMS-num-algo:LQ-ODE-Picard}}
\end{algorithm}

\subsection{Fictitious play}

Another procedure, which has been widely studied in algorithmic game theory, is the so-called fictitious play. It has also drawn interest in the MFG community~\cite{MR3608094,MR4030259,Elie2020OnTC,perrin2020fictitious}. Once again, it amounts to update in turn the mean-field term and the control, but here the control is computed as the best response to a weighted average of the mean-field terms in previous iterations. In the context of the LQ example studied here, one can recast this procedure as follows: $r^{(\mathtt{k}+1)}$ is computed given $\tilde{z}^{(\mathtt{k})}$ and 
$$
	\tilde{z}^{(\mathtt{k}+1)} = \frac{\mathtt{k}}{\mathtt{k}+1} \tilde{z}^{(\mathtt{k})} + \frac{1}{\mathtt{k}+1} z^{(\mathtt{k}+1)}
$$
where $z^{(\mathtt{k}+1)}$ solves~\eqref{AMS-num-eq:LQ-ODE-z} with $r$ replaced by $r^{(\mathtt{k}+1)}$. 
 This corresponds to Algorithm~\ref{AMS-num-algo:LQ-ODE-Picard} with $\delta(\mathtt{k}) = \frac{\mathtt{k}}{\mathtt{k}+1}$. 
 
\begin{figure}[h]
	\begin{subfigure}{.45\columnwidth}
		\centering
		\includegraphics[width=\columnwidth]{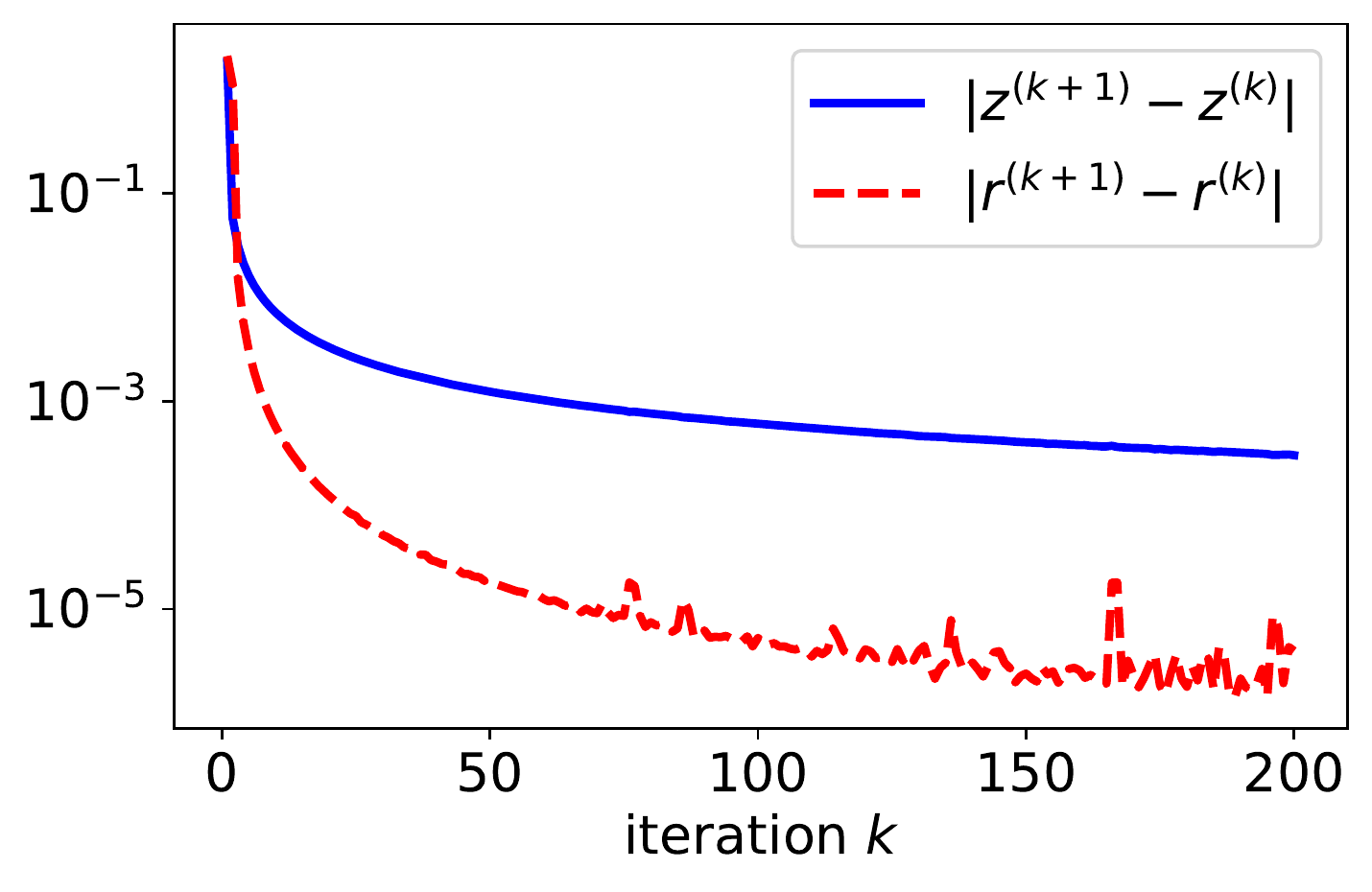}
	\end{subfigure}%
	\begin{subfigure}{.45\columnwidth}
		\centering 
		\includegraphics[width=\columnwidth]{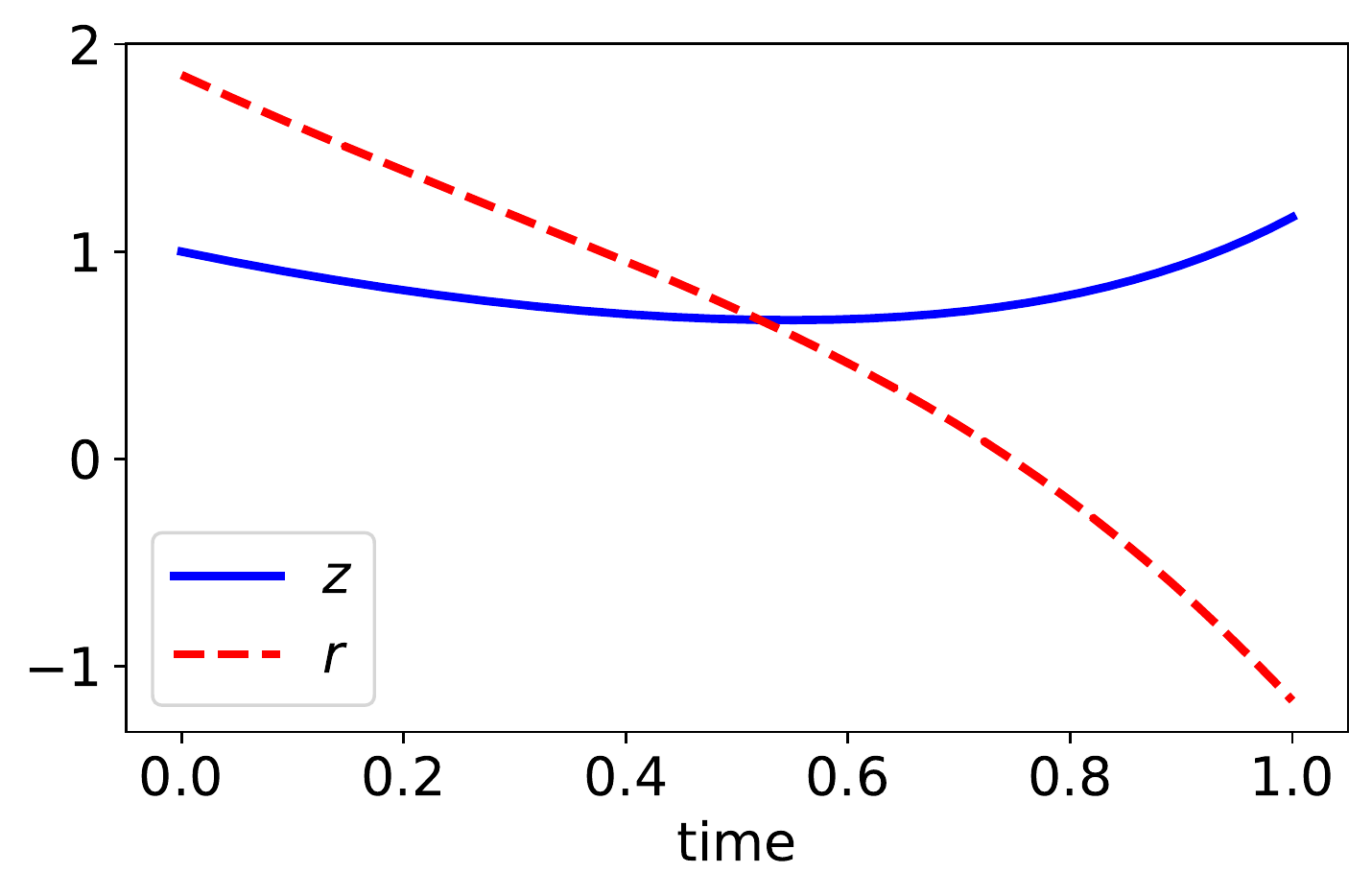}
	\end{subfigure}
	\caption{Fictitious play iterations for Test case~2 (see Table~\ref{AMS-num-tab:params-LQ-base}): $L^2$ difference between two successive iterates (left) and $z^{(200)}$ and $r^{(200)}$ (right). The curves related to $z$ and $r$ are respectively in blue (full line) and in red (dashed line). } 
 	\label{AMS-num-fig:LQ-ODE-FictitiousPlay-convergence}
\end{figure}

An important advantage of this method is that it can be proved to converge under less stringent conditions than Picard iterations, see~\cite{MR3608094}. However, the damping effect becomes stronger as the number of iterations increases, which means that convergence can be slower.  This is illustrated in Figure~\ref{AMS-num-fig:LQ-ODE-FictitiousPlay-convergence}. Heuristics to improve the convergence have been proposed in~\cite{MR4030259}.

\subsection{Newton iterations}

\begin{figure}[h]
	\begin{subfigure}{.45\columnwidth}
		\centering
		\includegraphics[width=\columnwidth]{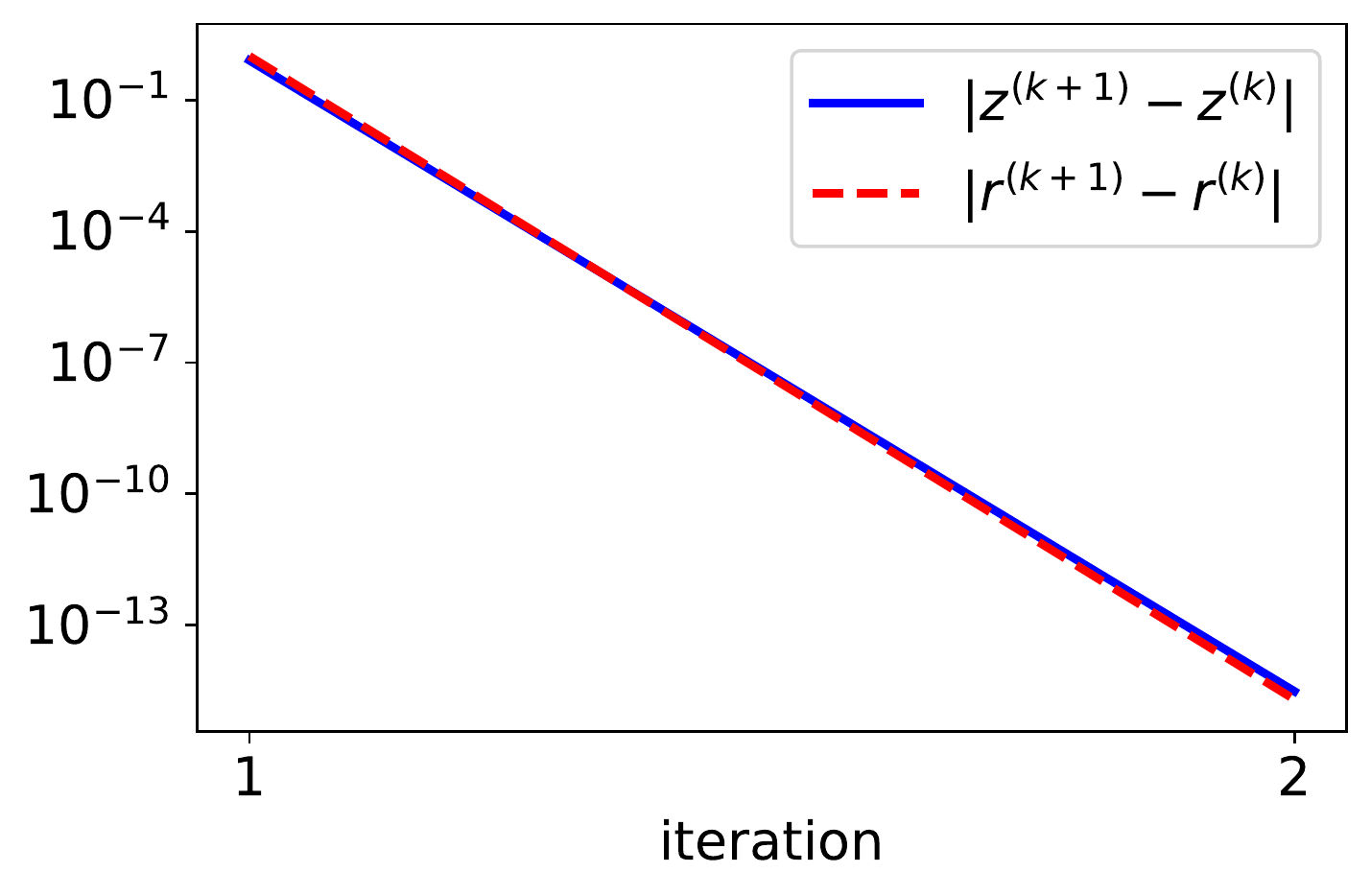}
	\end{subfigure}%
	\begin{subfigure}{.45\columnwidth}
		\centering 
		\includegraphics[width=\columnwidth]{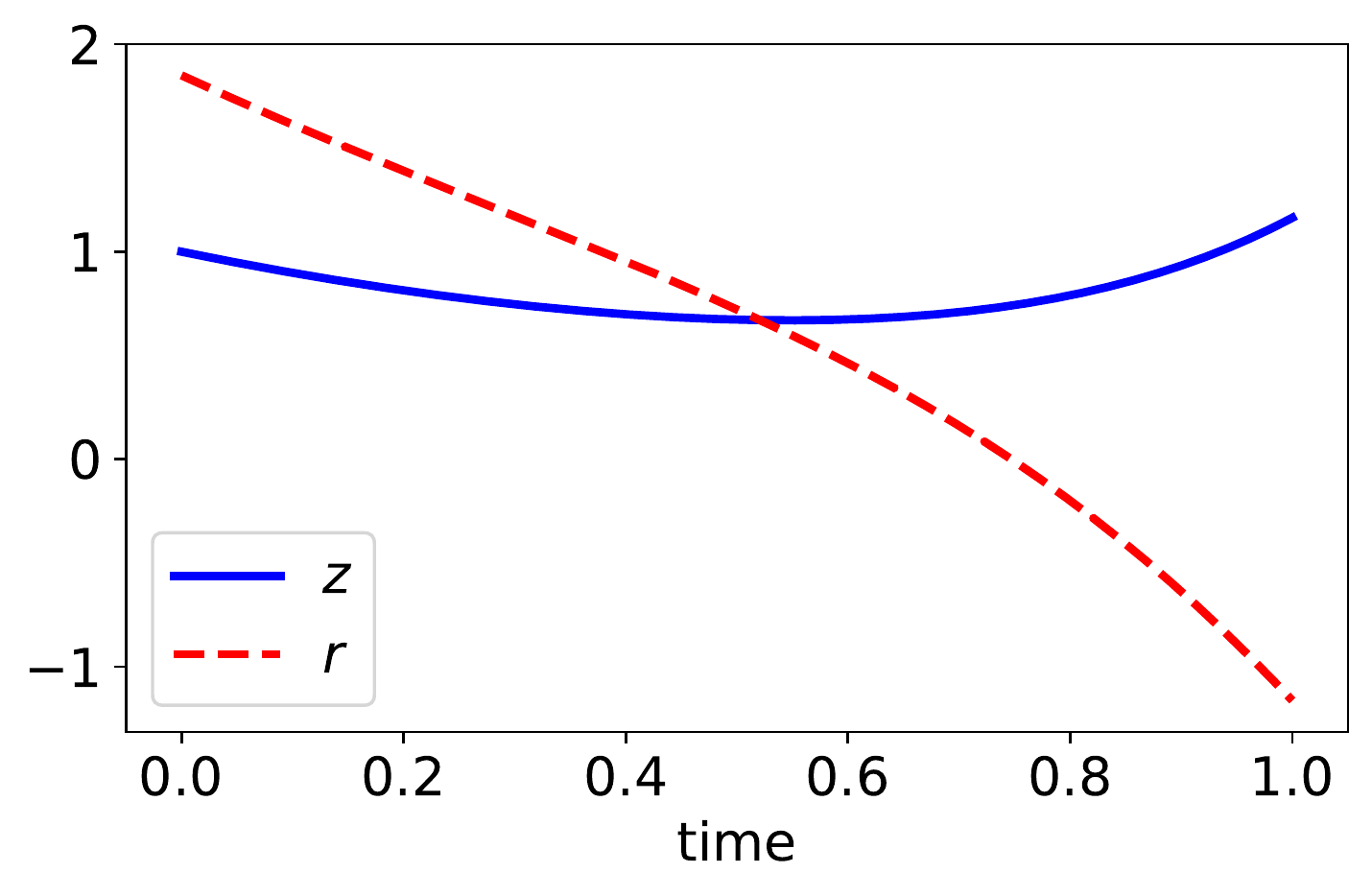}
	\end{subfigure}
	\caption{Newton iterations for Test case~2 (see Table~\ref{AMS-num-tab:params-LQ-base}): $L^2$ difference between two successive iterates (left) and $z^{(200)}$ and $r^{(200)}$ (right). The curves related to $z$ and $r$ are respectively in blue (full line) and in red (dashed line). } 
 	\label{AMS-num-fig:LQ-ODE-Newton}
\end{figure}

Instead of updating alternatively $z$ and $r$ by solving a forward equation and then a backward equation, another approach consists in solving the whole forward-backward system using Newton iterations. After discretizing time as in \S~\ref{AMS-num-sec:LQ-ODE-time-discr}, the idea is to view the system~\eqref{AMS-num-eq:LQ-ODE-z-discrete}\&\eqref{AMS-num-eq:LQ-ODE-r-discrete} as the zero of an operator $\mathcal{F}: \RR^{2(N_T+1)} \to \RR^{2(N_T+1)}$ defined as:
$$
	(Z,R) \hbox{ solve } \eqref{AMS-num-eq:LQ-ODE-z-discrete}-\eqref{AMS-num-eq:LQ-ODE-r-discrete} \Leftrightarrow \mathcal{F}(Z,R) = 0.
$$
Note that $\mathcal{F}$ must take into account the initial and terminal conditions. 
Denoting by $D \mathcal{F}$ the differential of this operator, Algorithm~\ref{AMS-num-algo:LQ-ODE-Newton} summarizes Newton iterations to find a zero of $\mathcal{F}$.  In the present LQ example, this operator is linear hence solving~\eqref{AMS-num-eq:LQ-Newton-step-k} is just as straightforward as solving directly the linear system~\eqref{AMS-num-eq:LQ-ODE-matrix-form} of interest. See Fig~\ref{AMS-num-fig:LQ-ODE-Newton}. However, this strategy based on Newton iterations will be particularly useful when solving forward-backward systems of PDEs which are not linear (see \S~\ref{sec:NewtonPDE}).

\begin{algorithm}[H]
\DontPrintSemicolon
  
  \KwInput{Initial guess $(\tilde Z, \tilde R)$; number of iterations $\mathtt{K}$}
  \KwOutput{Approximation of $(\hat z, \hat r)$ solving~\eqref{AMS-num-eq:LQ-ODE-z}\&\eqref{AMS-num-eq:LQ-ODE-r}}
  Initialize $(Z^{(0)}, R^{(0)}) = (\tilde Z, \tilde R)$ \;
        \For{$\mathtt{k}=0,1,2,\dots, \mathtt{K}-1$}    
        { 
        	Let $(\tilde Z^{(\mathtt{k}+1)},\tilde R^{(\mathtt{k}+1)})$ solve 
	\begin{equation}
	\label{AMS-num-eq:LQ-Newton-step-k}
		 D\mathcal{F}(Z^{(\mathtt{k})}, R^{(\mathtt{k})})(\tilde Z^{(\mathtt{k}+1)},\tilde R^{(\mathtt{k}+1)}) = \mathcal{F}(Z^{(\mathtt{k})}, R^{(\mathtt{k})})  
	\end{equation}\;
	Let $(Z^{(\mathtt{k}+1)}, R^{(\mathtt{k}+1)}) = (\tilde Z^{(\mathtt{k}+1)},\tilde R^{(\mathtt{k}+1)}) + (Z^{(\mathtt{k})}, R^{(\mathtt{k})})$ \;
	\Return{$(Z^{(\mathtt{K})}, R^{(\mathtt{K})})$}\;
        }

\caption{Newton Iterations \label{AMS-num-algo:LQ-ODE-Newton}}
\end{algorithm}

\subsection{Mean Field Type Control and the Price of Anarchy}
\label{AMS-num-sec:LQ-sub-MFC}

Using the linear dynamics and quadratic costs introduced \S~\ref{AMS-num-sec:LQ-pb}, we now turn our attention to the corresponding MFC problem. For this problem too, the solution can be characterized using a system of ODEs. More precisely, denoting by $\ctrl^*$ and $m^*$ respectively the optimal control and the induced flow of densities, we have
\begin{subequations}
     \begin{empheq}[left=\empheqlbrace]{align*}
		\int \xi m^*(t,\xi) d \xi &= \check{z}_t,
		\\
		\ctrl^*(t,x) &= -B(\check{p}_t x + \check{r}_t) / C,
		\\
		J^{MFC}(\ctrl^*) &= \frac{1}{2} \check{p}_0(\sigma_0^2 + \bar x_0^2) + \check{r}_0 \bar x_0 + \check{s}_0 + (1-S_T)\bar Q_T S_T \check{z}_T^2
		\\&\qquad - \int_0^T \big[(\check{p}_t \check{z}_t + \check{r}_t) \bar A \check{z}_t - (1-S_t)\bar Q S \check{z}_t^2\big]dt
     \end{empheq}
   \end{subequations}
where $(\check{z}, \check{p}, \check{r}, \check{s})$ solve the following system of ordinary differential equations (ODEs):
    \begin{subequations}
     \begin{empheq}[left=\empheqlbrace\,\,]{align*}
     \frac{d \check{z}}{dt} &= (A+\bar{A} - B^2 C^{-1}) \check{z}_t - B^2 C^{-1} \check{r}_t, 			&& \check{z}_0 = \bar x_0,
     \\
     -\frac{d \check{p}}{dt} &= 2A\check{p}_t - B^2C^{-1}\check{p}_t^2 + Q + \bar{Q}, 				&& \check{p}_T = Q_T + \bar{Q}_T,
     \\
     - \frac{d \check{r}}{dt} &= (A + \bar A - \check{p}_t B^2 C^{-1})\check{r}_t + (2\check{p}_t \bar{A} - 2\bar{Q} S + \bar Q S^2) \check{z}_t, 	&& \check{r}_T = - \bar{Q}_T S_T \check{z}_T,
     \\
     - \frac{d s}{dt} &= \nu \check{p}_t - \frac{1}{2} B^2 C^{-1} \check{r}_t^2 + \check{r}_t \bar{A} \check{z}_t + \frac{1}{2} S^2 \bar{Q} \check{z}_t^2, && \check{s}_T = \frac{1}{2} \bar{Q}_T S_T^2 \check{z}_T^2.
     \end{empheq}
   \end{subequations}
   Note that this system is very similar to~\eqref{AMS-num-eq:LQ-ODE-z}--\eqref{AMS-num-eq:LQ-ODE-s} except that the ODEs for $r$ and $\check r$ are different. Furthermore, the expressions for $J^{MFG}$ and $J^{MFC}$ are also different. This is due to the fact that, in the MFC, the function $\check u(t,x) = \frac{1}{2}\check{p}_t x^2 + \check{r}_t x + \check{s}_t$ is not the value function of the control problem but rather an adjoint state, as we shall see in Section~\ref{AMS-num-sec:PDE-scheme}.
     
\begin{figure}[h]
	\begin{subfigure}{.33\columnwidth}
		\centering
		\includegraphics[width=\columnwidth]{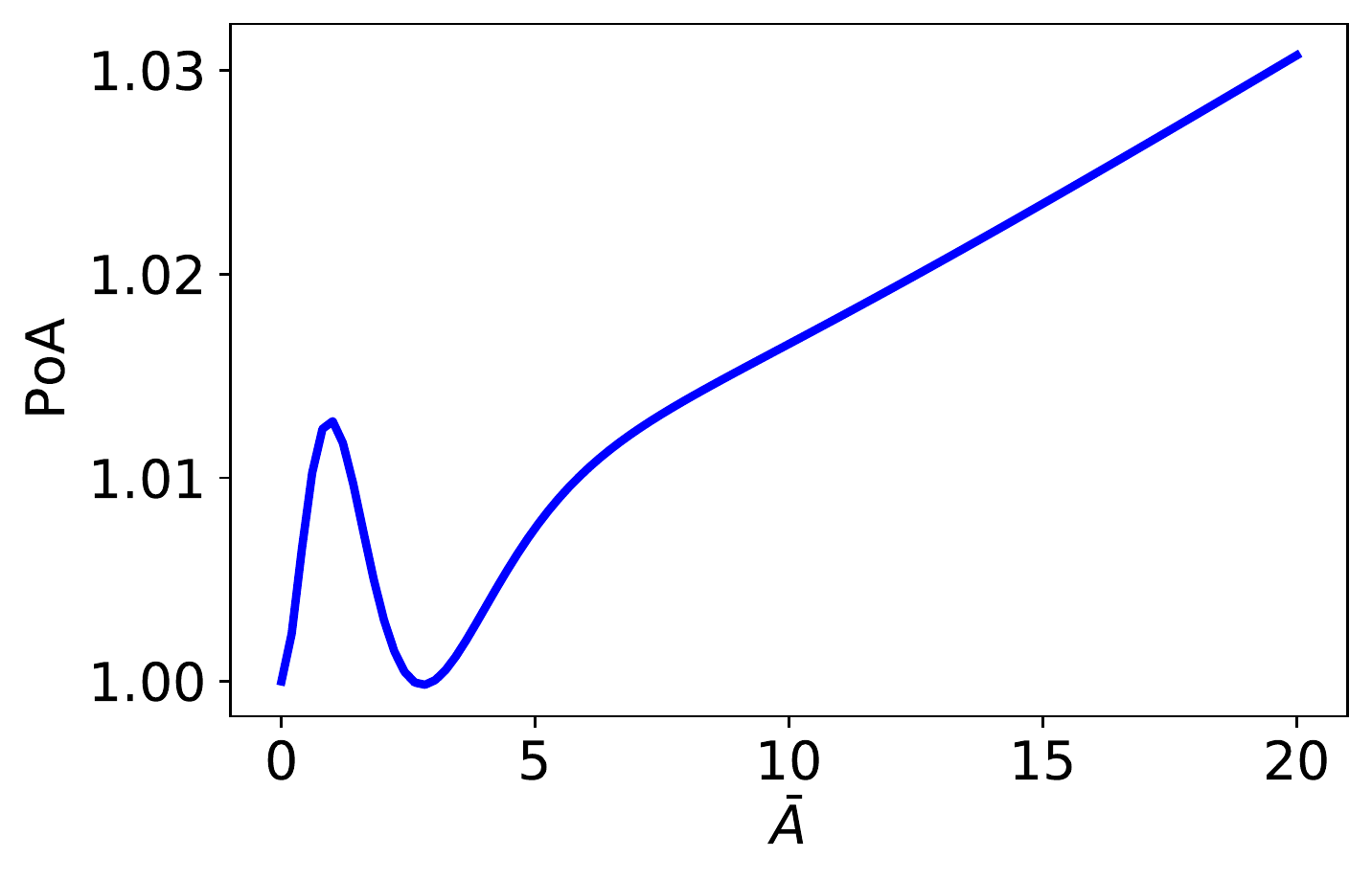}
	\end{subfigure}%
	\begin{subfigure}{.33\columnwidth}
		\centering 
		\includegraphics[width=\columnwidth]{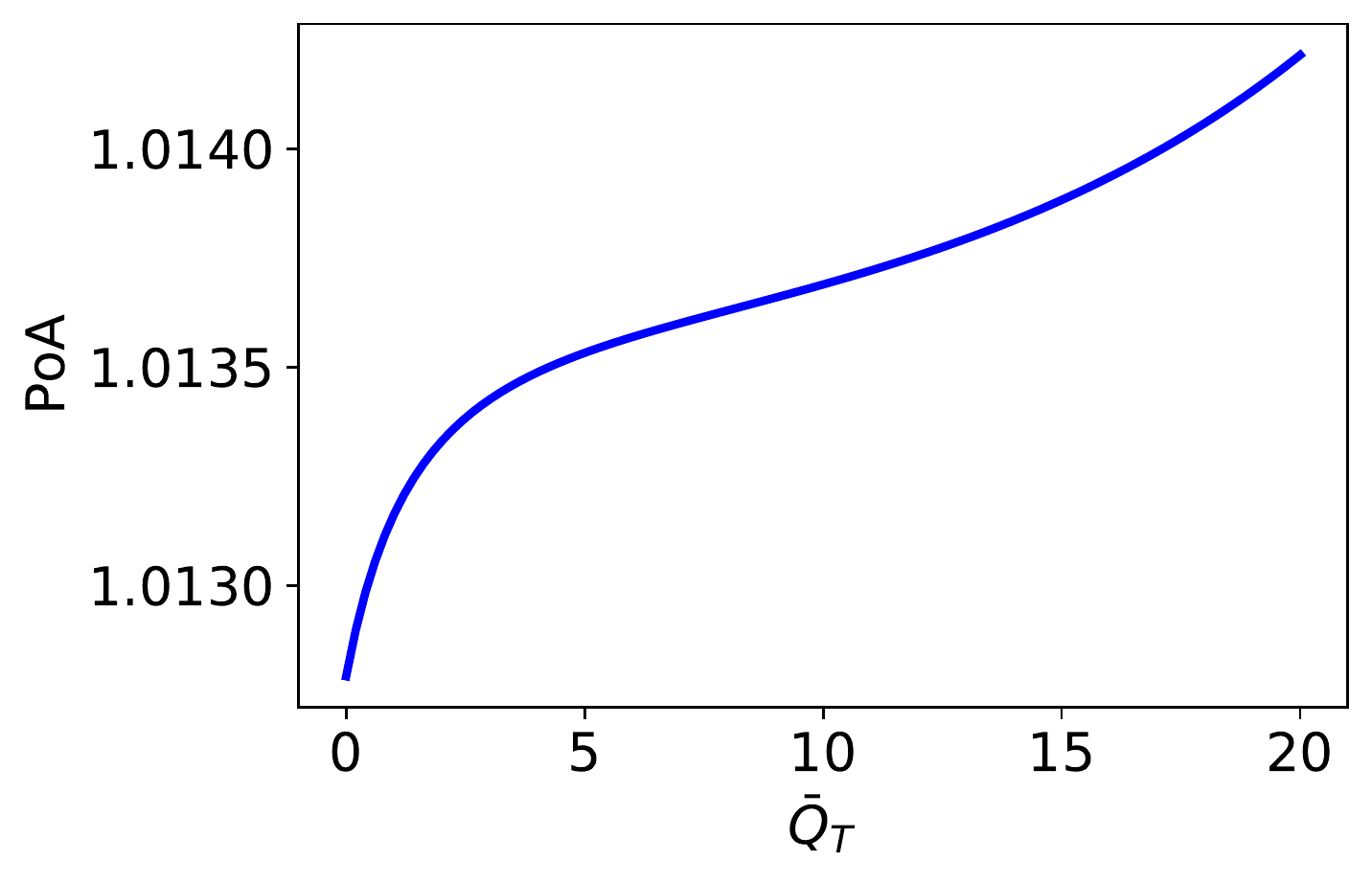}
	\end{subfigure}
	\begin{subfigure}{.33\columnwidth}
		\centering 
		\includegraphics[width=\columnwidth]{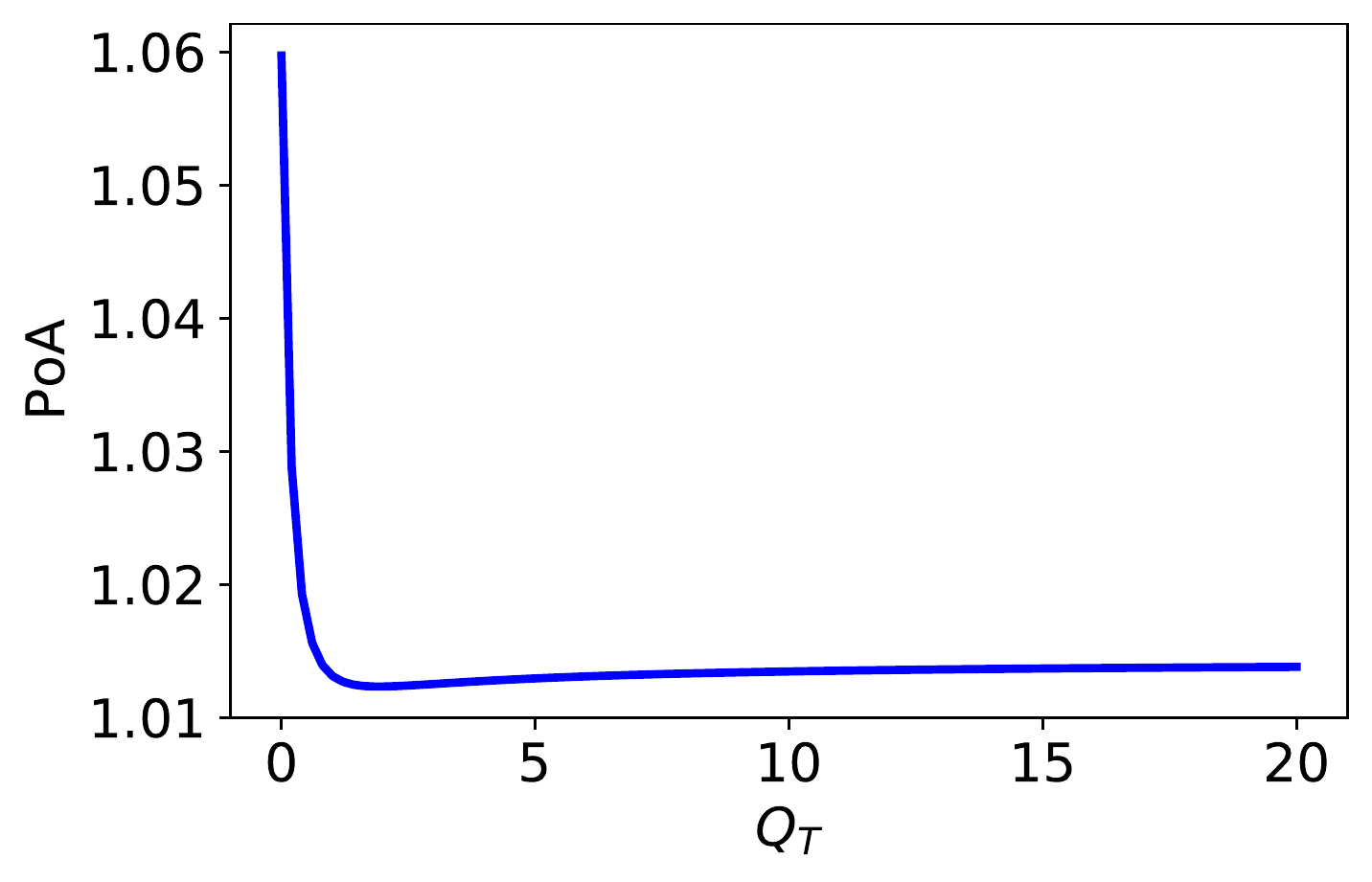}
	\end{subfigure}
	\caption{Price of anarchy for Test case~3 (left), Test case~4 (middle) and Test case~5 (right). See Table~\ref{AMS-num-tab:params-LQ-base} for the parameters values. } 
 	\label{AMS-num-fig:LQ-ODE-PoA}
\end{figure}

Since MFC corresponds to a social optimum, the price paid by a typical player can only be lower than in a MFG (with the same costs and dynamics). In other words, the price of anarchy (PoA), defined as:
\begin{equation}
\label{eq:def-MFPoA}
	PoA = \frac{J^{MFG}_{\hat{m}}(\hat{\ctrl})}{J^{MFC}(\ctrl^*)}, 
\end{equation}
is at least $1$. Figure~\ref{AMS-num-fig:LQ-ODE-PoA} illustrates the price of anarchy when $\bar A, \bar Q$ or $Q_T$ vary. In particular, in Test case~3, $\bar Q = \bar Q_T = 0$ and when $\bar A=0$ too, the price of anarchy is $1$, see Table~\ref{AMS-num-tab:params-LQ-base}. Indeed, in this scenario, the mean-field interactions completely disappear from the cost and the dynamics, hence MFG and MFC amount to the same problem. Test cases~4 and~5 illustrate the dependence of the price of anarchy with respect to $\bar Q_T$ and $Q_T$ respectively. For more details on the price of anarchy in LQ MFG and MFC, see~\cite{MR3968548}.

\clearpage

\section{\bf PDE Systems and Numerical Schemes}
\label{AMS-num-sec:PDE-scheme}

In general, besides the LQ case, optimality conditions for MFG and MFC can not be phrased in terms of ODEs. The two main approaches are through systems of partial differential equations (PDEs) or systems stochastic differential equations (SDEs). In both cases, these systems have a forward-backward structure. A third approach is via the so-called master equation, which we will discuss in~\S~\ref{AMS-num-sec:NN-finite-master}. For now, we focus on the first approach, which leads to a system composed of a forward Kolmogorov-Fokker-Planck PDE for the distribution of the population's state coupled with a backward Hamilton-Jacobi-Bellman equation for the value function of a typical player. To solve numerically this PDE system, the first question that arises is the approximation of functions of space and time. A possible approach is to consider a space-time mesh and approximate these functions by their values at the points of the mesh. The functions are approximated by vectors, which satisfy finite difference equations coming from a suitable numerical scheme of the PDEs.

\subsection{Optimality condition through PDE systems}
\label{AMS-num-sec:PDEsystem}

For the sake of simplicity, unless otherwise specified, we consider the $d$-dimensional torus for the spatial domain, \textit{i.e.}, $\dom = \TT^d$. We first consider the MFG and then the MFC. Despite their similarities, the two PDE systems have important differences reflecting the differences between the two problems.

\textbf{Mean field game.} We first note that for a given flow of densities $\hat{m}$ and a given feedback control $\ctrl$, the density $m^{\hat{m}, \ctrl}_t$ of the law of $X^{\hat{m}, \ctrl}_t$ solving~\eqref{AMS-num-eq:dyn-X-general-MFG} satisfies the Kolmogorov-Fokker-Planck (KFP) equation:
\begin{equation}
\label{AMS-num-eq:KFP-ctrl}
	\frac{\partial m^{\hat{m}, \ctrl}}{\partial t}(t,x)  - \nu \Delta m^{\hat{m}, \ctrl}(t,x) + \diver\left( m^{\hat{m}, \ctrl}(t,\cdot) b(\cdot, \hat m(t), \ctrl(t,\cdot)) \right)(x) = 0, 
	\hbox{ in } (0,T] \times \TT^d ,
\end{equation}
with the initial condition:
\begin{equation}
\label{AMS-num-eq:KFP-ctrl-init}
	m^{\hat{m}, \ctrl}(0,x) = m_0(x), 
	\hbox{ in } \TT^d.
\end{equation}

Let $H: \TT^d \times L^2(\TT^d) \times \RR^d  \to \in \RR$ be the Hamiltonian of the control problem faced by an infinitesimal agent in the first point above, which is defined by:
$$
	H: (x,m,p) \mapsto H(x,m,p) = \max_{\ctrl \in \RR^\ctrldim} \{ -L(x,m,\ctrl,p) \},
$$
where $L: \TT^d \times L^2(\TT^d) \times \RR^\ctrldim \times \RR^d \to  \RR$ is the Lagrangian, defined by
$$
L: (x,m,\ctrl,p) \mapsto L(x,m,\ctrl,p) =
f(x,m,\ctrl) + \langle b(x,m,\ctrl) , p \rangle.
$$
In the sequel, we will assume that the running cost $f$ and the drift $b$ are such that $H$ is well-defined, $\mathcal{C}^1$ with respect to $(x,p)$, and strictly convex with respect to $p$. We will write $\partial_p H$ for the partial derivative of $H$ with respect to the $p$ variable.

In the optimal control problem faced by an infinitesimal agent, the flow of population densities is given. From standard optimal control theory (for example using dynamic programming), the best strategy can be characterized through the value function $u$ of the above optimal control problem for a typical agent, which satisfies a Hamilton-Jacobi-Bellman (HJB) equation. Together with the equilibrium condition on the distribution, we obtain that the equilibrium best response $\hat{\ctrl}$ is characterized by
$$
	\hat{\ctrl}(t,x) = \argmax_{\ctrl \in \RR^\ctrldim} \big\{ -L(x, m(t,\cdot), \ctrl, \nabla u(t,x))  \big\},
$$

where $(u,m)$ solves the following forward-backward PDE system:
\begin{subequations}\label{AMS-num-eq:PDE-system-MFG}
     \begin{empheq}[left=\empheqlbrace]{alignat=2}
     	\displaystyle
	0&=-\frac{\partial u}{\partial t}(t,x)  - \nu \Delta u(t,x) + H(x, m(t,\cdot), \grad u(t,x)), 
	&&\hbox{ in } [0,T) \times \TT^d,
	\label{AMS-num-eq:PDE-system-MFG-HJB}
	\\
	\displaystyle
	0&=\frac{\partial m}{\partial t}(t,x)  - \nu \Delta m(t,x) 
	\notag
	\\
	&\qquad- \diver\left( m(t,\cdot) \partial_p H (\cdot, m(t), \grad u(t,\cdot)) \right)(x), 
	&&\hbox{ in } (0,T] \times \TT^d,
	\label{AMS-num-eq:PDE-system-MFG-KFP}
	\\
	&u(T,x) = g(x, m(T,\cdot)), \qquad m(0,x) = m_0(x), 
	&&\hbox{ in } \TT^d.
	\label{AMS-num-eq:PDE-system-MFG-IT}
     \end{empheq}
\end{subequations}
Recall that $\nu = \frac{1}{2} \sigma^2$.

\begin{example}[$f$ depends separately on $\ctrl$ and $m$]
\label{AMS-num-ex:special-case-sep}
Consider the case where the drift is the control, \textit{i.e.}, $b(x,m,\ctrl) = \ctrl$ with $\ctrldim = d$, and the running cost is of the form $f(x,m,\ctrl) = L_0(x,\ctrl) + f_0(x,m)$ where $L_0(x,\cdot): \RR^d \ni \ctrl \mapsto L_0(x,\ctrl) \in \RR$ is strictly convex and such that $\lim_{|\ctrl|\to \infty} \min_{x\in \TT^d} \frac {L_0(x,\ctrl)} {|\ctrl|}=+\infty$. We set $H_0(x,p)= \max_{\ctrl\in \RR^d} \langle -p,\ctrl\rangle - L_0(x,\ctrl)$, which is convex with respect to $p$. Then
	\begin{equation}
	\label{AMS-num-ex:special-case-sep-H}
	H(x,m,p)
	 = \max_{\ctrl \in \RR^d} \{- L_0(x,\ctrl) -  \langle \ctrl , p  \rangle  \} - f_0(x,m)
	 =  H_0(x, p)  - f_0(x,m).
	\end{equation}
\end{example}
\begin{example}[Quadratic Hamiltonian with separate dependence]
\label{AMS-num-ex:special-case-sep-quadra}
	Consider the setting of Example~\ref{AMS-num-ex:special-case-sep} and assume $L_0(x,\ctrl) = \tfrac{1}{2}|\ctrl|^2$, \textit{i.e.}, $H_0(x,\cdot) = H^*_0(x,\cdot) = \tfrac{1}{2}|\cdot|^2$, where $H^*_0$ denotes the convex conjugate of $H_0$ with respect to the second variable. Then the maximizer in~\eqref{AMS-num-ex:special-case-sep-H} is
        $ - p$,  the Hamiltonian reads $H(x,m,p) = \tfrac{1}{2} |p|^2 -f_0(x,m)$
        and the equilibrium best response is $\hat{\ctrl}(t,x) = - \nabla u(t,x)$ where $(u,m)$ solves the PDE system
	\begin{subequations}
	\label{AMS-num-ex:special-case-sep-quadra-PDE}
     \begin{empheq}[left=\empheqlbrace]{alignat=2}
     	\displaystyle
	&-\frac{\partial u}{\partial t}(t,x)  - \nu \Delta u(t,x) + \frac{1}{2} |\grad u(t,x)|^2 = f_0(x, m(t,\cdot)), 
	&&\hbox{ in } [0,T) \times \TT^d,
	\\
	\displaystyle
	&\frac{\partial m}{\partial t}(t,x)  - \nu \Delta m(t,x) - \diver\left( m(t,\cdot) \grad u(t,\cdot) \right)(x) = 0, 
	&&\hbox{ in } (0,T] \times \TT^d,
	\\
	&u(T,x) = g(x, m(T,\cdot)), \qquad m(0,x) = m_0(x), 
	&&\hbox{ in } \TT^d.
     \end{empheq}
\end{subequations}
\end{example}
Note that in the first equation, the coupling is only through the source term in the right hand side.

\textbf{Mean field type control problem.} In the MFC problem, the value function is not a function of $(t,x)$ only. Indeed, from the definition~\eqref{AMS-num-eq:def-J-MFC}--\eqref{AMS-num-eq:dyn-X-general-MFC} one can see that changing the control $\ctrl$ also changes the density $m^{\ctrl}$ (which is not the case in the optimal control problem faced by an infinitesimal agent in a MFG, see~\eqref{AMS-num-eq:def-J-MFG}). Hence the value function of the central planner is a function of the population's distribution and dynamic programming arguments have been developed using this point of view~\cite{MR3258261,MR3501391,MR3652408,MR3631380,MR3739204,djete2019mckean}. However, it is still possible to characterize the solution through a system of PDEs over the finite dimensional space $\domT = [0,T] \times \dom$. This system can be obtained either from the value function or via calculus of variation (see \textit{e.g.}~\cite[Chapter 4]{MR3134900} for more details on the latter approach).  A necessary condition for the existence of a smooth feedback function $\ctrl^* $ achieving
$J^{MFC}( \ctrl^*)= \min J^{MFC}(\ctrl)$ is that
$$
  \ctrl^*(t,x) = \argmax_{\ctrl \in \RR^\ctrldim} \big\{ -L(x, m(t,\cdot), \ctrl, \nabla u(t,x))  \big\},
$$
where $(u,m)$ solve the
following system of partial differential equations
\begin{subequations} 
     \begin{empheq}[left=\empheqlbrace]{alignat=2}
	0 &=
	\displaystyle
	-\frac{\partial u} {\partial t} (t,x) - \nu \Delta u(t,x) + H( x, m(t,\cdot), \grad u(t,x)) 
	\notag
	\\
	&\qquad + \int_{\TT^d} \frac{\partial H} {\partial m}(\xi, m(t,\cdot), \grad u(t, \xi))(x) m(t,\xi) d\xi  , 
	&& \hbox{ in } (0,T] \times \TT^d,
\\
	0 &= 
	\displaystyle \frac{\partial m} {\partial t} (t,x)  - \nu \Delta  m(t,x)
	\notag 
	\\
	&\qquad - \diver\Bigl( m(t,\cdot) \partial_p H(\cdot, m(t),\grad u(t,\cdot))\Bigr)(x) ,
	&& \hbox{ in } [0,T) \times \TT^d,
	\\
  	& u(T,x) = g (x, m(T,\cdot)) 
	\notag \\
	&\qquad\qquad\qquad   + \int_{\TT^d}
\frac{\partial g} {\partial m} (\xi, m(T,\cdot))(x) m(T,\xi) d\xi, 
	&& \hbox{ in }  \TT^d,
	\\
	& m(0,x) = m_0(x),  && \hbox{ in } \TT^d.
\end{empheq}
\end{subequations}
Compared with the MFG system~\eqref{AMS-num-eq:PDE-system-MFG}, there are extra terms involving the partial derivatives with respect to $m$, which should be understood in the following sense: if $\varphi: L^2(\RR^d) \to \RR$ is differentiable, 
$$
	\frac{d}{d\varepsilon} \varphi(m + \varepsilon \tilde m)(x)_{\big| \varepsilon = 0} = \int_{\RR^d} \frac{\partial \varphi}{\partial m}(m)(\xi) \tilde m(\xi) d\xi.
$$
 See~\cite[Chapter 4]{MR3134900} for more details. If the cost functions and the drift function depends on the density only locally (\textit{i.e.}, only on the density at the current position of the agent), $\frac{\partial}{\partial m}$ becomes a derivative in the usual sense.

\textbf{Existence and uniqueness. } The interested reader is referred to \textit{e.g.}~\cite{PLL-CDF,MR2295621,Cardaliaguet-2013-notes} for details on the question of existence and uniqueness of solutions for the MFG PDE system, and to \textit{e.g.}~\cite{MR3392611,MR3498932} for the corresponding MFC system. We simply mention here that the existence of solutions can be typically be obtained when the mean-field interactions are local and smooth or occurring through a regularizing kernel as explained in~\cite{MR2295621}.  As for uniqueness, in the MFG setting, a general sufficient condition for uniqueness is the so-called Lasry-Lions \index{monotonicity condition}monotonicity condition~\cite{MR2295621}. Intuitively, it holds when the cost function discourages the players from gathering, \textit{i.e.}, from having a density taking large values. More precisely, when considering a local dependence on the distribution and when the terminal cost $g$ does not depend on $m$, a sufficient condition for uniqueness of a classical solution the MFG PDE system is the positive definiteness of the matrix:
$$
  \begin{pmatrix}
 - \frac{2}{m} \frac {\partial H} {\partial m}    \left(x, m, p\right) &   \frac {\partial } {\partial m}  \partial_p ^T H(x, m, p)
\\
 \frac {\partial } {\partial m}  \partial_p H(x, m, p) & 2D^2_{p,p} H(x, m, p)
 \end{pmatrix}
$$
for all $x\in \TT^d$, $m>0$ and $p\in \RR^d$.  If $H$ depends separately on $p$ and $m$, then $\frac {\partial } {\partial m}  \partial_p H (x, m, p)=0$ and the condition for MFG becomes: $H$ is strictly convex with respect to $p$ for $m>0$ and non-increasing with respect to $m$, or $H$ is convex with respect to $p$ and strictly decreasing with respect to $m$.

As for the MFC problem, an analogous condition has been proposed in~\cite{MR3392611}, when the mean field interactions are of local type. Letting $\mathcal{H}(x, m, p) = m H(x, m, p)$, a sufficient condition is that: for every $x \in \TT^d$ and $m>0$, $p \mapsto \mathcal{H}(x, m, p)$ is strictly convex, and for every $x \in \TT^d$ and $p \in \RR^d$, $\RR_+  \ni m \mapsto \mathcal{H}(x, m, p)$ is strictly concave.

An interesting example in which the Lasry-Lions monotonicity condition can be understood quite intuitively is crowd motion with congestion. We come back to this point below in~\S~\ref{subsec:num-evacuation-crowd}.

\subsection{A Finite difference scheme}
\label{AMS-num-sec:fin-diff}

In this section, we present a finite-difference scheme first introduced in~\cite{MR2679575}; see also~\cite{MR3135339}. We consider the special case described in Example~\ref{AMS-num-ex:special-case-sep} and we focus on the case of local interactions, so that $H(x,m,p) = H_0(x,p) - f_0(x,m)$ for all $x \in \RR^d, m \in [0,+\infty), p \in \RR^d$. Similar methods have been  applied and at least partially analyzed in situations when the Hamiltonian does not depend separately on $m$ and $p$ (for example models addressing congestion, see \textit{e.g.}~\cite{YAJML}). 

To alleviate the notation, we present the scheme in the one-dimensional setting, \textit{i.e.}, $d=1$, so the domain is the one-dimensional unit torus, denoted by $\TT$.

\paragraph{\textbf{Discretization.} }

Let $N_T$ and $N_h$ be two positive integers. We consider $(N_T+1)$ and $(N_h+1)$ points in  time and space respectively. For any integers $i < j$, let $\llbracket i,j \rrbracket = \{i,\dots,j\}$ and $\llbracket j \rrbracket = \llbracket 0, j \rrbracket$. Let $\Delta t = T/N_T$ and $h = 1/N_h$, and $t_n = n \times \Delta t, x_i = i \times h$ for $(n,i) \in \llbracket N_T \rrbracket \times \llbracket N_h \rrbracket$. 

We approximate $u$ and $m$ respectively by vectors $U$ and $M \in \RR^{(N_T+1)\times(N_h+1)}$, that is, $u(t_n,x_i) \approx U^{n}_{i}$ and $m(t_n,x_i) \approx M^{n}_{i}$ for each $(n,i) \in \llbracket N_T \rrbracket \times \llbracket N_h \rrbracket$. We use a superscript and a subscript respectively for the time and space indices. Since we consider periodic boundary conditions, we slightly abuse notation and for any $W \in \RR^{N_h+1}$, we identify $W_{N_h+1}$ with $W_1$, and $W_{-1}$ with $W_{N_h-1}$. The periodic boundary condition will be translated into the constraint $W_{N_h} = W_0$.

We introduce the finite difference operators
\begin{align*}
	(D_t W)^n &= \frac{1}{\Delta t}(W^{n+1} - W^n), \qquad &&n \in \llbracket N_T-1 \rrbracket, \qquad W \in \RR^{N_T+1},
	\\
	(D W)_i &= \frac{1}{h} (W_{i+1} - W_{i}), \qquad &&i  \in \llbracket N_h \rrbracket, \qquad W \in \RR^{N_h+1},
	\\
	(\Delta_h W)_i &= -\frac{1}{h^2} \left(2 W_i - W_{i+1} - W_{i-1}\right), \qquad &&i \in \llbracket N_h \rrbracket, \qquad W \in \RR^{N_h+1},
	\\
	[\grad_h W]_i &= \left( (D W)_{i}, (D W)_{i-1} \right)^\top, \qquad &&i \in \llbracket N_h \rrbracket, \qquad W \in \RR^{N_h+1}.
\end{align*}

\paragraph{\textbf{Discrete Hamiltonian.} }
Let $\tilde{H}_0: \TT \times \RR \times \RR \to \RR, (x,p_1,p_2) \mapsto \tilde{H}_0(x,p_1,p_2)$ be  a discrete Hamiltonian, assumed to satisfy the following properties:
\begin{enumerate}[label=\normalfont{\textbf{($\mathrm{\tilde H}$\arabic{*})}}, ref=\textbf{($\mathrm{\tilde H}$\arabic{*})}]
	\it
	\item\label{AMS-num-hyp-discH-monotonicity} Monotonicity:  for every $x \in \TT$, $\tilde{H}_0$ is nonincreasing in $p_1$ and nondecreasing in $p_2$.
	\item\label{AMS-num-hyp-discH-consistency}  Consistency:  for every $x \in \TT, p \in \RR$, $\tilde{H}_0(x,p,p) = H_0(x,p)$.
	\item\label{AMS-num-hyp-discH-differentiability}  Differentiability:  for every $x \in \TT$, $\tilde{H}_0$ is of class $\mathcal C^1$ in $p_1,p_2$.
	\item\label{AMS-num-hyp-discH-convexity}  Convexity:   for every $x \in \TT$, $(p_1,p_2) \mapsto \tilde{H}_0(x,p_1,p_2)$ is convex. 
\end{enumerate}
We refer to~\cite{MR2679575,MR3135339} for more details.

\begin{example}
\label{AMS-num-eq-quadratic-discrete}
For instance, if $H_0(x,p) = \frac{1}{2} |p|^2$, then one can take $\tilde{H}_0(x, p_1, p_2) = \frac{1}{2}|P_K(p_1,p_2)|^2$ where $P_K$ denotes the projection on $K = \RR_- \times \RR_+$.
\end{example}

\begin{remark}
  Analogously, for $d$-dimensional problems, the discrete Hamiltonians that we consider are real valued functions defined on $\TT^d \times (\RR^d)^2 $.
\end{remark}

\paragraph{\textbf{Discrete HJB equation.} }
We consider the following discrete version of the HJB equation~\eqref{AMS-num-eq:PDE-system-MFG-HJB}:
\begin{subequations}\label{AMS-num-eq:discrete-HJB-finitediff}
     \begin{empheq}[left=\empheqlbrace]{align}
     	\displaystyle
	\,\, &- (D_t U_{i})^{n} - \nu (\Delta_h U^{n})_{i}
	+ \tilde{H}_0(x_i, [\grad_h U^n]_i) 
	\\
	& \qquad = f_0(x_i,M^{n+1}_i) \, , 
	&& (n,i) \in \llbracket N_T-1 \rrbracket \times \llbracket N_h \rrbracket \, ,
	\notag
	\\
	\,\, &U^{n}_{0} = U^{n}_{N_h} \, ,  && n \in \llbracket N_T-1 \rrbracket \, ,
	\\
	\,\, &U^{N_T}_{i} = g(x_i, M^{N_T}_i) \, ,  && i \in \llbracket N_h \rrbracket \, .
	\end{empheq}
\end{subequations}

Note that it is an implicit scheme since the equation is backward in time.

\paragraph{\textbf{Discrete KFP equation.} } \label{AMS-num-sec:textbfd-kfp-equat}
To define an appropriate discretization of the KFP equation~\eqref{AMS-num-eq:PDE-system-MFG-KFP}, we start by considering the weak form. For a smooth test function $w \in \mathcal{C}^{\infty}([0,T] \times \TT)$, it involves, among other terms,  the expression
\begin{align}
	&- \int_{\TT} \partial_x \big( \partial_p H_0(x, \partial_x u(t,x)) m(t,x)\big) w(t,x) dx
	\notag
	\\
	&=
	\int_{\TT} \partial_p H_0(x, \partial_x u(t,x)) m(t,x) \, \partial_x w(t,x) dx  \, ,
	\label{AMS-num-eq:HJB-weak}
\end{align}
where we used an integration by parts and the periodic boundary conditions. In view of what precedes, it is quite natural to propose the following discrete version of the right hand side of \eqref{AMS-num-eq:HJB-weak}:
\begin{displaymath}
 h \sum_{i=0}^{N_h-1} M_i^{n+1} \left(  \partial_{p_1} \tilde{H}_0(x_i, [\grad_h U^n]_i)  \frac { W^n_{i+1}-W^n_i} h+
    \partial_{p_2} \tilde{H}_0(x_i, [\grad_h U^n]_i)  \frac { W^n_{i}-W^n_{i-1}} h   \right).
\end{displaymath}
Performing a discrete integration by parts, we obtain the discrete counterpart  of the left hand side of \eqref{AMS-num-eq:HJB-weak} as follows:
$ \displaystyle -
  h \sum_{i=0}^{N_h-1}  	\mathcal{T}_i(U^n , M^{n+1}) W_i^n$,
where $\mathcal{T}_i$ is the following discrete transport operator:
\begin{align*}
	\mathcal{T}_i(U, M)
	= \frac{1}{h} 
	&\Big( M_i  \partial_{p_1} \tilde{H}_0(x_i, [\grad_h U]_i) 
		- M_{i-1} \partial_{p_1} \tilde{H}_0(x_{i-1}, [\grad_h U]_{i-1})
		\\
		&\quad +
		M_{i+1} \partial_{p_2} \tilde{H}_0(x_{i+1}, [\grad_h U]_{i+1})
		- M_i \partial_{p_2} \tilde{H}_0(x_i, [\grad_h U]_i)
	\Big) \, .
\end{align*}

Then, for the discrete version of~\eqref{AMS-num-eq:PDE-system-MFG-KFP}, we consider:
\begin{subequations}\label{AMS-num-eq:discrete-KFP-finitediff}
     \begin{empheq}[left=\empheqlbrace]{align}
	\,\, &(D_t M_{i})^{n} - \nu (\Delta_h M^{n+1})_{i}
	\label{AMS-num-eq:sysH-FP-scheme-discrete-edo}
	\\
	& \qquad - \mathcal{T}_i(U^{n}, M^{n+1})
	= 0 \, , 
	&& (n,i) \in \llbracket N_T-1 \rrbracket \times \llbracket N_h \rrbracket \, ,
	\notag 
	\\
	\,\, &M^{n}_{0} = M^{n}_{N_h} \, , &&n \in \llbracket 1, N_T \rrbracket \, , 
	\label{AMS-num-eq:sysH-FP-scheme-discrete-bc}
	\\
	\,\, &M^{0}_{i} = \bar m_0(x_i) \, , \,  && i \in \llbracket N_h \rrbracket \, ,
	\label{AMS-num-eq:sysH-FP-scheme-discrete-initialc}
\end{empheq}
\end{subequations}
where, 
\begin{equation}
  \label{AMS-num-eq:discrete-barm-0}
	\bar m_0(x_i) = \int_{|x - x_i| \le h/2} m_0(x) dx.
\end{equation}
Here again, the scheme is implicit since the equation is forward in time.

\begin{remark}[Structure of the discrete system]
The finite difference system~\eqref{AMS-num-eq:discrete-HJB-finitediff}--\eqref{AMS-num-eq:discrete-KFP-finitediff} preserves the structure of the PDE system~\eqref{AMS-num-eq:PDE-system-MFG} in the following sense: The operator $M \mapsto - \nu (\Delta_h M)_{i}
	- \mathcal{T}_i(U, M)$
	is the adjoint of the linearization of the operator
	$U \mapsto - \nu (\Delta_h U)_{i}
	+ \tilde{H}_0(x_i, [\grad_h U]_i)$. Indeed, 
	\begin{align*}
		& \sum_{i} \mathcal{T}_i(U, M) W_i
		= - \sum_{i} M_i  \left\langle \partial_p \tilde H_0(x_i, [\grad_h U]_i) , [\grad_h W]_i \right\rangle \, .
	\end{align*}
\end{remark}

\paragraph{\bf Convergence results.}

Existence and uniqueness for the discrete system have been proved in~\cite[Theorems 6 and 7]{MR2679575}.
The monotonicity properties ensure that the grid function $M$ is nonnegative. By construction of $\mathcal{T}$, the scheme preserves the total mass $h \sum_{i=0}^{N_h-1} M_i^n$. Note that there is no restriction on the time step since the scheme is implicit. Thus, this method may be used for long horizons and
the scheme can be very naturally adapted to ergodic MFGs, see~\cite{MR2679575}.

Furthermore, convergence results are available. A first type of convergence theorems see ~\cite{MR2679575,MR2888257,MR3097034} (in particular~\cite[Theorem 8]{MR2679575} for finite horizon problems) make the assumption that the MFG system of PDEs has a unique classical solution and strong versions of Lasry-Lions monotonicity assumptions, see~\cite{MR2269875,MR2271747,MR2295621}. Under such assumptions, 
the solution to the discrete system converges towards the classical solution as the grid parameters $h$ and $\Delta t$ tend to zero. Another type of results, obtained in~\cite{MR3452251}, is the  convergence of the solution of the discrete problem to weak solutions of the system of forward-backward PDEs. We refer to~\cite{MR3452251} for the precise statement but let us stress that these results have been proved without assuming the existence of a (weak) solution to the MFG PDE system, nor  Lasry-Lions monotonicity assumptions. This approach can thus be seen as an alternative proof of existence of weak solutions of the MFG PDE system. Besides the setting presented here, similar finite difference schemes have been developed for mean field games with interaction through the law of the controls~\cite{achdouKobeissi2020mfgcfinitediff} or in a time-fractional setting~\cite{camilli2020approximation-timefrac}.

\subsection{A Semi-Lagrangian scheme}

In this section, we present an alternative numerical scheme which relies on a Lagrangian viewpoint instead of a Eulerian point of view as in the aforementioned finite difference scheme. Intuitively, the Lagrangian approach corresponds to idea of following the dynamics of typical player. In~\cite{MR3148086,MR3392626,MR3828859} Carlini and Silva have developed a semi-Lagrangian scheme for MFG in which the diffusion term can be absent or degenerate.  For ease of presentation, take $d=1$, $\dom = \RR$, and let us focus on the case without viscosity (first order setting). The Lagrangian point of view is particularly relevant in this situation, because in the absence of noise, a trajectory is completely determined by the initial position and the control. More precisely, if $b(x,m,\ctrl) = \ctrl$ and $\nu = 0$, then the solution of the state equation~\eqref{AMS-num-eq:dyn-X-general-MFG} is given by
$$
	X_t^{\ctrl} = X_0^{\ctrl} + \int_{0}^t \ctrl(s, X_s^{\ctrl}) ds, \qquad t \ge 0.
$$
Taking, as in Example~\ref{AMS-num-ex:special-case-sep-quadra}, a running cost function of the form $f(x,m,\ctrl) = \frac{1}{2}|\ctrl|^2 + f_0(x,m)$ and a terminal cost function $g(x, m)$ (where $f_0$ and $g$ depend on $m \in L^2(\dom)$ in a potentially non-local way) leads to the following MFG PDE system:
\begin{subequations}\label{AMS-num-eq:PDE-system-MFG}
     \begin{empheq}[left=\empheqlbrace]{alignat=2}
     	\displaystyle
	&-\frac{\partial u}{\partial t}(t,x)  + \frac{1}{2} |\grad u(t,x)|^2 = f_0(x, m(t,\cdot)), 
	&&\hbox{ in } [0,T) \times \RR,
	\label{AMS-num-eq:PDE-system-MFG-HJB-CarliniSilva}
	\\
	\displaystyle
	&\frac{\partial m}{\partial t}(t,x) - \diver\left( m(t,\cdot) \grad u(t,\cdot) \right)(x) = 0, 
	&&\hbox{ in } (0,T] \times \RR,
	\label{AMS-num-eq:PDE-system-MFG-KFP-CarliniSilva}
	\\
	&u(T,x) = g(x, m(T,\cdot)), \qquad m(0,x) = m_0(x), 
	&&\hbox{ in } \RR.
	\label{AMS-num-eq:PDE-system-MFG-IT-CarliniSilva}
     \end{empheq}
\end{subequations}
This amounts to taking $\nu=0$ (and changing the domain) in system~\eqref{AMS-num-ex:special-case-sep-quadra-PDE}.

\paragraph{\textbf{Discrete HJB equation. }} Given a flow of densities $m = (m_t)_{t\in[0,T]}$, the corresponding value function $u$ admits the following representation formula:
$$
	u(t,x) = \inf_{\ctrl \in L^2([t,T]; \RR)} \int_t^T \left[ \frac{1}{2}|\ctrl(s)|^2 + f_0(X^{\ctrl,t,x}_s, m(s,\cdot))\right] ds + g(X^{\ctrl,t,x}_T, m(T,\cdot)),
$$
where $X^{\ctrl,t,x}$ starts from $x$ at time $t$ and is controlled by $\ctrl$.

Based on this intuition, let us consider the equation:
\begin{equation}
\label{AMS-num-eq:semiLagrangian-discrete-HJB}
\begin{cases}
	U^n_i = S_{\Delta t, h}[m](U^{n+1}, i, n), & (n,i) \in \llbracket N_T-1 \rrbracket \times  \ZZ,
	\\
	U^{N_T}_i = g(x_i, m(T, \cdot)), & i \in \ZZ,
\end{cases}
\end{equation}
where $S_{\Delta t, h}$ is defined as 
\begin{equation}
\label{AMS-num-eq:semiLagrangian-S_Deltat-h}
	S_{\Delta t, h}[m](W, n, i) = \inf_{\ctrl \in \RR^d} \left\{ I[W](x_i + \ctrl \, \Delta t) + \frac{1}{2} |\ctrl|^2 \, \Delta t  \right\} + f_0(x_i, m(t_n, \cdot)) \, \Delta t ,
\end{equation}
with $I: \cB(\ZZ) \to \mathcal{C}_b(\RR)$ denoting the interpolation operator defined as
$$ 
	I[W](\cdot) = \sum_{i \in \ZZ} W_i \beta_i(\cdot), 
$$
where $\cB(\ZZ)$ is the set of bounded functions from $\ZZ$ to $\RR$, and $\beta_i =  \left[ 1 - \frac{|x - x_i|}{h} \right]_+$ is the triangular function with support $[x_{i-1}, x_{i+1}]$ and such that $\beta_i(x_i)=1$.

From the solution $U = (U^n_i)_{n,i}$ of the discrete HJB equation~\eqref{AMS-num-eq:semiLagrangian-discrete-HJB} for a given density flow $m$, we interpolate it to construct the following function $u_{\Delta t, h}[m](x,t):  [0,T] \times \RR \to \RR$,
$$
	u_{\Delta t, h}[m](t,x) = I[U^{[\frac{t}{\Delta t}]}](x), \qquad (t,x) \in [0,T] \times \RR.
$$

\paragraph{\textbf{Discrete KFP equation. }}

In order to write a discrete version of the KFP equation, the solution of the discrete HJB equation is replaced by a regularized version. Let $\rho \in \mathcal{C}^\infty_c(\RR)$ with $\rho \ge 0$ and $\int \rho(x) dx = 1$. For $\epsilon>0$, let us consider the mollifier $\rho_\epsilon(x) = \frac{1}{\epsilon} \rho(x/\epsilon)$ and define 
\begin{equation}
\label{AMS-num-eq:SL-u-epsilon-Deltat-h}
	u^{\epsilon}_{\Delta t, h}[m](t,\cdot) = \rho_\epsilon * u_{\Delta t, h}[m](t,\cdot), \qquad t \in [0,T].	
\end{equation}
Then let us introduce, for $(t,x) \in [0,T] \times \RR$, the induced control
\begin{equation}
\label{AMS-num-eq:SL-hat-ctrl-epsilon-Deltat-h}
	\hat \ctrl^{\epsilon}_{\Delta t, h}[m](t,x) = - \nabla u^{\epsilon}_{\Delta t, h}[m](t,x),
\end{equation}
and its discrete counter part: for $(n,i) \in \llbracket N_T \rrbracket \times \ZZ$,
$$
	\hat \ctrl^{\epsilon}_{n,i} = \hat \ctrl^{\epsilon}_{\Delta t, h}[m](t_n,x_i).
$$
Then define the following discrete flow: for $(n,i) \in \llbracket N_T-1 \rrbracket \times \ZZ$, 
$$
	\Phi^{\epsilon}_{n,n+1,i}[m] = x_i + \hat \ctrl^{\epsilon}_{\Delta t, h}[m](t_n,x_i) \Delta t \, .
$$

\begin{figure}[h]
	\begin{subfigure}{.9\columnwidth}
		\centering
		\includegraphics[width=\columnwidth]{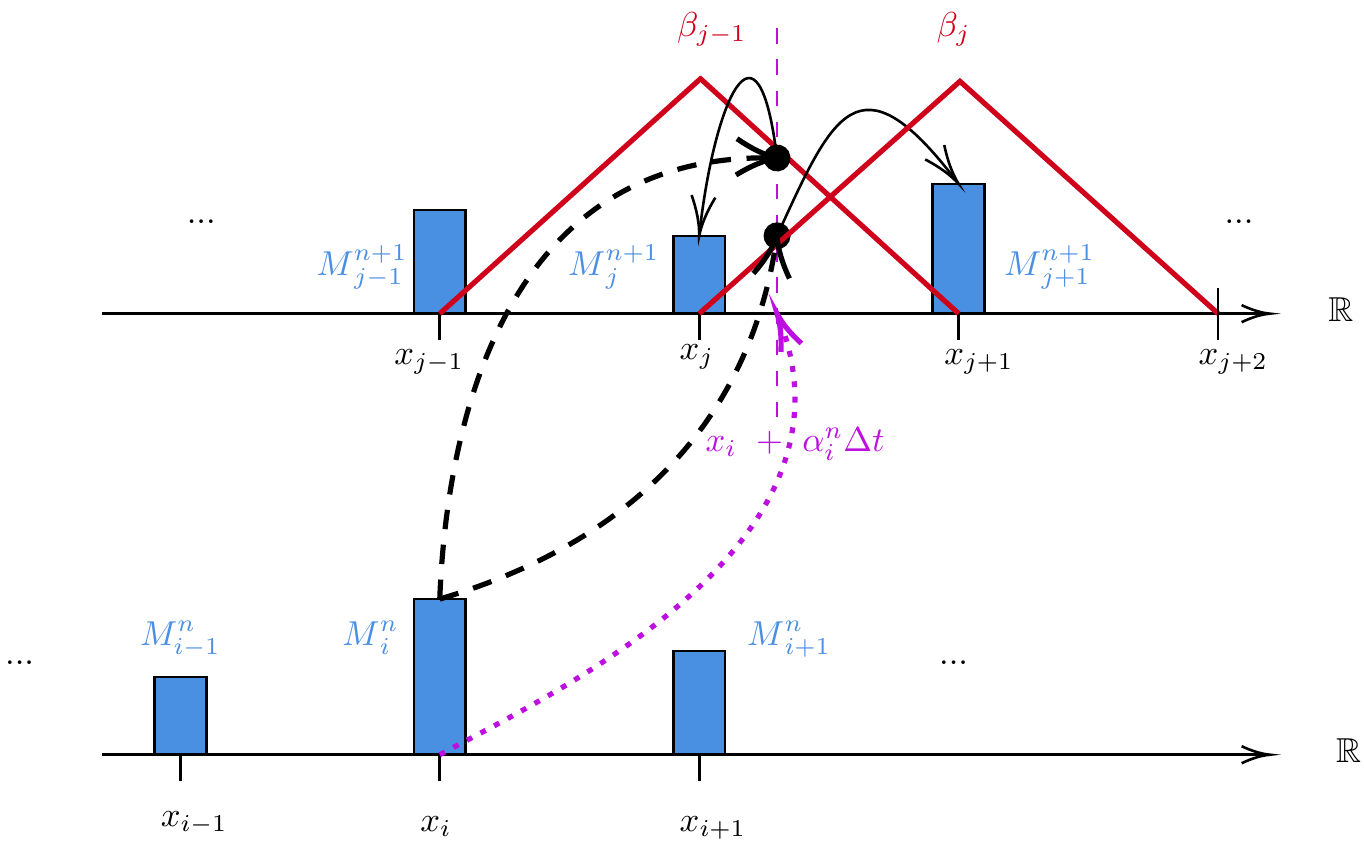}
	\end{subfigure}%
	\caption{Illustration of the evolution of the mass as given by~\eqref{AMS-num-eq:semiLagrangian-discrete-KFP} and explained in the text. The bottom row corresponds to time $n$ and the top row to time $n+1$.} 
 	\label{AMS-num-fig:SL-mass-evol}
\end{figure}

We can now introduce the discrete KFP equation for $M^{\epsilon}[m] = (M^{\epsilon,n}_i[m])_{n,i}$:
\begin{equation}
\label{AMS-num-eq:semiLagrangian-discrete-KFP}
\begin{cases}
	M^{\epsilon, n+1}_i[m] = \sum_{j} \beta_i \left( \Phi^{\epsilon}_{n, n+1, j}[m]\right) M^{\epsilon, n}_j[m]  , & (n,i) \in \llbracket N_T-1 \rrbracket \times  \ZZ,
	\\
	M^{\epsilon, 0}_i[m] = \int_{[x_i - h/2, x_i+h/2]} m_0(x) dx, & i \in \ZZ.
\end{cases}
\end{equation}
Intuitively, the idea is as depicted by Figure~\ref{AMS-num-fig:SL-mass-evol}: First, if control $(\ctrl^n_i)_{i,n}$, an agent starting at time $n$ at point $x_i$ is supposed to arrive at time $n+1$ at point $x_i + \alpha^n_i \Delta t$ (dotted arrow in fuschia). Supposing that this point is between the grid points $x_{j-1}$ and $x_j$, the mass $M^n_i$ (rectangle in blue) at $x_i$ is going to be split between $x_{j-1}$ and $x_j$. The proportion of mass moving to each point is proportional to the values of the hat functions $\beta_{j-1}$ and $\beta_{j}$ (in red) at the arrival point $x_i + \alpha^n_i \Delta t$.

From here, we can recover a function by defining $m^{\epsilon}_{\Delta t, h}[m]: [0,T] \times \RR \to \RR$ as: for $n \in \llbracket N_T-1 \rrbracket$, for $t \in [t_n, t_{n+1})$,
\begin{align}
	\label{AMS-num-eq:m-epsilon-Deltat-h}
	m^{\epsilon}_{\Delta t, h}[m](t,x) 
	&= \frac{1}{\rho} \left[ \frac{t_{n+1}-t}{\Delta t} \sum_{i \in \ZZ} M^{\epsilon, n}_{i}[m] \indic_{[x_i - h/2, x_i+h/2]}(x)  \right.
	\\
	&\qquad\qquad\qquad 
	\left. + \frac{t-t_{n}}{\Delta t} \sum_{i \in \ZZ} M^{\epsilon,n+1}_{i}[m] \indic_{[x_i - h/2, x_i+h/2]}(x) \right] \, .
	\notag
\end{align}
The goal is then to solve the following fixed-point problem: Find $\hat M$ such that
$$
	\hat M^{n}_i = M^n_i\big[m^{\epsilon}_{\Delta t, h}[\hat M]\big].
$$
Recall that $(M^n_i)_{i,n}$ depends on $\hat \ctrl^{\epsilon}_{\Delta t, h}$, which itself depends on $(U^n_i)_{i,n}$. Hence, as in the finite difference scheme, the equations for the distribution and the value function are coupled.   
Convergence of the scheme towards the continuous solution of the PDE system has been proved under suitable conditions even in the first order ($\nu=0$) as in the present discussion or degenerate case, see~\cite{MR3148086,MR3392626,MR3828859} for the details.

\clearpage

\section{\bf Algorithms to solve the discrete schemes}
\label{AMS-num-sec:sol-strat}

In this section, we review several algorithms to solve the numerical schemes introduced in the previous section. Their forward-backward structure is the main challenge. We revisit the methods presented in Section~\ref{AMS-num-sec:LQ}, namely, fixed point iterations and Newton iterations.

\subsection{Picard (fixed point) iterations}

Probably the most straightforward method is to iteratively solve the (discrete) HJB equation and the (discrete) FP equation in order to update respectively the estimate of the value function and the density, using in each equation the most recent estimate of the other function. As presented in Section~\ref{AMS-num-sec:LQ}, a damping coefficient can be introduced, which can even account for fictitious play type updates.  
This approach can be used with the finite difference scheme or the semi-Lagrangian scheme presented earlier in Section~\ref{AMS-num-sec:PDE-scheme}. We illustrate it here with the semi-Lagrangian scheme in the spirit of~\cite{MR3148086}. In this context, the iterations take the form described in Algorithm~\ref{AMS-num-algo:semi-Lagrangian-Picard}. Note that the step corresponding to~\eqref{AMS-num-eq:semiLagrangian-S_Deltat-h} requires computing an infimum over the controls. In the implementation, we can replace $\RR^d$ by a bounded set, which is then discretized, so that the infimum is taken over a finite number of values. Instead of fixing a priori the number of iterations, we can use a stopping criterion of the form:
$$
	\|U^{(\epsilon, \mathtt{k}+1)} - U^{(\epsilon, \mathtt{k})}\| < \varepsilon, \quad \hbox{ and}  \quad \|M^{(\epsilon, \mathtt{k}+1)} - M^{(\epsilon, \mathtt{k})}\| < \varepsilon,
$$
for some threshold $\varepsilon>0$.

\begin{algorithm}[H]
\DontPrintSemicolon
  
  \KwInput{Initial guess $\tilde M$; damping $\delta(\cdot)$; number of iterations $\mathtt{K}$}
  \KwOutput{Approximation of $(\hat U^{\epsilon}, \hat M^{\epsilon})$ solving~\eqref{AMS-num-eq:semiLagrangian-discrete-KFP}\&\eqref{AMS-num-eq:semiLagrangian-discrete-HJB}}
  Initialize $\tilde{M}^{(\epsilon, 0)} = \tilde M$ \;
        \For{$\mathtt{k}=0,1,2,\dots, \mathtt{K}-1$}    
        { 
        	Let $U^{(\mathtt{k}+1)}$ be a solution of~\eqref{AMS-num-eq:semiLagrangian-discrete-HJB} with $m$ replaced by $m^{\epsilon}_{\Delta t, h}[\tilde{M}^{(\epsilon, \mathtt{k})}]$ defined by~\eqref{AMS-num-eq:m-epsilon-Deltat-h}\;
	Let $U^{(\epsilon, \mathtt{k}+1)}$ be defined by~\eqref{AMS-num-eq:SL-u-epsilon-Deltat-h} with the convolution replaced by a discrete convolution\;
	Let $M^{(\epsilon, \mathtt{k}+1)}$ be solution of~\eqref{AMS-num-eq:semiLagrangian-discrete-KFP} where~\eqref{AMS-num-eq:SL-hat-ctrl-epsilon-Deltat-h} is replaced with a discrete version obtained using $U^{(\epsilon, \mathtt{k}+1)}$ and the gradient is approximated by a finite difference operator \;
	Let $\tilde{M}^{(\epsilon, \mathtt{k}+1)} = \delta(\mathtt{k}) \tilde{M}^{(\epsilon, \mathtt{k})} + (1-\delta(\mathtt{k})) M^{(\epsilon, \mathtt{k}+1)}$\;
	}
	\Return{$(U^{(\epsilon, \mathtt{K})}, M^{(\epsilon, \mathtt{K})})$} 
\caption{Fixed-point iterations with semi-Lagrangian scheme \label{AMS-num-algo:semi-Lagrangian-Picard}}
\end{algorithm}

\subsection{Numerical illustration: an example with concentration of mass}

\begin{figure}%
	\centering
	\begin{subfigure}[t]{0.45\linewidth}%
		\centering
		\includegraphics[width=\textwidth]{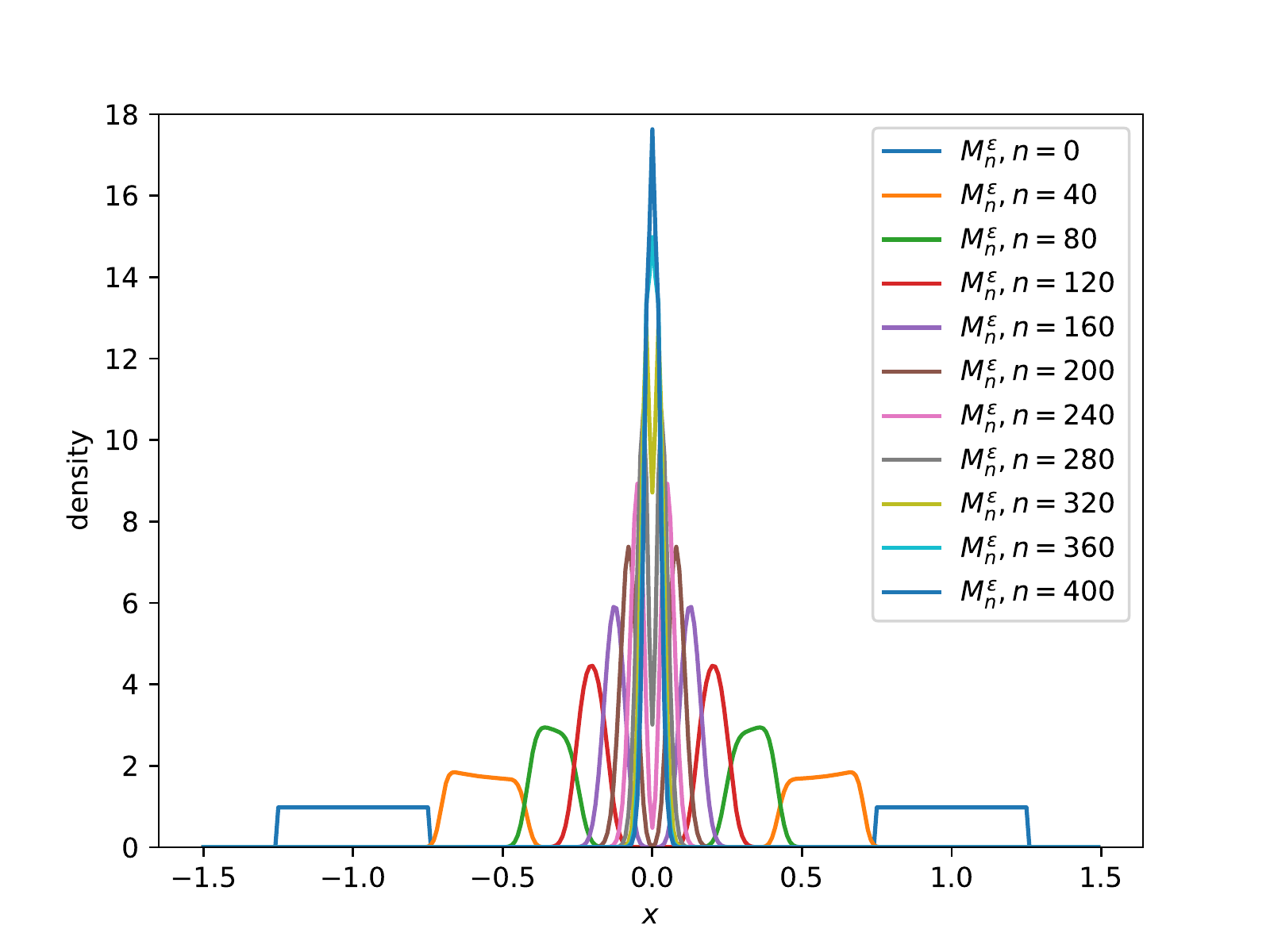}  
	\end{subfigure}
	\quad
	\begin{subfigure}[t]{0.45\linewidth}
		\centering
	  	\includegraphics[width=\textwidth]{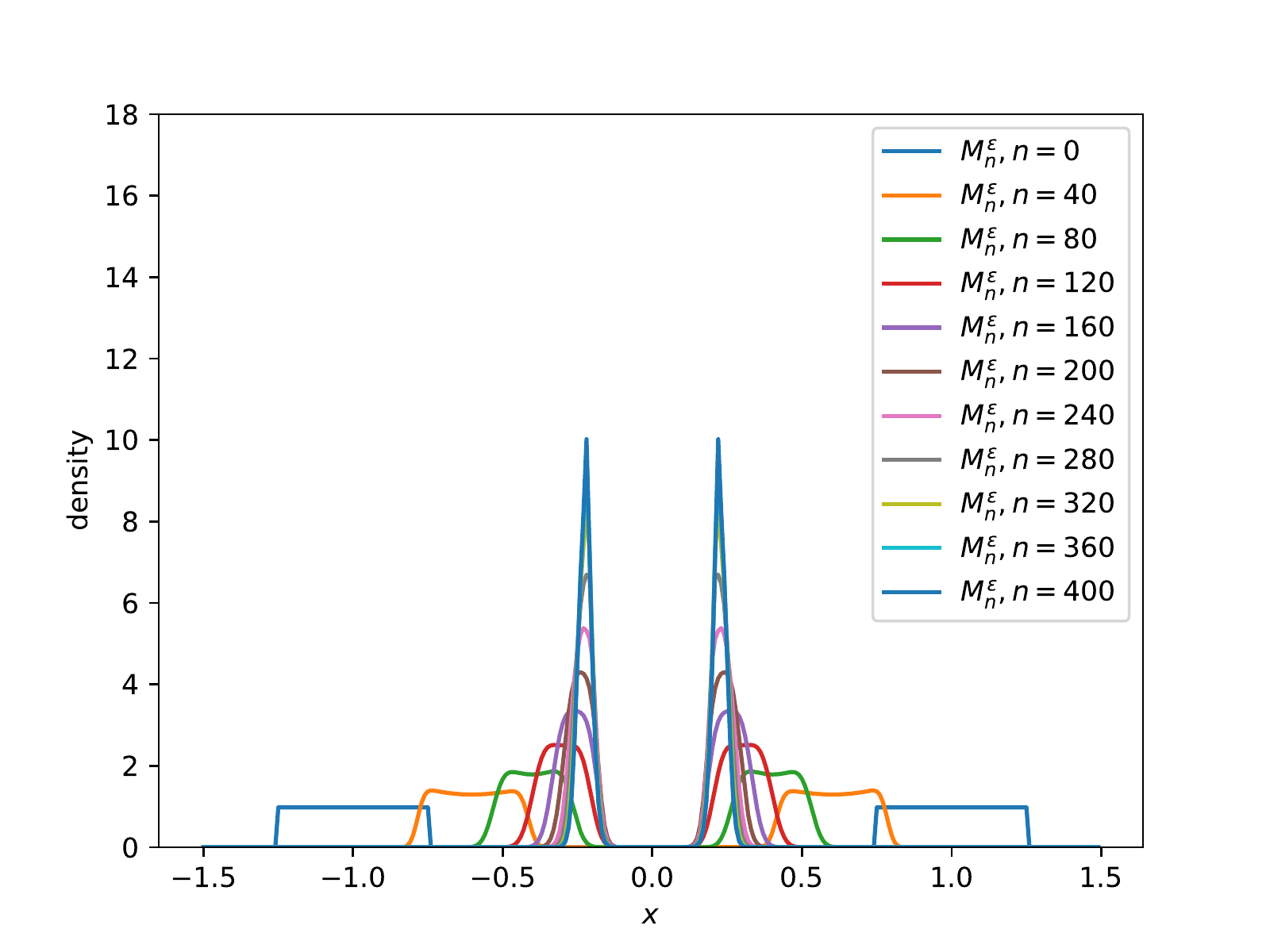}
	\end{subfigure}
	\caption{\label{fig:SEMILAG-numres-m} Evolution of the distribution for the example solved with the semi-Lagrangian scheme. Left: $\kappa_{MF} = 0.5$; right: $\kappa_{MF} = 0.9$.}
\end{figure}

We borrow an example from~\cite{MR3148086}. 
There is no terminal cost, \textit{i.e.}, $g \equiv 0$, and the running cost is:
$$
	f(x,m,\ctrl)
	= \frac{1}{2} |\ctrl|^2 + (x - c^*)^2 + \kappa_{MF} V(x,m),
$$
where $c^* \in \RR$ plays the role of a target position, $\kappa_{MF} \ge 0$ is a parameter for the strength of the mean-field interactions. The first term penalizes large velocity, the second term penalizes being far from the target, and the last term involves mean-field interactions. $V$ is a non-local interaction term which penalizes having a high density around the position $x$, defined by:
$$
	V(x,m) = \rho_{\sigma_V} * \big(\rho_{\sigma_V} * m\big)(x),
$$
where we recall that $\rho_{\sigma_V}$ is a mollifier and $*$ denotes the convolution. For the numerical implementation, we used $c^*=0$, $\sigma_V = 0.25$. The initial distribution is composed of two uniform parts: on $[-1.25, -0.75]$ and on $[0.75, 1.25]$. In Figure~\ref{fig:SEMILAG-numres-m}, we compare the evolution of the density for two values of $\kappa_{MF}$: $0.5$ or $0.9$. In the first case, the mass concentrates around $c^*=0$. In the second case, due to the higher penalty, the mass does not concentrate at $0$ and the two parts of the density remain separate, with lower peaks. We used $T = 4$ with $400$ time steps and the domain $[-2.5,2.5]$, although for the sake of clarity we show the results for the truncated domain $[-1.5, 1.5]$. It can be noticed that at each time step the density remains zero on most of the domain, which is due to the fact that there is no diffusion in this model (\textit{i.e.}, $\nu=0$ in our notation).

\clearpage

\subsection{Newton method}
\label{sec:NewtonPDE}

Another solution strategy, proposed by Achdou \textit{et al.}~\cite{MR2888257}, consists in using Newton method to solve the finite-difference system~\eqref{AMS-num-eq:discrete-HJB-finitediff}\&\eqref{AMS-num-eq:discrete-KFP-finitediff}. The main idea is to directly look for a zero of the function 
$$
	\varphi = (\varphi_{\cU}, \varphi_{\cM})^\top
$$
 with $\varphi_{\cU}$ and $\varphi_{\cM}$ defined such that, for $U$ and $M \in \RR^{(N_T+1)\times(N_h+1)}$,
$$
\begin{cases}
	\varphi_{\cU}(U,M) = 0
	&\Leftrightarrow
	\hbox{ $(U,M)$ solves discrete HJB equation~\eqref{AMS-num-eq:discrete-HJB-finitediff}}
	\\
	\varphi_{\cM}(U,M) = 0
	&\Leftrightarrow
	\hbox{ $(U,M)$ solves discrete KFP equation~\eqref{AMS-num-eq:discrete-KFP-finitediff}.}
\end{cases}
$$

Let $J_\varphi$ denote the Jacobian of $\varphi$. Newton method then consists in starting from an initial guess $(U^0,M^0)$ and iteratively computing $(U^{(\mathtt{k}+1)}, M^{(\mathtt{k}+1)})$ given by:
$$
	(U^{(\mathtt{k}+1)}, M^{(\mathtt{k}+1)})^\top = (U^{(\mathtt{k})}, M^{(\mathtt{k})})^\top  - J_\varphi(U^{(\mathtt{k})}, M^{(\mathtt{k})})^{-1} \varphi(U^{(\mathtt{k})}, M^{(\mathtt{k})})
$$
or rather solving for $(\tilde U^{(\mathtt{k}+1)}, \tilde M^{(\mathtt{k}+1)})$:
$$
	J_\varphi(U^{(\mathtt{k})}, M^{(\mathtt{k})})  (\tilde U^{(\mathtt{k}+1)}, \tilde M^{(\mathtt{k}+1)}) = - \varphi(U^{(\mathtt{k})}, M^{(\mathtt{k})}),
$$
and then setting $(U^{(\mathtt{k}+1)}, M^{(\mathtt{k}+1)})^\top = (\tilde U^{(\mathtt{k}+1)}, \tilde M^{(\mathtt{k}+1)}) + (U^{(\mathtt{k})}, M^{(\mathtt{k})})^\top$.

Hence each step amounts to solving a linear system of the form:
\begin{equation}
\label{AMS-num-eq:Newton-system-sol}
	\begin{pmatrix}
		A_{\cU,\cU} & A_{\cU,\cM}
		\\
		A_{\cM,\cU} & A_{\cM,\cM}
	\end{pmatrix}
	\begin{pmatrix}
		U
		\\
		M
	\end{pmatrix}
	=
	\begin{pmatrix}
		G_{\cU}
		\\
		G_{\cM}
	\end{pmatrix},
\end{equation}
where
\begin{align*}
	&A_{\cU,\cM}(U,M) = \nabla_{U} \varphi_{\cM}(U,M), 
	\quad &&A_{\cU,\cU}(U,M) = \nabla_{U} \varphi_{\cU}(U,M), 
	\\ 
	&A_{\cM,\cU}(U,M) = \nabla_{M} \varphi_{\cU}(U,M), 
	\quad &&A_{\cM,\cM}(U,M) = \nabla_{M} \varphi_{\cM}(U,M). 
\end{align*}

For the numerical implementation, we note that $A_{\cU, \cM}$ and $A_{\cM, \cU}$ are block-diagonal, $A_{\cU, \cU} = A_{\cM, \cM}^\top$, and:
$$
	A_{\cU, \cU}
	=
	\begin{pmatrix}
	D_1 & -\frac{1}{\Delta t}\Id_{N_h} & \dots & \dots & 0 
	\\
	0 & D_2 & \ddots & 0 & \vdots 
	\\
	0 & \ddots & \ddots & \ddots & \vdots
	\\
	\vdots & \ddots & \ddots & \ddots & -\frac{1}{\Delta t}\Id_{N_h}
	\\
	0 & \ddots & 0 & 0 & D_{N_T}
	\end{pmatrix},
$$
where $D_n$ is the matrix corresponding to the discrete operator:
$$
	Z = (Z_{i,j})_{i,j}
	\mapsto
	\left( \frac{1}{\Delta t} Z_{i,j} - \nu (\Delta_h Z)_{i,j} + [\nabla_h Z]_{i,j} \cdot \nabla_q H_h(x_{i,j}, [\nabla_h U^{(\mathtt{k}), n}]_{i,j})
	\right)_{i,j},
$$
which comes from the linearization of the discrete HJB equation~\eqref{AMS-num-eq:discrete-HJB-finitediff}. 

We refer to~\cite[Section 4]{MR3135339} and~\cite{MR2928376} for more details on possible strategies to solve~\eqref{AMS-num-eq:Newton-system-sol}. We simply stress that, as usual with Newton method, the choice of the initial guess $(U^{(0)}, M^{(0)})$ is important. A possible choice is to use the solution to the corresponding ergodic problem (if it is known or if it can be computed easily). Another possibility is to exploit the fact that the method converges more easily when the viscosity coefficient  $\nu$ is large. It is thus possible to use a continuation method: we start by solving the problem with a large $\nu$, then use the solution as an initial guess for the problem with a smaller $\nu$, and so on until the desired viscosity coefficient is reached.

\subsection{Numerical illustration: evacuation of a room with congestion effects}
\label{subsec:num-evacuation-crowd}

\begin{figure}%
	\centering
	\begin{subfigure}[t]{0.4\linewidth}%
		\centering
		\includegraphics[height=4cm]{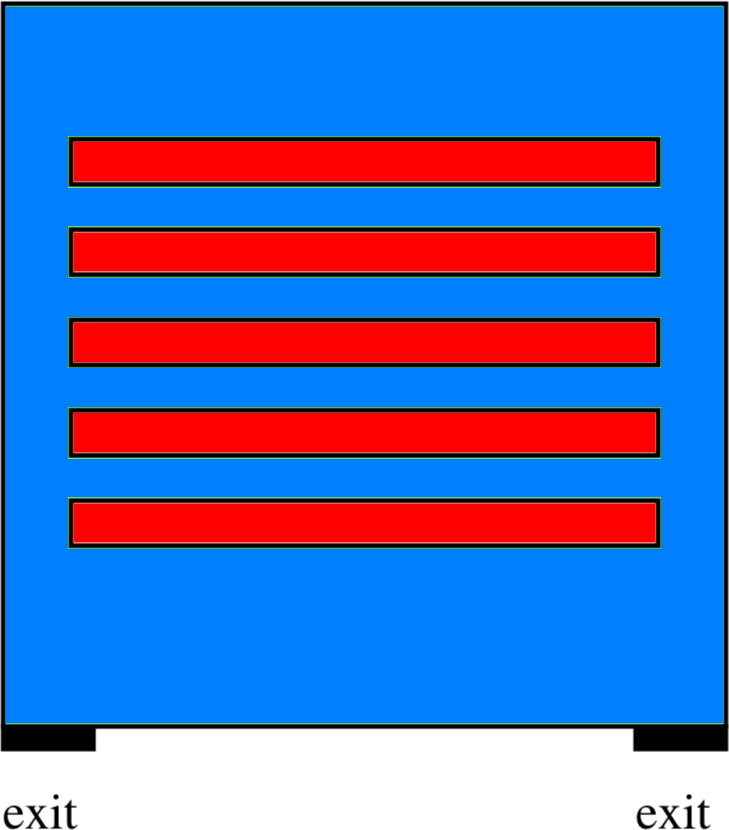}
	\end{subfigure}
	\begin{subfigure}[t]{0.4\linewidth}
		\centering
	  	\includegraphics[height=5cm]{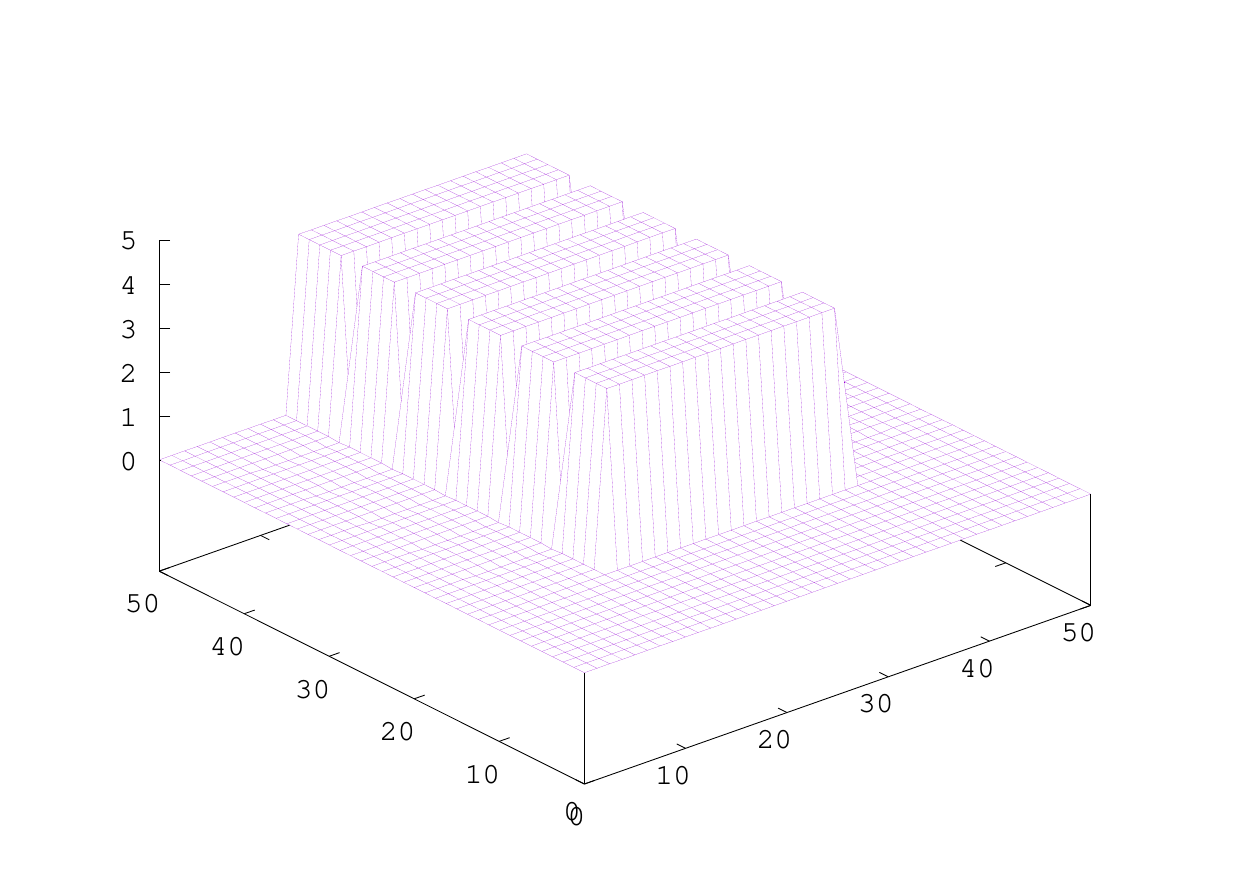}
	\end{subfigure}
	\caption{\label{fig:MFTC-geom-m0} Left: the geometry (obstacles are in red). Right: the density at $t=0$.}
\end{figure}

We provide an example, borrowed from~\cite{MR3392611}, which we solve using the finite difference scheme and Newton method. The model represents a crowd of pedestrians who want to leave a room represented by a square hall (whose side is 50 meters long) containing rectangular obstacles. There are two doors. The chosen geometry and the initial distribution are represented on Figure~\ref{fig:MFTC-geom-m0}.

For simplicity we focus here on a model with local interactions; see \textit{e.g.}~\cite{LachapelleWolfram-2011-MFG-congestion-aversion,YAJML,MR3765549,MR3498932} for other models of this type and~\cite{MR3763083,MR3392611} for crowd motion models with non-local interactions. Here, the congestion is taken into account through the cost (leading to so-called soft congestion, see Remark~\ref{rem:congestion} below): the higher the density at the current position, the higher the price to move. In particular, the running cost does not have a separate dependence in $\ctrl$ and $m$ as in Example~\ref{AMS-num-ex:special-case-sep}. We consider the following  Hamiltonian, which depends locally  on $m$ and captures congestion effects:
  \begin{equation}
  \label{eq:YAMLcongestion-def-H}
    H(x, m, p)= \frac {8 |p|^2} {(1+m)^{\frac 3 4}} - \frac 1 {3200} \,.
  \end{equation} 
 We compare the evolution of the density in the non-cooperative and the cooperative situation:
\begin{enumerate}
\item Mean field games: the MFG PDE system~\eqref{AMS-num-eq:PDE-system-MFG} becomes
\begin{subequations}
     \begin{empheq}[left=\empheqlbrace]{alignat=2}
     	\displaystyle
	- \frac {\partial  u}{\partial t}  -   0.05 \;  \Delta  u
     + \frac {8} {(1+m)^{\frac 3 4}} \;   |\nabla u|^2 &=      \frac  1  {3200} \,, 
	\\
	\displaystyle
	\frac {\partial  m}{\partial t}  -   0.05\;  \Delta m   - 16 \, \diver\left(   \frac { m \nabla u } {(1+m)^{\frac 3 4}}  \right) &= 0 \,.
	\label{eq:MFG-KFP-congestion}
     \end{empheq}
\end{subequations}
\item Mean field control:  the KFP equation~\eqref{eq:MFG-KFP-congestion} is the same but the HJB equation becomes
  \begin{equation*}
   -\frac {\partial  u}{\partial t}  -   0.05 \;  \Delta  u   
   + \left( \frac {2} {(1+m)^{\frac 3 4}}    +   \frac {6} {(1+m)^{\frac 7 4}}\right) \;   |\nabla u|^2=      \frac  1  {3200}.
  \end{equation*}
\end{enumerate}
We choose $\nu=0.05$. The time horizon corresponds to $T=50$ minutes and we do not put any terminal cost, meaning $g \equiv 0$. %

The boundary consists of several parts. On the part corresponding to the doors, for $u$, we impose a Dirichlet condition $u=0$, which corresponds to an exit cost.  For $m$, we assume that $m =0$ outside the domain, so we also get the Dirichlet condition  $m=0$ on this part of the boundary. On the part of the boundary corresponding to the solid walls, for $u$,  we impose a homogeneous Neumann boundary condition: $ \frac{\partial u}{\partial n}=0$, which means that the velocity of the pedestrians is tangential to the walls. For $m$, we choose the boundary condition: 
$\nu \frac {\partial  m}{\partial n}+m     \frac {\partial H} {\partial p} (\cdot, m, \nabla  u )\cdot n=0$, therefore $ \frac {\partial  m}{\partial n}=0$
on this part of the boundary.

The initial density $m_0$ is piecewise constant and takes two values $0$ and $4$ people/m$^2$, see Figure \ref{fig:MFTC-geom-m0}.
At $t=0$, there are 3300 people in the hall.%

We use  Newton iterations with the finite difference scheme discussed in~\S~\ref{AMS-num-sec:fin-diff} and originally proposed in \cite{MR2679575}.

Figure~\ref{fig:MFTC-MFG-congestion-evol-m} displays the density $m$ for the two models, at $t=1$, $2$, $5$ and $15$ minutes. 
In both cases, the pedestrians move between the obstacles towards the narrow corridors leading to the exits, at the left and right sides of the hall. The density thus reaches high values at the intersections of corridors. Then, due to congestion effects, the velocity is lower in the regions where the density is higher. As shown on Figure~\ref{fig:MFTC-MFG-congestion-evol-m}, MFC leads to
 lower values of density. This is consistent with Figure~\ref{fig:MFTC-MFG-congestion-evol-totalmass} (left), which shows that MFC leads to a faster exit of the room. Furthermore, the price of anarchy can be estimated numerically. As shown on Figure~\ref{fig:MFTC-MFG-congestion-evol-totalmass} (right), it is found to increase with respect to the sensitivity to the crowd, captured by the exponent $\frac{3}{4}$ in~\eqref{eq:YAMLcongestion-def-H}. This coefficient can be replaced by a value $\beta\ge0$, which parameterizes the congestion effects (\textit{i.e.}, the strength of the mean-field interactions). We conclude that the non-cooperative scenario (MFG) is less and less efficient (relatively to the MFC) when the congestion effects increase. In the extreme case where $\beta=0$, there is no more mean-field interactions, which explains why the PoA is $1$.

 \begin{remark}
 \label{rem:congestion}
 	The form of soft congestion through the cost function considered here is in contrast with a form of hard congestion taken into account through a density constraint, see \textit{e.g.}~\cite{MR3420414,MR3556062}. In the latter case, the density cannot exceed a given threshold, which leads to flat regions in the density surface at points where the constraint is binding, see \textit{e.g.}~\cite{MR3772008} for numerical examples. Moreover, a congestion cost is different from a crowd aversion cost represented by an increasing function of $m$ which is in general independent of the control. The latter cost is a penalty term which discourages the agents from being in a crowded region independently of whether they are moving or not. This type of cost typically makes the population spread although not exactly in the same way as a diffusion, see \S~\ref{sec:auglag-crowd} for more details.
 \end{remark}

 \begin{center}
   \includegraphics[width=0.45\linewidth]{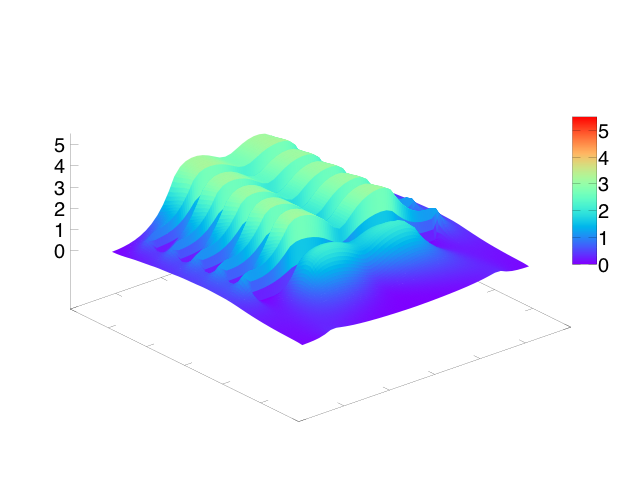} 
   \quad \includegraphics[width=0.45\linewidth]{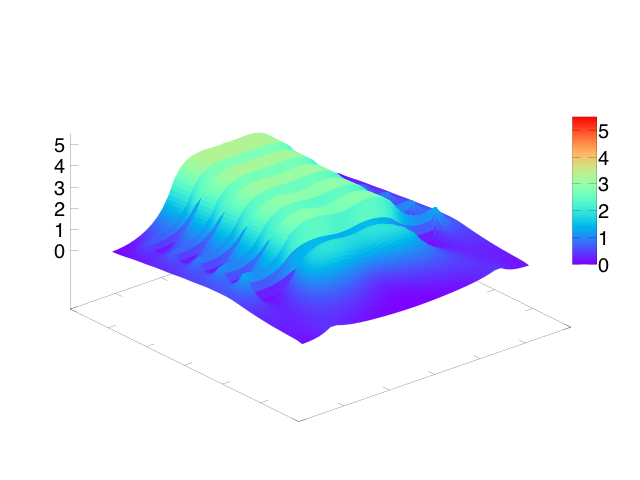}
   \\    \vskip -0.65cm
    \includegraphics[width=0.45\linewidth]{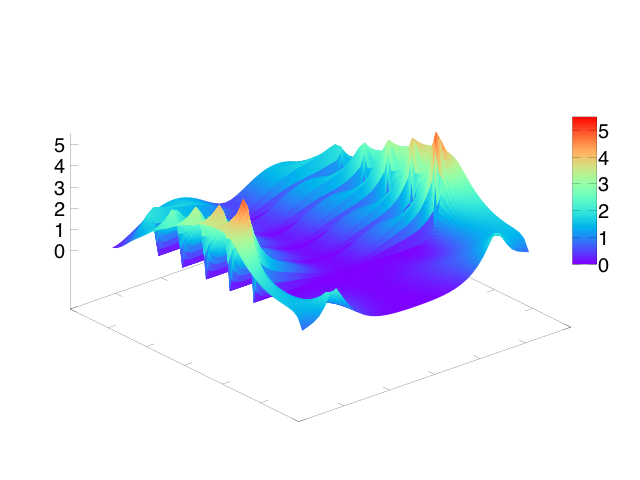} 
  \quad  
   \includegraphics[width=0.45\linewidth]{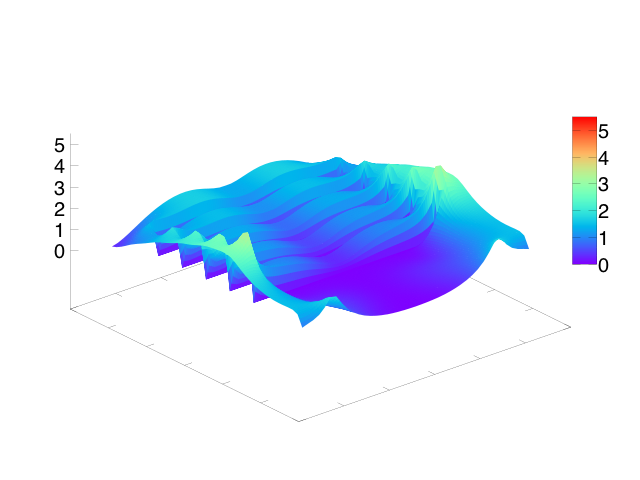}
   \\    \vskip -0.65cm
   \includegraphics[width=0.45\linewidth]{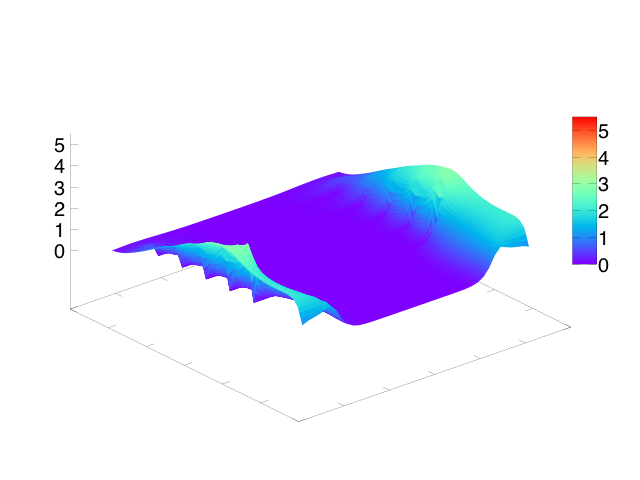}
  \quad  
  \includegraphics[width=0.45\linewidth]{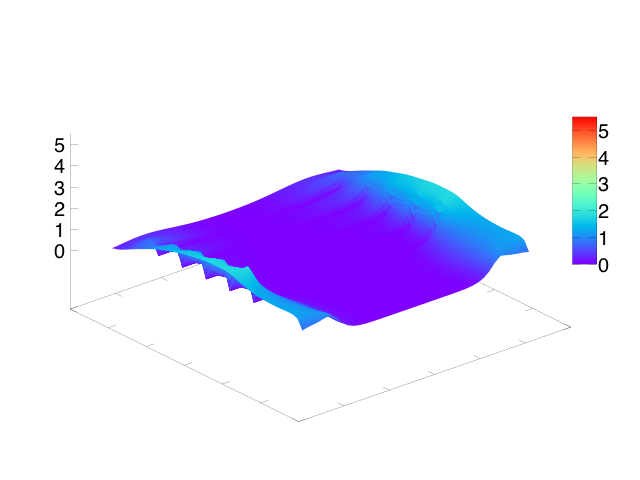}
     \vskip -0.75cm
  \captionof{figure}{ \label{fig:MFTC-MFG-congestion-evol-m} The density computed with the two models at different dates, $t=1,5$ and $15$ minutes (from top to bottom). Left: Mean field game. Right: Mean field type control.}
 \end{center}

\begin{center} %
		\includegraphics[width=0.45\linewidth]{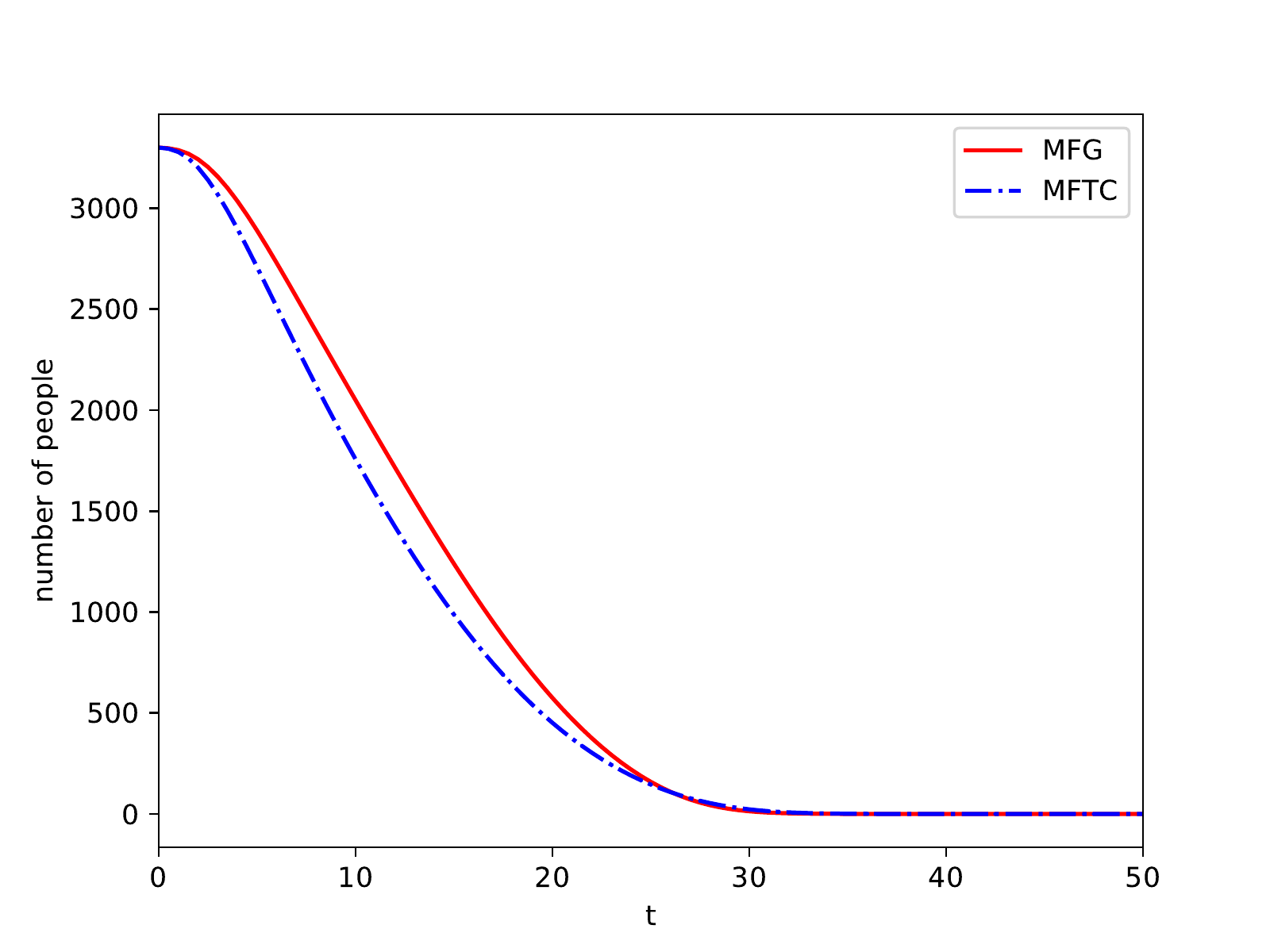}
		\includegraphics[width=0.45\linewidth]{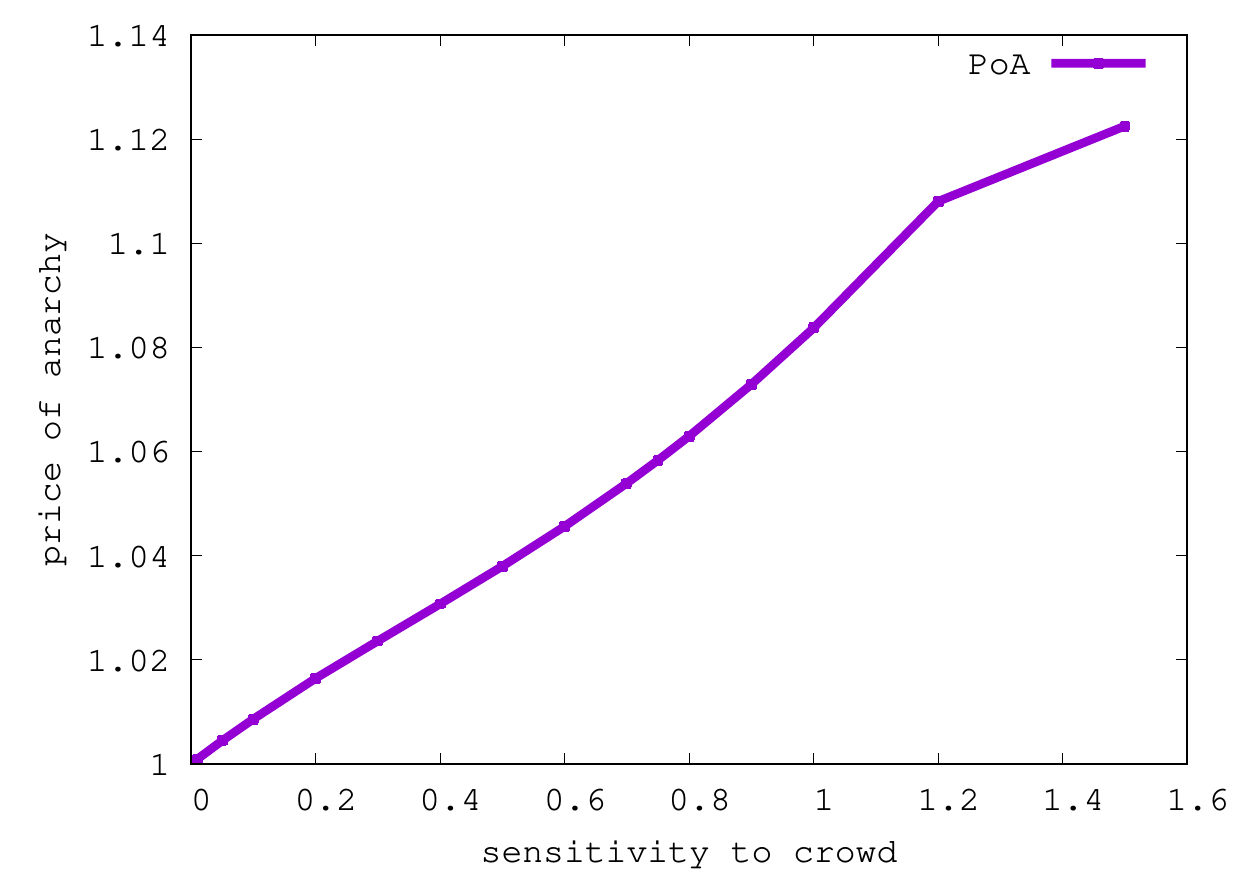}
	\captionof{figure}{\label{fig:MFTC-MFG-congestion-evol-totalmass} Left: Evolution of the total remaining number of people in the room for the mean field game (red line) and the mean field type control (dashed blue line). Right: Price of Anarchy as a function of the parameter $\beta$.}
\end{center}

\clearpage

\clearpage

\section{\bf Optimization methods for MFC and variational MFGs}
\label{AMS-num-sec:optim-variational}

We now turn our attention to mean-field problems for which the PDE system can be interpreted as the optimality condition of a minimization problem for an energy functional subject to a constraint given by a continuity equation. While mean-field control problems can always be formulated in this way, this is not true for all mean-field games. Mean-field games with this property are called variational MFGs.  From a theoretical and a numerical viewpoint, the aforementioned formulation as a minimization problem is important because it allows us to solve the problem with optimization techniques, which is not the case in general for fixed-point problems. In this section, after describing the main ideas of variational MFGs, we focus on the numerical aspects and describe two methods.

\subsection{Variational formulation of MFC and MFG}

Let us start with the MFC problem~\eqref{AMS-num-eq:def-J-MFC}--\eqref{AMS-num-eq:dyn-X-general-MFC}. The expectation as an integral against $m^{\ctrl}(t, \cdot)$, the probability density of the law of $X_t^{\ctrl}$, which leads to:
\begin{align*}
	J^{MFC}(\ctrl) = \int_0^T \int_{\dom} f(x, m^{\ctrl}(t,\cdot), \ctrl(t,x) ) m^{\ctrl}(t,x) dx dt + \int_{\dom} g(x, m^{\ctrl}(T,\cdot)) m^{\ctrl}(T,x) dx. 
\end{align*}
The cost is thus formulated in a deterministic way and is to be minimized under the constraint given by the KFP equation for $m^{\ctrl}$:
\begin{equation}
\label{AMS-num-eq:KFP-ctrl-Q}
	\frac{\partial m^{\ctrl}}{\partial t}(t,x)  - \nu \Delta m^{\ctrl}(t,x) + \diver\left( m^{\ctrl}(t,\cdot) b(\cdot, m^{\ctrl}(t), \ctrl(t,\cdot)) \right)(x) = 0, 
	\hbox{ in } (0,T] \times \domT ,
\end{equation}
with the initial condition:
\begin{equation}
\label{AMS-num-eq:KFP-ctrl-Q-init}
	m^{\ctrl}(0,x) = m_0(x), 
	\hbox{ in } \dom.
\end{equation}

For a MFG, it is not always possible to characterize the equilibrium flow of densities as the minimizer of a functional. A special class of games for which it is possible is the class of so-called potential mean-field games, namely, games in which the costs $f$ and $g$ derive from a potential. For instance, let us consider the setting where $b$ and $f$ are as in Example~\ref{AMS-num-ex:special-case-sep-quadra}, and suppose in addition that there exists $F_0$ and $G$ such that
$$
	f_0(x,m) = \frac{\delta F_0}{\delta m}(x,m), 
	\qquad g (x,m) = \frac{\delta G}{\delta m}(x,m).
$$
Here the derivative should be understood in the sense of measures and we implicitly identify the probability density $m$ with the corresponding probability measure. If we are concerned with square-integrable densities, we can assume that there exist $\bF_0, \bG: \cP_2(\dom) \to \RR$, where $\cP_2(\dom)$ denotes the set of probability measures on $\dom$ with a second moment, such that:
$$
	\bF_0(m) - \bF_0(m') = \int_0^1 \int_{\dom} f_0(x, (1-\theta) m + \theta m') (m-m')(dx)  d\theta, \qquad \forall m,m' \in \cP_2(\dom),
$$ 
and likewise for $\bG$ and $g$. We can then consider the cost functional
\begin{align*}
	\bJ(\ctrl) = \int_0^T \left[ \int_{\dom} \frac{1}{2}|\ctrl(t,x)|^2 m^{\ctrl}(t,x) dx + \bF_0(m^{\ctrl}(t,\cdot)) \right]  dt + \bG(m^{\ctrl}(T,\cdot)). 
\end{align*}
where, as before, $m^{\ctrl}$ solves~\eqref{AMS-num-eq:KFP-ctrl-Q}--\eqref{AMS-num-eq:KFP-ctrl-Q-init} with our choice of $b$, namely,
\begin{equation*}
	\frac{\partial m^{\ctrl}}{\partial t}(t,x)  - \nu \Delta m^{\ctrl}(t,x) + \diver\left( m^{\ctrl}(t,\cdot) \ctrl(t,\cdot) \right)(x) = 0, 
	\hbox{ in } (0,T] \times \domT ,
\end{equation*}
with $m^{\ctrl}(0,x) = m_0(x)$ in $\dom$. Under suitable conditions, this problem admits a minimizer $\ctrl^*$ and a corresponding flow of densities $m^{\ctrl^*}$, from which it is possible to recover the MFG Nash equilibrium $(\hat m, \hat u)$ satisfying the MFG PDE system~\eqref{AMS-num-eq:PDE-system-MFG}. 

In the above situation, the variational problem directly stems from the fact that the cost can be interpreted as a potential. More generally, MFG for which the equilibrium can be characterized via a critical point of a variational problem are called variational MFG. In this case, the PDE system can be viewed either as the equilibrium condition of the MFG or as the critical point condition of the variational problem. In the latter problem, the energy functional to be minimized does not necessarily correspond to the original cost functional of an infinitesimal player. We refer the interested reader to \textit{e.g.}~\cite{Cardaliaguet-2013-notes,MR3644590} for more details on variational and potential MFGs.

\subsection{A PDE driven optimal control problem}

\label{AMS-num-subsec:PDE-optc}

To fix ideas, let us consider the following problem: Minimize the function
\begin{equation}
\label{AMS-num-eq:variational-primal-B}
	\cB(m,w) = \int_0^T \int_{\TT^d} \left(   \mathfrak{L}(x, m(t,x), w(t,x)) +   \mathfrak{F}(x, m(t,x)) \right) dx dt + \int_{\TT^d}   \mathfrak{G}(x, m(T,x)) dx
\end{equation}
on pairs $(m,w)$
such that $m\ge 0$ and 
\begin{equation}
\label{AMS-num-eq:variational-form-linearcons}
\begin{cases}
	\displaystyle
	\,\, \frac{\partial m}{\partial t}  - \nu \Delta m + \diver w = 0, 
	&\hbox{ in } (0,T] \times \TT^d ,
	\\
	\,\, m|_{t=0} = m_0,  &\hbox{ in }  \TT^d ,
\end{cases}
\end{equation}
where we use the notation:
$$
	  \mathfrak{F}(x,m) = 
	\begin{cases}
		\int_0^{m} \tilde f(x, s) ds, & \hbox{ if } m \ge 0,
		\\
		+\infty, 	& \hbox{otherwise,}
	\end{cases}
	\qquad
	  \mathfrak{G}(x,m) = 
	\begin{cases}
		\int_0^m \, \tilde g(x, s) ds, & \hbox{ if } m \ge 0,
		\\
		+\infty, 	& \hbox{otherwise,}
	\end{cases}
$$ 
and
$$
	 \mathfrak{L}(x, m, w) = 
	\begin{cases}
		m \tilde \ell \left(x, \frac{w}{m}\right), & \hbox{ if } m > 0,
		\\
		0, 		& \hbox{ if } m = 0 \hbox{ and } w = 0,
		\\
		+\infty, 	& \hbox{otherwise}.
	\end{cases}
$$

For simplicity, we assume that for every $x$, $\tilde \ell(x,\cdot)$ is a power-like function with exponent greater than $1$ (\textit{i.e.}, it is bounded above and below by $p \mapsto p^{r}$ for some $r>1$, up to multiplicative and additive constants; see \textit{e.g.}~\cite{MR3358627} for detailed assumptions). If $\tilde f(x, \cdot)$ and $\tilde g(x, \cdot)$ are the antiderivatives of $m \mapsto  m f_0(x, m)$ and $m \mapsto  m g(x, m)$ respectively, and if $\tilde \ell  = L_0$, then this problem corresponds to the MFC problem in the setting of Example~\ref{AMS-num-ex:special-case-sep} when the cost functions depend only locally on $m$. Problems of the above type can also stem from the variational formulation of some MFG, as originally explained by Lasry and Lions in~\cite{MR2295621}. We refer the interested read to~\cite{MR3399179,MR3358627,MR3498932} and the references therein for a rigorous development of these ideas.

Here $w$ plays the role of the product $m \ctrl$. Using this new variable, the continuity equation becomes~\eqref{AMS-num-eq:variational-form-linearcons}, which is linear in $(m,w)$. Under suitable conditions, $\cB$ is convex and hence this is a convex minimization problem under a linear constraint.  For instance, we typically assume that $\tilde f$ and $\tilde g$ are nondecreasing, which implies that $\mathfrak{F}$ and $\mathfrak{G}$ are convex, and that $\tilde \ell(x,\cdot)$ is a convex power-like function, which implies that $\mathfrak{L}$ is convex with respect to $(m,w)$.

In this case, the problem admits the following dual formulation: Maximize over $u$ such that $u(T,x) = g(x)$, the function
 \begin{equation}
 \label{AMS-num-eq:variational-dual-A}
 	\cA(u) = \inf_{m} \cA(u, m)
 \end{equation}
 with:
\begin{align*}
	\cA(u,m)
	&=
	\int_0^T \int_{\TT}
	m(t,x) \Big(\partial_t u (t,x) + \nu \Delta u(t,x) - \mathfrak{H}(x, m(t,x), \nabla u(t,x))  \Big)
	dx dt
	\\
	&\qquad+
	\int_{\TT} m_0(x) u(0,x) dx,
\end{align*}
where
$$
	\mathfrak{H}(x,m,p) = \sup_{\ctrl} \{ - \mathfrak{L}(x,m,\ctrl) - \langle \ctrl, p \rangle \}.
$$
It can be shown that $\cA$ and $\cB$ are indeed in duality, \textit{i.e.} (with suitable spaces for $u, m$ and $w$)
$$
	\textbf{(A)} = \sup_{u} \cA(u) = \inf_{(m,w)} \cB(m,w) = \textbf{(B)}.
$$

This property relies on Fenchel-Rockafellar duality theorem~\cite{MR1451876} and the following observation:
$$
	\textbf{(A)}
	=
	- \inf_{u} \Big\{ \cF(u) + \cG(\Lambda(u)) \Big\},
	\qquad
	\textbf{(B)} = \min_{(m,w)} \Big\{ \cF^*(\Lambda^*(m,w)) + \cG^*(-m,-w)\Big\}
$$
where $ \Lambda(u)= \left(\frac{\partial u} {\partial t}   +   \nu \Delta u, \grad u \right)$, 
  \begin{equation*}
    \cF(u)= \chi_T(u)- \int_{\TT^d} m_0(x)u(0,x)dx, 
    \qquad
    \chi_T(u)=
    \begin{cases}
    0, & \hbox{ if } u|_{t=T} = g \\
    +\infty, & \hbox{ otherwise,}
    \end{cases}
  \end{equation*}
\begin{equation*}
\cG(\varphi_1, \varphi_2)= -\inf_{ m\ge 0    } \int_0^T \int_{\TT^d}    m(t,x) \left( \varphi_1(t,x) -  \mathfrak{H} (x, m(t,x), \varphi_2(t,x) )  \right)    dx dt   .
\end{equation*}
$\cF^*, \cG^*$ are the convex conjugates of $\cF,\cG$, and $\Lambda^*$ is the adjoint operator of $\Lambda$.

\subsection{Discrete version of the PDE driven optimal control problem}
\label{AMS-num-subsec:PDE-optc-discrete}

In order to implement optimization methods, we need to discretize the problem introduced above in~\S~\ref{AMS-num-subsec:PDE-optc}.  To alleviate notation, we consider the one dimensional case, \textit{i.e.}, $d=1$. We have in mind a MFC problem in the setting of Example~\ref{AMS-num-ex:special-case-sep}, in which case $\tilde \ell(x,p)$ corresponds to $H_0^*(x, -p)$. For simplicity, we assume that for every $x$, $\tilde \ell(x,\cdot)$ is a power-like function with exponent greater than $1$. 

 We introduce the following spaces, respectively for the discrete counterpart of $m,w,$ and $u$:
$$
	\cM = \RR^{(N_T+1) \times N_h}, 
	\qquad 
	\cW = (\RR^2)^{N_T \times N_h}, 
	\qquad
	\cU = \RR^{N_T \times N_h}.
$$
Note that at each space-time grid point, we will have two variables for the value of $w$, which will be useful to define an upwind scheme.

We denote by $\tilde H_0$ a discrete Hamiltonian satisfying  the properties~\ref{AMS-num-hyp-discH-monotonicity}--\ref{AMS-num-hyp-discH-convexity}. Let $\tilde H^*_0: \TT^d \times (\RR^2)^d \to \RR \cup \{+\infty\}$ be its convex conjugate w.r.t. the $p$ variable:
\begin{equation}
\label{eq:var:1}  
	\tilde H^*_0: 
	(x,\gamma) \mapsto 
	\tilde H^*_0(x, \gamma)
	=
	\max_{p \in \RR^2} \left\{ \langle \gamma , p \rangle - \tilde H_0(x, p) \right\}.
\end{equation}

\begin{example}\label{sec:discr-vers-vari-2}
	In the setting of Example~\ref{AMS-num-eq-quadratic-discrete}, $\tilde H_0(x, p_1, p_2) = \tfrac{1}{2} |P_K(p_1,p_2)|^2$ with $K = \RR_- \times \RR_+$, so
	$$
		\tilde H^*_0(x, \gamma_1, \gamma_2)
		=
		\begin{cases}
			\displaystyle
			\,\,  \tfrac{1}{2}(|\gamma_1|^2 + |\gamma_2|^2), &\hbox{ if } (\gamma_1,\gamma_2) \in K,
			\\
			\,\, +\infty, &\hbox{ otherwise.}
		\end{cases}
	$$
\end{example}
A discrete counterpart $\tilde\cB: \cM \times \cW \to \RR$ to the functional $\cB$ introduced in~\eqref{AMS-num-eq:variational-primal-B} can be defined as
\begin{equation}
\label{eq:var:2}  
	\tilde\cB:
	(M,W) \mapsto
	\sum_{n=1}^{N_T} \sum_{i=0}^{N_h-1} \left[ \tilde{\mathfrak{L}}(x_i, M^n_i, W^{n-1}_i) + \mathfrak{F}(x_i, M^n_i) \right] + \frac{1}{\Delta t}  \sum_{i=0}^{N_h-1}  \mathfrak{G}(x_i, M^{N_T}_i),
      \end{equation}
      where $\tilde{\mathfrak{L}}: \TT \times \RR \times \RR^2 \to \RR \cup \{+\infty\}$ is a discrete version of $\mathfrak{L}$ defined as:
\begin{equation}
\label{eq:var:3}  
	\tilde{\mathfrak{L}}:
	 (x,m,w) \mapsto \tilde{\mathfrak{L}}(x,m,w)
	 =
	\begin{cases}
		\displaystyle
		\,\, m \tilde H^*_0\left(x,-\frac{w}{m}\right), &\hbox{ if } m > 0 \hbox{ and } w \in K,
		\\
		\,\, 0, &\hbox{ if } m=0, \hbox{ and } w= 0,
		\\
		\,\, +\infty, &\hbox{ otherwise.}
	\end{cases}
      \end{equation}
      Furthermore, a discrete version of the linear constraint~\eqref{AMS-num-eq:variational-form-linearcons} can be written as
      \begin{equation}
        \label{eq:var:5}  
	\Sigma(M,W) = (0_\cU, \bar M^0)
      \end{equation}
      where $0_\cU \in \cU$ is the vector with $0$ on all coordinates, $\bar M^0 = (\bar m_0(x_0), \dots, \bar m_0(x_{N_h-1})) \in \RR^{N_h}$, see~\eqref{AMS-num-eq:discrete-barm-0} for the definition of $\bar m_0$, and 
      \begin{equation}
\label{eq:var:6}  
	\Sigma: \cM \times \cW \to \cU \times \RR^{N_h},
	\qquad
	(M,W) \mapsto \Sigma(M,W) = (\Lambda(M,W), M^0),
      \end{equation}
      with 
      $$
      	\Lambda(M,W) = AM + BW,
	$$
	with $A$ and $B$ being discrete versions of respectively the heat operator and the divergence operator, defined as follows:
\begin{equation}
\label{eq:def-op-A}
	A: \cM \to \cU,
	\qquad
	(AM)^{n}_i = \frac{M^{n+1}_i - M^{n}_i}{\Delta t} - \nu (\Delta_h M^{n+1})_i, \qquad 0 \le n < N_T, 0 \le i < N_h,
\end{equation}
and
\begin{equation}
\label{eq:def-op-B}
	B: \cW \to \cU,
	\qquad
	(BW)^{n}_i = \frac{W^{n}_{i+1,2} - W^{n}_{i,2}}{h} + \frac{W^{n}_{i,1} - W^{n}_{i-1,1}}{h} 
              \qquad 0 \le n \le N_T, 0 \le i < N_h.
\end{equation}
The discrete counterpart of the (primal) variational problem~\eqref{AMS-num-eq:variational-primal-B}--\eqref{AMS-num-eq:variational-form-linearcons} is therefore
\begin{equation}
\label{eq:discrete-pb-short}
	\inf_{\substack{(M,W) \in \cM \times \cW}} \tilde\cB(M,W),  \quad \hbox {subject to } \Sigma(M,W) = (0_\cU, \bar M^0).
\end{equation}
It can be shown that this problem admits a unique solution, and this solution satisfies $M^n_i >0$ for all $i \in \{0, \dots, N_h\}$ and all $n>0$; see \textit{e.g.}~\cite[Theorem 3.1]{BricenoAriasetalCEMRACS2017} for more details on a special case (see also~\cite[Section 6]{MR3135339} and~\cite[Theorem 2.1]{MR3772008}  for similar results respectively in the context of the planning problem and in the context of an ergodic MFG).

Analogously to the continuous setting, the discrete problem admits a dual formulation. Furthermore, a key feature of the discretization we chose is that the optimality condition of the discrete problem coincides with the finite-difference scheme presented in Section~\ref{AMS-num-sec:fin-diff}, adapted to the current setting. We refer to~\cite[Section 3]{achdoulauriere-2020-mfg-numerical} and the references therein for more details.

\subsection{ADMM and Primal-dual method}

In order to draw a connection with numerical methods for optimization problems, let us reformulate the problem from a more abstract perspective. We consider the following primal problem:
\begin{equation}
\label{AMS-num-eq:generic-primal-constraint}
	\inf_{\xi \in \RR^{\mathrm{d}}} \varphi(\xi), \hbox{ subject to } \Xi \xi = 0, 
\end{equation}
where $\varphi: \RR^{\mathrm{d}} \to (-\infty, +\infty]$ is a lower semi-continuous convex proper function, $\Xi: \RR^{\mathrm{d}} \to \RR^{\mathrm{d}'}$ is a linear operator,  and $0 \in \RR^{\mathrm{d}'}$. In the setting of the variational mean-field problem described above, problem~\eqref{AMS-num-eq:generic-primal-constraint} can cover~\eqref{eq:discrete-pb-short} with:
$$
	\mathrm{d} = 2(N_T+1) N_T (N_h)^2, \quad \xi = (M,W), \quad \varphi = \tilde\cB, \quad \Xi = \Sigma.
$$

Equivalently, the above primal problem can be written as follows:
\begin{equation}
\label{AMS-num-eq:generic-primal-sum}
	\inf_{\xi \in \RR^{\mathrm{d}}} \varphi(\xi) + \psi(\widetilde \Xi \xi), 
\end{equation}
where, for $y \in \RR^{\mathrm{d}'}$, 
$$
	\psi(y) = \iota_{\{0\}}(y) 
	= \begin{cases}
	0, &\hbox{ if } y = 0,
	\\
	+\infty, &\hbox{ otherwise},
	\end{cases}
	\quad \hbox{ and } 
	\widetilde \Xi = \Xi,
$$
or
$$
	\psi(y) = \iota_{\{0\}}(\Xi y) 
	= \begin{cases}
	0, &\hbox{ if }  \Xi y = 0,
	\\
	+\infty, &\hbox{ otherwise},
	\end{cases}
	\quad \hbox{ and } 
	\widetilde \Xi = \mathrm{id}.
$$

Problems of the form~\eqref{AMS-num-eq:generic-primal-sum} have been extensively studied and various numerical methods have been introduced. Here, we will focus on two of them: the Alternating Direction Method of Multipliers for the Augmented Lagrangian of the dual problem, and a primal-dual method proposed by Chambolle and Pock~\cite{MR2782122}. Some of the main advantages of these methods are that proofs of convergence are readily available, and they can be applied to first order problems (\textit{i.e.}, when there is no diffusion and $\nu=0$).

\textbf{Augmented Lagrangian and ADMM.} We first note that~\eqref{AMS-num-eq:generic-primal-sum} admits as a dual formulation the following problem:
\begin{equation}
\label{AMS-num-eq:generic-dual-sum}
	\inf_{\zeta \in \RR^{\mathrm{d}'}} \varphi^*(- \widetilde{\Xi}^* \zeta) + \psi^*(\zeta), 
\end{equation}
where $\varphi^*$ and $\psi^*$ are the convex conjugates of $\varphi$ and $\psi$ respectively and $\widetilde{\Xi}^*$ is the adjoint operator of $\widetilde\Xi$, \textit{i.e.}, $\langle \widetilde{\Xi} \xi, \zeta \rangle =  \langle \xi, \widetilde{\Xi}^* \zeta \rangle$. In the setting of \S~\ref{AMS-num-subsec:PDE-optc} and \S~\ref{AMS-num-subsec:PDE-optc-discrete}, $\zeta$ plays the role of a vector $U\in\cU$ approximating the function $u$ appearing in the dual problem~\eqref{AMS-num-eq:variational-dual-A}.

For numerical purposes, it is useful to exploit the additive structure of the objective function and introduce an extra variable  to rewrite the problem as:
\begin{equation*}
	\inf_{\zeta, \nu \in \RR^{\mathrm{d}'}} \varphi^*(\nu) + \psi^*(\zeta), \hbox{ subject to } \nu = - \widetilde{\Xi}^* \zeta.
\end{equation*}
The Lagrangian associated to this constrained optimization problem is defined, introducing a Lagrange multiplier $\lambda$, as:
\begin{equation*}
	\cL(\zeta, \nu, \lambda) = \varphi^*(\nu) + \psi^*(\zeta) + \langle \lambda , \nu + \widetilde{\Xi}^* \zeta \rangle .
\end{equation*}
Solving the dual problem~\eqref{AMS-num-eq:generic-dual-sum} thus amounts to finding a saddle point of $\cL$, \textit{i.e.}, $(\zeta^*, \nu^*, \lambda^*)$ achieving
\begin{equation}
\label{AMS-num-eq:generic-LagSaddlePoint}
	\cL(\zeta^*, \nu^*, \lambda^*) = \sup_{\lambda} \inf_{\zeta, \nu} \cL(\zeta, \nu, \lambda).
\end{equation}
This suggests to use a steepest descent method in order to approximate the maximizer of $\lambda \mapsto \inf_{\zeta, \nu} \cL(\zeta, \nu, \lambda)$, where the $\inf$ can be split into two separate optimization sub-problems. For more details on this method, see the algorithm ALG1 in~\cite[Chapter 3]{MR724072}. A variation of this approach consists in adding an extra penalty term to obtain the so-called Augmented Lagrangian defined, for $r>0$, as:
$$
	\cL_r(\zeta, \nu, \lambda) = \varphi^*(\nu) + \psi^*(\zeta) + \langle \lambda , \nu + \widetilde{\Xi}^* \zeta \rangle  + \frac{r}{2} \|\nu + \widetilde{\Xi}^* \zeta \|_2^2.
$$
This new function can be interpreted as the Lagrangian of a modified version of the dual problem~\eqref{AMS-num-eq:generic-dual-sum} in which the penalty term $\frac{r}{2} \|\nu + \widetilde{\Xi}^* \zeta \|_2^2$ is added to the objective. The primal problem corresponding to this penalized dual problem can be interpreted as a regularized version of~\eqref{AMS-num-eq:generic-primal-constraint}.  Note that $\cL$ and $\cL_r$ have the same saddle point, for any $r>0$.

The principle of the Alternating Direction Method of Multipliers applied to this Augmented Lagrangian $\cL_r$ is to iteratively update in turn $\zeta, \nu$ and $\lambda$, as summarized in Algorithm~\ref{AMS-num-admm-algo-generic}. See the algorithm ALG2 in~\cite[Chapter 3]{MR724072} for more details. Here $\prox$ denotes the proximal operator defined, for a lower semicontinuous convex proper function $\varphi: \RR^{\mathrm{d}} \to (-\infty, +\infty]$ as:
$$
	\prox_\varphi(y) = \argmin_{y'} \left\{ \varphi(y') + \frac{1}{2} \|y' - y\|^2\right\}.
$$
It generalizes the notion of orthogonal projection in the sense that if $K \subseteq \RR^{\mathrm{d}}$ is a non-empty, closed, convex set and $\iota_K: \RR^{\mathrm{d}} \to \{0,+\infty\}$ denotes its characteristic function, then $\prox_{\iota_K}$ is the orthogonal projection on $K$. In Algorithm~\ref{AMS-num-admm-algo-generic}, the first two steps (optimization with respect to $\zeta$ and proximal step for $\nu$) are the most costly in terms of computation. 

The convergence of this method can be ensured under suitable conditions, see~\cite[Theorem 8]{MR1168183} for more details.

\begin{algorithm}[H]
\DontPrintSemicolon
  
  \KwInput{Initial guess $(\zeta_0, \nu_0, \lambda_0)$;  number of iterations $\mathtt{K}$}
  \KwOutput{Approximation of a saddle point $(\zeta^*, \nu^*, \lambda^*)$ achieving~\eqref{AMS-num-eq:generic-LagSaddlePoint}}
  Initialize $(\zeta^{(0)}, \nu^{(0)}, \lambda^{(0)}) = (\zeta_0, \nu_0, \lambda_0)$ \;
        \For{$\mathtt{k}=0,1,2,\dots, \mathtt{K}-1$}    
        { 
        Let $\zeta^{(\mathtt{k}+1)} = \argmin_{\zeta}\left\{ \psi^*(\zeta) + \langle \lambda^{(\mathtt{k})} , \widetilde{\Xi}^* \zeta \rangle  + \frac{r}{2} \|\nu^{(\mathtt{k})} + \widetilde{\Xi}^* \zeta \|_2^2 \right\}$\;
        Let 
        \begin{align*}
        \nu^{(\mathtt{k}+1)} 
        &= \argmin_{\nu} \left\{ \varphi^*(\nu)  + \langle \lambda^{(\mathtt{k})} , \nu  \rangle  + \frac{r}{2} \|\nu + \widetilde{\Xi}^* \zeta^{(\mathtt{k}+1)} \|_2^2 \right\} 
        \\
        &= \prox_{\varphi^*/r}\left( - \frac{1}{r}\lambda^{(\mathtt{k})} - \widetilde{\Xi}^* \zeta^{(\mathtt{k}+1)} \right)
	\end{align*} \;
	Let $\lambda^{(\mathtt{k}+1)} = \lambda^{(\mathtt{k})} - r (\nu^{(\mathtt{k}+1)} +\widetilde{\Xi}^* \zeta^{(\mathtt{k}+1)} )$\;
	}
	\Return{$(\zeta^{(\mathtt{K})}, \nu^{(\mathtt{K})}, \lambda^{(\mathtt{K})})$}
\caption{Alternating Direction Method of Multipliers  \label{AMS-num-admm-algo-generic}}
\end{algorithm}

The ADMM was made popular  by Benamou and Brenier  for optimal transport,  see~\cite{MR1738163}, and  first used in the context of MFGs by Benamou and Carlier in~\cite{MR3395203}; see also~\cite{MR3731033} for an application to second order MFG using multilevel preconditioners and~\cite{MR3575615} for an application to mean field type control.

When applied to the dual of the variational problem~\eqref{eq:discrete-pb-short}, in the first step, the first order optimality condition for the minimization yields that $U^{(\mathtt{k}+1)}$ is the solution to a finite difference equation which in the general case $\nu>0$ corresponds to a PDE with a fourth order elliptic operator. Then a preconditioner is needed, see \textit{e.g.} ~\cite{MR2928376,BricenoAriasetalCEMRACS2017}, except if $\nu=0$, in which case a direct solver can be used.  In the second step, the minimization problem can be done separately at each point of the grid, which allows parallelization.

\textbf{Chambolle and Pock's primal-dual algorithm.} We now turn our attention to an algorithm proposed in~\cite{MR2782122} which relies on both the primal and the dual problems~\eqref{AMS-num-eq:generic-primal-sum}--\eqref{AMS-num-eq:generic-dual-sum}. It is based on the Lagrangian:
\begin{equation}
\label{AMS-num-eq:generic-LagSaddlePoint-PD}
	\cL_{PD} (\xi, \zeta) = \varphi(\xi) - \psi^*(\zeta) + \langle \zeta, \widetilde{\Xi} \xi \rangle,
\end{equation}
for which the optimality conditions in each variable are:
\begin{align*}
	\begin{cases}
	\,\, \widetilde{\Xi}^*\hat\zeta &\in -\partial \varphi(\hat\xi) 
	\\
	\,\,  \widetilde{\Xi}\hat\xi &\in \partial \psi^*(\hat \zeta)
	\end{cases}
	&\quad
	\Leftrightarrow
	\quad
	\begin{cases}
	\,\, \hat\xi - \tau \widetilde{\Xi}^*\hat\zeta &\in \hat\xi + \tau \partial \varphi(\hat\xi) 
	\\
	\,\, \hat\zeta + \gamma \widetilde{\Xi}\hat\xi &\in \hat\zeta + \gamma \partial \psi^*(\hat \zeta)
	\end{cases}
	\\
	&\quad
	\Leftrightarrow
	\quad
	\begin{cases}
	\,\, \hat \zeta &\in \argmin_\zeta \left\{ \psi^*(\zeta) + \frac{1}{2 \gamma}\| \zeta - (\hat \zeta + \gamma \widetilde{\Xi} \hat \xi )\|^2 \right\}
	\\
	\,\, \hat \xi &\in \argmin_\xi \left\{ \varphi(\xi) + \frac{1}{2 \tau} \| \xi - (\hat \xi - \tau \widetilde{\Xi}^* \hat \zeta) \|^2 \right\}
	\end{cases}
	\\
	&\quad
	\Leftrightarrow
	\quad
	\begin{cases}
	\,\, \hat \zeta &\in \prox_{\gamma \psi^*} (\hat \zeta + \gamma \widetilde{\Xi} \hat \xi ) 
	\\
	\,\, \hat \xi &\in \prox_{\tau \varphi} (\hat \xi - \tau \widetilde{\Xi}^* \hat \zeta).
	\end{cases}
\end{align*}
The method proposed by Chambolle and Pock is basically to iterate over the two proximal steps, with an additional extrapolation step, as summarized in Algorithm~\ref{AMS-num-chambollepock-algo-generic}. The method has been proved to converged when $\gamma \tau < 1$, see~\cite{MR2782122}.

\begin{algorithm}[H]
\DontPrintSemicolon
  
  \KwInput{Initial guess $(\xi_0, \zeta_0, \tilde \xi_0)$;  number of iterations $\mathtt{K}$}
  \KwOutput{Approximation of a saddle point $(\xi^*, \zeta^*)$ of $\cL_{PD}$ defined by~\eqref{AMS-num-eq:generic-LagSaddlePoint-PD}}
  Initialize $(\xi^{(0)}, \zeta^{(0)}, \tilde \xi^{(0)}) = (\xi_0, \zeta_0, \tilde \xi_0)$ \;
        \For{$\mathtt{k}=0,1,2,\dots, \mathtt{K}-1$}    
        { 
        Let $\zeta^{(\mathtt{k}+1)} = \prox_{\gamma \psi^*} (\zeta^{(\mathtt{k})} + \gamma \widetilde{\Xi} \tilde\xi^{(\mathtt{k})} ) $\;
        Let $\xi^{(\mathtt{k}+1)} = \prox_{\tau \varphi} (\xi^{(\mathtt{k})} - \tau \widetilde{\Xi}^* \zeta^{(\mathtt{k}+1)})$\;
	Let $\tilde\xi^{(\mathtt{k}+1)} = 2 \xi^{(\mathtt{k}+1)} - \xi^{(\mathtt{k})}$\;
	}
	\Return{$(\xi^{(\mathtt{K})}, \zeta^{(\mathtt{K})})$}
\caption{Chambolle-Pock's method \label{AMS-num-chambollepock-algo-generic}}
\end{algorithm}
 
The application of this method to stationary MFGs was first investigated by Brice{\~n}o-Arias, Kalise and Silva in~\cite{MR3772008}; see also~\cite{BricenoAriasetalCEMRACS2017} for an extension to the dynamic setting. When applied to the problem~\eqref{eq:discrete-pb-short}, the first step is similar to the first step in the ADMM method described above and amounts to solving a linear fourth order PDE. The second step is easier thanks to the choice of $\cB$ and the form of its prox.

\subsection{Numerical illustration: crowd motion without diffusion}
\label{sec:auglag-crowd}
We present an example, borrowed from~\cite{MR3575615}, of a mean field model for crowd motion with congestion. The interactions are local and, in contrast with the example of~\S~\ref{subsec:num-evacuation-crowd}, there is no diffusion, \textit{i.e.}, $\nu=0$. Due to the fact that the running cost is not separable in $\ctrl$ and $m$, the MFG does not have a variational structure. We focus here on the MFC problem.

On top of the congestion effects (meaning that moving quickly in a dense region is expensive), we incorporate aversion effects (\textit{i.e.}, being in a crowded region is uncomfortable). The latter aspect is modeled by a function $\ell: \dom \times \RR \to \RR, (x,m) \mapsto \ell(x,m)$ which is increasing with respect to $m$. The domain is  $\dom = [0,1]^2 \backslash [0.4,0.6]^2$, which is a square with obstacle at the center. The Lagrangian (corresponding to the running cost) is:
	$$
		L(x, m,\ctrl) = (b-1)^{b^*}  m^{\frac{a}{b-1}}|\ctrl|^{b^*} + \ell(x, m), 
		\qquad 1<b \leq 2, 0\leq a < 1,
	$$
	where $b^* = b/(b-1)$ is the conjugate exponent of $b$.  
	We define the Hamiltonian by duality (such that on $\partial \Omega$ the vector speed is towards the interior):
	\begin{align*}
		 H(x, m, p) &=
		 \begin{cases}
		 	\displaystyle \inf_{\ctrl \in \RR^2} \big\{\ctrl \cdot p + L(x, m,\ctrl)\big\} = - m^{-a} |p|^{b} + \ell(x, m), &\mbox{ if } x \in \Omega, \\
			\displaystyle \inf_{\ctrl \in \RR^2 \,:\, \ctrl \cdot \mathbf{n} \leq 0} \big\{\ctrl \cdot p + L(x, m,\ctrl)\big\}, &\mbox{ if } x \in \partial \Omega,
		\end{cases}
	\end{align*}
	where $\mathbf{n}$ denotes the outward normal.

\begin{figure}
\centering
\begin{subfigure}{.45\textwidth}
  \centering
  \includegraphics[width=\linewidth]{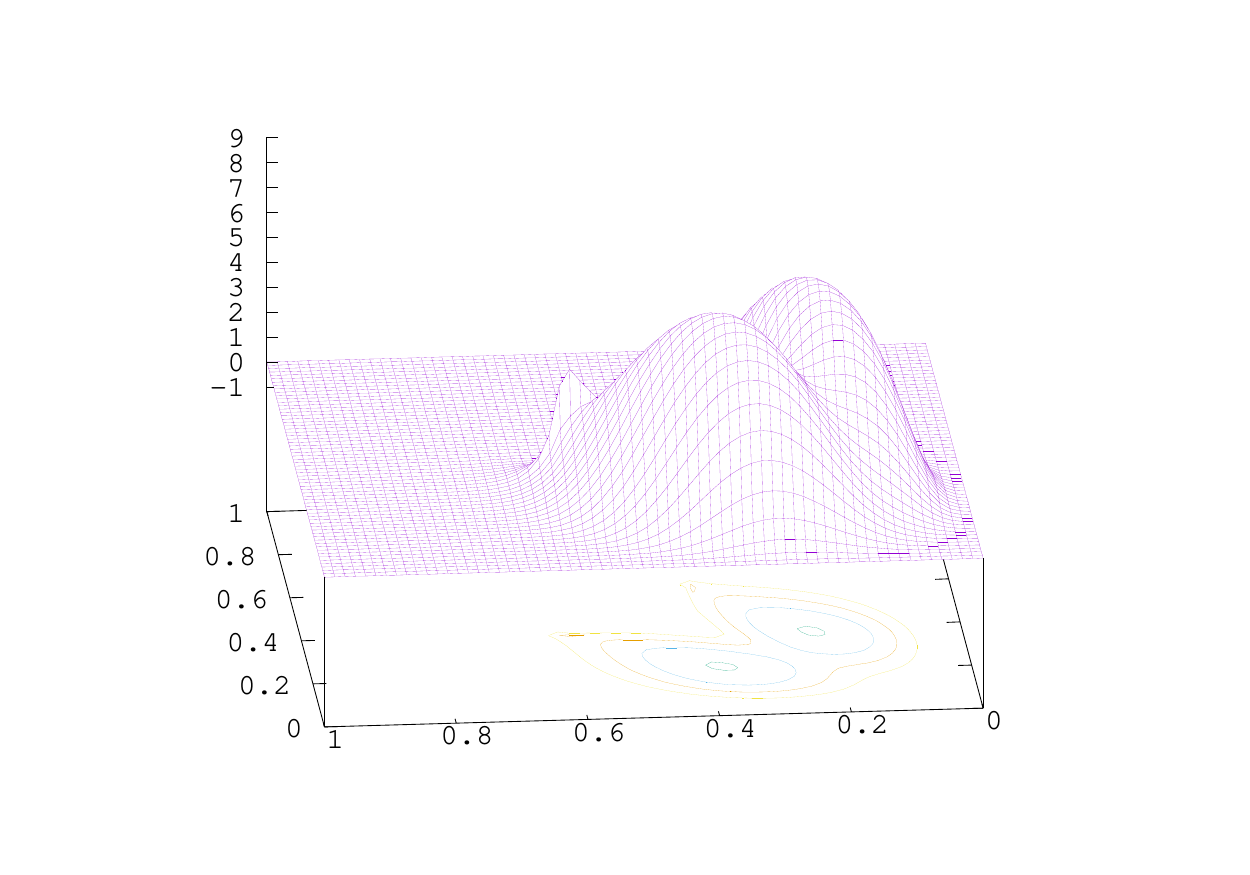}
  \caption*{$t = T/4$}
\end{subfigure}%
\begin{subfigure}{.45\textwidth}
  \centering
  \includegraphics[width=\linewidth]{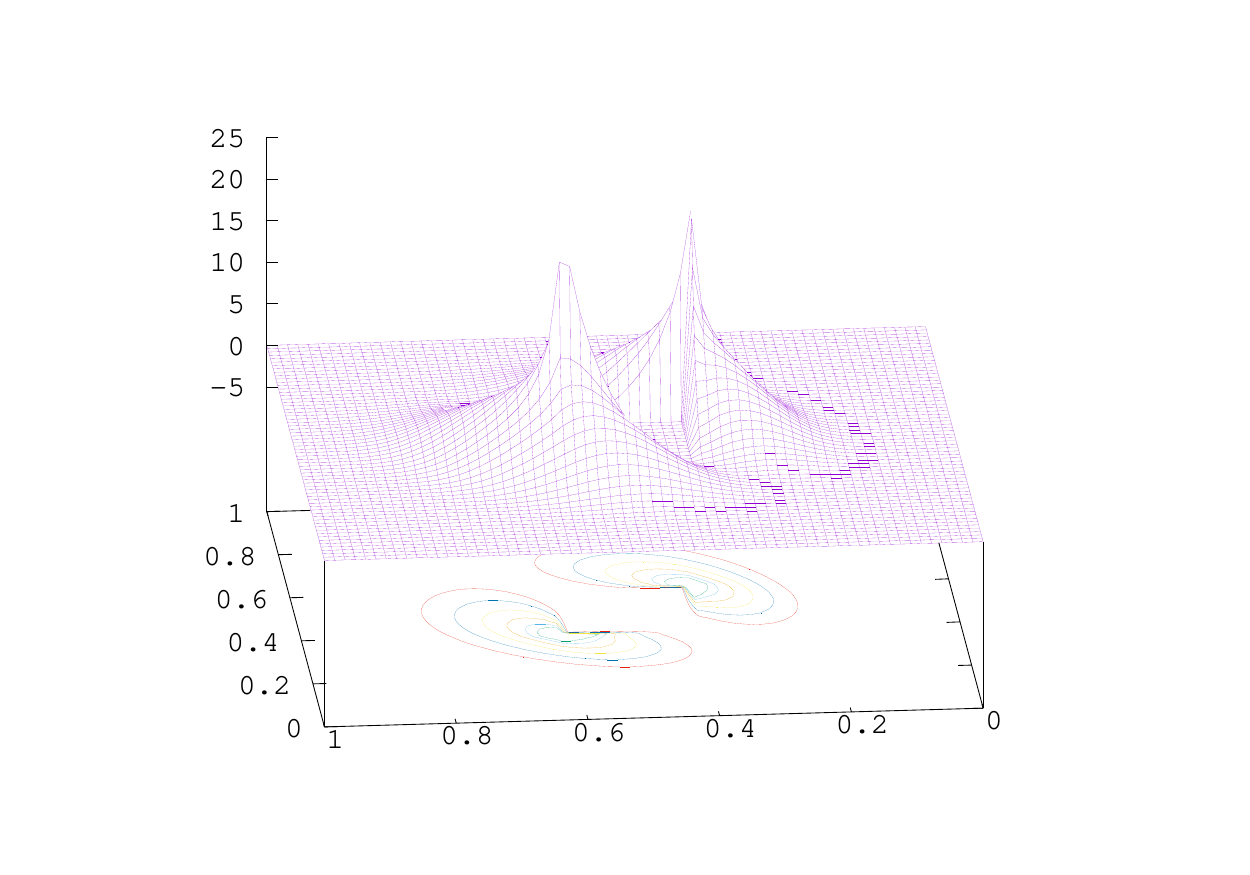}
  \caption*{$t = T/2$}
\end{subfigure}
\\
\begin{subfigure}{.45\textwidth}
  \centering
  \includegraphics[width=\linewidth]{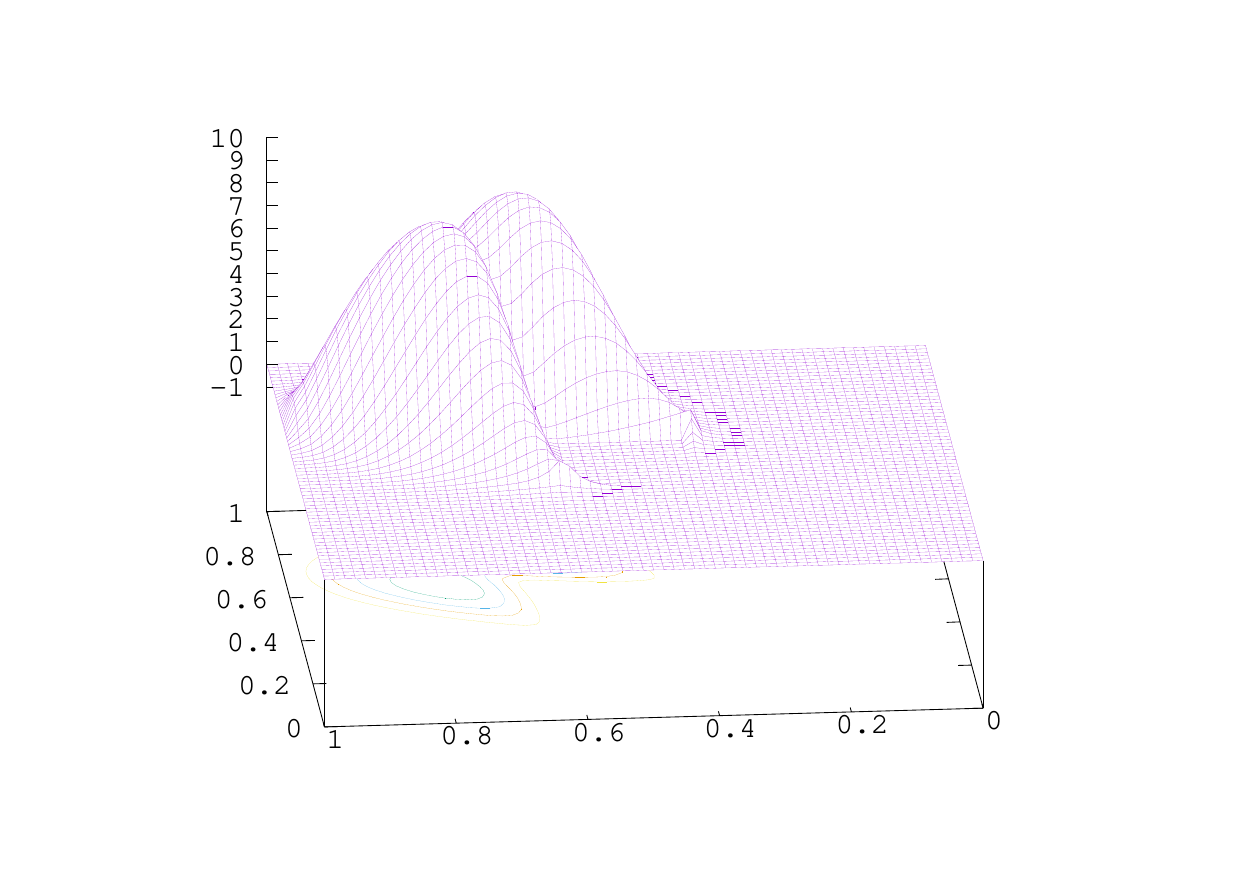}
  \caption*{$t = 3T/4$}
\end{subfigure}%
\begin{subfigure}{.45\textwidth}
  \centering
  \includegraphics[width=\linewidth]{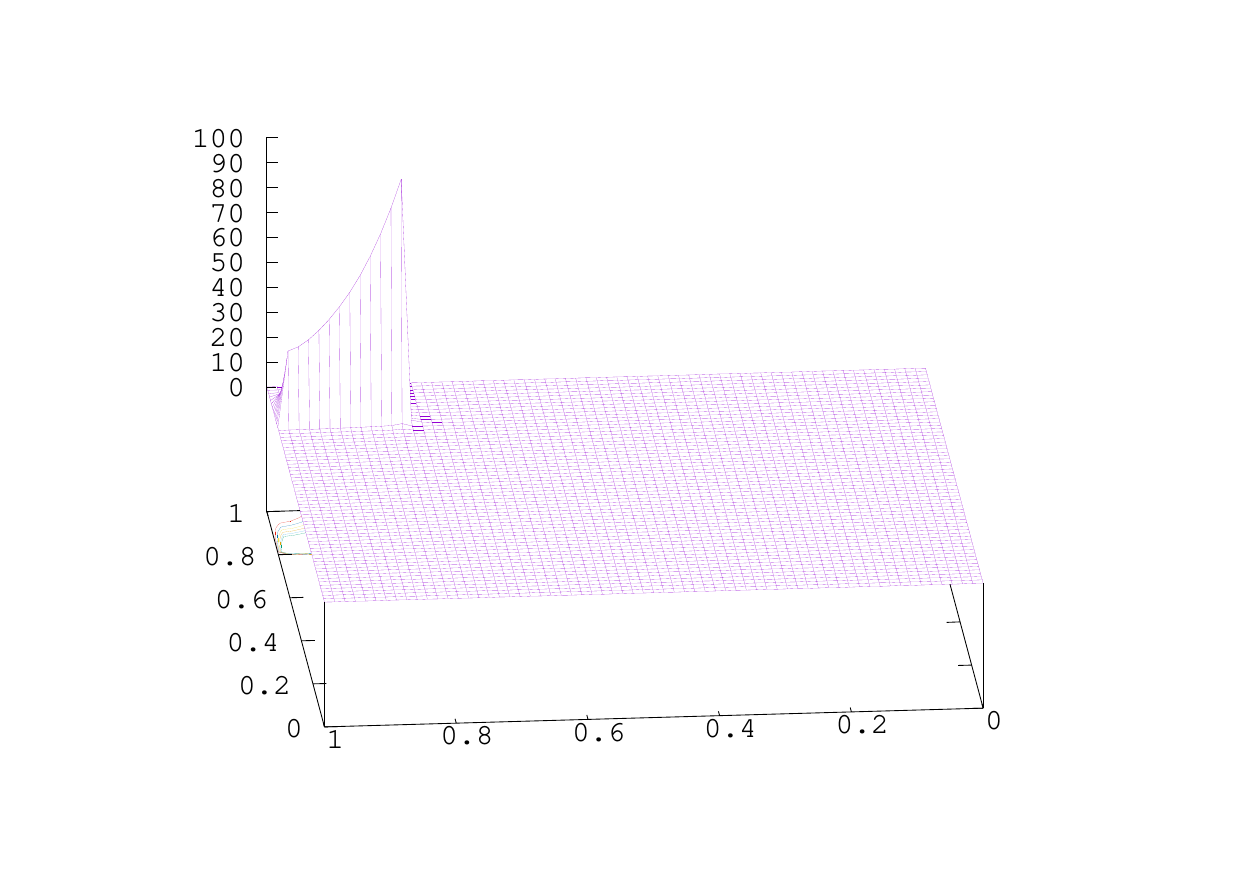}
  \caption*{$t = T$}
\end{subfigure}
\caption{Evolution of the density.}
\label{fig:fig:AugLag-MFTC-congestion-evol-m}
\end{figure}

	As shown in Figure~\ref{fig:AugLag-MFTC-congestion-m0-uT}, we take $m_0$ uniform in a corner and the terminal cost $g$ minimal at the opposite corner. For the numerical experiments reported here, we took $a = 0.01, b = 2, \ell(x, m) = 0.01m$. The evolution of the density is displayed in Figure~\ref{fig:fig:AugLag-MFTC-congestion-evol-m}. We see that it avoids the obstacle and ends up in the opposite corner. Since high velocity is prohibitive (and there is no terminal cost penalizing high density), most of the mass concentrates near the arrival point. Due to the absence of diffusion, at every time step the density remains zero on a large part of the domain. Here, we used the ADMM on the Augmented Lagrangian for the discrete problem as discussed above in this section.

\begin{center} %
		\includegraphics[width=0.45\linewidth]{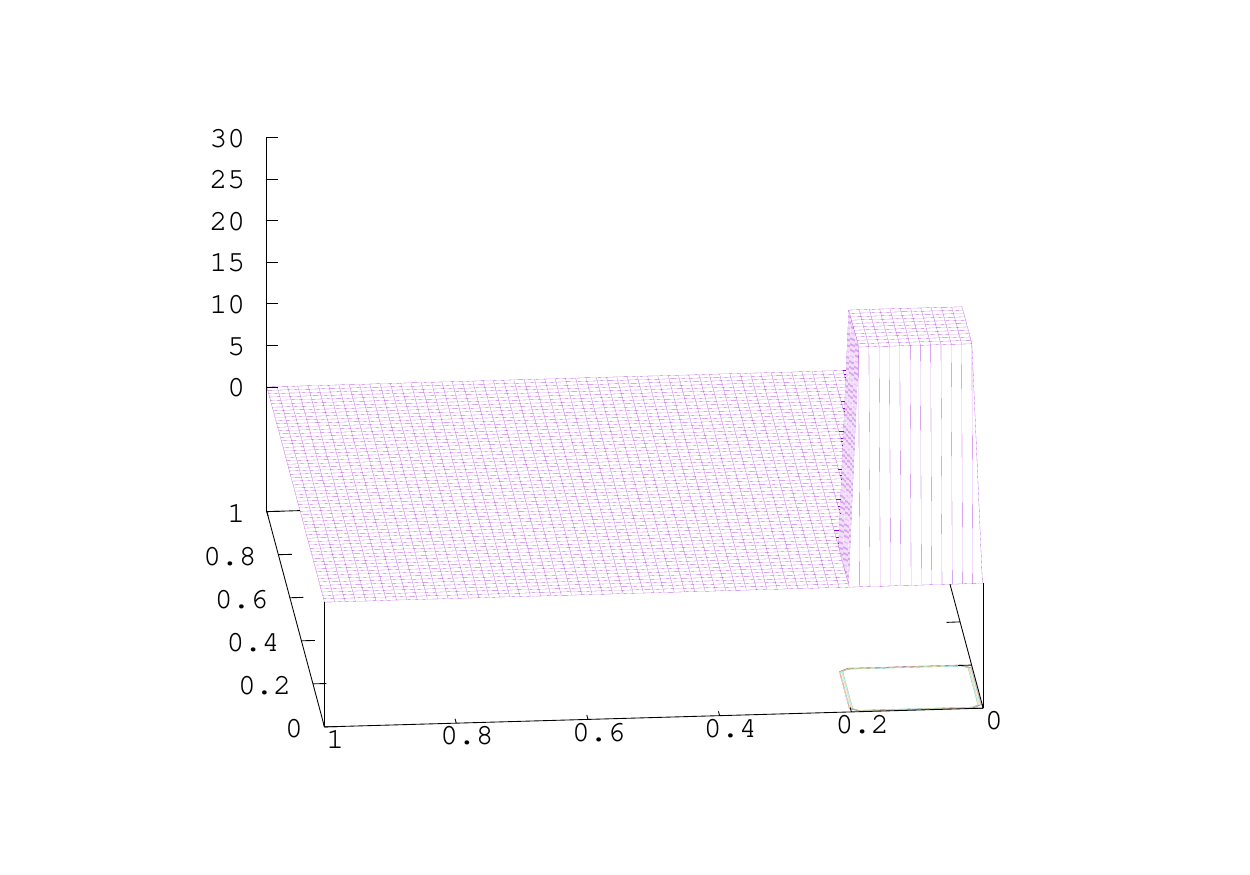} %
		\includegraphics[width=0.45\linewidth]{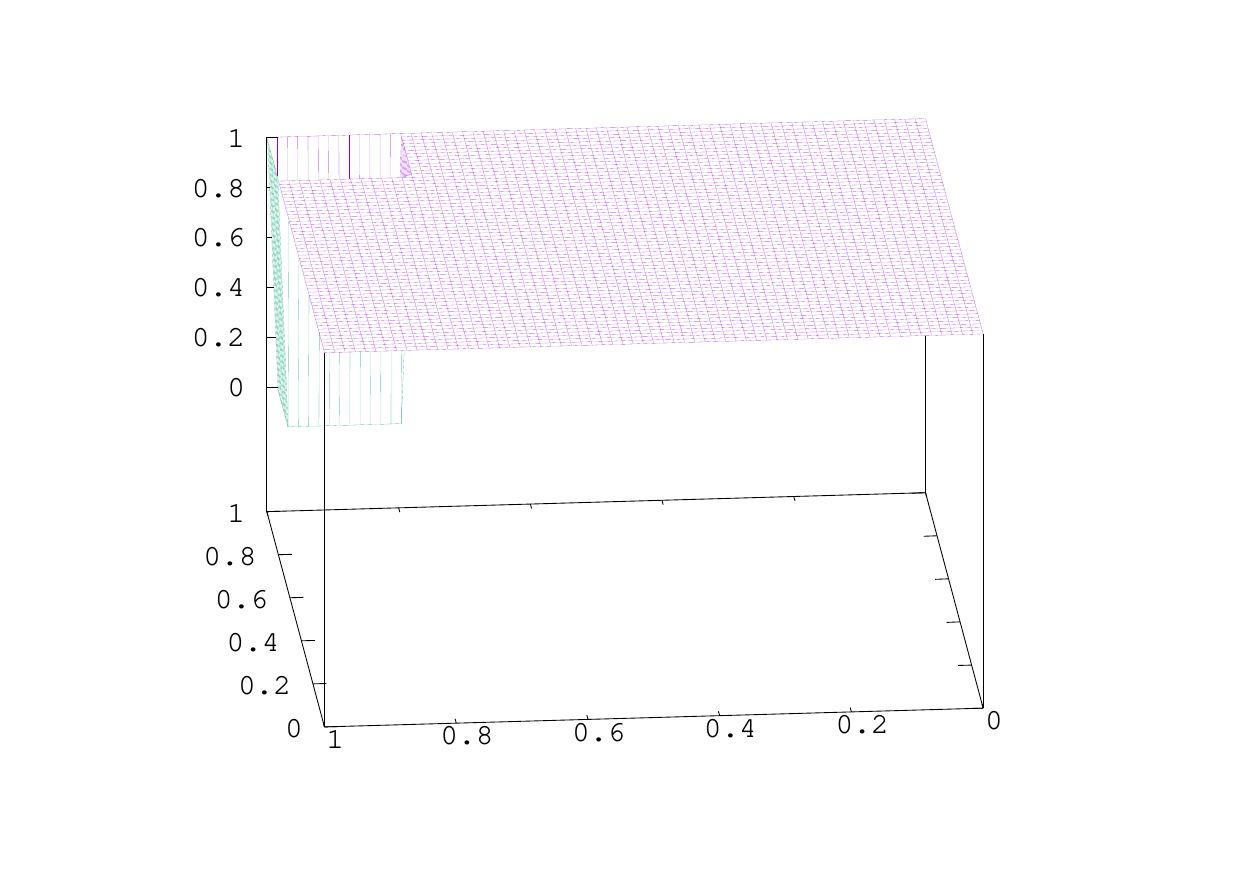}
	\captionof{figure}{\label{fig:AugLag-MFTC-congestion-m0-uT} Left: Initial distribution $m_0$. Right: terminal cost $g$.}
\end{center}

\clearpage

\section{\bf A Method based on monotone operators}
\label{AMS-num-sec:monotone-op}

Almulla, Ferreira and Gomes proposed in~\cite{MR3698446} two numerical methods for mean field games in infinite  horizon, in which the solution is stationary.
 These methods are further studied in~\cite{GomesSaude2018,GomesYang2018hessian}. The first method can be applied to variational MFG (or more generally to MFC) for which the PDE system corresponds to 
the Euler-Lagrange conditions arising  in the minimization of some energy functional.
 The idea is then to follow the gradient-flow generated by this energy until reaching a zero of the gradient.

 The second method can be applied more broadly to infinite horizon MFGs, provided a  monotonicity condition is satisfied.
 The main idea is that, under some assumptions, the solution to the system of PDEs can be recast as a zero of a monotone operator.
 The strategy is then to follow the flow generated by this operator until finding a zero.

\subsection{A Monotonic flow}
We focus here on the second method. Let us consider the case of Example~\ref{AMS-num-ex:special-case-sep} when, in addition, the dependence on the distribution is local and more specifically $f_0(x,m) = \log(m(x))$. The ergodic MFG (see \textit{e.g.}~\cite{MR2295621}) leads to the following system of PDEs: 
\begin{subequations}\label{AMS-num-eq:PDE-system-MFG-ergodic}
     \begin{empheq}[left=\empheqlbrace\,\,\,]{align}
     	\displaystyle
	& 0 = \lambda - \nu \Delta u(x) + H_0(x, \grad u(x)) - \log(m(x)), 
	&&\hbox{ in }  \TT^d,
	\\
	\displaystyle
	& 0 = -\nu \Delta m(x) - \diver\left( m \partial_pH_0 (x, \grad u(x)) \right), 
	&&\hbox{ in } \TT^d,
	\\
	& \int_{\TT^d} m = 1, \qquad \int_{\TT^d} u = 0, \qquad m>0 \hbox{ in } \TT^d,
     \end{empheq}
\end{subequations}
where the unknowns are the functions $u$ and $m$, and the ergodic constant $\lambda$. 
Let $A$ be the operator defined on the domain 
$$
	\mathcal{D} = \{(u,m) \in W^{2,2}(\TT^d) \times W^{1,2}(\TT^d) \,:\, \inf_{\TT^d} m > 0 \} \subseteq L^2(\TT^d) \times L^2(\TT^d)
$$ 
by
$$
	A \begin{pmatrix}
	u
	\\
	m
	\end{pmatrix} 
	= 
	\begin{pmatrix}
	 - \nu \Delta m- \diver\left( m(\cdot) \partial_p H_0 (\cdot, \grad u(\cdot)) \right)
\\
	 \nu \Delta u - H_0(\cdot, \grad u) + \log(m)
	\end{pmatrix}.
$$
It can be checked that this operator,
 seen as an operator from $L^2(\TT^d) \times L^2(\TT^d)$ to $L^2(\TT^d) \times L^2(\TT^d)$, is monotone on its domain, namely,
 $$
 	\langle A(u,m) - A(u',m'), (u,m) - (u',m') \rangle \ge 0, \qquad (u,m), (u',m') \in \mathcal{D}.
 $$

Then, we can consider the following flow $(u_\tau, m_\tau)_{\tau \ge 0}$, where  $(u_\tau, m_\tau) \in L^2([0,T] \times \TT^d) \times L^2([0,T] \times \TT^d)$ for all $\tau \ge 0$:
$$
	\frac{d}{d \tau}
	\begin{pmatrix}
	u_\tau
	\\
	m_\tau
	\end{pmatrix} 
	=
	- A \begin{pmatrix}
	u_\tau
	\\
	m_\tau
	\end{pmatrix} - \begin{pmatrix}
	0
	\\
	\lambda_\tau
	\end{pmatrix} ,
$$
with $\lambda_\tau$ such that $\int_{\TT^d} m_\tau = 1$. We stress that here $\tau$ denotes the time in the evolution of the flow (and recall that we focus on a stationary MFG here). Under suitable conditions, if the initial point $(u_0, m_0)$ is in $\mathcal{D}$, the flow remains in $\mathcal{D}$ until it reaches a point where the right-hand side is $0$, \textit{i.e.}, the solution to the ergodic MFG PDE system.

This evolution can be approximated using the finite difference discretization introduced in \S~\ref{AMS-num-sec:fin-diff}. More precisely, if the state space is of dimension $d=1$, we consider the following system for $n=0,1,\dots$:
$$
	\begin{pmatrix}
	(D_t U)^n
	\\
	(D_t M)^n
	\end{pmatrix} 
	=
	- \begin{pmatrix}
	 -\nu (\Delta_h M^{n+1})_{i}
	+ \mathcal{T}_i(U^{n+1}, M^{n+1})
	\\
	\lambda^{n+1} + \nu (\Delta_h U^{n+1})_{i}
	- \tilde{H}_0(x_i, [\grad_h U^{n+1}]_i) 
	+ \log(M^{n+1}_i)
	\end{pmatrix},
$$
where $\lambda^{n+1}$ is such that $\sum_i M^{n+1}_i = 1$.  For any $n \ge 0$, assuming that $\sum_i M^{n}_i = 1$ and summing over $i$ the second half of the above system of equations yields
\begin{align*}
	& \hbox{ $\lambda^{n+1}$ is such that $\sum_i M^{n+1} = 1$} 
	\\
	\Leftrightarrow \quad
	&0 = N_h \lambda^{n+1} + \sum_{i}  \left[ \nu (\Delta_h U^{n+1})_{i}
	- \tilde{H}_0(x_i, [\grad_h U^{n+1}]_i) 
	+ \log(M^{n+1}_i) \right]
	\\
	\Leftrightarrow \quad 
	&\lambda^{n+1} = \psi(U^{n+1}, M^{n+1}),
\end{align*}
where $\psi(U, M) = - h \sum_{i}  \left[ \nu (\Delta_h U)_{i}
	- \tilde{H}_0(x_i, [\grad_h U]_i) 
	+ \log(M_i) \right]$, for $U, M \in \RR^{N_h+1}$.
	
This leads directly to an iterative procedure, by starting from an initial guess $(U^0, M^0)$ and repeating the following step for increasing values of $n=0,1,\dots$:
\begin{equation}
\label{AMS-num-eq:monotoneop-discrete-flow}
	\begin{pmatrix}
	(D_t U)^n
	\\
	(D_t M)^n
	\end{pmatrix} 
	=
	- \begin{pmatrix}
	 -\nu (\Delta_h M^{n+1})_{i}
	+ \mathcal{T}_i(U^{n+1}, M^{n+1})
	\\
	\psi(U^{n+1}, M^{n+1}) + \nu (\Delta_h U^{n+1})_{i}
	- \tilde{H}_0(x_i, [\grad_h U^{n+1}]_i) 
	+ \log(M^{n+1}_i)
	\end{pmatrix},
\end{equation}
At each step, this computation involves solving a non-linear system to obtain $(U^{n+1}, M^{n+1})$, which can be done with Newton method for example. In order to tackle the fact that $M_i$ must remain positive for every $i$ (in order for the $\log$ to make sense), one can for instance use a damped version of Newton iterations.

\subsection{Numerical illustration}

 We illustrate this method on two examples. To the best of our knowledge, the convergence has been proved only for first order (\textit{i.e.}, $\nu=0$) stationary problems with a logarithmic coupling, see~\cite{MR3698446}. For other coupling costs, a projection step might be needed in order to ensure non-negativity of the distribution.

\textbf{Test case 1 (first order).} The following example is borrowed from~\cite{MR3698446}. We consider the case where $d=1, \nu = 0$ and $H_0$ is of the form
$$
	H_0(x,p) = \frac{1}{2} |p|^2 + b(x) p + V(x).
$$
with $b(x) = 2 c \pi \cos(2 \pi x)$ where $c \in [0,1]$, and 
$$
	V(x) = \sin(2 \pi x).
$$
In this case, the system~\eqref{AMS-num-eq:PDE-system-MFG-ergodic} admits a unique solution, given explicitly by
$$
	u^*(x) = c \sin(2 \pi x) \, , 
	\qquad
	m^*(x) = \frac{e^{V(x) - \frac{b(x)^2}{2}}}{ \int_{\TT} e^{V(y) - \frac{b(y)^2}{2}} dy } \,, 
	\qquad 
	\lambda =  \int_{\TT} e^{V(y) - \frac{b(y)^2}{2}} dy \,.
$$

\begin{figure}	
	\centering
	\begin{subfigure}[t]{0.45\linewidth}%
		\centering
		\includegraphics[width=\linewidth]{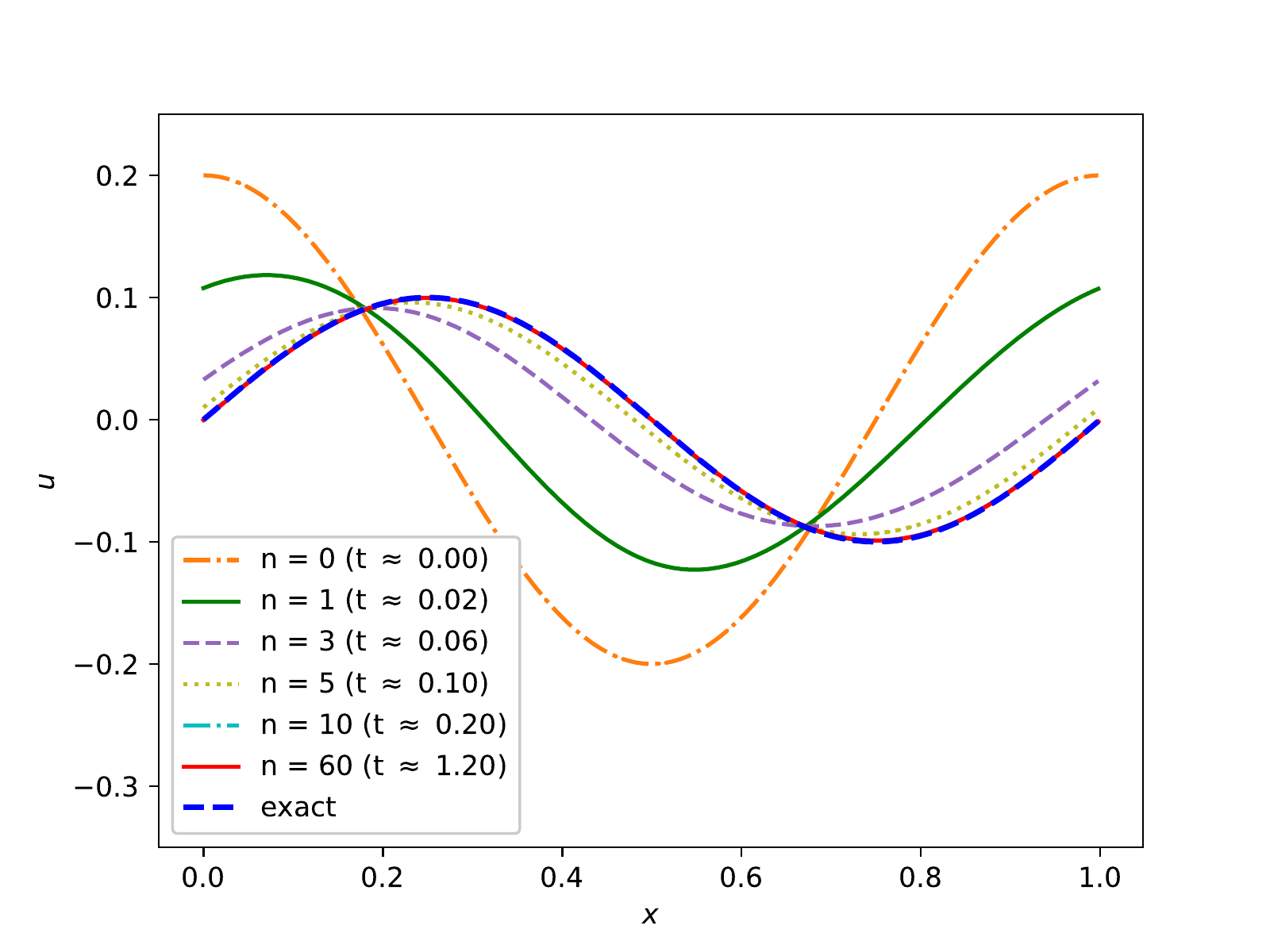}
	\end{subfigure}
	\quad
	\begin{subfigure}[t]{0.45\linewidth}
		\centering
	  	\includegraphics[width=\linewidth]{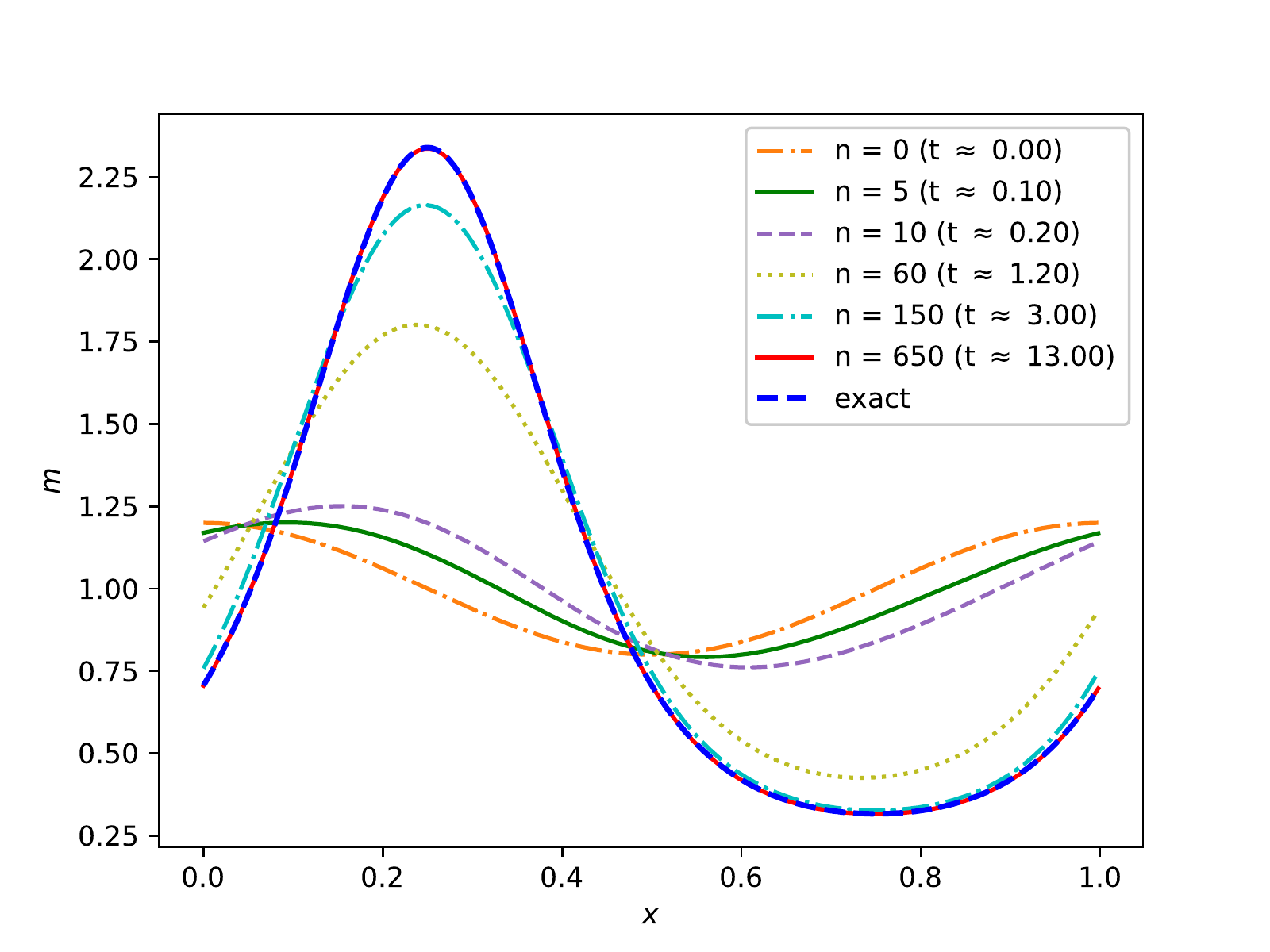}
	\end{subfigure}
	\caption{\label{AMS-num-fig:monotone-flow-evol-firstorder} Test case 1: Value function $U^n$ (left) and density $M^n$ (right) for several values of $n$, following the evolution~\eqref{AMS-num-eq:monotoneop-discrete-flow}. }
\end{figure}

\begin{figure}	
	\centering
	\begin{subfigure}[t]{0.45\linewidth}%
		\centering
		\includegraphics[width=\linewidth]{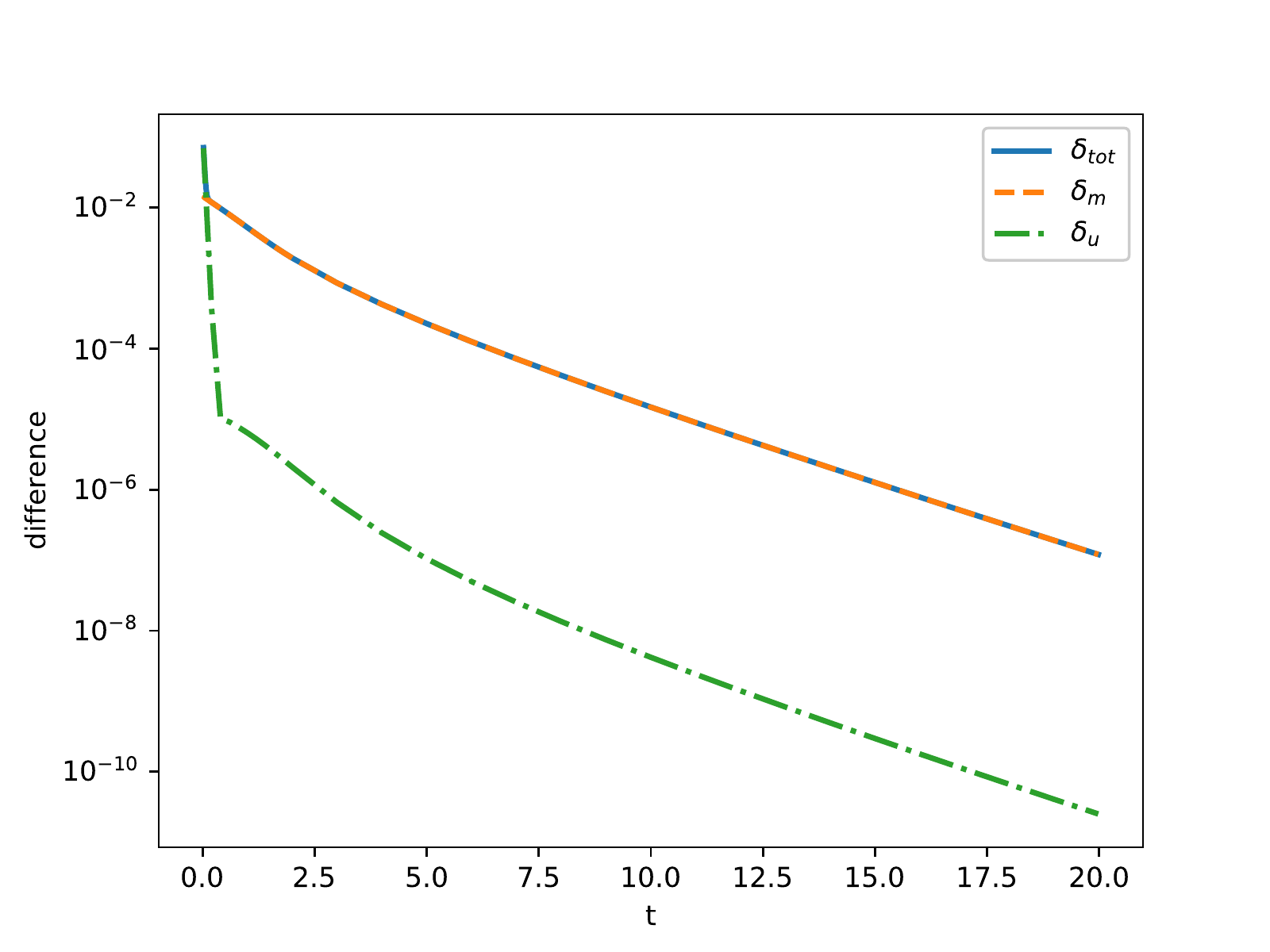}
	\end{subfigure}
	\quad
	\begin{subfigure}[t]{0.45\linewidth}
		\centering
	  	\includegraphics[width=\linewidth]{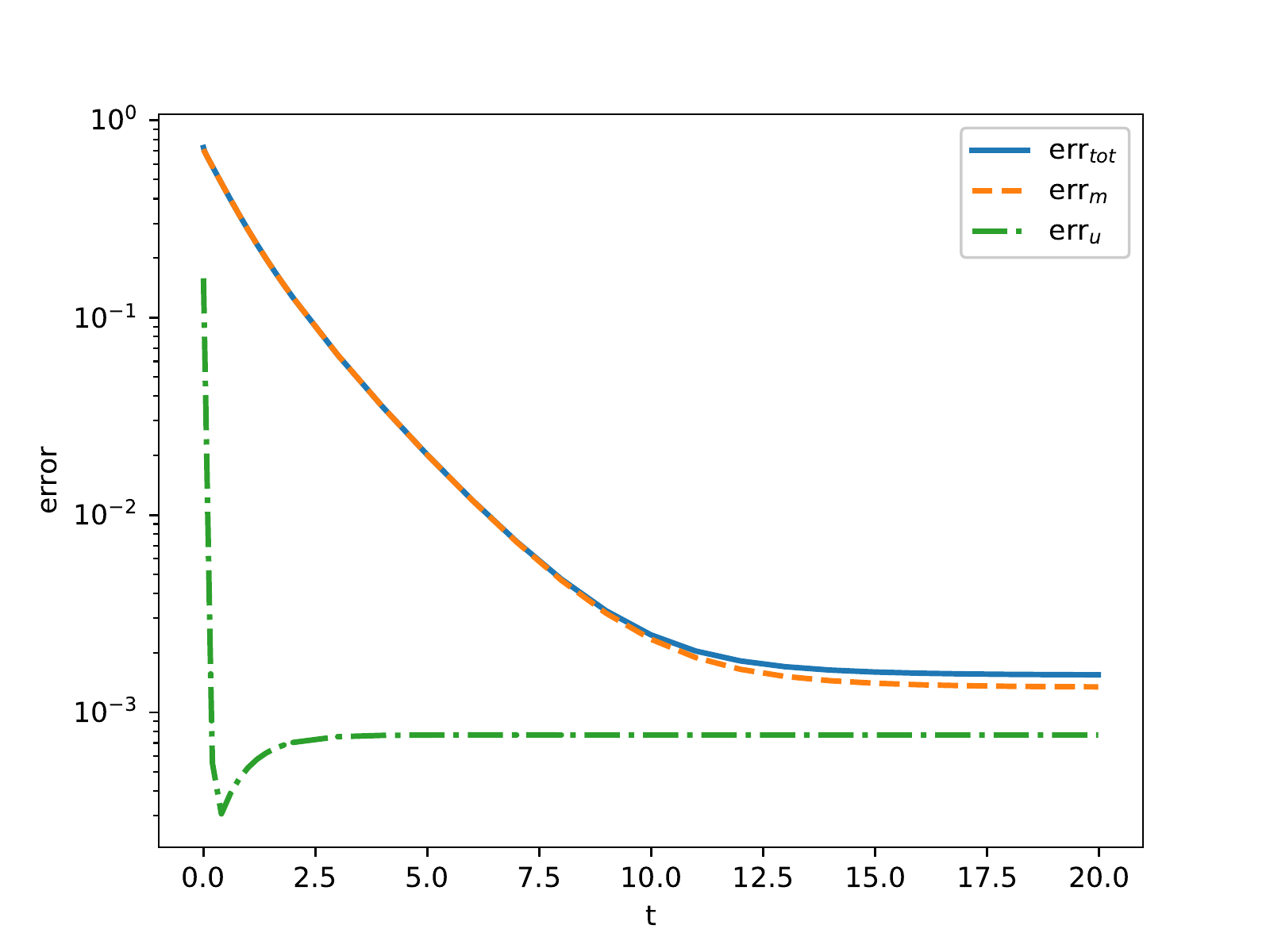}
	\end{subfigure}
	\caption{\label{AMS-num-fig:monotone-flow-error-firstorder} Test case 1: Difference between two iterations of the monotonic flow~\eqref{AMS-num-eq:monotoneop-discrete-flow} (left) and error with respect to the exact solution (right), computed as indicated in the text, see~\eqref{AMS-num-eq:monotoneop-def-deltas}-\eqref{AMS-num-eq:monotoneop-def-errors}.}
\end{figure}

We then consider the following discrete Hamiltonian (which is a modification of the one discussed in Example~\ref{AMS-num-eq-quadratic-discrete} to incorporate the effect of the extra terms with $b$ and $V$):
$$
	\tilde{H}_0(x, p_1, p_2) = \tilde{H}_0^{(1)}(p_1, p_2) + \tilde{H}_0^{(2)}(x, p_1, p_2) + V(x) ,
$$ 
with
$$
	\tilde{H}_0^{(1)}(p_1, p_2) = \frac{1}{2}|P_K(p_1,p_2)|^2,
	\qquad
	\tilde{H}_0^{(2)}(x, p_1, p_2) = 
	\begin{cases}
		b(x) p_1, & \hbox{ if } b(x) \le 0,
		\\
		b(x) p_2, & \hbox{ otherwise, }
	\end{cases}
$$
where $P_K$ denotes the projection on $K = \RR_- \times \RR_+$. Figure~\ref{AMS-num-fig:monotone-flow-evol-firstorder} displays $U^n$ and $M^n$ for several values of $n$ (note that the values of $n$ are not the same in each figure). Figure~\ref{AMS-num-fig:monotone-flow-error-firstorder} displays the $L^2$ difference between two iterations and the $L^2$ error with respect to the exact solution, defined respectively by:
\begin{equation}
\label{AMS-num-eq:monotoneop-def-deltas}
	\delta^n_u = \sqrt{ h \sum_i | U^{n}_i - U^{n-1}_i|^2 }, \quad 
	\delta^n_m = \sqrt{ h \sum_i | M^{n}_i - M^{n-1}_i|^2 }, \quad   
	\delta^n_{tot} = \sqrt{|\delta^n_u|^2 + |\delta^n_m|^2},
\end{equation}
and
\begin{equation}
\label{AMS-num-eq:monotoneop-def-errors}
	\text{err}^n_u = \sqrt{ h \sum_i | U^{n}_i - U^{*}_i|^2 }, \quad 
	\text{err}^n_m = \sqrt{ h \sum_i | M^{n}_i - M^{*}_i|^2 }, \quad   
	\text{err}^n_{tot} = \sqrt{|\text{err}^n_u|^2 + |\text{err}^n_m|^2}, 
\end{equation}
where 
$$
	U^{*}_i =  u^*(x_i),
	\qquad
	M^{*}_i = m^*(x_i),
$$
are given by the explicit solution. In the numerical computations to obtain these figures, we used $c = 0.1$, $N_h = 500$, $T = 20$, $N_T = 1000$.

\begin{figure}	
	\centering
	\begin{subfigure}[t]{0.45\linewidth}%
		\centering
		\includegraphics[width=\linewidth]{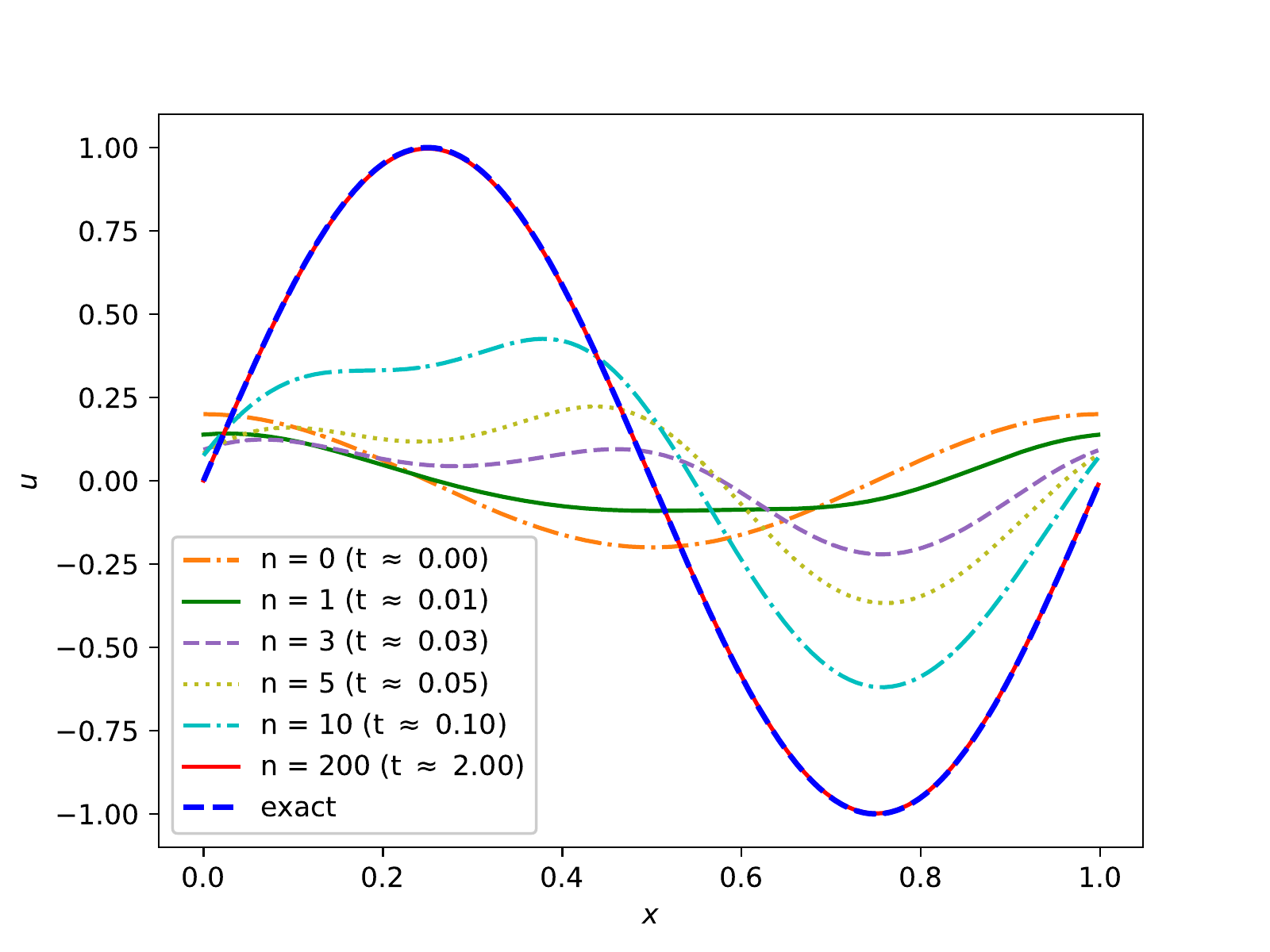}
	\end{subfigure}
	\quad
	\begin{subfigure}[t]{0.45\linewidth}
		\centering
	  	\includegraphics[width=\linewidth]{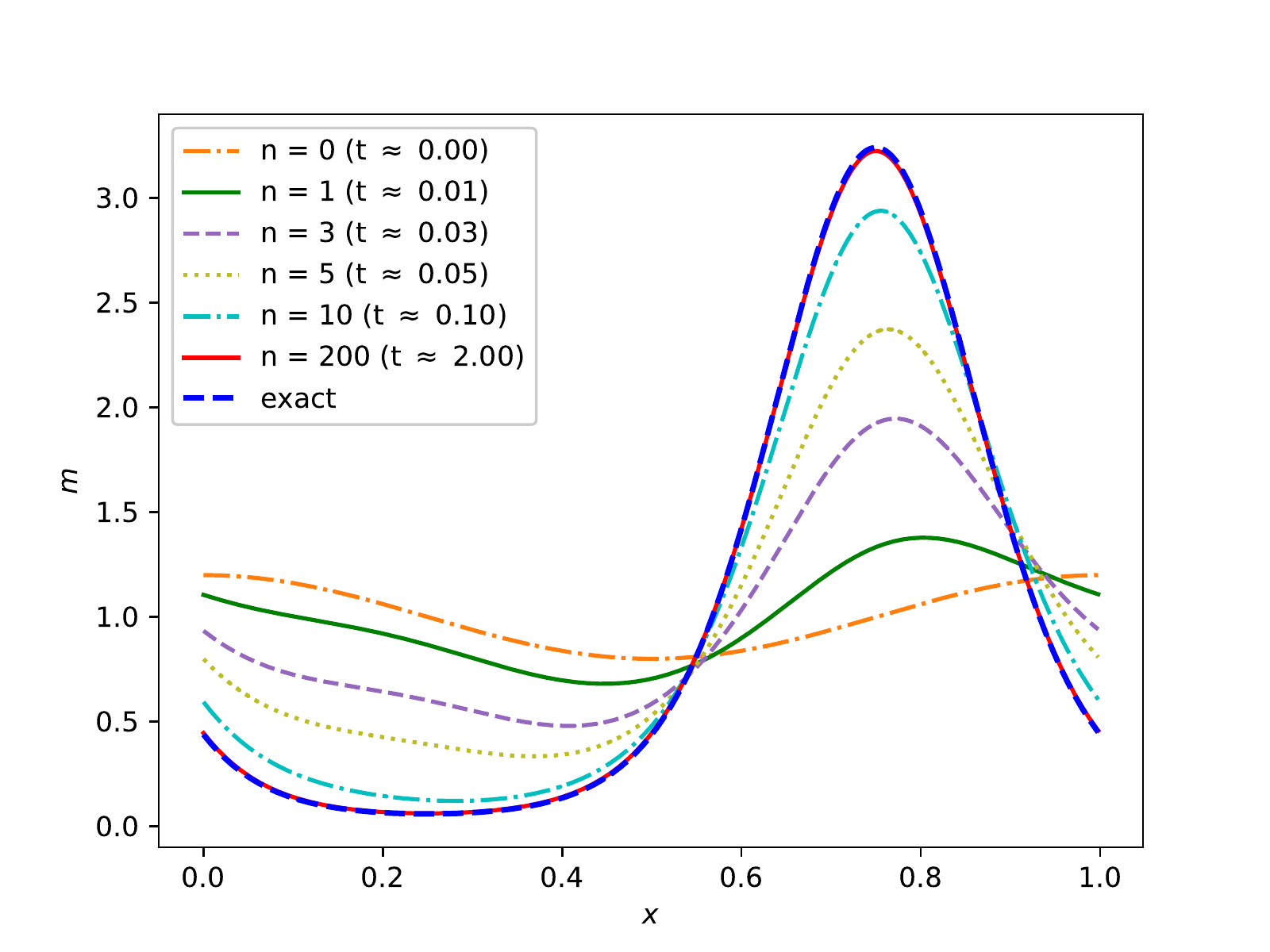}
	\end{subfigure}
	\caption{\label{AMS-num-fig:monotone-flow-evol-secondorder} Test case 2: Value function $U^n$ (left) and density $M^n$ (right) for several values of $n$, following the evolution~\eqref{AMS-num-eq:monotoneop-discrete-flow}. }
\end{figure}

\begin{figure}
	\centering
	\begin{subfigure}[t]{0.45\linewidth}%
		\centering
		\includegraphics[width=\linewidth]{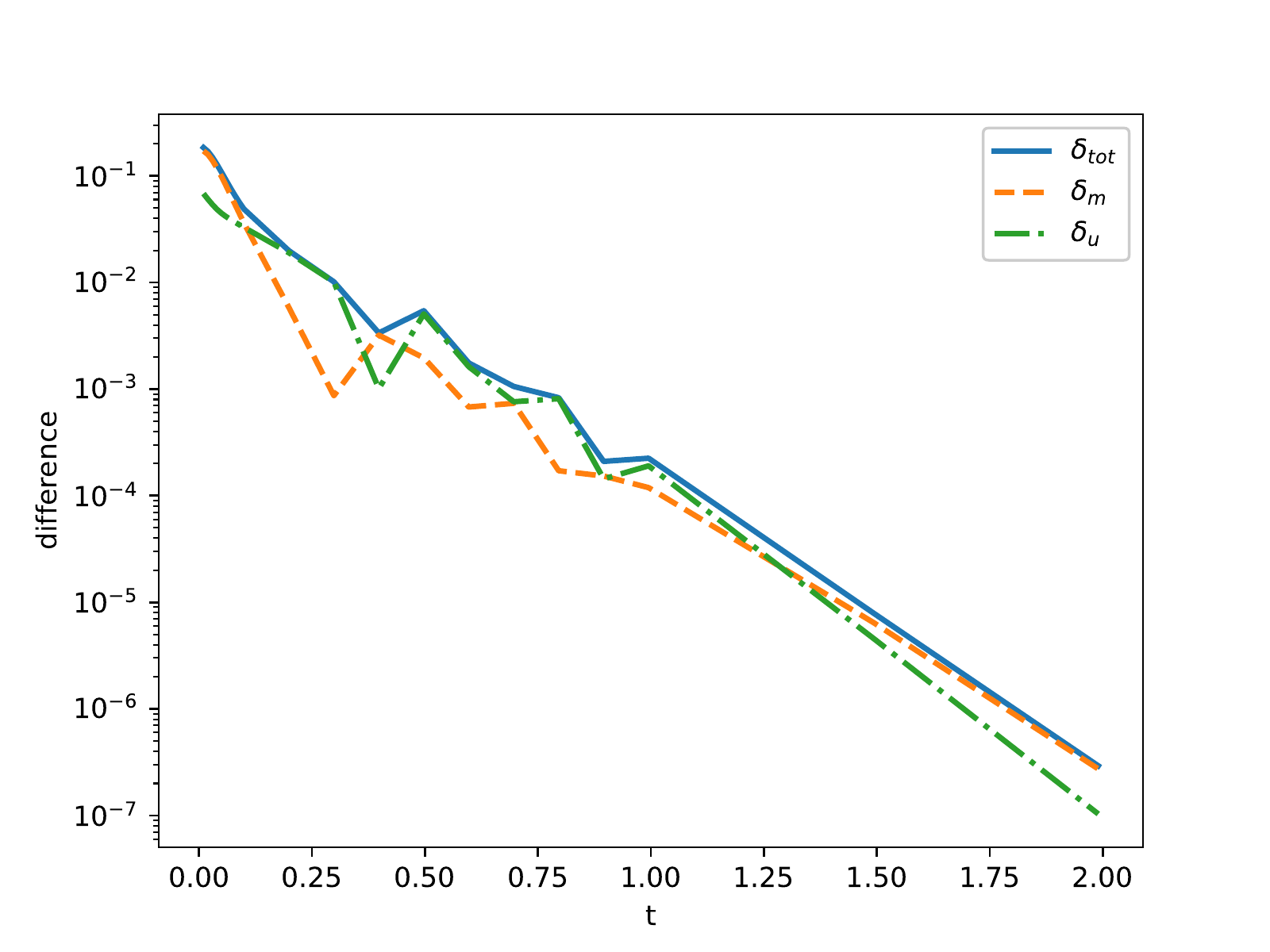}
	\end{subfigure}
	\quad
	\begin{subfigure}[t]{0.45\linewidth}
		\centering
	  	\includegraphics[width=\linewidth]{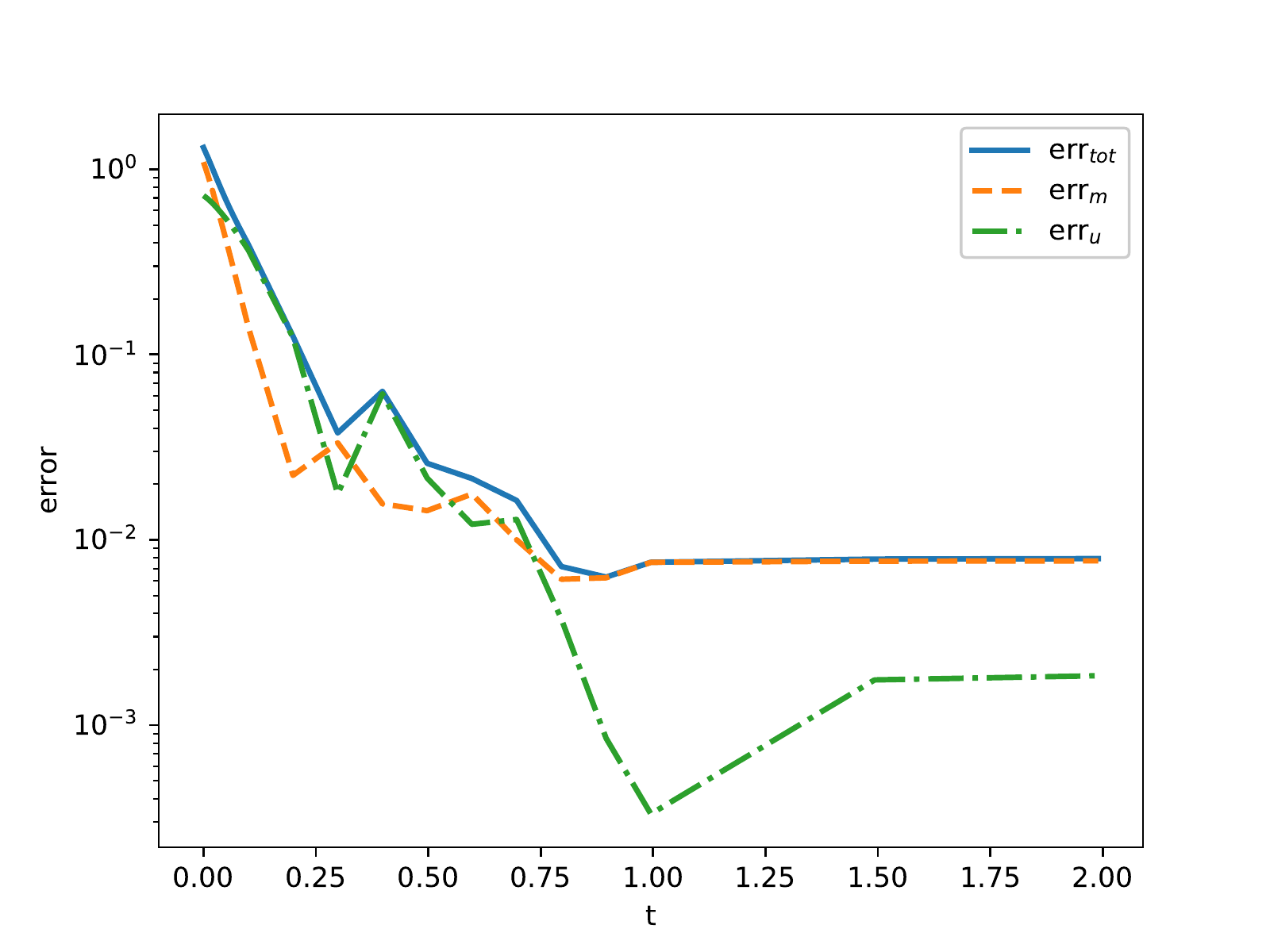}
	\end{subfigure}
	\caption{\label{AMS-num-fig:monotone-flow-error-secondorder} Test case 2: Difference between two iterations of the monotonic flow~\eqref{AMS-num-eq:monotoneop-discrete-flow} (left) and error with respect to the exact solution (right), computed as indicated in the text, see~\eqref{AMS-num-eq:monotoneop-def-deltas}-\eqref{AMS-num-eq:monotoneop-def-errors}.}
\end{figure}

\textbf{Test case 2 (second order). } We now consider an example of second-order MFG (\textit{i.e.}, with $\nu>0$) admitting an explicit solution. Although, to the best of our knowledge, the scheme has not been proved to converge in this setting, the numerical results seems to indicate that convergence holds. Let $b=0, \nu=0.5, V(x)=2 \pi^2  (-\kappa \sin(2 \pi x) - \kappa^2 \cos^2(2 \pi x)) - 2 \kappa \sin(2 \pi x)$, where $\kappa $ is a constant. Then the explicit solution is given by
$$
	u(x) = \kappa \sin(2 \pi x), \qquad m(x) = \kappa \frac{e^{2u(x)}}{\int e^{2u}}.
$$
We show the results for $\kappa = 1.0$. See Figures~\ref{AMS-num-fig:monotone-flow-evol-secondorder} and~\ref{AMS-num-fig:monotone-flow-error-secondorder}.

\clearpage

\section{\bf Methods based on approximation by neural networks}
\label{AMS-num-sec:NNapprox}

In this section, we present methods introduced recently which are based on tools borrowed from machine learning. To wit, we will use (artificial) neural networks to approximate functions of interest, and stochastic gradient descent to optimize over the weights of these neural networks. The first method we present is for MFC (or variational MFG) and directly based on the stochastic formulation. The second method tackles the master equation for finite-state mean-field problems.\footnote{For the sake of brevity, we do not discuss here other neural network based methods which have been developed for instance for the systems of partial differential equations or stochastic differential equations of MFGs. See the references in this section and in the conclusion.}

\subsection{Direct optimization for MFC}

Since MFC is an optimization problem, we present a first method which uses stochastic optimization tools from machine learning and which is directly based on the definition~\eqref{AMS-num-eq:def-J-MFC}--\eqref{AMS-num-eq:dyn-X-general-MFC} of the MFC problem. Starting from this formulation, we make three approximations in order to obtain a new problem that is more amenable to numerical treatment. 

\textbf{Approximation steps.} First, we restrict the set of (feedback) controls to be the set of neural networks with a given architecture. We introduce new notation to define this class of controls. 
We denote by:
\begin{align*}
	\mathbf{L}^\psi_{d_1, d_2} = 
	&\Big\{ \phi: \RR^{d_1} \to \RR^{d_2} \,\Big|\,  \exists (\beta, w) \in  \RR^{d_2} \times \RR^{d_2 \times d_1}, \forall  i \in \{1,\dots,d_2\}, 
	\\
	&\qquad \;\phi(x)_i = \psi\Big(\beta_i + \sum_{j=1}^{d_1} w_{i,j} x_j\Big) \Big\} 
\end{align*}
the set of layer functions with input dimension $d_1$, output dimension $d_2$, and activation function $\psi: \RR \to \RR$. Typical choices for $\psi$ are ReLU function $\psi(x) = x^+$ or the sigmoid function $\psi(x) = 1/(1+e^{-x})$. 
Building on this notation and denoting by $\circ$ the composition of functions, we define:
\begin{align}
\label{eq:NNarchi-FFC}
	\bN^\psi_{d_0, \dots, d_{\ell+1}} 
	= 
	&\Big\{ \varphi: \RR^{d_0} \to \RR^{d_{\ell+1}} \,\Big|\, \exists (\phi_i)_{i=0, \dots, \ell-1} \in \bigtimes_{i=0}^{i=\ell-1} \mathbf{L}^\psi_{d_i, d_{i+1}}, 
	\\
	&\qquad \exists \phi_\ell \in \mathbf{L}^{\mathrm{id}}_{d_{\ell}, d_{\ell+1}}, \varphi = \phi_\ell \circ \phi_{\ell-1} \circ \dots \circ \phi_0  \Big\} \, 
	\notag
\end{align}
 the set of regression neural networks with $\ell$ hidden layers and one output layer, the activation function of the output layer being the identity. The number $\ell$ of hidden layers, the numbers $d_0$, $d_1$, $\cdots$ , $d_{\ell+1}$ of units per layer, and the activation functions (one single function $\psi$ in the present situation), are what is usually called the architecture of the network. Once it is fixed, the actual network function $\varphi\in \bN^\psi_{d_0, \dots, d_{\ell+1}} $ is determined by the remaining real-valued parameters:
 $$
 \theta=(\beta^{(0)}, w^{(0)},\beta^{(1)}, w^{(1)},\cdots\cdots,\beta^{(\ell-1)}, w^{(\ell-1)},\beta^{(\ell)}, w^{(\ell)})
 $$
defining the functions $\phi_0$, $\phi_1$, $\cdots$ , $\phi_{\ell-1}$ and $\phi_\ell$ respectively. The set of such parameters is denoted by $\Theta$. For each $\theta\in\Theta$, the function $\varphi$ computed by the network will be denoted by $\ctrl_\theta \in \bN^\psi_{d_0, \dots, d_{\ell+1}}$.
As it should be clear from the discussion of the previous section, here, we are interested in the  case where $d_0 = d+1$ (since the inputs are time and state) and $d_{\ell+1} = \ctrldim$ (\textit{i.e.}, the control dimension). The problem becomes to minimize $J^{MFC}$ defined by~\eqref{AMS-num-eq:def-J-MFC}--\eqref{AMS-num-eq:dyn-X-general-MFC} over $\ctrl \in \bN^\psi_{d+1, d_1, \dots, d_{\ell}, \ctrldim} $, or equivalently, to minimize over $\theta \in \Theta$ the function:
\begin{align*}
	\bJ: \theta \mapsto \EE \left[\int_0^T f(X_t^{\ctrl_\theta}, m^{\ctrl_\theta}(t,\cdot), \ctrl_\theta(t,X_t^{\ctrl_\theta}) ) dt + g(X_T^{\ctrl_\theta}, m^{\ctrl_\theta}(T,\cdot)) \right]
\end{align*}
where $m^{\ctrl_\theta}(t,\cdot)$ is the probability density of the law of $X_t^{\ctrl_\theta}$, under the constraint that the process $X^{\ctrl_\theta}$ solves the SDE~\eqref{AMS-num-eq:dyn-X-general-MFC} with control $\ctrl_\theta$.

Second, we also need to approximate the density. A (computationally) simple option is to replace it by the empirical distribution of a system of $N$ interacting particles. Given a feedback control $\ctrl$, we denote by $(\underline X_t^{\ctrl})_{t} = (X_t^{1, \ctrl}, \dots, X_t^{N, \ctrl})_{t}$ the solution of:
\begin{equation}
\label{AMS-num-eq:NN-MFC-Nparticles}
	d X_t^{i, \ctrl} = b(X_t^{i, \ctrl}, m^{N, \ctrl}_t, \ctrl(t, X_t^{i, \ctrl})) dt + \sigma d W^i_t, \qquad t \ge 0,
\end{equation}
where 
$$
	m^{N, \ctrl}_t = \frac{1}{N} \sum_{j=1}^N \delta_{X_t^{j, \ctrl}}, 
$$ 
$(W^i)_{i=1,\dots,N}$ is a family of $N$ independent $d$-dimensional Brownian motions, and the initial positions $(X_0^{i, \ctrl})_{i=1,\dots,N}$ are i.i.d. with distribution given by the density $m_0$. Note that in~\eqref{AMS-num-eq:NN-MFC-Nparticles} the particles interact only through $m^{N, \ctrl}_t$ in the drift function. The controls are distributed in the sense that the control used in the dynamics of $X^{i, \ctrl}$ is a function of $t$ and $X_t^{i, \ctrl}$ itself only (and not of the position of the other particles). This choice is motivated by the fact that MFC and $N$-agent control problems with distributed controls are tightly connected, see~\cite{MR3753660}. We thus obtain the following new problem: Minimize over $\theta \in \Theta$ the function
\begin{align*}
	\bJ^{N}: \theta \mapsto \frac{1}{N} \sum_{i=1}^N \EE \left[\int_0^T f(X_t^{i, \ctrl_\theta}, m^{N, \ctrl_\theta}_t, \ctrl_\theta(t,X_t^{i, \ctrl_\theta}) ) dt + g(X_T^{i, \ctrl_\theta}, m^{N,\ctrl_\theta}_T) \right],
\end{align*}
under the dynamics~\eqref{AMS-num-eq:NN-MFC-Nparticles} with control $\ctrl_\theta$. The average over $i$ is here for the sake of analogy with $N$-player games, but each expectation in the sum has the same value since the agents are identically distributed. We interpret $\ctrl_\theta(t,X_t^{i, \ctrl_\theta})$ as the control used by player $i$ at time $t$. Note that it can be viewed as a feedback control, function of time and player $i$'s state. Since it is not a function of the empirical distribution $m^{N, \ctrl_\theta}_t$, one may wonder how the solution to this problem could be close to the solution of the corresponding MFC. This is where a mean-field approximation is helpful: player $i$'s control can implicitly depend on the MFC distribution  through its input $t$. Intuitively, the MFC distribution $m^{\ctrl_\theta}(t,\cdot)$ is a deterministic proxy for the $N$-agent empirical (stochastic) distribution $m^{N, \ctrl_\theta}_t$. This is because the initial distribution is given and the idiosyncratic noises disappear in the limiting distribution evolution.   

Last, the time interval is discretized. Let $N_T$ be a positive number, let $\Delta t = T/N_T$ and $t_n = n \Delta t$, $n =0,\dots, N_T$. We consider the problem: Minimize over $\theta \in \Theta$ the function
\begin{equation}
\label{AMS-num-eq:NN-MFC-cost-totalapprox}
	\bJ^{N, \Delta t}: \theta \mapsto \frac{1}{N} \sum_{i=1}^N \EE \left[\sum_{n=0}^{N_T-1} f(\check X_{t_n}^{i, \ctrl_\theta}, \check m^{N, \ctrl_\theta}_{t_n}, \ctrl_\theta(t_n,\check X_{t_n}^{i, \ctrl_\theta}) ) \Delta t + g(\check X_T^{i, \ctrl_\theta}, \check m^{N,\ctrl_\theta}_T) \right],
\end{equation}
under the dynamic constraint:
\begin{equation}
\label{AMS-num-eq:NN-MFC-Nparticles-Deltat}
	\check X_{t_{n+1}}^{i, \ctrl_\theta} = \check X_{t_{n}}^{i, \ctrl_\theta} + b( \check X_{t_{n}}^{i, \ctrl_\theta}, \check m^{N, \ctrl_\theta}_{t_{n}}, \ctrl_\theta(t, \check X_{t_{n}}^{i, \ctrl_\theta})) \Delta t + \sigma \Delta \check W^i_n, \qquad n = 0, \dots, N_T-1,
\end{equation}
and the initial positions $(\check X_0^{i, \ctrl_\theta})_{i=1,\dots,N}$ are i.i.d. with distribution given by the density $m_0$, where 
$$
	\check m^{N, \ctrl_\theta}_{t_{n}} = \frac{1}{N} \sum_{j=1}^N \delta_{\check X_{t_{n}}^{j, \ctrl_\theta}}, 
$$ 
and the $(\Delta \check W_n^i)_{i=1,\dots,N,n=0,\dots,N_T-1}$ are i.i.d. random variables with Gaussian distribution $\mathcal N(0, \Delta t)$.

Under suitable assumptions on the problem and the neural network architecture, the difference between $\inf_\theta\bJ^{N, \Delta t}(\theta)$ and $\inf_\ctrl J^{MFC(\ctrl)}$ goes to $0$ as $N_t, N$ and the number of parameters in the neural network go to infinity. See~\cite{carmona2019convergence2} for more details.

\textbf{Optimization procedure. } Note that the cost function~\eqref{AMS-num-eq:NN-MFC-cost-totalapprox} is in general non-convex due to the complicated way in which $\theta$ is involved in the cost (particularly since we consider neural networks). Moreover, $\theta$ is typically in (finite but) high dimension. In order to compute an (approximate) optimizer $\theta^*$, it is possible to run a stochastic gradient descent (SGD) by exploiting the fact that the cost~\eqref{AMS-num-eq:NN-MFC-cost-totalapprox} is written as an expectation. The randomness in this problem comes from the initial positions $\underline {\check X}_0 = (\check X_0^{i, \ctrl_\theta})_{i}$ and the noise increments $(\Delta \underline {\check W}_n)_{n=0,\dots,N_T} = (\Delta \check W_n^i)_{i,n}$. Hence $S = (\underline {\check X}_0, (\Delta \underline {\check W}_n)_{n})$ is going to play the role of a random sample in SGD. Given  a realization of $S$ and a choice of parameter $\theta$, we can construct the trajectory $(\check X_{t_n}^{i, \ctrl_\theta,S})_{i=1,\dots,N, n=0,\dots,N_T}$ by following~\eqref{AMS-num-eq:NN-MFC-Nparticles-Deltat} and we can compute the induced cost:
\begin{equation}
\label{AMS-num-eq:NN-MFC-cost-totalapprox-oneS}
	\bJ^{N, \Delta t}_S(\theta) = \frac{1}{N} \sum_{i=1}^N \left[\sum_{n=0}^{N_T-1} f(\check X_{t_n}^{i, \ctrl_\theta, S}, \check m^{N, \ctrl_\theta, S}_{t_n}, \ctrl_\theta(t_n,\check X_{t_n}^{i, \ctrl_\theta, S}) ) \Delta t + g(\check X_T^{i, \ctrl_\theta, S}, \check m^{N,\ctrl_\theta, S}_T) \right].
\end{equation}

The SGD procedure in this context is summarized in Algorithm~\ref{AMS-num-algo:SGD-MFC}. We refer to, \textit{e.g.},~\cite{MR3797719} for more details on SGD. The most costly step is the computation of the gradient $\nabla_\theta \bJ^{N, \Delta t}_S(\theta^{(\mathtt{k})})$  with respect to $\theta$. However, modern programming libraries (such as TensorFlow or PyTorch) allow us to perform this computation automatically using backpropagation, and to adjust the learning rate in an efficient way. The present method is thus very straightforward to implement: contrary to the methods presented in the previous sections, there is no need to derive by hand any PDE, any FBSDE, or any gradient, and one works directly with the definition of the MFC.

Besides this aspect, the main reasons behind the success of this method are the expressive power of neural networks (meaning that complex functions can be well approximated with relatively few parameters) and the fact that there is a priori no limitation on the number of iterations $\mathtt{K}$ (because the samples $S$ come from Monte-Carlo simulation and not from a finite set of data).

\begin{algorithm}[H]
\DontPrintSemicolon
\KwData{An initial parameter $\theta_0\in\Theta$; a number of steps $\mathtt{K}$;  a sequence $(\beta^{(\mathtt{k})})_{\mathtt{k}=0,\dots,\mathtt{K}-1}$ of learning rates.}
\KwResult{A parameter $\theta$ such that $\ctrl_{\theta}$ approximately minimizes $J^{MFC}$}
\Begin{
  \For{$\mathtt{k} = 0, 1, 2, \dots, \mathtt{K}-1 $}{
    Pick $S  = (\underline {\check X}_0, (\Delta \underline {\check W}_n)_{n})$\;
    Compute the gradient $\nabla_\theta \bJ^{N, \Delta t}_S(\theta^{(\mathtt{k})})$, see~\eqref{AMS-num-eq:NN-MFC-cost-totalapprox-oneS}\;
    Set $\theta^{(\mathtt{k}+1)} = \theta^{(\mathtt{k})} - \beta^{(\mathtt{k})} \nabla_\theta \bJ^{N, \Delta t}_S(\theta^{(\mathtt{k})})$ \;
    }
  \KwRet{ $\theta^{(\mathtt{K})}$}
  }
\caption{SGD for MFC\label{AMS-num-algo:SGD-MFC}}
\end{algorithm}

\textbf{Numerical illustration: LQ MFC. } We revisit the linear-quadratic mean field control problem, see~\S~\ref{AMS-num-sec:LQ-sub-MFC}. We consider an example in dimension $d=10$. For simplicity, we take a setting where the solution is the same in every dimension. More precisely, we take:
$$
	A = 0, \bar A = 0, 
	\quad B = Q= \bar Q = S = C = Q_T = \bar Q_T = S_T = \mathrm{Id}, 
	\quad \sigma = 0.1.
$$
The initial distribution is $\mathcal{N}(x_0, \sigma_0)$ with $x_0 = 2$ and $\sigma_0 = 0.2$.
In the implementation, we used a feedforward fully connected neural network, as introduced in~\eqref{eq:NNarchi-FFC}, with sigmoid activation function. Figure~\ref{fig:NN-MFC-LQ-L2errctrl} shows, as a function of $k$, the $L^2$ error between the control $\ctrl_{\theta^{(\mathtt{k})}}$ learnt by the neural network and the true optimal control $\ctrl^*$. We see that the error decreases when the number of time steps or the number of particles increase. Although the optimal control is simply a linear function of the state at each time $t$, the neural network does not know this information a priori and nevertheless manages to approximate the optimal solution.

\begin{center} %
		\includegraphics[width=0.45\linewidth]{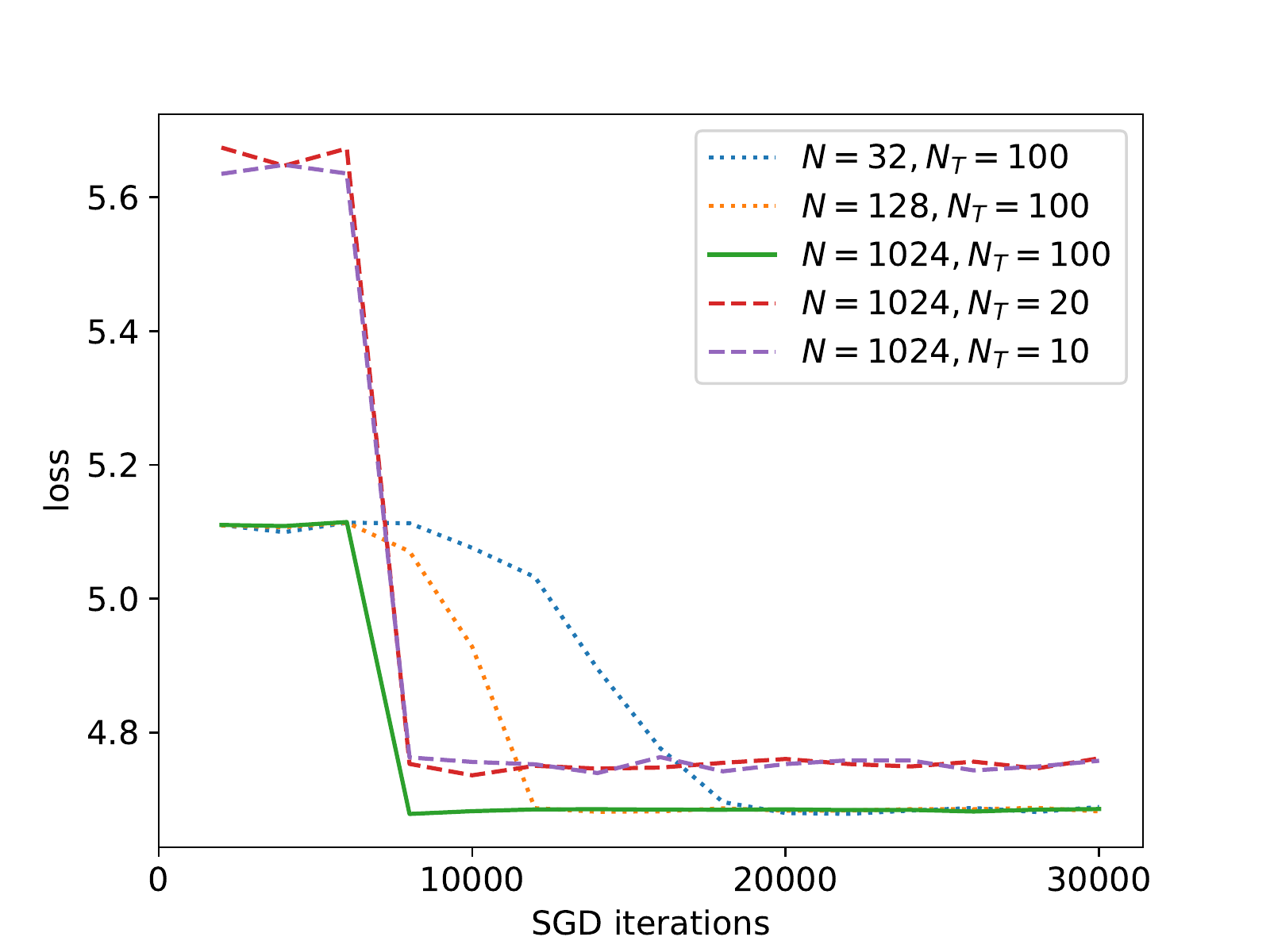}
	\quad
		\includegraphics[width=0.45\linewidth]{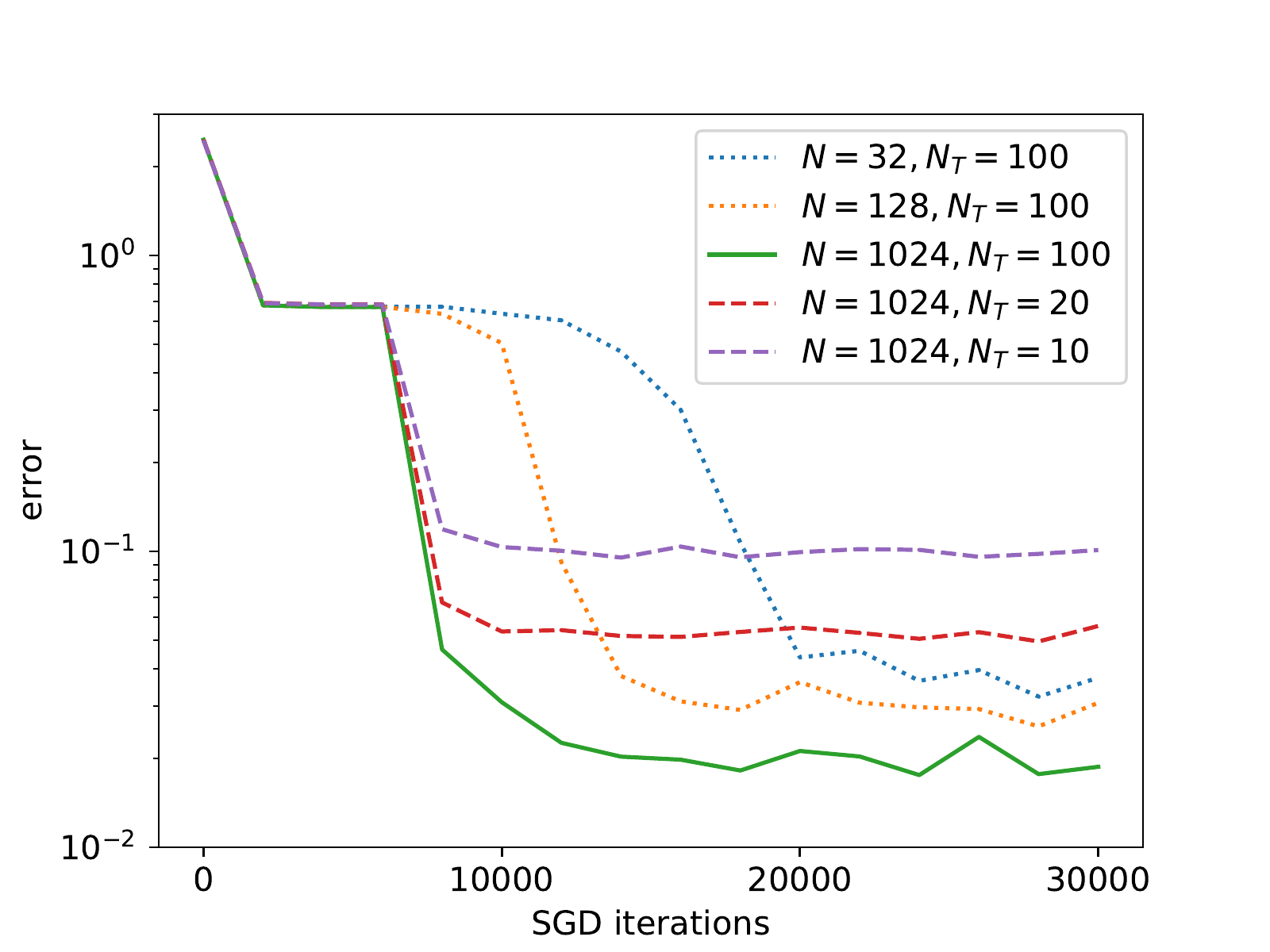}
	\captionof{figure}{\label{fig:NN-MFC-LQ-L2errctrl} LQ MFC with neural network approximation: total cost (left) and $L^2$ error on the optimal control (right) in dimension $d=10$ with $N_T$ time steps and $N$ particles as shown in the figures.}
\end{center}

\textbf{Extensions.} This technique has already been used for a while for standard stochastic optimal control problems, see \textit{e.g.},~\cite{parisini1996neural,MR2137498,han2016deep-googlecitations}. Its flexibility has lead to various applications such as, in the context of MFC, to problems with delay~\cite{fouque2019deep} or (biological) neural networks~\cite{agram2020deep}. Although it is very well-suited to MFC, this method is not immediately applicable to MFG which do not have a variational structure. However, it is possible to extend this technique to solve forward-backward stochastic differential equations (FBSDEs) of McKean-Vlasov type by rewriting the FBSDE as a mean-field optimal control of two forward SDEs with a terminal penalty, see \textit{e.g.}~\cite{carmona2019convergence2} for more details and numerical results. This line of research extends to the mean-field framework techniques that have been investigated for (non-MKV) FBSDE by E, Han and Jentzen in~\cite{MR3736669}.

\subsection{Solving the finite state master equation}
\label{AMS-num-sec:NN-finite-master}

All the methods presented so far to solve time-dependent MFC or MFG assume that the initial distribution is fixed. The solution is thus computed to obtain a flow of distribution which starts from this prescribed initial condition. However, in some applications, one does not know for sure the initial distribution. A possible strategy is to repeatedly use one of the methods presented earlier for various initial distributions. However, the computational cost would be prohibitive. It is thus interesting to be able solve the problem for any initial condition at once. This provides a motivation for solving the \index{master equation}master equation, which is an equation posed on the space of probability distributions, whose system of characteristics is the forward-backward PDE system introduced in \S~\ref{AMS-num-sec:PDEsystem}. Furthermore, the master equation is also useful to prove convergence of $N$-player games to MFG~\cite{MR3967062} or to obtain large deviation principles~\cite{MR4079435}.

\textbf{Master equation for finite state MFG.} We focus here on a setting in which the state space is finite. Following the discussion in~\cite[Section 7.2]{MR3752669}, let $\cE = \{e_1, \dots, e_d\}$ be a finite set, and let $\cA \subseteq \RR^\ctrldim$ be a Borel set, corresponding respectively to the state space and the action space. We will view $\cE$ as a subset of $\RR^d$ by identifying its elements with the canonical basis, and we will identify $\cP(\cE)$ with the simplex $\{m \in \RR^d \,| \, \sum_{i=1}^d m_i = 1\}$. Let $f: \cE \times \cP(\cE) \times \cA \to \RR$ and  $g: \cE \times \cP(\cE) \to \RR$ be respectively a running cost and a terminal cost functions. Let $\lambda : \cE \times \cP(\cE) \times \cA \to \RR$ be a jump rate function. To alleviate the notation, for every $m \in \cP(\cE)$ we write $m(x)$ instead of $m(\{x\})$. We will write $\RR^\cE$ for the set of functions from $\cE$ to $\RR$.

We then consider the following MFG equilibrium problem: Find a flow of probability densities $\hat{m}: [0,T] \times \cE \to \RR$ and a feedback control $\hat{\ctrl}: [0,T] \times \cE \to \cA$ satisfying  the following two conditions:
\begin{enumerate}
	\item $\hat{\ctrl}$ minimizes
\begin{align*}
	J^{MFG}_{\hat{m}}: \ctrl \mapsto  \EE \left[\int_0^T f(X_t^{\hat{m}, \ctrl}, \hat{m}(t,\cdot), \ctrl(t,X_t^{\hat{m}, \ctrl}) ) dt + g(X_T^{\hat{m}, \ctrl}, \hat{m}(T,\cdot)) \right]
\end{align*}
under the constraint that the  process $X^{\hat{m}, \ctrl} = (X_t^{\hat{m}, \ctrl})_{t \ge 0}$ is a nonhomogeneous $\cE$-valued Markov chain with transition probabilities determined by the $Q$-matrix of rates $q^{\hat{m}, \ctrl}: [0,T] \times \cE \times \cE \to \RR$ given by:
\begin{equation}
\label{AMS-num-eq:q-finite-MFG}
	q^{\hat{m}, \ctrl}(t, x,x') = \lambda(x, x', \hat{m}(t,\cdot), \ctrl(t,\cdot)), \qquad (t,x,x') \in [0,T] \times \cE \times \cE,
\end{equation}
and $X_0^{\hat{m}, \ctrl}$ has distribution with density $m_0$;
	\item For all $t \in [0,T]$, $\hat{m}(t,\cdot)$ is the law of $X_t^{\hat{m}, \hat{\ctrl}}$.
\end{enumerate}

To formulate an optimality condition, we introduce the Lagrangian $L: \cE \times \cP(\cE) \times \RR^\cE \times \cA \to \RR$ defined by:
$$
	L(x, m, h, \ctrl) = \sum_{x' \in \cE} \lambda(x, x', m, \ctrl) h(x') + f(x, m, \ctrl) ,
$$
and the Hamiltonian:
$$
	H(x, m, h) = \sup_{\ctrl \in \cA} -L(x, m, h, \ctrl).
$$ 
Assuming there is a unique maximizer for every $(x, m, h) \in \cE \times \cP(\cE) \times \RR^\cE$, we denote:
\begin{equation}
\label{AMS-num-eq:finiteMFG-ctrl-argmax}
	\ctrl^*(x, m, h) = \argmax_{\ctrl \in \cA} -L(x, m, h, \ctrl).
\end{equation}
It will also be useful to introduce the function $q^*: \cE \times \cE \times \cP(\cE) \times \RR^\cE \to \RR$ defined by
$$
	q^*(x,x',m,h) = \lambda\big(x, x', m, \ctrl^*(x,m, h)\big).
$$

In the spirit of the forward-backward PDE system, the solution of a finite state MFG can be characterized through a system of  ODEs: a forward ODE for the distribution $m: [0,T] \times \cE \to \RR$ and a backward ODE for the value function $u: [0,T] \times \cE \to \RR$ of an infinitesimal player. To wit, under suitable conditions (see \textit{e.g.}~\cite[section 7.2]{MR3752669}), there is a unique MFG equilibrium $(\hat{m}, \hat{\ctrl})$ and it is given by:
$$
	\hat{\ctrl}(t,x) = \ctrl^*(x, m(t,\cdot), u(t,\cdot)), 
$$  
where $\ctrl^*$ is defined by~\eqref{AMS-num-eq:finiteMFG-ctrl-argmax} and $(u,m)$ solves the forward-backward system:
\begin{subequations}\label{AMS-num-eq:ODE-system-finiteMFG}
     \begin{empheq}[left=\empheqlbrace]{alignat=2}
     	\displaystyle
	0 
	&= - \partial_t u(t, x) + H(x, m(t,\cdot), u(t,\cdot)),
	&&\quad (t,x) \in [0,T) \times \cE,
	\label{AMS-num-eq:ODE-system-finiteMFG-HJB}
	\\
	\displaystyle
	0 
	&= \partial_t m(t, x) - \sum_{x' \in \cE} m(t, x') q^*(x', x, m(t,\cdot), u(t,\cdot)), 
	&&\quad (t,x) \in (0,T] \times \cE,
	\label{AMS-num-eq:ODE-system-finiteMFG-KFP}
	\\
	&u(T,x) = g(x, m(T,\cdot)), \qquad m(0,x) = m_0(x), 
	&&\quad x \in \cE.
	\label{AMS-num-eq:ODE-system-finiteMFG-IT}
     \end{empheq}
\end{subequations}

Since $\cE$ is finite, $m(t,\cdot)$ and $u(t,\cdot)$ can be identified with vectors and each equation can be viewed as an ODE. This forward-backward system can be tackled using techniques similar to the ones described in the previous sections for the PDE system arising in the continuous state space case. For instance, one can first discretize time with a semi-implicit scheme and then use fixed point iterations, alternating between the HJB and the FP equation, or use Newton iterations; see \S~\ref{AMS-num-sec:sol-strat}.

 Due to the coupling between the two equations, the value function $u$ depends implicitly on the distribution $m$. The solution to the master equation allows us to make this dependence explicit. In the present setting, this equation takes the following form (see \textit{e.g.}~\cite[section 7.2]{MR3752669} for more details):
 \begin{equation}
 \label{AMS-num-eq:master-finiteMFG}
 	-\partial_t \cU(t,x,m) 
	+ H(x,m,\cU(t, \cdot, m))
	- \sum_{x' \in \cE} h^*(m ,\cU(t, \cdot, m))(x') \frac{\partial \cU(t, x ,m)}{\partial m(x')} = 0, 
 \end{equation}
for $(t,x,m) \in [0,T] \times \cE \times \cP(\cE)$, with the terminal condition $\cU(T,x,m) = g(x,m)$, for $(x,m) \in \cE \times \cP(\cE)$. Here, $h^*: [0,T] \times \cP(\cE) \times \RR^\cE \times \cE \to \RR$ is defined as:
$$
	h^*(m, u)(x') = \sum_{x \in \cE} \lambda(x, x', m, \ctrl^*(x, m, u)) m(x).
$$
The notation $\frac{\partial \cU(t, x ,m)}{\partial m(x')}$ represents the (classical) partial derivative of $\RR^d \ni m \mapsto \cU(t,x,m)$ with respect to the coordinate corresponding to $x'$ when $m$ is viewed as a vector of dimension $d$. The link between the master equation and the forward-backward system is the following. For every initial distribution $m_0 \in \cP(\cE)$, 
\begin{equation}
\label{AMS-num-eq:finiteMFG-link-U-u}
	\cU(t,x, m^{m_0}(t, \cdot)) = u^{m_0}(t,x), \qquad (t,x) \in [0,T] \times \cE,
\end{equation}
where $(u^{m_0}, m^{m_0})$ is the solution to~\eqref{AMS-num-eq:ODE-system-finiteMFG} starting with $m(0,\cdot) = m_0$. In words, the solution to the forward-backward system plays the role of characteristics for the master equation, and $\cU$ explicitly captures the (implicit) dependence of the infinitesimal player's value function $u$ on the population's distribution $m$.

\textbf{Principle of the Deep Galerkin method. } Note that the master equation~\eqref{AMS-num-eq:master-finiteMFG} is posed on the space $[0,T] \times \cE \times \cP(\cE) \subset [0,T] \times \cE \times \RR^d$, which is high dimensional as soon as the number of states, $d$, is large. Solving high-dimensional PDEs using for instance finite difference schemes is in general computationally expensive  due to the fact that the number of grid points increases exponentially with the dimension. Recently, several methods based on machine learning tools have been developed for this purpose. Here, we propose to use the Deep Galerkin Method (DGM) introduced in~\cite{MR3874585} to solve the master equation~\eqref{AMS-num-eq:master-finiteMFG}. The main idea is to first replace the unknown function, namely $\cU$ in our case, by a neural network, say $\cU_\theta$, with parameters $\theta$, and then to optimize over $\theta$ in order to minimize the residual of the PDE~\eqref{AMS-num-eq:master-finiteMFG}.  This optimization is achieved by using SGD in the following way: at each iteration, a point in the domain is picked according to a chosen distribution, then the gradient of the PDE residual with respect to the neural network parameters is computed and a gradient step is performed. The procedure is summarized in Algorithm~\ref{AMS-num-algo:DGM-master-finiteMFG}. Here, we fix a neural network architecture and denote by $\Theta$ the set of possible parameters for neural networks with this architecture. For the PDE residual we use the notation: for $\theta \in \Theta$ and $S = (t,x,m) \in [0,T] \times \cE \times \cP(\cE)$,
\begin{equation}
\label{AMS-num-eq:master-finiteMFG-residual}
	Res_S(\theta) = \left|\partial_t \cU_\theta(S) - H(x,m,\cU_\theta(S))
	+ \sum_{x' \in \cE} h^*(m ,\cU_\theta(S))(x') \frac{\partial \cU_\theta(S)}{\partial m(x')}\right|^2.
\end{equation}
In the implementation, the neural network $\cU_\theta$ is a function of $S$ for which partial derivatives can be computed using automatic differentiation.  The computation of the residual can thus be implemented without approximating the derivatives. Then, the gradient $\nabla_\theta Res_S(\theta)$ of the residual with respect to $\theta$ can be computed using backpropagation.  The DGM amounts to try to minimize the loss function:
$$
	L(\theta) = \EE_{S \sim \nu}\left[ Res_S(\theta) \right].
$$

\begin{algorithm}[H]
\DontPrintSemicolon
\KwData{An initial parameter $\theta_0\in\Theta$; a number of steps $\mathtt{K}$;  a sequence $(\beta^{(\mathtt{k})})_{\mathtt{k}=0,\dots,\mathtt{K}-1}$ of learning rates; a probability distribution $\nu$ on $[0,T] \times \cE \times \cP(\cE)$.}
\KwResult{A parameter $\theta$ such that $\cU_{\theta}$ approximately solves~\eqref{AMS-num-eq:master-finiteMFG}}
\Begin{
  \For{$\mathtt{k} = 0, 1, 2, \dots, \mathtt{K}-1 $}{
    Pick $S  = (t,x,m) \sim \nu$\;
    Compute the gradient $\nabla_\theta Res_S(\theta^{(\mathtt{k})})$ (see~\eqref{AMS-num-eq:master-finiteMFG-residual})\;
    Set $\theta^{(\mathtt{k}+1)} = \theta^{(\mathtt{k})} - \beta^{(\mathtt{k})} \nabla_\theta Res_S(\theta^{(\mathtt{k})})$ \;
    }
  \KwRet{ $\theta^{(\mathtt{K})}$}
  }
\caption{DGM for the master equation\label{AMS-num-algo:DGM-master-finiteMFG}}
\end{algorithm}

\textbf{Example 1: A Cybersecurity model. } We consider a model introduced in~\cite{MR3575619} and revisited in~\cite[Section 7.2.3]{MR3752669}. In this model, each player owns a computer which can be either defended (D) or undefended (U), and either infected (I) or susceptible (S) of infection. Hence the set $\cE$ has four elements corresponding to the four possible combinations: $\cE = \{DI, DS, UI, US\}$. The action set is $\cA = \{0,1\}$, where $0$ is interpreted as the fact that the player is satisfied with the current level of protection (D or U) of its computer, whereas $1$ means that she wants to change this level of protection (\textit{i.e.}, she wants to go from D to U or vice versa). In the latter case, the update occurs at a (fixed) rate $\rho >0$. At each of the four states, all the computers are indistinguishable. When infected, each computer may recover at rate $q_{rec}^D$ or $q_{rec}^U$ depending on whether it is defended or not. On the other hand, a computer may be infected directly by a hacker, at rate $v_H q_{inf}^D$ (resp. $v_H q_{inf}^U$) if it is defended (resp. undefended), or it may be infected by undefended infected computers, at rate $\beta_{UU}\mu(\{UI\})$ (resp. $\beta_{UD}\mu(\{UI\})$) if it is undefended (resp. defended), or by defended infected computers, at rate $\beta_{DU}\mu(\{DI\})$ (resp. $\beta_{DD}\mu(\{DI\})$) if it is undefended (resp. defended). In short, the matrix of transition rates is given by: for $m \in \cP(\cE), a \in \cA$, 
$$
	\lambda(\cdot, \cdot, m, a) 
	= \left( \lambda(x, x', m, a) \right)_{x,x' \in \cE}
	= \begin{pmatrix}
	\dots 	& 		P^{m,a}_{DS \rightarrow DI}	&	 \rho a 	&	0
	\\
	q_{rec}^D 	& 	\dots 		&	 0	&	\rho a
	\\
	\rho a 	& 	0 		&	 \dots	&	P^{m,a}_{US \rightarrow UI}
	\\
	0	&	\rho a	&	q_{rec}^U	& \dots
	\end{pmatrix}
$$
where 
\begin{align*}
	&P^{m,a}_{DS \rightarrow DI} = v_H q_{inf}^D + \beta_{DD} m(\{DI\})  + \beta_{UD} m(\{UI\}) ,
	\\
	&P^{m,a}_{US \rightarrow UI} = v_H q_{inf}^U + \beta_{UU} m(\{UI\}) + \beta_{DU} m(\{DI\}),
\end{align*}
and all the instances of $\dots$ should be replaced by the negative of the sum of the entries of the row in which $\dots$ appears on the diagonal. At each time, the player pays a protection cost $k_D>0$ if its computer is defended, and a penalty $k_I>0$ if it is infected. 
There is no terminal cost and, given a mean-field flow $m$ and a control $\ctrl$, the instantaneous cost at time $t$ is hence 
$$
	f(X^{m,\ctrl}_t, m(t,\cdot), \ctrl(t, X^{m,\ctrl}_t)) = -\left[ k_D \indic_{\{DI, DS\}}(X^{m,\ctrl}_t) + k_I \indic_{\{DI, UI\}}(X^{m,\ctrl}_t)\right].
$$
Note that in this example, the cost is independent of the population's distribution and the mean-field interactions are only in the dynamics. 

\begin{figure}[h]
	\begin{subfigure}{.45\columnwidth}
		\centering
		\includegraphics[width=\columnwidth]{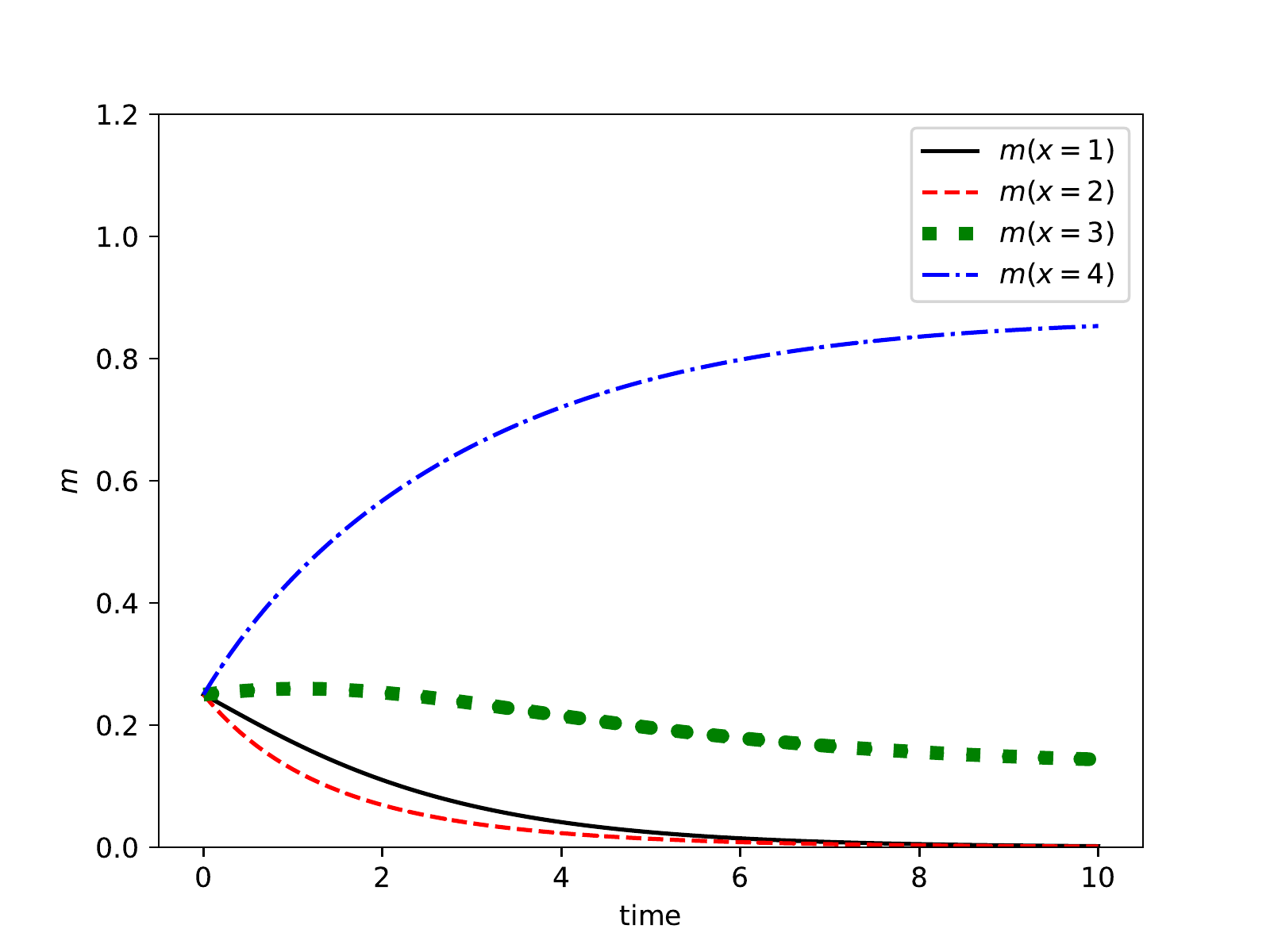}
	\end{subfigure}%
	\begin{subfigure}{.45\columnwidth}
		\centering 
		\includegraphics[width=\columnwidth]{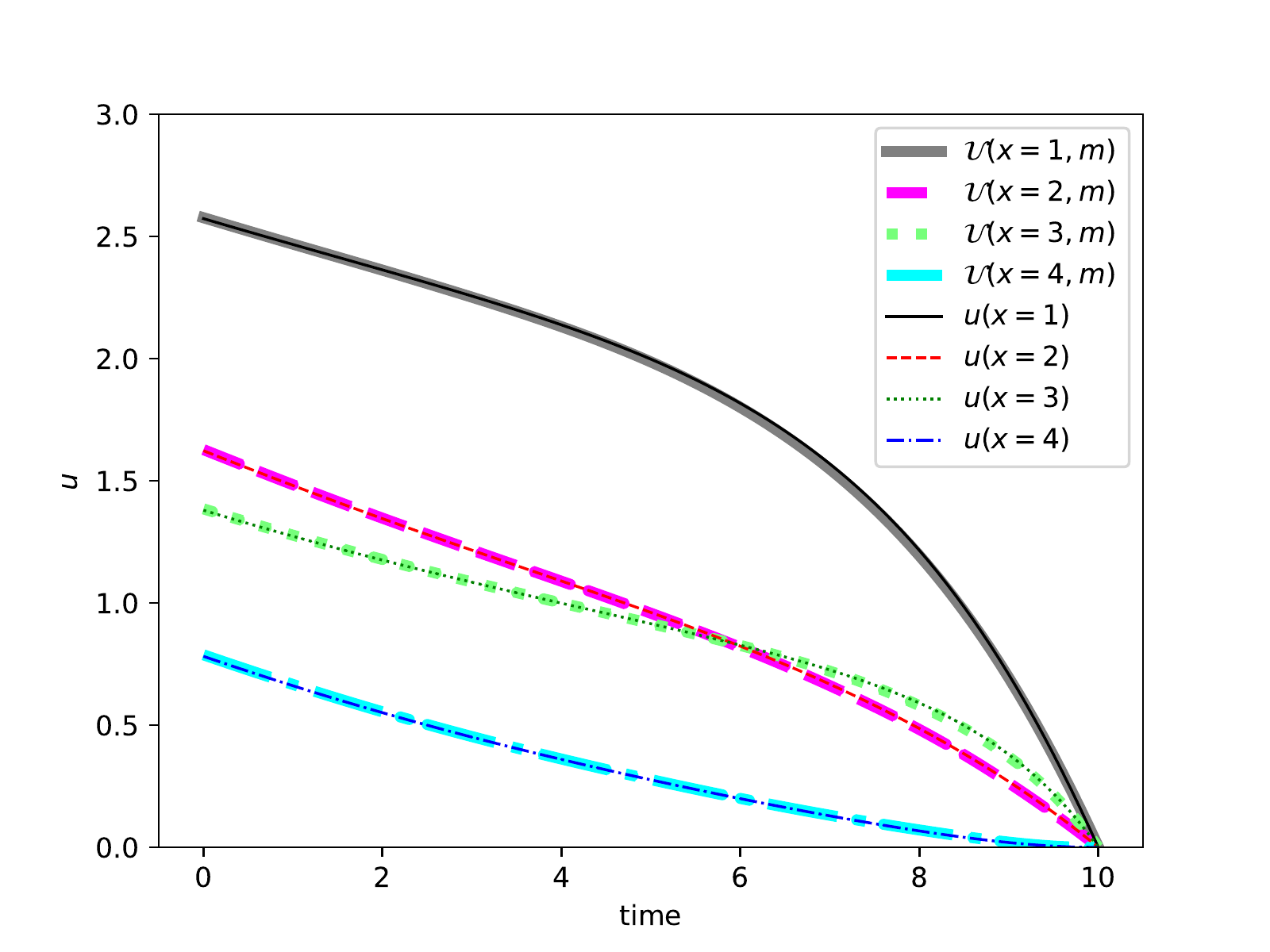}
	\end{subfigure}
	\caption{MFG Cyber-security example, test case 1: Evolution of the distribution $m^{m_0}$ (left) and the value function $u^{m_0}$ and $\cU(\cdot, \cdot, m^{m_0}(\cdot))$ (right) for $m_0 = (1/4, 1/4, 1/4, 1/4)$. } 
 	\label{AMS-num-fig:finiteMFG-cyber-Master-m0-1}
\end{figure}

\begin{figure}[h]
	\begin{subfigure}{.45\columnwidth}
		\centering
		\includegraphics[width=\columnwidth]{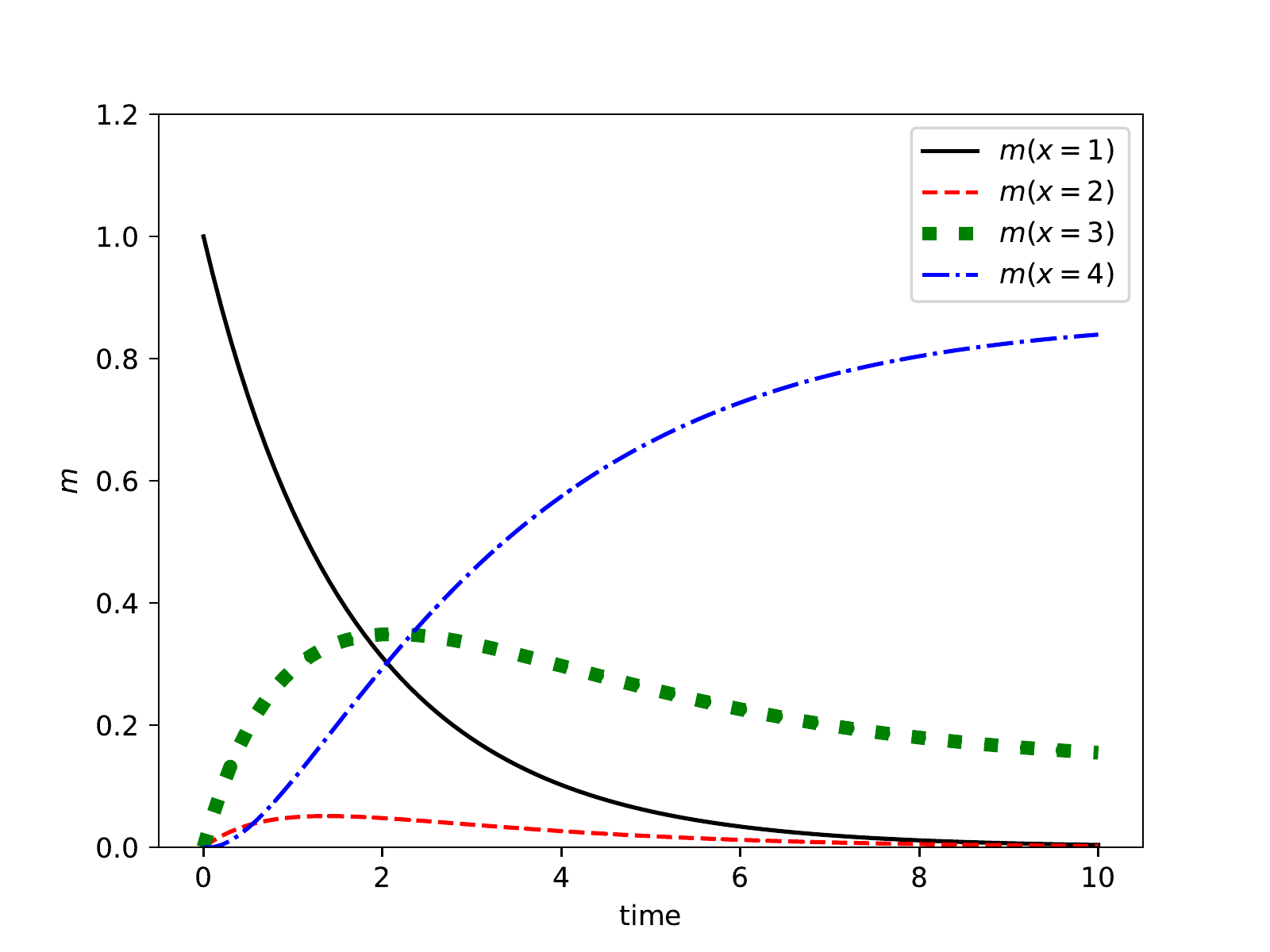}
	\end{subfigure}%
	\begin{subfigure}{.45\columnwidth}
		\centering 
		\includegraphics[width=\columnwidth]{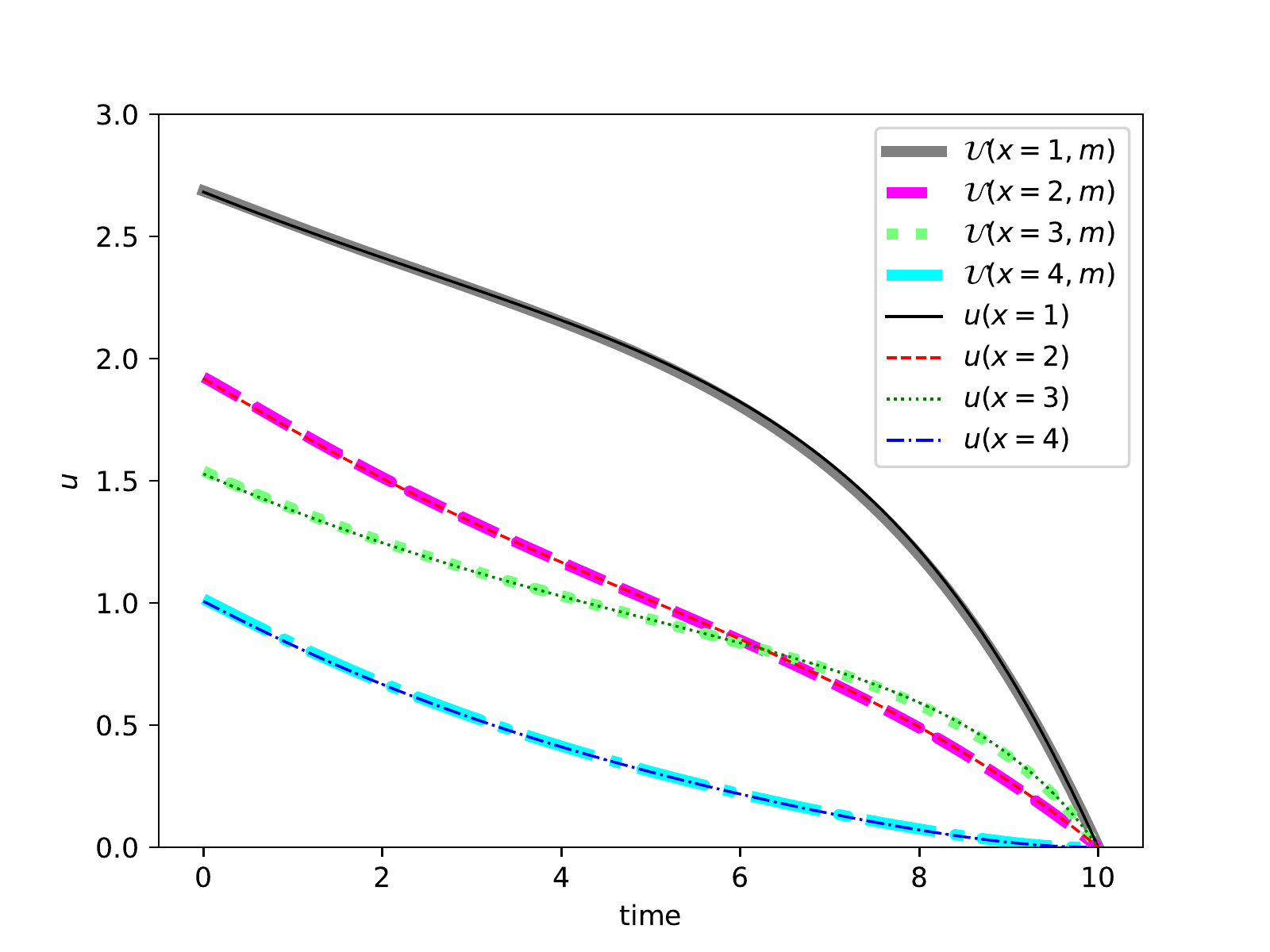}
	\end{subfigure}
	\caption{MFG Cyber-security example, test case 2: Evolution of the distribution $m^{m_0}$ (left) and the value function $u^{m_0}$ and $\cU(\cdot, \cdot, m^{m_0}(\cdot))$ (right) for  $m_0 = (1, 0, 0, 0)$. } 
 	\label{AMS-num-fig:finiteMFG-cyber-Master-m0-2}
\end{figure}

\begin{figure}[h]
	\begin{subfigure}{.45\columnwidth}
		\centering
		\includegraphics[width=\columnwidth]{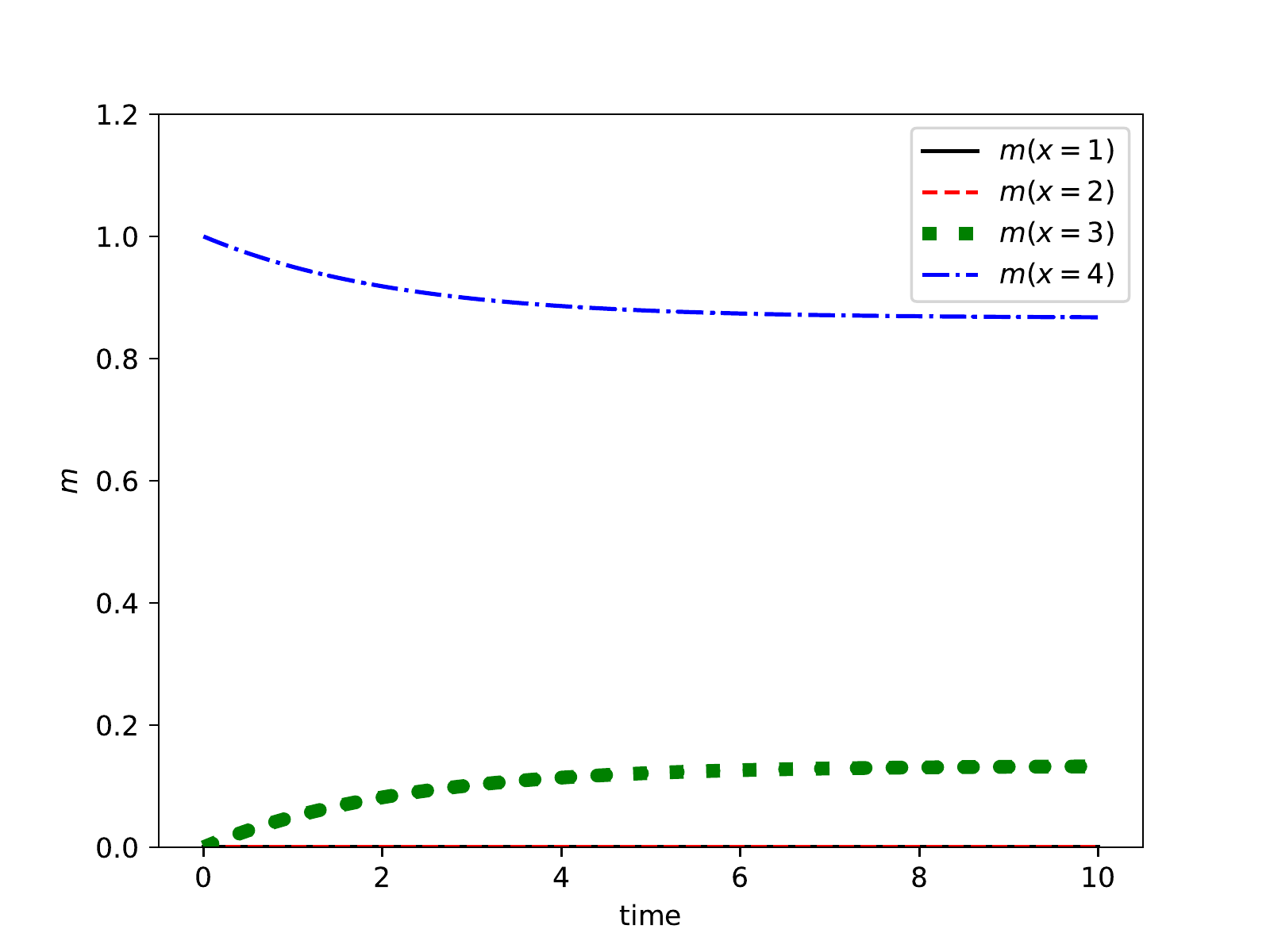}
	\end{subfigure}%
	\begin{subfigure}{.45\columnwidth}
		\centering 
		\includegraphics[width=\columnwidth]{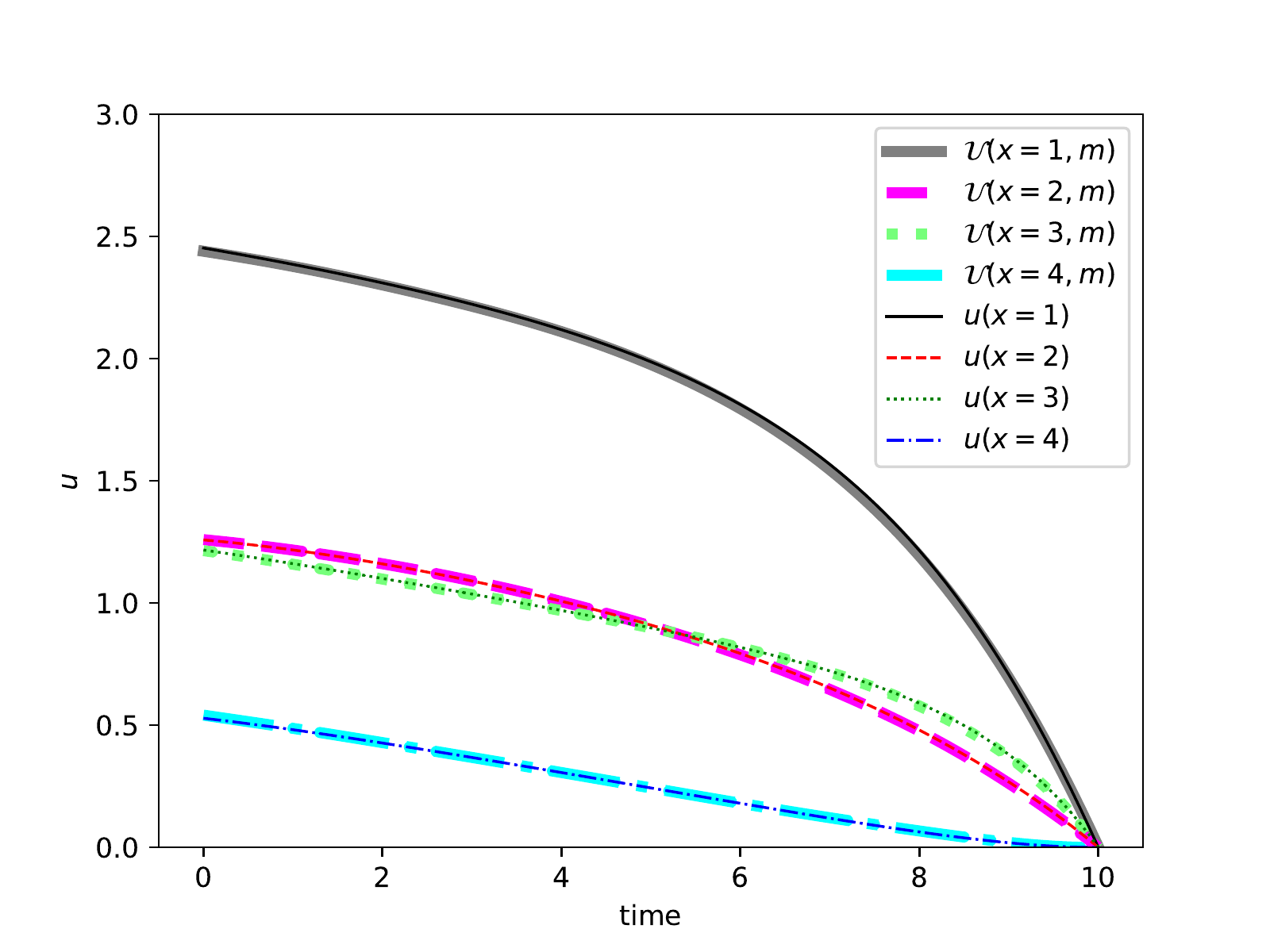}
	\end{subfigure}
	\caption{MFG Cyber-security example, test case 3: Evolution of the distribution $m^{m_0}$ (left) and the value function $u^{m_0}$ and $\cU(\cdot, \cdot, m^{m_0}(\cdot))$ (right) for  $m_0 = (0, 0, 0, 1)$. } 
 	\label{AMS-num-fig:finiteMFG-cyber-Master-m0-3}
\end{figure}

We apply the DGM method described above to the master equation~\eqref{AMS-num-eq:master-finiteMFG} in this cyber-security example. We obtain a neural network $\cU_\theta$ which is an approximation of $\cU$. For the sake of comparison, we (separately) solve the forward-backward system~\eqref{AMS-num-eq:ODE-system-finiteMFG} for various initial distributions $m_0$ and obtain solutions $(u^{m_0}, m^{m_0})$. We then compare $\cU_\theta(t,x,m^{m_0}(t,\cdot))$ and $u^{m_0}(t,x)$ which gives two curves for each $x \in \cE$. According to the relation~\eqref{AMS-num-eq:finiteMFG-link-U-u}, we know that these two curves should coincide, and this is verified in our numerical experiments, see Figures~\ref{AMS-num-fig:finiteMFG-cyber-Master-m0-1}--\ref{AMS-num-fig:finiteMFG-cyber-Master-m0-3} where we consider three test cases corresponding to three different initial conditions.  Here we used the following values for the parameters:
\begin{equation*}
\left\{
\begin{split}
&\beta_{UU} = 0.3,
\beta_{UD} = 0.4,
\beta_{DU} = 0.3,
\beta_{DD} = 0.4,
\\
&v_H = 0.2,
\lambda = 0.5,
\\
&q_{rec}^D = 0.1, 
q_{rec}^U = 0.65, 
q_{inf}^D = 0.4, 
q_{inf}^U = 0.3, 
\\
&k_D = 0.3, 
k_I = 0.5.
\end{split}
\right.
\end{equation*}
In the implementation, we used a feedforward fully connected neural network, as introduced in~\eqref{eq:NNarchi-FFC}, with sigmoid activation function. For problems in higher dimension, other architectures are sometimes more suitable.

\textbf{Example 2: An Example without uniqueness. } In the above example, the solution to the master equation is very smooth. We now turn our attention to an example proposed by Cecchin \textit{et al.} in~\cite{MR3981375} for which the master equation has multiple solutions but there is only one satisfying an entropy condition. Numerical results show that the neural network manages to approximate this latter solution despite its lack of smoothness.

The state space is $\cE = \{-1, 1\}$. Since there are only two states, every element of $\cP(\cE)$ can be characterized by its mean $\bar{m} = m(1) - m(-1)$. The transition rate is the control. In other words, at each time $t$, given its state $X_t$, an infinitesimal player can choose the rate $\ctrl(t, X_t) \in [0,+\infty)$ at which she wants to flip her state. The running and terminal costs are given by:
$$
	f(x, m, a) = \frac{a^2}{2}, \qquad g(x,m) = - x \bar{m}.
$$

The master equation takes the following form:  For $(t,x,\bar{m}) \in [0,T] \times \{-1,1\} \times [-1,1]$,
\begin{align*}
	&-\partial_t \cU(t,x,\bar{m})
	+
	\frac{1}{2} \left[ (\Delta^x \cU(t,x,\bar{m}))^- \right]^2 
	\\
	&\qquad - \sum_{y \in\{-1,1\}} \frac{1+ y \bar{m}}{2} D^{\bar{m}} \cU(t,x,\bar{m}, y) (\Delta^x \cU(t, y, \bar{m}))^- 
	= 0,
\end{align*}
with the terminal condition $\cU(T,x,\bar{m}) = - x \bar{m}$, for $(x,\bar{m}) \in \{-1,1\} \times [-1,1]$. Here $\Delta^x$ denotes the first difference, defined as:
$$
	\Delta^x \cU(t,x,\bar{m}) = \cU(t,-x,\bar{m}) - \cU(t,x,\bar{m}),
$$
and the derivative with respect to $\bar{m}$ is given by:
$$
	D^{\bar{m}} \cU(t,x,\bar{m},1) = 2 \partial_{\bar{m}} \cU(t,x,\bar{m}) = -D^{\bar{m}} \cU(t,x,\bar{m},-1),
$$ 
where $\partial_{\bar{m}}$ denotes a partial derivative with respect to the real variable $\bar{m}$ in the usual sense.

Introducing the new variable:
\begin{equation}
\label{AMS-num-eq:entropy-U-to-Z}
	\cZ(t,\bar{m}) = \cU(T-t,-1,\bar{m}) - \cU(T-t,1,\bar{m}), \qquad (t, \bar m) \in [0,T] \times [-1,1],
\end{equation}
we can check that $\cZ$ solves: 
\begin{equation*}
\left\{
\begin{split}
	& \partial_t \cZ(t,\bar{m}) + \partial_{\bar{m}} \mathfrak{g}(\bar{m}, \cZ(t,\bar{m})) = 0,  &&(t,\bar{m}) \in [-1,1],
	\\
	& \cZ(0,\bar{m}) = \mathfrak{f}(\bar{m}), &&\bar{m} \in [-1,1],
\end{split}
\right.
\end{equation*}
where
$$
	\mathfrak{f}(\bar{m}) = 2 \bar{m},
	\quad \mathfrak{g}(\bar{m}, z) = \bar{m} \frac{z |z|}{2} - \frac{z^2}{2},
	\qquad (\bar m, z) \in [-1,1] \times \RR.
$$
This equation is a scalar conservation law admitting three solutions, among which only one is an entropy solution (see~\cite[Proposition 3]{MR3981375}), which is given by:
\begin{equation}
\label{AMS-num-eq:entropy-Z-expli}
	\cZ(t, \bar{m}) = \frac{2 M(t,\bar{m})}{t |M(t, \bar{m})| + 1},
\end{equation}
where $M(t, 0) = 0$ and,  if $\bar{m} \neq 0$, $M(t, \bar{m})$ denotes the unique solution to: 
$$
	t^2M^3 + t(2-t)M |M| + (1-2t)M - \bar{m} = 0, 
$$
with the same sign as $\bar{m}$.

\begin{figure}[h]
	\begin{subfigure}{.45\columnwidth}
		\centering
		\includegraphics[width=\columnwidth]{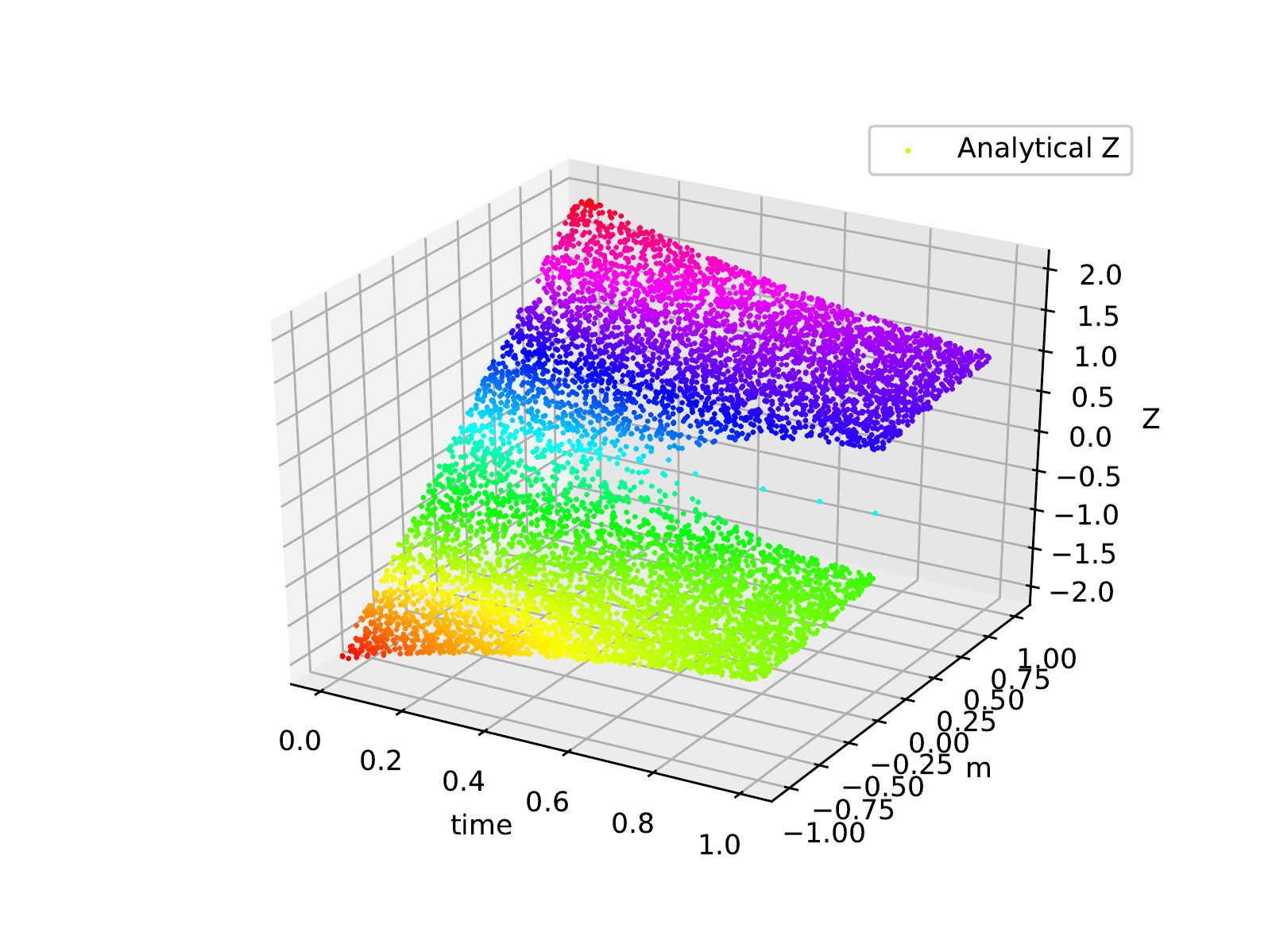}
	\end{subfigure}%
	\begin{subfigure}{.45\columnwidth}
		\centering 
		\includegraphics[width=\columnwidth]{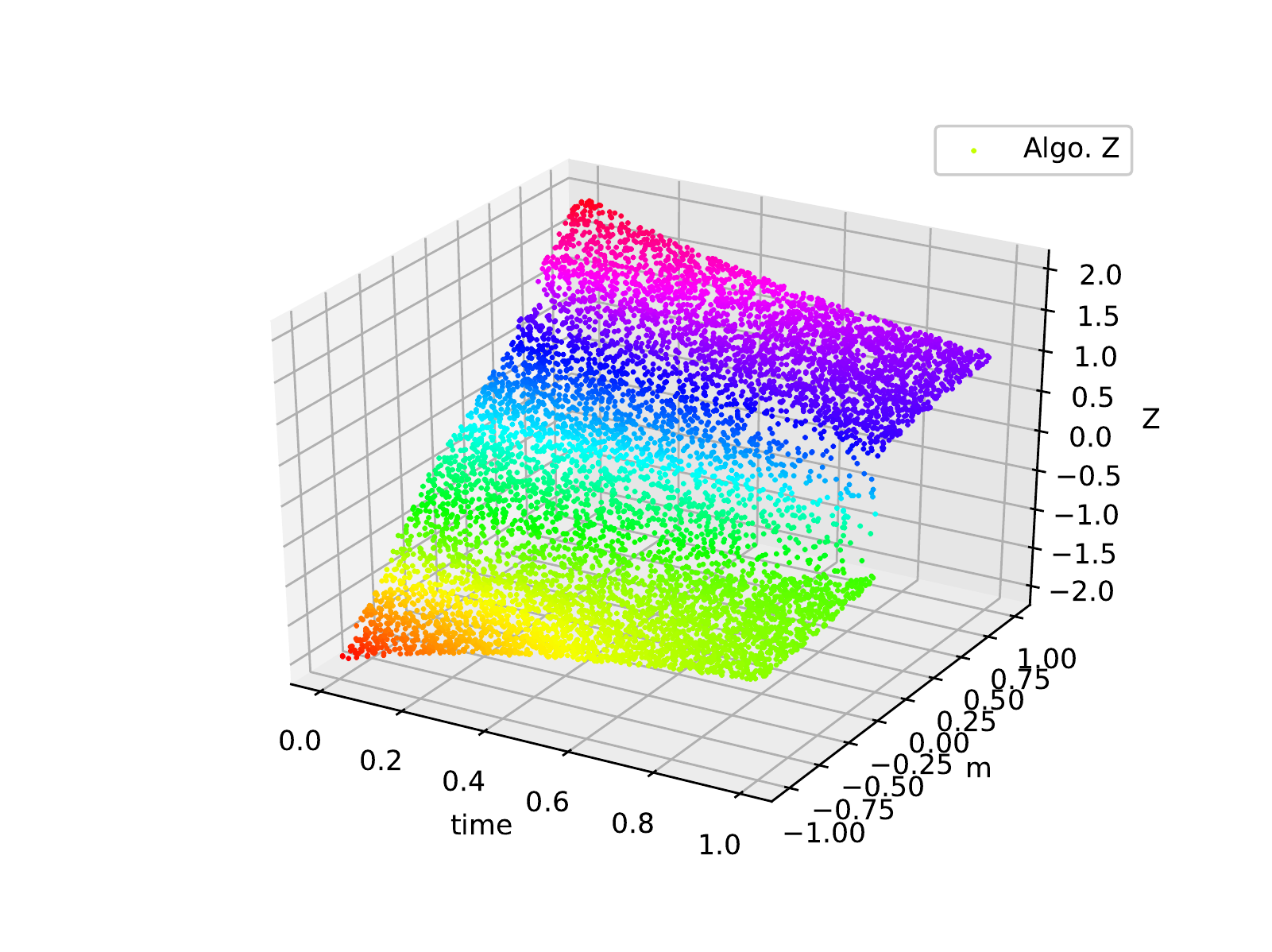}
	\end{subfigure}
	\caption{Example with non-uniqueness: $\cZ$ from analytical formula~\eqref{AMS-num-eq:entropy-Z-expli} (left) and from Deep Galerkin Method applied to compute $\cU$, from which $\cZ$ is deduced by~\eqref{AMS-num-eq:entropy-U-to-Z} (right). \label{AMS-num-fig:finiteMFG-entropy-3D} } 
\end{figure}

\begin{figure}[h]
	\begin{subfigure}{.45\columnwidth}
		\centering
		\includegraphics[width=\columnwidth]{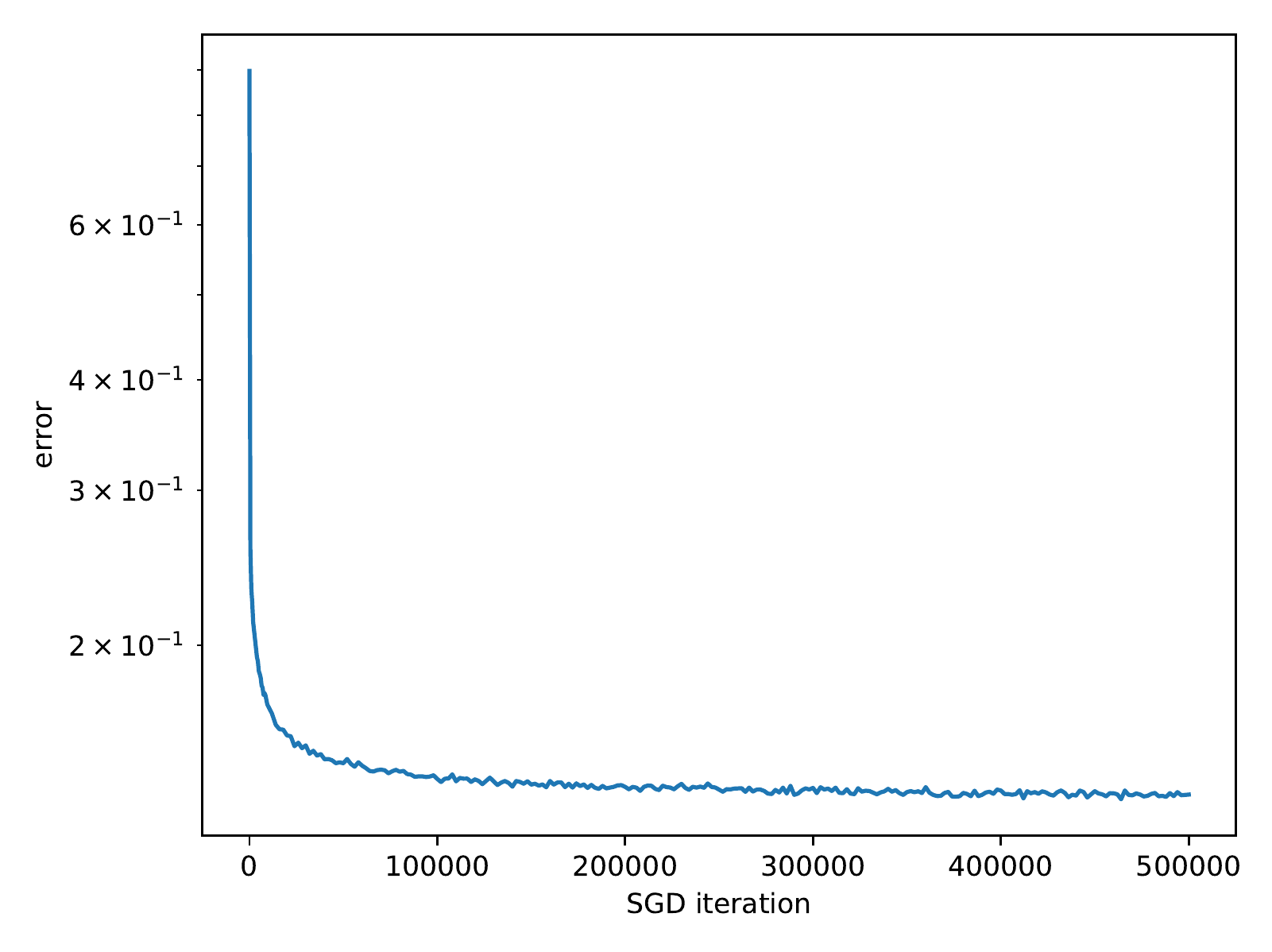}
	\end{subfigure}%
	\begin{subfigure}{.45\columnwidth}
		\centering 
		\includegraphics[width=\columnwidth]{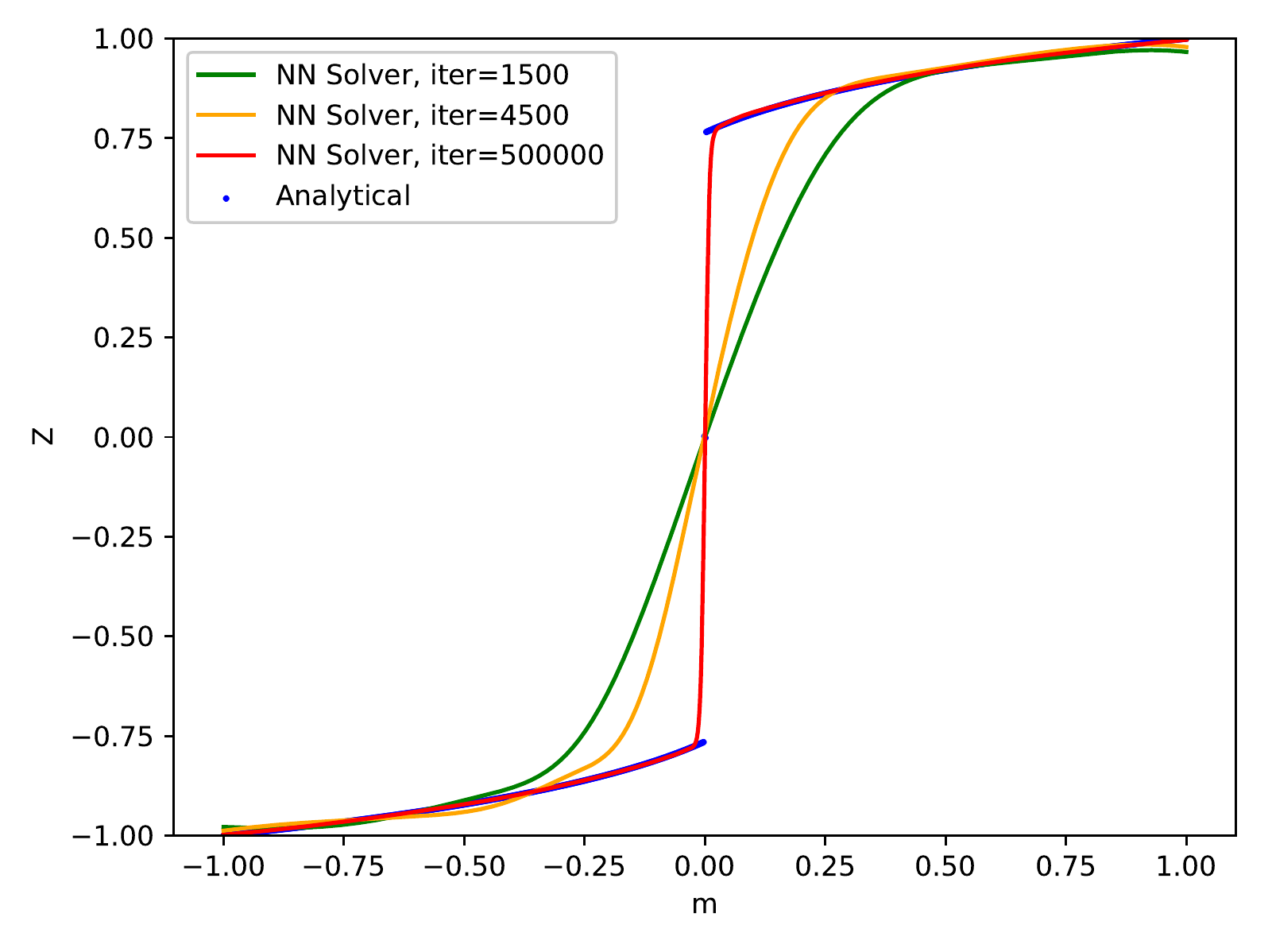}
	\end{subfigure}
	\caption{Example with non-uniqueness: error versus number of SGD iterations (left), and $\cZ(T,\cdot)$  (right), where the true solution is in blue, and the solution learnt by the neural network is in green, yellow and red for $1500$, $4500$ and $5\times10^5$ iterations respectively.
 	\label{AMS-num-fig:finiteMFG-entropy-2D}} 
\end{figure}

Figure~\ref{AMS-num-fig:finiteMFG-entropy-3D} shows the true $\cZ$ given by~\eqref{AMS-num-eq:entropy-Z-expli} and the one obtained by the change of variable~\eqref{AMS-num-eq:entropy-U-to-Z} after the neural network has been trained to approximate $\cU$. We see that the true  $\cZ$ exhibits a discontinuity at $\bar{m}$ after time $t = 0.5$, whereas the neural network is continuous because we used the sigmoid function as an activation function. However, Figure~\ref{AMS-num-fig:finiteMFG-entropy-2D} shows that the $L^2$ error decays with the number of SGD iterations, and after a large enough number of iterations, the neural network manages to approximate the discontinuity as shown for the terminal time $T=1$. 

To conclude this section, let us mention that the DGM method can also be used to tackle the system of HJB-KFP PDEs characterizing the MFG or MFC solutions; see ~\cite{al2018solving,carmona2019convergence1} for more details.

\clearpage

\section{\bf A Glance at model-free methods}
\label{AMS-num-sec:modelfree}

All the previous methods rely, in one way or another, on the fact that the cost functions $f$ and $g$ as well as the drift $b$ and the volatility $\sigma$ are known. However, in many applications, coming up with a realistic and accurate model is a daunting task. It is sometimes impossible to guess the form of the dynamics, or the way the costs are incurred. This provides a motivation to study so-called model-free methods. The theory of reinforcement learning (RL) has formalized this framework and numerous algorithms have been developed. Intuitively, an agent evolving in an environment can take actions and observe the consequences of her actions: the state of the environment (or her own state) changes, and a cost is incurred to the agent. The agent does not know how the new state and the cost are computed. The goal for the agent is then to learn an optimal behavior (\textit{i.e.}, which minimizes the sum of future costs) by trial and error. The problem is even more complex if multiple agents try to learn simultaneously and their actions influence each other's costs. We briefly describe some recent progress on the connection between reinforcement learning and mean field problems. We start with MFC, which, as a control problem, can be recast as a RL problem. We then consider MFG, which requires learning a Nash equilibrium.

\subsection{Reinforcement learning for mean field control}

As discussed in Section~\ref{AMS-num-sec:optim-variational}, MFC problems can be viewed as optimal control problems driven by the evolution of a probability distribution. This corresponds to the viewpoint of a central planner trying to find the optimal control that the agents should use in order to minimize the social cost. To fit more easily in the framework of reinforcement learning, we revisit the setting with finite state and action spaces considered in~\S~\ref{AMS-num-sec:NN-finite-master}. Instead of the finite horizon setting, we consider here the infinite horizon discounted case. We modify the class of controls and instead of considering functions of $t$ (time) and $x$ (individual state), we consider here functions of $m$ (population distribution) and $x$ (individual state). Intuitively, the reason behind this choice is that in many cases, the optimal control depends on time only through the population's distribution. The MFC counterpart to the MFG is to find a feedback control $\ctrl^*: \cP(\cE) \times \cE \to \cA$ minimizing:
\begin{align*}
	J^{MFC}: \ctrl \mapsto  \EE \left[\int_0^{+\infty} e^{-\beta t} f(X_t^{\ctrl}, m^\ctrl(t,\cdot), \ctrl(m^\ctrl(t,\cdot),X_t^{\ctrl}) ) dt  \right]
\end{align*}
under the constraint that the  process $X^{\ctrl} = (X_t^{\ctrl})_{t \ge 0}$ is a nonhomogeneous $\cE$-valued Markov chain with transition probabilities determined by the $Q$-matrix of rates $q^{m^\ctrl, \ctrl}: [0,T] \times \cE \times \cE \to \RR$ given by~\eqref{AMS-num-eq:q-finite-MFG} and $m^\ctrl$ is the flow of distributions corresponding to $X^{\ctrl}$.

Notice that the cost can be rewritten as:
$$
	J^{MFC}(\ctrl) = \int_0^{+\infty} e^{-\beta t} \sum_{x\in\cE} f(x, m^\ctrl(t,\cdot), \ctrl(m^\ctrl(t,\cdot),x) ) m^\ctrl(t,x) dt ,
$$
subject to the constraint that $m^\ctrl$ solves an evolution equation analogous to~\eqref{AMS-num-eq:ODE-system-finiteMFG-KFP} but controlled by $\ctrl$:
$$
	\partial_t m^\ctrl(t,x) = \sum_{x' \in \cE} m^\ctrl(t,x') \lambda(x', x, m^\ctrl(t,\cdot), \ctrl(m^\ctrl(t,\cdot), x')), \qquad (t,x) \in [0,T] \times \cE.
$$

Discretizing time with a grid $\{t_n = n \Delta t, n = 0,1,2,\dots\}$ with $\Delta t>0$, the cost is approximately equal to:
\begin{align*}
	J^{MFC, \Delta t}(\ctrl) 
	&= \sum_{n=0}^{+\infty} e^{-\beta t_n} \sum_{x\in\cE} f(x, m^\ctrl_n, \ctrl(m^\ctrl_n,x) ) m^\ctrl_n(x) \Delta t 
	\\
	&=  \sum_{n=0}^{+\infty} \gamma^n \tilde f(m^\ctrl_n, \ctrl(m^\ctrl_n,\cdot)),
\end{align*}
subject to:
$$
	m^\ctrl_{n+1} 
	= (m^\ctrl_{n})^\top (I + P^{\ctrl(m^\ctrl_{n},\cdot), m^\ctrl_{n}} \Delta t)
	=: \Phi(m^\ctrl_{n}, \ctrl(m^\ctrl_{n},\cdot)), \qquad n = 0, 1, 2, \dots
$$
with a given initial condition $m^\ctrl_0$, where $\gamma = e^{-\beta \Delta t}$, $\tilde f: \cP(\cE) \times \cA^{\cE} \to \RR$ is defined by:
$$
	\tilde f(m, \tilde\ctrl) = \sum_{x\in\cE} f(x, m, \tilde\ctrl(x) ) m(x) \Delta t, \qquad (m, \tilde\ctrl) \in \cP(\cE) \times \cA^\cE, 
$$
and
$$
	P^{\tilde\ctrl, m}(x', x) = \lambda(x', x, m, \tilde\ctrl(x')), \qquad (x',x,m,\tilde\ctrl) \in \cE \times \cE \times \cP(\cE) \times \cA^\cE.
$$
The problem thus fits in the framework of Markov decision processes (MDP for short) by considering the distribution $m^\ctrl$ as the state (which is consistent with the point of view of the central planner). Instead of considering the dynamics of $X^\ctrl$ evolving in $\cE$, the MDP is for the ``lifted'' problem driven by the dynamics of $m^\ctrl$ evolving in $\cP(\cE)$.  
The value function $V$ associated to this optimal control problem represents the optimal (infinite horizon discounted) cost-to-go starting from a given state. Here, a state is a population distribution and hence the value function is a function of $m \in \cP(\cE)$:
$$
	V(m) = 
	\inf_{\ctrl : \cP(\cE) \times \cE \to \cA} \sum_{n=0}^{+\infty} \gamma^n \tilde f(m^\ctrl_n, \ctrl(m^\ctrl_n,\cdot)), \quad \hbox{ with } m^\ctrl_0 = m.
$$
 However, even if we manage to compute $V$, it is not clear how to recover the optimal control if the cost function and the dynamics are not known. From this perspective, a more useful function is the so-called $Q$-function. It is a function of a state and an action, and it represents the optimal cost-to-go given that we start with the prescribed state-action pair (and then we behave optimally). In our context, it takes the following form: 
 $$
 	Q(m,\tilde\ctrl) = \tilde f(m,\tilde\ctrl) + \gamma \inf_{\ctrl : \cP(\cE) \times \cE \to \cA} \sum_{n=0}^{+\infty} \gamma^n \tilde f(m^\ctrl_n, \ctrl(m^\ctrl_n,\cdot)), \quad (m, \tilde\ctrl) \in \cP(\cE) \times \cA^{\cE},
 $$
 where $m^\ctrl_0 = \Phi(m, \tilde\ctrl)$ is the distribution obtained by starting from $m$ and moving forward by one time step using action $\tilde\ctrl$, and for $n \ge 0$, $m^\ctrl_{n+1} = \Phi(m^\ctrl_{n}, \ctrl(m^\ctrl_{n},\cdot))$. We have
 \begin{equation}
 \label{eq:MFCQ-link-VQ}
 	V(m) = \inf_{\tilde\ctrl \in \cA^{\cE}} Q(m, \tilde\ctrl), \qquad m \in \cP(\cE),
 \end{equation}
 and the optimal control is readily recovered from $Q$ as: 
\begin{equation}
 \label{eq:MFCQ-link-ctrloptQ}
 	\ctrl^*(m,\cdot) = \arginf_{\tilde\ctrl \in \cA^{\cE}} Q(m, \tilde\ctrl).
\end{equation}
 This is one of the main advantages of the state-action value function $Q$ over the (state only) value function $V$. 
 Furthermore, it can be shown that $Q$ satisfies a dynamic programming equation (or Bellman equation):
$$
	 Q(m,\tilde\ctrl) = \tilde f(m,\tilde\ctrl) + \gamma \inf_{\tilde\ctrl' \in \cA^{\cE}} Q(\Phi(m, \tilde\ctrl), \tilde\ctrl'), \quad (m,\tilde\ctrl) \in \cP(\cE) \times \cA^{\cE}.
$$
Thus, a straightforward strategy to compute $Q$ is to iteratively update an approximation by plugging in the right hand side above the previous estimate. This leads to the so-called (synchronous) $Q$-learning updates:
\begin{equation}
\label{eq:MFC-Q-updates}
	 Q^{(\mathtt{k+1})}(m,\tilde\ctrl) = \tilde f(m,\tilde\ctrl) + \gamma \inf_{\tilde\ctrl' \in \cA^{\cE}} Q^{(\mathtt{k})}(\Phi(m, \tilde\ctrl), \tilde\ctrl'), \quad (m,\tilde\ctrl) \in \cP(\cE) \times \cA^{\cE}.
\end{equation}
Notice that in order to perform this kind of updates, all we need is to have access to the cost $\tilde f(m,\tilde\ctrl)$ and to the next state $\Phi(m, \tilde\ctrl)$, which can be provided by some environment whose inner functioning is not revealed to the learning agent. In other words, to perform the $Q$-learning updates, we do not need to know the functions $\tilde f$ and $\Phi$. In this sense, the algorithm is model-free.

To implement this method, we need a finite dimensional approximation of the $Q$ function. We do not discuss here function approximation techniques. For simplicity, we replace $\cP(\cE)$ by a discrete version, denoted by $\tilde\cP_\cE$, which contains a finite number of points of $\cP(\cE)$. Since $\cA^\cE$ is also a finite set, the $Q$-function can be approximated by its values on $\tilde\cP_\cE \times \cA^\cE$, \textit{i.e.}, by a matrix $\tilde Q \in \RR^{|\tilde\cP_\cE| \times |\cA^\cE|}$. The update rule become:
\begin{equation}
\label{eq:MFC-Q-updates-proj}
	 \tilde Q^{(\mathtt{k+1})}(\tilde m,\tilde\ctrl) = \tilde f(\tilde m,\tilde\ctrl) + \gamma \inf_{\tilde\ctrl' \in \cA^{\cE}} \tilde Q^{(\mathtt{k})}(\Pi_{\tilde \cP_\cE} \Phi(\tilde m, \tilde\ctrl), \tilde\ctrl'), \quad (\tilde m, \tilde\ctrl) \in \tilde\cP(\cE) \times \cA^{\cE}.
\end{equation}
where $\Pi_{\tilde P_\cE}: \cP_\cE \to \tilde \cP_\cE$ denotes a projection on $\tilde\cP_\cE$ (note that $\Phi(\tilde m, \tilde\ctrl)$ is not necessarily an element of $\tilde P_\cE$). 

Note that in the current setting, in fact it is enough to interact once with the environment and query $(\tilde f(\tilde m,\tilde \ctrl), \Phi(\tilde m, \tilde \ctrl))$ for every pair $(\tilde m, \tilde \ctrl) \in \tilde \cP_\cE \times \cA^\cE$ in order to know everything we need about the cost and the dynamics to compute the $\tilde Q$ function. However, this is no longer true when the dynamics is stochastic, \textit{e.g.} due to the presence of common noise. 
 Moreover, when the number of states $|\tilde\cP_\cE|$ or the number of actions $|\cA^\cE|$ is large, updating $\tilde Q$ for every pair $(\tilde m,\tilde \ctrl) \in \tilde\cP_\cE \times \cA^\cE$ at each iteration is prohibitive. One can rely on asynchronous updates or function approximation to improve efficiency. 

\begin{figure}[h]
	\begin{subfigure}{.45\columnwidth}
		\centering
		\includegraphics[width=\columnwidth]{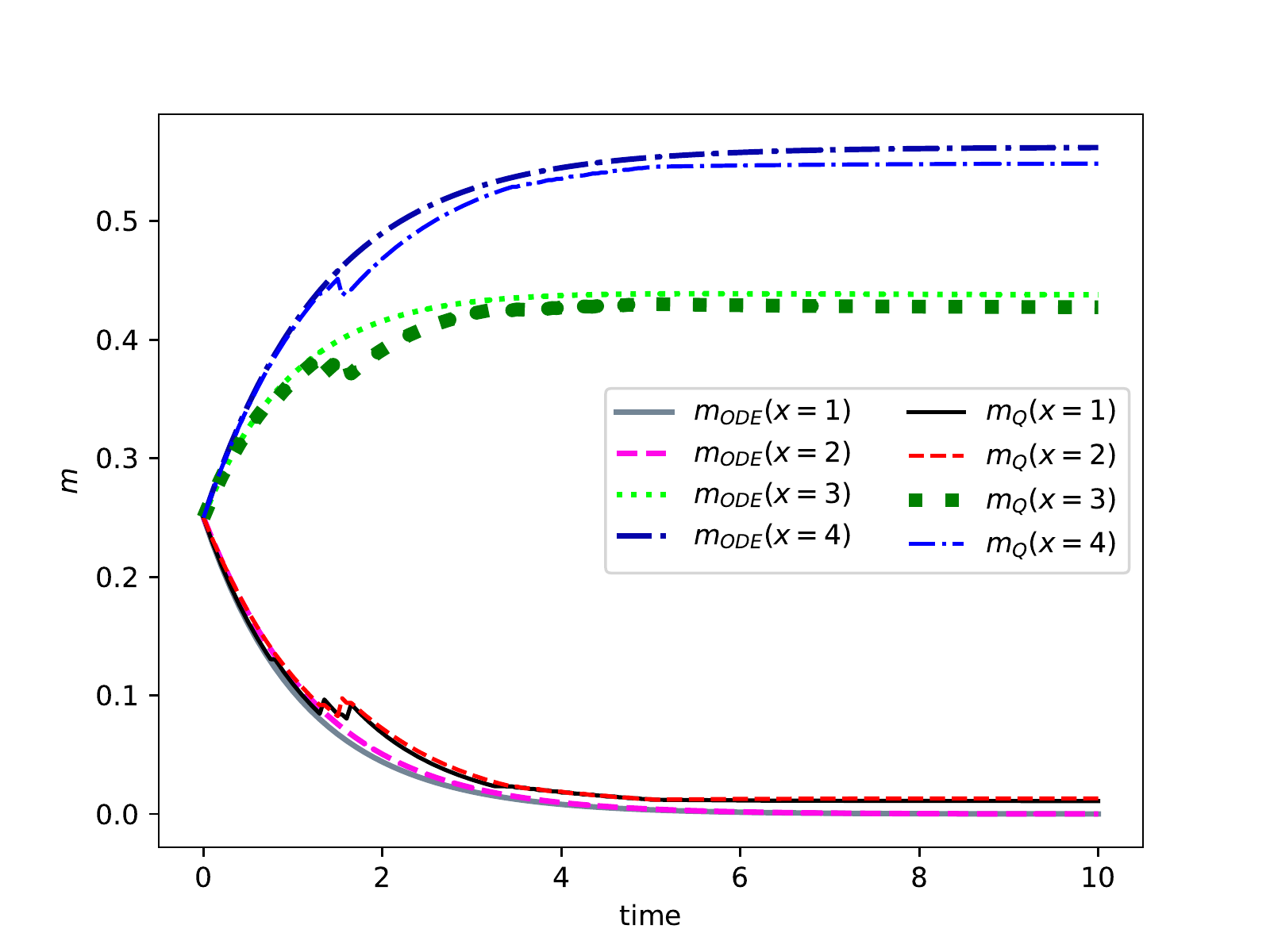}
	\end{subfigure}%
	\begin{subfigure}{.45\columnwidth}
		\centering 
		\includegraphics[width=\columnwidth]{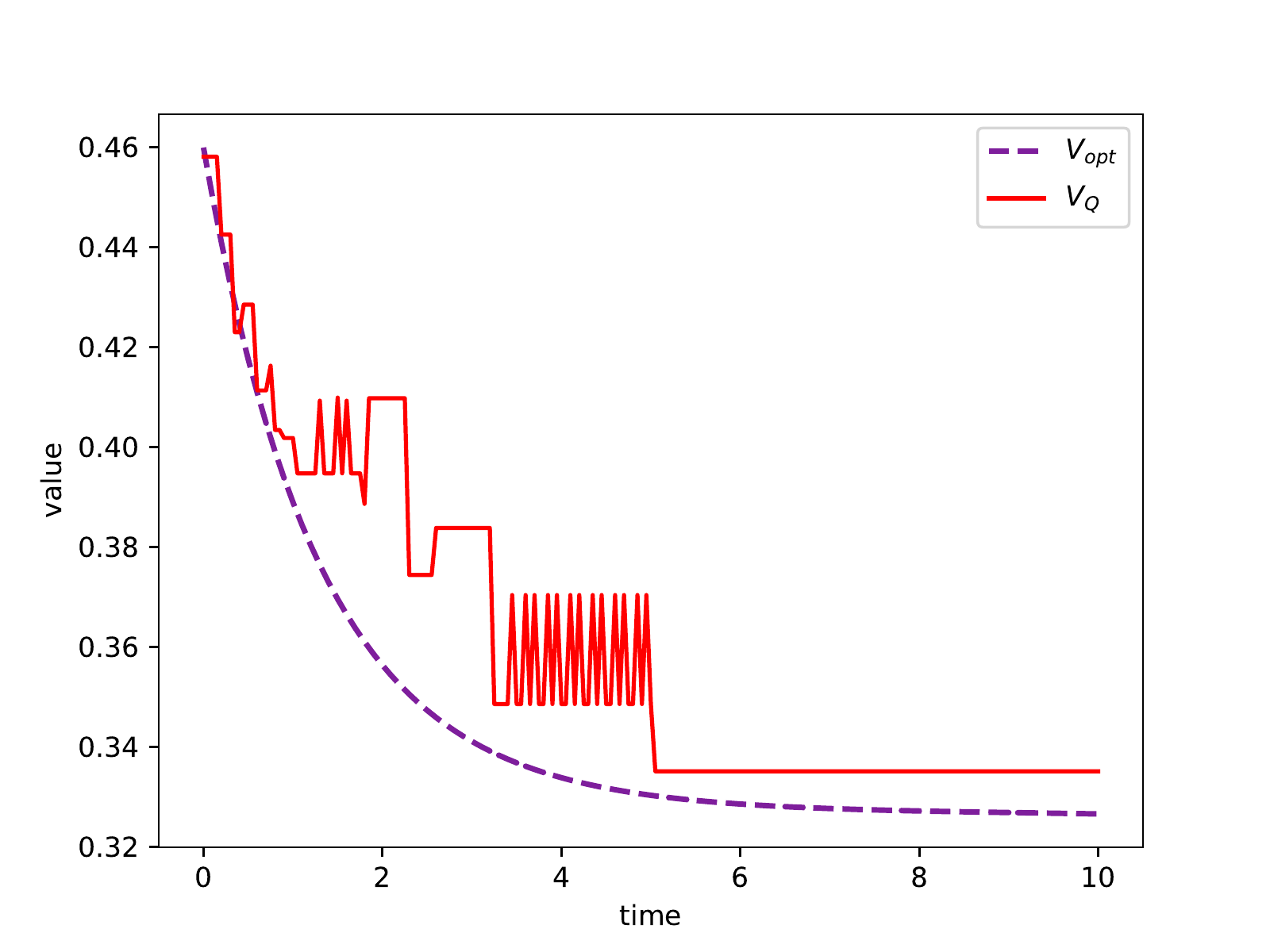}
	\end{subfigure}
	\caption{MFC Cyber-security example, test case 1: $m_0 = (1/4, 1/4, 1/4, 1/4)$. Left: Evolution of the distribution $m^{m_0}$ optimally controlled ($m_{ODE}$) or controlled using the approximate $Q$-function ($m_Q$). Right: $V$ function ($V_{opt}$) and approximate $Q$-function ($V_{Q}$) along the optimal flow. }    
 	\label{AMS-num-fig:finiteMFCQ-cyber-Master-m0-1}
\end{figure}

\begin{figure}[h]
	\begin{subfigure}{.45\columnwidth}
		\centering
		\includegraphics[width=\columnwidth]{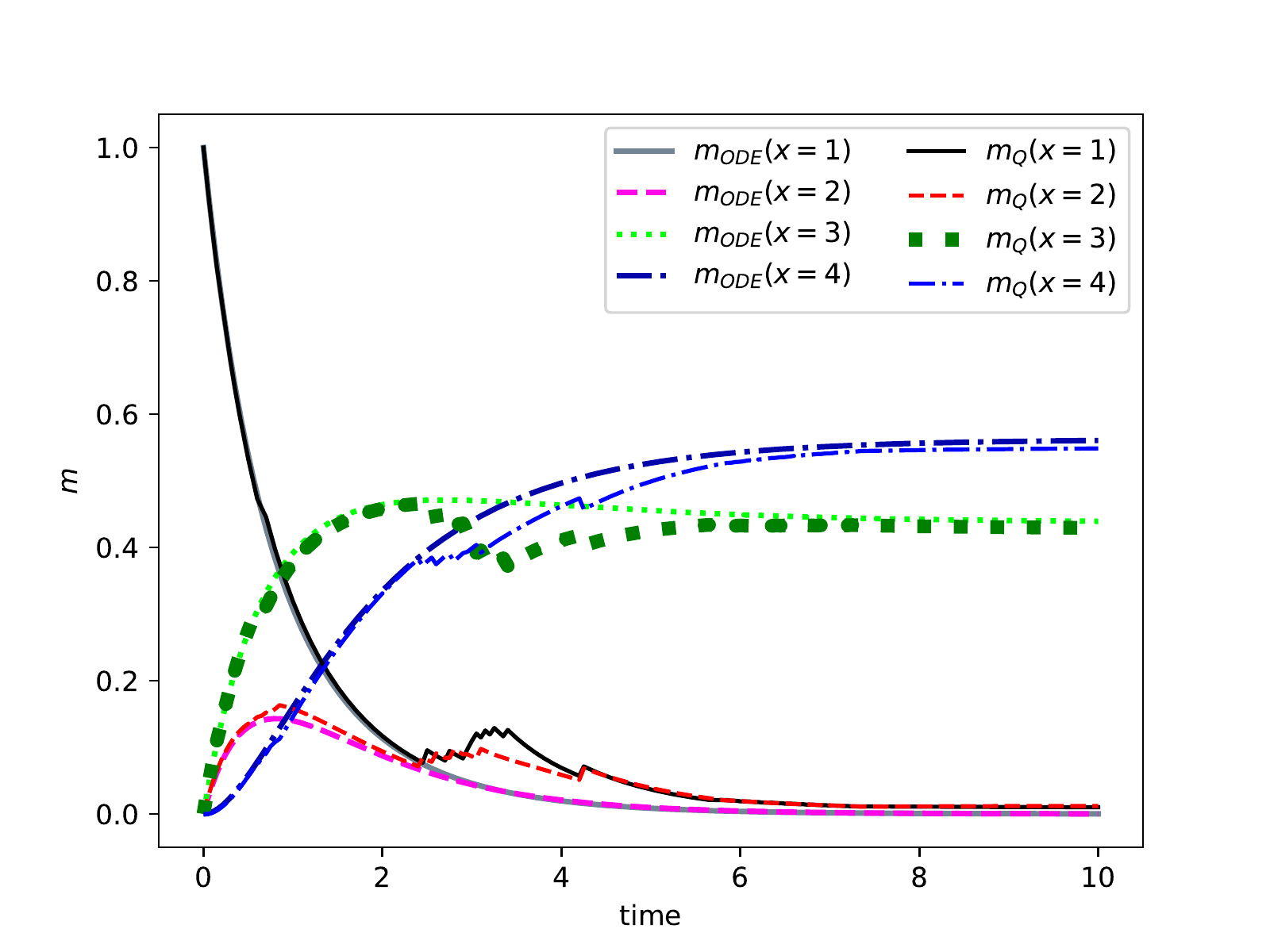}
	\end{subfigure}%
	\begin{subfigure}{.45\columnwidth}
		\centering 
		\includegraphics[width=\columnwidth]{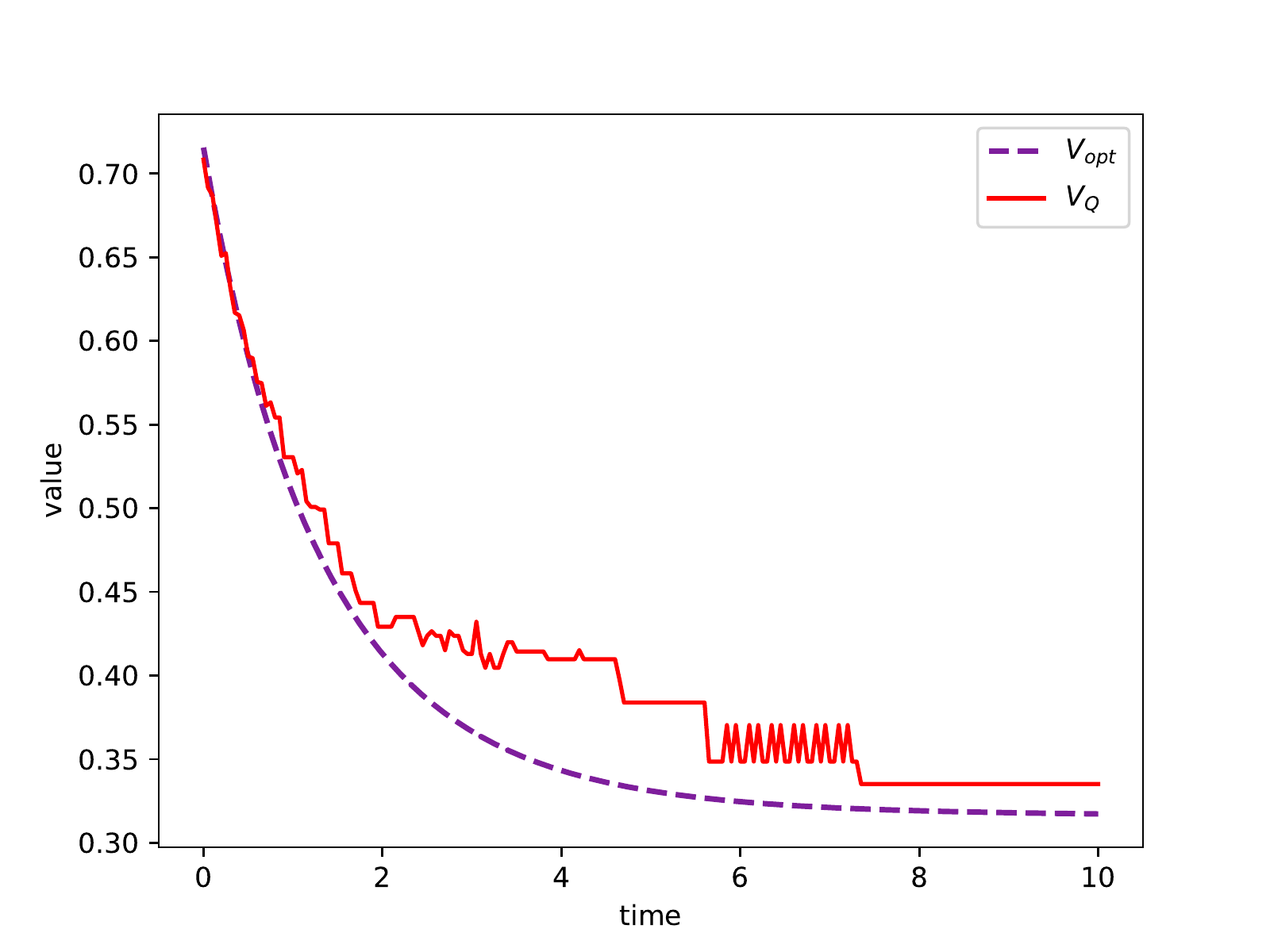}
	\end{subfigure}
	\caption{MFC Cyber-security example, test case 2: $m_0 = (1, 0, 0, 0)$. Left: Evolution of the distribution $m^{m_0}$ optimally controlled ($m_{ODE}$) or controlled using the approximate $Q$-function ($m_Q$). Right: $V$ function ($V_{opt}$) and approximate $Q$-function ($V_{Q}$) along the optimal flow. }   
 	\label{AMS-num-fig:finiteMFCQ-cyber-Master-m0-2}
\end{figure}

\begin{figure}[h]
	\begin{subfigure}{.45\columnwidth}
		\centering
		\includegraphics[width=\columnwidth]{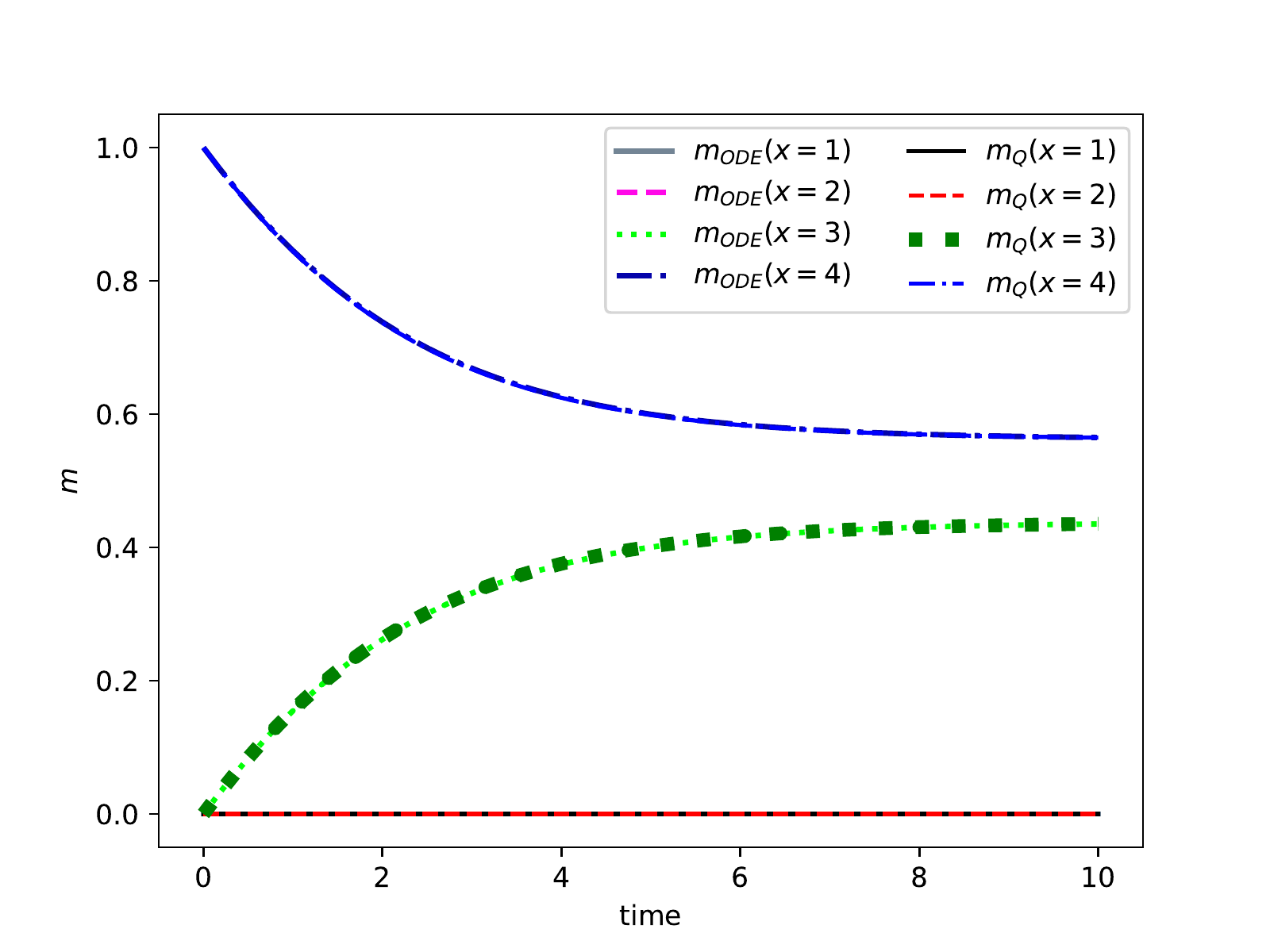}
	\end{subfigure}%
	\begin{subfigure}{.45\columnwidth}
		\centering 
		\includegraphics[width=\columnwidth]{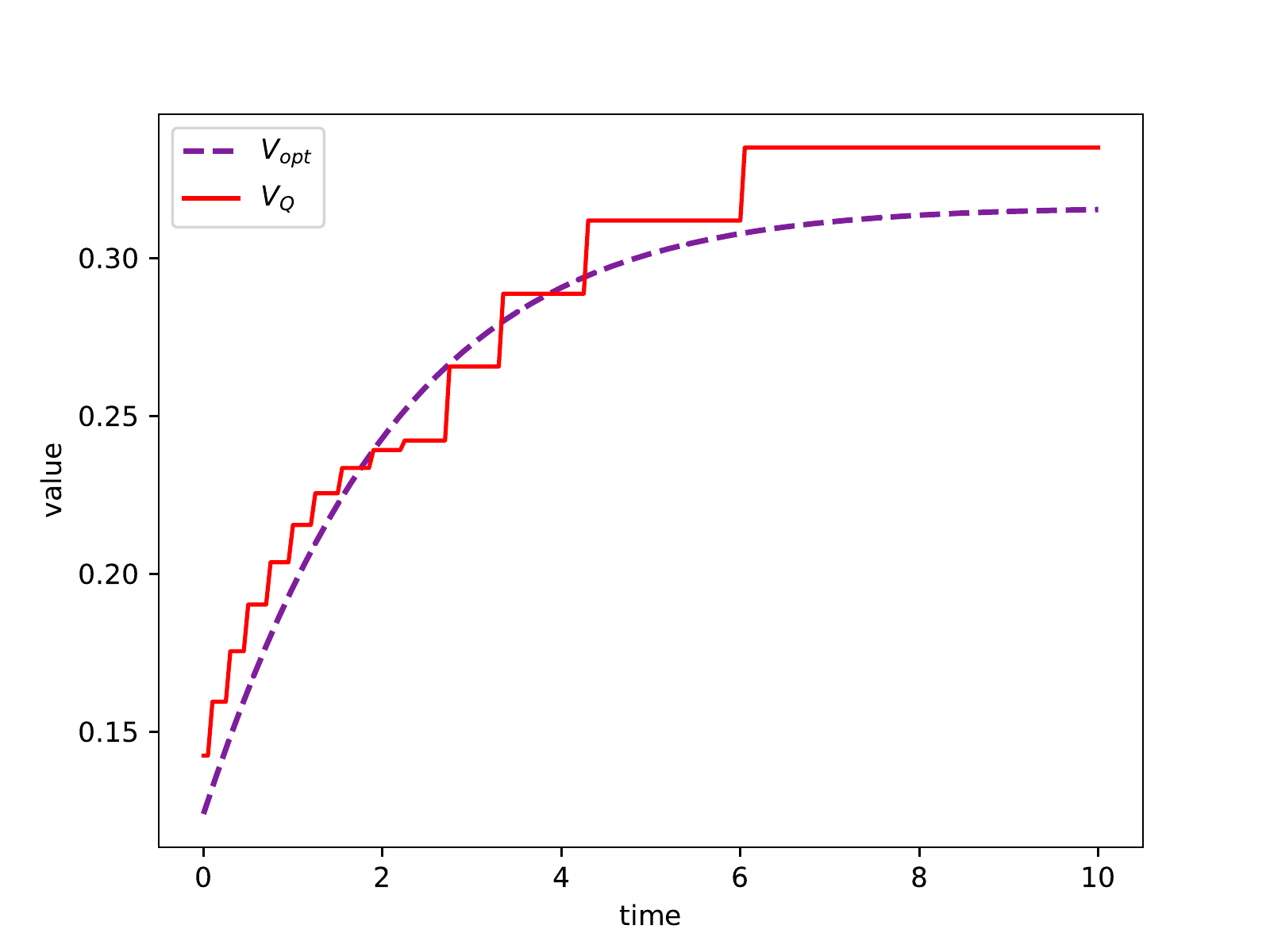}
	\end{subfigure}
	\caption{MFC Cyber-security example, test case 3: $m_0 = (0, 0, 0, 1)$. Left: Evolution of the distribution $m^{m_0}$ optimally controlled ($m_{ODE}$) or controlled using the approximate $Q$-function ($m_Q$). Right: $V$ function ($V_{opt}$) and approximate $Q$-function ($V_{Q}$) along the optimal flow. }   
 	\label{AMS-num-fig:finiteMFCQ-cyber-Master-m0-3}
\end{figure}

For the purpose of illustration, we present the behavior of the basic mean-field $Q$-learning method~\eqref{eq:MFC-Q-updates-proj}  on  
the cyber-security model presented in~\S~\ref{AMS-num-sec:NN-finite-master} but now from the point of view of a central planner controlling a large group of computers. Recall that $\cE = \{DI, DS, UI, US\}$ has four elements. 
We replace $\cP(\cE)$ by:
$$
	\tilde\cP_\cE = \left\{ (m_{DI}, m_{DS}, m_{UI}, m_{US} ) \in \cG_{N_m}^4 \,\big|\, m_{DI} + m_{DS} + m_{UI} + m_{US} = 1 \right\},
$$
where $\cG_{N_m}$ is a uniform grid over $[0,1]$ with $N_m \ge 2$ points, which contains $0$ and $1$. After $\mathtt{K}$ steps of $Q$-learning, we obtain an approximation $\tilde Q^{(\mathtt{K})}$ of the $Q$-function, from which we can recover an approximation $\tilde\ctrl^{\mathtt{K}}$ of the optimal control using~\eqref{eq:MFCQ-link-ctrloptQ},
namely: $\tilde\ctrl^{\mathtt{K}}(m, \cdot) = \argmax_{\tilde \ctrl \in \cA^\cE} \tilde Q^{(\mathtt{K})}(\Pi_{\tilde P_\cE} m, \tilde \ctrl)$. 
 For three different initial conditions, we compare the flow of distributions induced by this control $\tilde\ctrl^{\mathtt{K}}$ with the optimally controlled flow. In line with~\eqref{eq:MFCQ-link-VQ}, we also compare the value function $V(m^{\ctrl^*}(t))$ with $\max_{\cA^\cE}\tilde Q^{(\mathtt{K})}(\Pi_{\tilde P_\cE} m^{\ctrl^*}(t), \cdot)$ along the optimal flow $(m^{\ctrl^*}(t))_{t\in[0,T]}$. Figures~\ref{AMS-num-fig:finiteMFCQ-cyber-Master-m0-1}--\ref{AMS-num-fig:finiteMFCQ-cyber-Master-m0-3} show the results for three initial conditions. We see that the learnt $Q$-function approximately matches the $V$ value function.   
For these simulations, we used $N_m = 30$, $\gamma = 0.5$, and the following parameters:
\begin{equation*}
\left\{
\begin{split}
&\beta_{UU} = 0.3,
\beta_{UD} = 0.4,
\beta_{DU} = 0.3,
\beta_{DD} = 0.4,
\\
&v_H = 0.6,
\lambda = 0.8,
\\
&q_{rec}^D = 0.5, 
q_{rec}^U = 0.4, 
q_{inf}^D = 0.4, 
q_{inf}^U = 0.3, 
\\
&k_D = 0.3, 
k_I = 0.5.
\end{split}
\right.
\end{equation*}

As already mentioned, the update rule~\eqref{eq:MFC-Q-updates-proj} is provided here simply for the purpose of explaining the basic idea of $Q$-learning in a mean-field setup. 
 We refer to~\cite{MR2968782,carmona2019model,gu2019dynamic,gu2020qlearning,motte2019mean} for more details on mean-field MDPs and $Q$-learning for mean field control. Other methods have also been investigated, such as policy gradient, see~\cite{subramanianpolicy}, which can be proved to converge with a linear rate for linear-quadratic MFC problems, see~\cite{carmona2019linear}.

\subsection{Reinforcement learning in mean field games}

Since MFGs are fixed-point problems which can not, in general, be written as optimization problems, Nash equilibria can not be immediately recast as MDPs. When the model is not known, it is not clear how to use optimality conditions such as PDE of FBSDE systems. Instead, coming back to the definition~\eqref{AMS-num-eq:def-J-MFG}--\eqref{AMS-num-eq:dyn-X-general-MFG} suggests to iteratively update the distribution and, at each iteration, solve the optimal control problem faced by a typical player using reinforcement learning. This corresponds to the following scenario: an infinitesimal agent can observe a given population distribution and aims at learning her best response by trial and error (\textit{i.e.}, without knowing how her cost or her evolution, which depend on the mean-field term, are computed). This step can be recast as a (standard) MDP parameterized by the mean-field term. The key difference with MFC discussed above, is that the population's distribution does not change while the infinitesimal player learns her optimal control. Then, the mean-field term is updated using this (approximate) best response, and so on. Once again, the update of the distribution can be for instance of Picard type or of fictitious play type. We refer respectively to~\cite{guo2019learning} and~\cite{Elie2020OnTC,perrin2020fictitious} for more details on these methods; see also~\cite{subramanianpolicy} for a different method based on a two-timescale approach. To conclude this section, let us mention that despite the crucial differences between MFG and MFC, reinforcement learning algorithms for these problems can be unified through the aforementioned two-timescale approach: in fact, solutions to MFG and MFC can be computed by learning simultaneously the value function and the distribution, with different learning rates; see~\cite{angiulifouquelauriere2020unified} for more details. For the sake of brevity we refrain from developing here RL for MFGs.

\clearpage

\section{\bf Conclusion and perspectives}
\label{sec:conclusion}

We have presented several aspects of numerical methods for mean field games and mean field control problems. Besides the approaches presented above, several other directions have been investigated, such as (to cite just a few):
\begin{itemize}
	\item methods for forward-backward system of McKean-Vlasov SDEs: Chaudru de Raynal and Garcia Trillos have proposed in~\cite{MR3322862} a cubature method for FBSDE of MKV type when the system is decoupled; Chassagneux \textit{et al.} have proposed in~\cite{MR3914553,Angiulietal-2019} a recursive procedure, while other methods based on neural network approximation are studied in~\cite{fouque2019deep,carmona2019convergence2,germain2019numerical} by reformulating the FBSDE system as an optimal control of two forward SDEs where the cost function is a penalty for not matching the terminal condition;
	\item algorithms for ergodic MFGs: in~\cite{MR3601001,MR3882530}, Cacace \textit{et al.} proposed an approach relying on a finite difference scheme and a least square formulation, which is then solved using Gauss-Newton iterations; 
	\item gradient-type methods for MFC and variational MFG: see~\cite{LachapelleWolfram-2011-MFG-congestion-aversion, MR3501391,MR3619691} for applications to crowd motion or systemic risk;
	\item computational methods aiming specifically at MFGs with a non-local dependence on the population distribution, for instance using Fourier expansion techniques~\cite{nurbekyan2019fourier};
	\item methods based on machine learning tools besides the ones discussed above: for instance methods based on neural networks for the PDE system of MFG have been investigated in~\cite{al2018solving,carmona2019convergence1}, where the value function and the distribution are replaced by two neural networks and the optimization follows the DGM strategy; the link between this idea and Generative Adversarial Networks (GANs) is investigated in~\cite{cao2020connecting}; on the other hand,~\cite{ruthotto2020machine} proposes a method based on neural network for MFC and variational MFG;
	\item learning or reinforcement learning methods besides the ones discussed above, such as policy gradient~\cite{subramanianpolicy}, sequential decomposition~\cite{mishra2020model}, actor-critic method~\cite{fu2019actor}, online mirror descent~\cite{Hadikhanloo-phdthesis}, value iteration~\cite{anahtarci2019value} or policy iteration~\cite{cacacecamilligoffi2020policy}; see also~\cite{tembine2012meanfieldlearning};
	\item purely probabilistic methods such as quantization or regression based methods for a class of MFC~\cite{balata2019class}, or the approach based Markov chain approximation for MFG studied by Bayraktar \textit{et al.} in~\cite{MR3873029}.
\end{itemize}

For the sake of the presentation, the framework used in these notes is often much more restrictive and much simpler than what is covered in the literature. The interested reader is referred to the references cited throughout the text. However, some aspects such as general boundary conditions or general diffusion terms (which could depend on the control and be degenerate) remain challenging and will hopefully be addressed in future works.   Besides the development of numerical methods per se, we believe that these tools will unlock many applications for instance in economics, robotics, biology or epidemiology.

\vskip 6pt
\noindent \textbf{Acknowledgements. } M. Lauri{\`e}re is very grateful to Yves Achdou, Luis Brice\~{n}o-Arias, Ren\'e Carmona, Fran{\c c}ois Delarue, Francisco Silva and Xianjin Yang for fruitful discussions, and to the participants of the AMS Short Course for their questions and comments. M. Lauri{\`e}re would like to thank the anonymous referee for their insightful comments and suggestions that helped to improve the quality of the manuscript. 
The research of M. Lauri\`ere is supported by NSF grant DMS--1716673 and ARO grant W911NF--17--1--0578.

\bibliographystyle{amsplain}
\bibliography{mfg-num-bib2}

\providecommand{\bysame}{\leavevmode\hbox to3em{\hrulefill}\thinspace}
\providecommand{\MR}{\relax\ifhmode\unskip\space\fi MR }
\providecommand{\MRhref}[2]{%
  \href{http://www.ams.org/mathscinet-getitem?mr=#1}{#2}
}
\providecommand{\href}[2]{#2}
\begin{thebibliography}{100}

\bibitem{MR3135339}
Yves Achdou, \emph{Finite difference methods for mean field games},
  Hamilton-{J}acobi equations: approximations, numerical analysis and
  applications, Lecture Notes in Math., vol. 2074, Springer, Heidelberg, 2013,
  pp.~1--47. \MR{3135339}

\bibitem{MR2888257}
Yves Achdou, Fabio Camilli, and Italo Capuzzo-Dolcetta, \emph{Mean field games:
  numerical methods for the planning problem}, SIAM J. Control Optim.
  \textbf{50} (2012), no.~1, 77--109. \MR{2888257}

\bibitem{MR3097034}
\bysame, \emph{Mean field games: convergence of a finite difference method},
  SIAM J. Numer. Anal. \textbf{51} (2013), no.~5, 2585--2612. \MR{3097034}

\bibitem{MR2679575}
Yves Achdou and Italo Capuzzo-Dolcetta, \emph{Mean field games: numerical
  methods}, SIAM J. Numer. Anal. \textbf{48} (2010), no.~3, 1136--1162.
  \MR{2679575}

\bibitem{achdouKobeissi2020mfgcfinitediff}
Yves Achdou and Ziad Kobeissi, \emph{Mean field games of controls: Finite
  difference approximations}, arXiv preprint arXiv:2003.03968 (2020).

\bibitem{YAJML}
Yves Achdou and Jean-Michel Lasry, \emph{Mean field games for modeling crowd
  motion}, Contributions to Partial Differential Equations and Applications,
  Springer International Publishing, 2019, pp.~17--42.

\bibitem{MR3392611}
Yves Achdou and Mathieu Lauri\`ere, \emph{On the system of partial differential
  equations arising in mean field type control}, Discrete Contin. Dyn. Syst.
  \textbf{35} (2015), no.~9, 3879--3900. \MR{3392611}

\bibitem{MR3498932}
Yves Achdou and Mathieu Lauri{\`e}re, \emph{Mean {F}ield {T}ype {C}ontrol with
  {C}ongestion}, Appl. Math. Optim. \textbf{73} (2016), no.~3, 393--418.
  \MR{3498932}

\bibitem{MR3575615}
\bysame, \emph{Mean {F}ield {T}ype {C}ontrol with {C}ongestion ({II}): {A}n
  augmented {L}agrangian method}, Appl. Math. Optim. \textbf{74} (2016), no.~3,
  535--578. \MR{3575615}

\bibitem{achdoulauriere-2020-mfg-numerical}
\bysame, \emph{Mean field games and applications: Numerical aspects}, Mean
  Field Games, C.I.M.E. Foundation Subseries, vol. 2281, Springer International
  Publishing, 2020.

\bibitem{achdoulaurierelions2020optimal}
Yves Achdou, Mathieu Lauri{\`e}re, and Pierre-Louis Lions, \emph{Optimal
  control of conditioned processes with feedback controls}, Journal de
  Math{\'e}matiques Pures et Appliqu{\'e}es (2020).

\bibitem{MR2928376}
Yves Achdou and Victor Perez, \emph{Iterative strategies for solving linearized
  discrete mean field games systems}, Netw. Heterog. Media \textbf{7} (2012),
  no.~2, 197--217. \MR{2928376}

\bibitem{MR3452251}
Yves Achdou and Alessio Porretta, \emph{Convergence of a finite difference
  scheme to weak solutions of the system of partial differential equations
  arising in mean field games}, SIAM J. Numer. Anal. \textbf{54} (2016), no.~1,
  161--186. \MR{3452251}

\bibitem{MR3765549}
\bysame, \emph{Mean field games with congestion}, Ann. Inst. H. Poincar\'{e}
  Anal. Non Lin\'{e}aire \textbf{35} (2018), no.~2, 443--480. \MR{3765549}

\bibitem{agram2020deep}
Nacira Agram, Azzeddine Bakdi, and Bernt Oksendal, \emph{Deep learning and
  stochastic mean-field control for a neural network model}, Available at SSRN
  3639022 (2020).

\bibitem{al2018solving}
Ali Al-Aradi, Adolfo Correia, Danilo Naiff, Gabriel Jardim, and Yuri Saporito,
  \emph{Solving nonlinear and high-dimensional partial differential equations
  via deep learning}, arXiv preprint arXiv:1811.08782 (2018).

\bibitem{MR3698446}
Noha Almulla, Rita Ferreira, and Diogo Gomes, \emph{Two numerical approaches to
  stationary mean-field games}, Dyn. Games Appl. \textbf{7} (2017), no.~4,
  657--682. \MR{3698446}

\bibitem{anahtarci2019value}
Berkay Anahtarci, Can~Deha Kariksiz, and Naci Saldi, \emph{Value iteration
  algorithm for mean-field games}, arXiv preprint arXiv:1909.01758 (2019).

\bibitem{MR2784835}
Daniel Andersson and Boualem Djehiche, \emph{A maximum principle for {SDE}s of
  mean-field type}, Appl. Math. Optim. \textbf{63} (2011), no.~3, 341--356.
  \MR{2784835}

\bibitem{MR3731033}
Roman Andreev, \emph{Preconditioning the augmented {L}agrangian method for
  instationary mean field games with diffusion}, SIAM J. Sci. Comput.
  \textbf{39} (2017), no.~6, A2763--A2783. \MR{3731033}

\bibitem{angiulifouquelauriere2020unified}
Andrea Angiuli, Jean-Pierre Fouque, and Mathieu Lauri{\`e}re, \emph{Unified
  reinforcement {Q}-learning for mean field game and control problems}, arXiv
  preprint arXiv:2006.13912 (2020).

\bibitem{Angiulietal-2019}
Andrea Angiuli, Christy~V. Graves, Houzhi Li, Jean-Fran\c{c}ois Chassagneux,
  Fran\c{c}ois Delarue, and Ren\'e Carmona, \emph{Cemracs 2017: numerical
  probabilistic approach to mfg}, ESAIM: ProcS \textbf{65} (2019), 84--113.

\bibitem{MR3763083}
Alexander Aurell and Boualem Djehiche, \emph{Mean-field type modeling of
  nonlocal crowd aversion in pedestrian crowd dynamics}, SIAM J. Control Optim.
  \textbf{56} (2018), no.~1, 434--455. \MR{3763083}

\bibitem{balata2019class}
Alessandro Balata, C{\^o}me Hur{\'e}, Mathieu Lauri{\`e}re, Huy{\^e}n Pham, and
  Isaque Pimentel, \emph{A class of finite-dimensional numerically solvable
  {M}c{K}ean-{V}lasov control problems}, ESAIM: Proceedings and Surveys
  \textbf{65} (2019), 114--144.

\bibitem{MR3873029}
Erhan Bayraktar, Amarjit Budhiraja, and Asaf Cohen, \emph{A numerical scheme
  for a mean field game in some queueing systems based on {M}arkov chain
  approximation method}, SIAM J. Control Optim. \textbf{56} (2018), no.~6,
  4017--4044. \MR{3873029}

\bibitem{MR3739204}
Erhan Bayraktar, Andrea Cosso, and Huy\^{e}n Pham, \emph{Randomized dynamic
  programming principle and {F}eynman-{K}ac representation for optimal control
  of {M}c{K}ean-{V}lasov dynamics}, Trans. Amer. Math. Soc. \textbf{370}
  (2018), no.~3, 2115--2160. \MR{3739204}

\bibitem{MR1738163}
Jean-David Benamou and Yann Brenier, \emph{A computational fluid mechanics
  solution to the {M}onge-{K}antorovich mass transfer problem}, Numer. Math.
  \textbf{84} (2000), no.~3, 375--393. \MR{1738163}

\bibitem{MR3395203}
Jean-David Benamou and Guillaume Carlier, \emph{Augmented {L}agrangian methods
  for transport optimization, mean field games and degenerate elliptic
  equations}, J. Optim. Theory Appl. \textbf{167} (2015), no.~1, 1--26.
  \MR{3395203}

\bibitem{MR3644590}
Jean-David Benamou, Guillaume Carlier, and Filippo Santambrogio,
  \emph{Variational mean field games}, Active particles. {V}ol. 1. {A}dvances
  in theory, models, and applications, Model. Simul. Sci. Eng. Technol.,
  Birkh\"{a}user/Springer, Cham, 2017, pp.~141--171. \MR{3644590}

\bibitem{MR3134900}
Alain Bensoussan, Jens Frehse, and Sheung Chi~Phillip Yam, \emph{Mean field
  games and mean field type control theory}, Springer Briefs in Mathematics,
  Springer, New York, 2013.

\bibitem{MR3652408}
\bysame, \emph{On the interpretation of the {M}aster {E}quation}, Stochastic
  Process. Appl. \textbf{127} (2017), no.~7, 2093--2137. \MR{3652408}

\bibitem{MR3797719}
L\'{e}on Bottou, Frank~E. Curtis, and Jorge Nocedal, \emph{Optimization methods
  for large-scale machine learning}, SIAM Rev. \textbf{60} (2018), no.~2,
  223--311. \MR{3797719}

\bibitem{BricenoAriasetalCEMRACS2017}
Luis~M. Brice\~no Arias, Dante Kalise, Ziad Kobeissi, Mathieu Lauri\`ere,
  \'Alvaro Mateos~Gonz\'alez, and Francisco~J. Silva, \emph{On the
  implementation of a primal-dual algorithm for second order time-dependent
  mean field games with local couplings}, ESAIM: ProcS \textbf{65} (2019),
  330--348.

\bibitem{MR3772008}
Luis~M. Brice\~{n}o Arias, Dante Kalise, and Francisco~J. Silva, \emph{Proximal
  methods for stationary mean field games with local couplings}, SIAM J.
  Control Optim. \textbf{56} (2018), no.~2, 801--836. \MR{3772008}

\bibitem{MR3882530}
Simone Cacace, Fabio Camilli, Annalisa Cesaroni, and Claudio Marchi, \emph{An
  ergodic problem for mean field games: qualitative properties and numerical
  simulations}, Minimax Theory Appl. \textbf{3} (2018), no.~2, 211--226.
  \MR{3882530}

\bibitem{cacacecamilligoffi2020policy}
Simone Cacace, Fabio Camilli, and Alessandro Goffi, \emph{A policy iteration
  method for mean field games}, arXiv preprint arXiv:2007.04818 (2020).

\bibitem{MR3601001}
Simone Cacace, Fabio Camilli, and Claudio Marchi, \emph{A numerical method for
  mean field games on networks}, ESAIM Math. Model. Numer. Anal. \textbf{51}
  (2017), no.~1, 63--88. \MR{3601001}

\bibitem{CHM-2017-survey}
Peter~E. Caines, Minyi Huang, and Roland~P. Malham{\'e}, \emph{Handbook of
  dynamic game theory, t. basar and g. zaccour (eds.)}, Springer, Cham, 2017.

\bibitem{camilli2020approximation-timefrac}
Fabio Camilli, Serikbolsyn Duisembay, and Qing Tang, \emph{Approximation of an
  optimal control problem for the time-fractional fokker-planck equation},
  arXiv preprint arXiv:2006.03518 (2020).

\bibitem{cao2020connecting}
Haoyang Cao, Xin Guo, and Mathieu Lauri{\`e}re, \emph{Connecting gans and
  mfgs}, arXiv preprint arXiv:2002.04112 (2020).

\bibitem{Cardaliaguet-2013-notes}
Pierre Cardaliaguet, \emph{Notes on mean field games}, 2013.

\bibitem{MR3967062}
Pierre Cardaliaguet, Fran\c{c}ois Delarue, Jean-Michel Lasry, and Pierre-Louis
  Lions, \emph{The master equation and the convergence problem in mean field
  games}, Annals of Mathematics Studies, vol. 201, Princeton University Press,
  Princeton, NJ, 2019. \MR{3967062}

\bibitem{MR3358627}
Pierre Cardaliaguet and P.~Jameson Graber, \emph{Mean field games systems of
  first order}, ESAIM Control Optim. Calc. Var. \textbf{21} (2015), no.~3,
  690--722. \MR{3358627}

\bibitem{MR3399179}
Pierre Cardaliaguet, P.~Jameson Graber, Alessio Porretta, and Daniela Tonon,
  \emph{Second order mean field games with degenerate diffusion and local
  coupling}, NoDEA Nonlinear Differential Equations Appl. \textbf{22} (2015),
  no.~5, 1287--1317. \MR{3399179}

\bibitem{MR3608094}
Pierre Cardaliaguet and Saeed Hadikhanloo, \emph{Learning in mean field games:
  the fictitious play}, ESAIM Control Optim. Calc. Var. \textbf{23} (2017),
  no.~2, 569--591. \MR{3608094}

\bibitem{MR3556062}
Pierre Cardaliaguet, Alp\'{a}r~R. M\'{e}sz\'{a}ros, and Filippo Santambrogio,
  \emph{First order mean field games with density constraints: pressure equals
  price}, SIAM J. Control Optim. \textbf{54} (2016), no.~5, 2672--2709.
  \MR{3556062}

\bibitem{MR3148086}
Elisabetta Carlini and Francisco~J. Silva, \emph{A fully discrete
  semi-{L}agrangian scheme for a first order mean field game problem}, SIAM J.
  Numer. Anal. \textbf{52} (2014), no.~1, 45--67. \MR{3148086}

\bibitem{MR3392626}
\bysame, \emph{A semi-{L}agrangian scheme for a degenerate second order mean
  field game system}, Discrete Contin. Dyn. Syst. \textbf{35} (2015), no.~9,
  4269--4292. \MR{3392626}

\bibitem{MR3828859}
\bysame, \emph{On the discretization of some nonlinear
  {F}okker-{P}lanck-{K}olmogorov equations and applications}, SIAM J. Numer.
  Anal. \textbf{56} (2018), no.~4, 2148--2177. \MR{3828859}

\bibitem{carmona-AMS-tmp}
Ren\'e Carmona, \emph{Applications of mean field games to economic theory},
  Proc. AMS Short Course, arXiv preprint arXiv:2012.05237 (2020).

\bibitem{MR3395471}
Ren\'{e} Carmona and Fran\c{c}ois Delarue, \emph{Forward-backward stochastic
  differential equations and controlled {M}c{K}ean-{V}lasov dynamics}, Ann.
  Probab. \textbf{43} (2015), no.~5, 2647--2700. \MR{3395471}

\bibitem{MR3752669}
\bysame, \emph{Probabilistic theory of mean field games with applications.
  {I}}, Probability Theory and Stochastic Modelling, vol.~83, Springer, Cham,
  2018, Mean field FBSDEs, control, and games. \MR{3752669}

\bibitem{MR3753660}
\bysame, \emph{Probabilistic theory of mean field games with applications.
  {II}}, Probability Theory and Stochastic Modelling, vol.~84, Springer, Cham,
  2018, Mean field games with common noise and master equations. \MR{3753660}

\bibitem{MR3968548}
Ren\'{e} Carmona, Christy~V. Graves, and Zongjun Tan, \emph{Price of anarchy
  for mean field games}, C{EMRACS} 2017---numerical methods for stochastic
  models: control, uncertainty quantification, mean-field, ESAIM Proc. Surveys,
  vol.~65, EDP Sci., Les Ulis, 2019, pp.~349--383. \MR{3968548}

\bibitem{carmona2019convergence2}
Ren{\'e} Carmona and Mathieu Lauri{\`e}re, \emph{Convergence analysis of
  machine learning algorithms for the numerical solution of mean field control
  and games: Ii--the finite horizon case}, arXiv preprint arXiv:1908.01613
  (2019).

\bibitem{carmona2019convergence1}
\bysame, \emph{Convergence analysis of machine learning algorithms for the
  numerical solution of mean field control and games: {I}--the ergodic case},
  To appear in SIAM Journal on Numerical Analysis (2021).

\bibitem{carmona2019linear}
Ren{\'e} Carmona, Mathieu Lauri{\`e}re, and Zongjun Tan, \emph{Linear-quadratic
  mean-field reinforcement learning: convergence of policy gradient methods},
  arXiv preprint arXiv:1910.04295 (2019).

\bibitem{carmona2019model}
\bysame, \emph{Model-free mean-field reinforcement learning: mean-field {MDP}
  and mean-field {Q}-learning}, arXiv preprint arXiv:1910.12802 (2019).

\bibitem{MR3981375}
Alekos Cecchin, Paolo Dai~Pra, Markus Fischer, and Guglielmo Pelino, \emph{On
  the convergence problem in mean field games: a two state model without
  uniqueness}, SIAM J. Control Optim. \textbf{57} (2019), no.~4, 2443--2466.
  \MR{3981375}

\bibitem{MR2782122}
Antonin Chambolle and Thomas Pock, \emph{A first-order primal-dual algorithm
  for convex problems with applications to imaging}, J. Math. Imaging Vision
  \textbf{40} (2011), no.~1, 120--145. \MR{2782122}

\bibitem{MR3914553}
Jean-Fran\c{c}ois Chassagneux, Dan Crisan, and Fran\c{c}ois Delarue,
  \emph{Numerical method for {FBSDE}s of {M}c{K}ean-{V}lasov type}, Ann. Appl.
  Probab. \textbf{29} (2019), no.~3, 1640--1684. \MR{3914553}

\bibitem{MR3322862}
Paul-Eric Chaudru~de Raynal and Camilo~A. Garcia~Trillos, \emph{A cubature
  based algorithm to solve decoupled {M}c{K}ean-{V}lasov forward-backward
  stochastic differential equations}, Stochastic Process. Appl. \textbf{125}
  (2015), no.~6, 2206--2255. \MR{3322862}

\bibitem{delarue-AMS-tmp}
Fran\c{c}ois Delarue, \emph{Master equation and mean field games}, Proc. AMS
  Short Course (2020).

\bibitem{MR4079435}
Fran\c{c}ois Delarue, Daniel Lacker, and Kavita Ramanan, \emph{From the master
  equation to mean field game limit theory: large deviations and concentration
  of measure}, Ann. Probab. \textbf{48} (2020), no.~1, 211--263. \MR{4079435}

\bibitem{djete2019mckean}
Mao~Fabrice Djete, Dylan Possama{\"\i}, and Xiaolu Tan, \emph{{McKean-Vlasov}
  optimal control: the dynamic programming principle}, arXiv preprint
  arXiv:1907.08860 (2019).

\bibitem{MR3736669}
Weinan E, Jiequn Han, and Arnulf Jentzen, \emph{Deep learning-based numerical
  methods for high-dimensional parabolic partial differential equations and
  backward stochastic differential equations}, Commun. Math. Stat. \textbf{5}
  (2017), no.~4, 349--380. \MR{3736669}

\bibitem{MR1168183}
Jonathan Eckstein and Dimitri~P. Bertsekas, \emph{On the {D}ouglas-{R}achford
  splitting method and the proximal point algorithm for maximal monotone
  operators}, Math. Programming \textbf{55} (1992), no.~3, Ser. A, 293--318.
  \MR{1168183}

\bibitem{Elie2020OnTC}
Romuald Elie, Julien P{\'e}rolat, Mathieu Lauri{\`e}re, Matthieu Geist, and
  Olivier Pietquin, \emph{On the convergence of model free learning in mean
  field games}, AAAI, 2020.

\bibitem{MR724072}
Michel Fortin and Roland Glowinski, \emph{Augmented {L}agrangian methods},
  Studies in Mathematics and its Applications, vol.~15, North-Holland
  Publishing Co., Amsterdam, 1983, Applications to the numerical solution of
  boundary value problems, Translated from the French by B. Hunt and D. C.
  Spicer. \MR{724072}

\bibitem{fouque2019deep}
Jean-Pierre Fouque and Zhaoyu Zhang, \emph{Deep learning methods for mean field
  control problems with delay}, Frontiers in Applied Mathematics and Statistics
  \textbf{6(11)} (2020).

\bibitem{fu2019actor}
Zuyue Fu, Zhuoran Yang, Yongxin Chen, and Zhaoran Wang, \emph{Actor-critic
  provably finds nash equilibria of linear-quadratic mean-field games}, arXiv
  preprint arXiv:1910.07498 (2019).

\bibitem{MR2968782}
Nicolas Gast, Bruno Gaujal, and Jean-Yves Le~Boudec, \emph{Mean field for
  {M}arkov decision processes: from discrete to continuous optimization}, IEEE
  Trans. Automat. Control \textbf{57} (2012), no.~9, 2266--2280. \MR{2968782}

\bibitem{germain2019numerical}
Maximilien Germain, Joseph Mikael, and Xavier Warin, \emph{Numerical resolution
  of mckean-vlasov fbsdes using neural networks}, arXiv preprint
  arXiv:1909.12678 (2019).

\bibitem{MR2137498}
Emmanuel Gobet and R\'{e}mi Munos, \emph{Sensitivity analysis using
  {I}t\^{o}-{M}alliavin calculus and martingales, and application to stochastic
  optimal control}, SIAM J. Control Optim. \textbf{43} (2005), no.~5,
  1676--1713. \MR{2137498}

\bibitem{MR3559742}
Diogo~A. Gomes, Edgard~A. Pimentel, and Vardan Voskanyan, \emph{Regularity
  theory for mean-field game systems}, SpringerBriefs in Mathematics, Springer,
  [Cham], 2016. \MR{3559742}

\bibitem{MR3195844}
Diogo~A. Gomes and Jo\~{a}o Sa\'{u}de, \emph{Mean field games models---a brief
  survey}, Dyn. Games Appl. \textbf{4} (2014), no.~2, 110--154. \MR{3195844}

\bibitem{GomesSaude2018}
Diogo~A. Gomes and Jo{\~a}o Sa{\'u}de, \emph{Numerical methods for finite-state
  mean-field games satisfying a monotonicity condition}, Applied Mathematics
  {\&} Optimization (2018).

\bibitem{GomesYang2018hessian}
{Gomes, Diogo A.} and {Yang, Xianjin}, \emph{The hessian riemannian flow and
  newton\'{}s method for effective hamiltonians and mather measures}, ESAIM:
  M2AN \textbf{54} (2020), no.~6, 1883--1915.

\bibitem{gravesmalhame-AMS-tmp}
Christy Graves and Roland~P. Malham{\'e}, \emph{Mean field games : A paradigm
  for individual-mass interactions}, Proc. AMS Short Course (2020).

\bibitem{gu2019dynamic}
Haotian Gu, Xin Guo, Xiaoli Wei, and Renyuan Xu, \emph{Dynamic programming
  principles for learning {MFCs}}, arXiv preprint arXiv:1911.07314 (2019).

\bibitem{gu2020qlearning}
\bysame, \emph{Mean-field controls with {Q}-learning for cooperative {MARL}:
  Convergence and complexity analysis}, arXiv preprint arXiv:2002.04131 (2020).

\bibitem{guo2019learning}
Xin Guo, Anran Hu, Renyuan Xu, and Junzi Zhang, \emph{Learning mean-field
  games}, in proc. of NeurIPS, 2019.

\bibitem{Hadikhanloo-phdthesis}
Saeed Hadikhanloo, \emph{Learning in mean field games}, Ph.D. thesis,
  University Paris-Dauphine, 2018.

\bibitem{MR4030259}
Saeed Hadikhanloo and Francisco~J. Silva, \emph{Finite mean field games:
  fictitious play and convergence to a first order continuous mean field game},
  J. Math. Pures Appl. (9) \textbf{132} (2019), 369--397. \MR{4030259}

\bibitem{han2016deep-googlecitations}
Jiequn Han and Weinan E, \emph{Deep learning approximation for stochastic
  control problems}, Deep Reinforcement Learning Workshop, NIPS, arXiv preprint
  arXiv:1611.07422 (2016).

\bibitem{MR2352434}
Minyi Huang, Peter~E. Caines, and Roland~P. Malham\'e, \emph{Large-population
  cost-coupled {LQG} problems with nonuniform agents: individual-mass behavior
  and decentralized {$\epsilon$}-{N}ash equilibria}, IEEE Trans. Automat.
  Control \textbf{52} (2007), no.~9, 1560--1571. \MR{2352434}

\bibitem{MR2346927}
Minyi Huang, Roland~P. Malham{\'e}, and Peter~E. Caines, \emph{Large population
  stochastic dynamic games: closed-loop {M}c{K}ean-{V}lasov systems and the
  {N}ash certainty equivalence principle}, Commun. Inf. Syst. \textbf{6}
  (2006), no.~3, 221--251. \MR{2346927}

\bibitem{MR3575619}
Vassili~N. Kolokoltsov and Alain Bensoussan, \emph{Mean-field-game model for
  botnet defense in cyber-security}, Appl. Math. Optim. \textbf{74} (2016),
  no.~3, 669--692. \MR{3575619}

\bibitem{LachapelleWolfram-2011-MFG-congestion-aversion}
Aim{\'e} Lachapelle and Marie-Therese Wolfram, \emph{On a mean field game
  approach modeling congestion and aversion in pedestrian crowds},
  Transportation research part B: methodological \textbf{45} (2011), no.~10,
  1572--1589.

\bibitem{lacker-AMS-tmp}
Daniel Lacker, \emph{Notes for the {AMS} short course on mean field games: The
  convergence problem}, Proc. AMS Short Course (2020).

\bibitem{MR2269875}
Jean-Michel Lasry and Pierre-Louis Lions, \emph{Jeux \`a champ moyen. {I}. {L}e
  cas stationnaire}, C. R. Math. Acad. Sci. Paris \textbf{343} (2006), no.~9,
  619--625. \MR{2269875}

\bibitem{MR2271747}
\bysame, \emph{Jeux \`a champ moyen. {II}. {H}orizon fini et contr\^ole
  optimal}, C. R. Math. Acad. Sci. Paris \textbf{343} (2006), no.~10, 679--684.
  \MR{2271747}

\bibitem{MR2295621}
\bysame, \emph{Mean field games}, Jpn. J. Math. \textbf{2} (2007), no.~1,
  229--260. \MR{2295621}

\bibitem{MR3258261}
Mathieu Lauri{\`e}re and Olivier Pironneau, \emph{Dynamic programming for
  mean-field type control}, C. R. Math. Acad. Sci. Paris \textbf{352} (2014),
  no.~9, 707--713. \MR{3258261}

\bibitem{MR3501391}
\bysame, \emph{Dynamic programming for mean-field type control}, J. Optim.
  Theory Appl. \textbf{169} (2016), no.~3, 902--924. \MR{3501391}

\bibitem{PLL-CDF}
Pierre-Louis Lions, \emph{Cours du {C}oll{\`e}ge de {F}rance},
  \url{https://www.college-de-france.fr/site/en-pierre-louis-lions/_course.htm},
  2007-2011.

\bibitem{MR3420414}
Alp\'{a}r~Rich\'{a}rd M\'{e}sz\'{a}ros and Francisco~J. Silva, \emph{A
  variational approach to second order mean field games with density
  constraints: the stationary case}, J. Math. Pures Appl. (9) \textbf{104}
  (2015), no.~6, 1135--1159. \MR{3420414}

\bibitem{mishra2020model}
Rajesh~K Mishra, Deepanshu Vasal, and Sriram Vishwanath, \emph{Model-free
  reinforcement learning for non-stationary mean field games}, 2020 59th IEEE
  Conference on Decision and Control (CDC), IEEE, 2020, pp.~1032--1037.

\bibitem{motte2019mean}
M{\'e}d{\'e}ric Motte and Huy{\^e}n Pham, \emph{Mean-field markov decision
  processes with common noise and open-loop controls}, arXiv preprint
  arXiv:1912.07883 (2019).

\bibitem{nurbekyan2019fourier}
Levon Nurbekyan et~al., \emph{Fourier approximation methods for first-order
  nonlocal mean-field games}, Portugaliae Mathematica \textbf{75} (2019),
  no.~3, 367--396.

\bibitem{MR4133380}
Marcel Nutz and Yuchong Zhang, \emph{Conditional optimal stopping: a
  time-inconsistent optimization}, Ann. Appl. Probab. \textbf{30} (2020),
  no.~4, 1669--1692. \MR{4133380}

\bibitem{parisini1996neural}
Thomas Parisini and R~Zoppoli, \emph{Neural approximations for multistage
  optimal control of nonlinear stochastic systems}, IEEE Transactions on
  Automatic Control \textbf{41} (1996), no.~6, 889--895.

\bibitem{perrin2020fictitious}
Sarah Perrin, Julien P{\'e}rolat, Mathieu Lauri{\`e}re, Matthieu Geist, Romuald
  Elie, and Olivier Pietquin, \emph{Fictitious play for mean field games:
  Continuous time analysis and applications}, in proc. of NeurIPS, 2020.

\bibitem{MR3619691}
Laurent Pfeiffer, \emph{Numerical methods for mean-field type optimal control
  problems}, Pure Appl. Funct. Anal. \textbf{1} (2016), no.~4, 629--655.
  \MR{3619691}

\bibitem{MR3631380}
Huy\^en Pham and Xiaoli Wei, \emph{Dynamic programming for optimal control of
  stochastic {M}c{K}ean-{V}lasov dynamics}, SIAM J. Control Optim. \textbf{55}
  (2017), no.~2, 1069--1101. \MR{3631380}

\bibitem{ramanan-AMS-tmp}
Kavita Ramanan, \emph{Deviations and fluctuations for mean field games}, Proc.
  AMS Short Course (2020).

\bibitem{MR1451876}
R.~Tyrrell Rockafellar, \emph{Convex analysis}, Princeton Landmarks in
  Mathematics, Princeton University Press, Princeton, NJ, 1997, Reprint of the
  1970 original, Princeton Paperbacks.

\bibitem{ruthotto2020machine}
Lars Ruthotto, Stanley~J Osher, Wuchen Li, Levon Nurbekyan, and Samy~Wu Fung,
  \emph{A machine learning framework for solving high-dimensional mean field
  game and mean field control problems}, Proceedings of the National Academy of
  Sciences \textbf{117} (2020), no.~17, 9183--9193.

\bibitem{MR3874585}
Justin Sirignano and Konstantinos Spiliopoulos, \emph{D{GM}: a deep learning
  algorithm for solving partial differential equations}, J. Comput. Phys.
  \textbf{375} (2018), 1339--1364. \MR{3874585}

\bibitem{subramanianpolicy}
Jayakumar Subramanian and Aditya Mahajan, \emph{Reinforcement learning in
  stationary mean-field games}, in proc. of AAMAS, 2019.

\bibitem{tembine2012meanfieldlearning}
Hamidou Tembine, Raul Tempone, and Pedro Vilanova, \emph{Mean-field learning: a
  survey}, arXiv preprint arXiv:1210.4657 (2012).

\end{thebibliography}

\printindex

\end{document}